\documentclass[11pt]{amsart}
\usepackage{amscd,amssymb,amsthm}
\usepackage{amsxtra}
\usepackage{amssymb, amsmath}
\usepackage{amsthm}
\usepackage{multicol}
\usepackage {hyperref}
\usepackage{latexsym}

\usepackage[condensed,math]{iwona} 

\begin{document}
\title[Single Digit Representations - Inder J. Taneja]{}
\begin{center}
\textbf{\Huge{Single Digit Representations of Natural Numbers}}
\end{center}

\bigskip
\begin{center}
\textbf{\Large{Inder J. Taneja}}\footnote{\textit{Formerly, Professor of Mathematics, Universidade Federal de Santa Catarina, 88.040-900 Florian\'{o}polis, SC, Brazil. e-mail: ijtaneja@gmail.com.}}
\end{center}

\begin{abstract}
 In this work, we established symmetric representation of numbers where one can use any of 9 digits giving the same number. The representations of natural numbers from 0 to 1000 are given using only single digit in all the nine cases, i.e., 1, 2, 3, 4, 5, 6, 7, 8 and 9. This is done only using basic operations: \textit{addition, subtraction, multiplication, potentiation, division}.
  \end{abstract}

\maketitle

\section{\textbf{Introduction}}

Let $a$ be a single digit positive natural numbers, i.e, $a\in \{1,\;2,\;3,\;4,\;5,\;6,\;7,\;8,\;9\}.$ We can always write

\begin{align*}
0& =a-a;&\\
1& =\frac{a}{a};&\\
2& =\frac{a+a}{a};&\\
3& =\frac{a+a+a}{a};&\\
4& =\frac{a+a+a+a}{a};&\\
5& =\frac{a+a+a+a+a}{a},&\\
6& =\frac{a+a+a+a+a+a}{a};&
\end{align*}

\textellipsis

\bigskip
We observe that as number increases, one need more digits to write. But it is not true, for example to write 10, we can write as $10=11-1=2\times 2\times 2+2$. Here we need only 3 digits for 1 and four digits for 2.

\bigskip
Author  \cite{ta1} studied representations of natural numbers using the digits from 1 to 9 in increasing and decreasing orders. For comments see \cite{ab1, ab2, ne1, ne2}. Historical study of numbers and their properties can be found in \cite{du, ma}.  Study of numbers in little different way calling \textit{"selfie numbers"} is given by author \cite{ta2}.

\bigskip
The aim of this work is to write natural numbers from 0 to 1000 in terms of each digits 1, 2, 3, 4, 5, 6, 7, 8 and 9, with \textit{"as less as possible digits''}, using only the basic operations:
\bigskip
\begin{center}
[\textit{addition, subtraction, multiplication, division, potentiation}].
\end{center}

\bigskip
Before proceeding further, here below ares some numbers that can be represented in a symmetric way using any of 9 digits from 1 to 9.

\section{\textbf{Symmetrical Relations}}

Let us consider a series:

\[
f^{n}(10)=10^{n}+10^{n-1}+...+10^{2}+10+10^{0},
\]

then we can write

\[
af^{n}(10)=\underbrace {aaa...a}_{(n+1)-times},
\]

where $a\in \{1,\;2,\;3,\;4,\;5,\;6,\;7,\;8,\;9\}$. 

In particular for $n=3$, we have

\[
f^{3}(10)=a10^{3}+a10^{2}+a10+a=aaaa.
\]

Then,
\[
5=\frac{af^{1}(10)-a}{a+a}=\frac{aa-a}{a+a},
\]

\[
55=\frac{af^{2}(10)-a}{a+a}=\frac{aaa-a}{a+a}
\]

And
\[
6=\frac{af^{1}(10)+a}{a+a}=\frac{aa+a}{a+a},
\]

\[
56=\frac{af^{2}(10)-a}{a+a}=\frac{aaa+a}{a+a}.
\]

\bigskip
Equivalently, one can write the following symmetric representations

\begin{flalign*}
\qquad 5& =\frac{11-1}{1+1}=\frac{22-2}{2+2}=\frac{33-3}{3+3}=\frac{44-4}{4+4}=\frac{55-5}{5+5}
=\frac{66-6}{6+6}=\frac{77-7}{7+7}=\frac{88-8}{8+8}=\frac{99-9}{9+9};&
\end{flalign*}

\begin{flalign*}
\qquad 55& =\frac{111-1}{1+1}=\frac{222-2}{2+2}=\frac{333-3}{3+3}=\frac{444-4}{4+4}=\frac{555-5}{5+5}
=\frac{666-6}{6+6}=\frac{777-7}{7+7}=\frac{888-8}{8+8}=\frac{999-9}{9+9};&
\end{flalign*}

\textellipsis

\begin{flalign*}
\qquad 6& =\frac{11+1}{1+1}=\frac{22+2}{2+2}=\frac{33+3}{3+3}=\frac{44+4}{4+4}=\frac{55+5}{5+5}
=\frac{66+6}{6+6}=\frac{77+7}{7+7}=\frac{88+8}{8+8}=\frac{99+9}{9+9};&
\end{flalign*}

\begin{flalign*}
\qquad 56&=\frac{111+1}{1+1}=\frac{222+2}{2+2}=\frac{333+3}{3+3}=\frac{444+4}{4+4}=\frac{555+5}{5+5}
=\frac{666+6}{6+6}=\frac{777+7}{7+7}=\frac{888+8}{8+8}=\frac{999+9}{9+9};&
\end{flalign*}

\textellipsis

\bigskip
Following the same procedure, we have more symmetries, such as,

\begin{flalign*}
\qquad 11& =\frac{22}{2}=\frac{33}{3}=\frac{44}{4}=\frac{55}{5}=\frac{66}{6}=\frac{77}{7}=
\frac{88}{8}=\frac{99}{9};&
\end{flalign*}

\begin{flalign*}
\qquad 11&=\frac{22+22}{2+2}=\frac{33+33}{3+3}=\frac{44+44}{4+4}=\frac{55+55}{5+5}
=\frac{66+66}{6+6}=\frac{77+77}{7+7}=\frac{88+88}{8+8}=\frac{99+99}{9+9};&
\end{flalign*}

\begin{flalign*}
\qquad 37&=\frac{111}{1+1+1}=\frac{222}{2+2+2}=\frac{333}{3+3+3}=\frac{444}{4+4+4}=
\frac{555}{5+5+5}&\\\\
\qquad  &=\frac{666}{6+6+6}=\frac{777}{7+7+7}=\frac{888}{8+8+8}=\frac{999}{9+9+9};&
\end{flalign*}

\begin{flalign*}
\qquad 100&=\frac{111-11}{1+1}=\frac{222-22}{2+2}=\frac{333-33}{3+3}=\frac{444-44}{4+4}=
\frac{555-55}{5+5}&\\\\
&=\frac{666-66}{6+6}=\frac{777-77}{7+7}=\frac{888-88}{8+8}=\frac{999-99}{9+9};&
\end{flalign*}

\begin{flalign*}
\qquad 101&=\frac{1111}{11}=\frac{2222}{22}=\frac{3333}{33}=\frac{4444}{44}=
\frac{5555}{55}=\frac{6666}{66}=\frac{7777}{77}=\frac{8888}{88}=\frac{9999}{99};&
\end{flalign*}

\begin{flalign*}
\qquad 111& =\frac{222}{2}=\frac{333}{3}=\frac{444}{4}=\frac{555}{5}=\frac{666}{6}=
\frac{777}{7}=\frac{888}{8}=\frac{999}{9};&
\end{flalign*}

\begin{flalign*}
\qquad 925& =\frac{11111-11}{11+1}=\frac{22222-22}{22+2}=\frac{33333-33}{33+3}=\frac{44444-44}{44+4}
=\frac{55555-55}{55+5}&\\\\
\qquad  &=\frac{66666-66}{66+6}=\frac{77777-77}{77+7}=\frac{88888-88}{88+8}=\frac{99999-99}{99+9};&
\end{flalign*}

\begin{flalign*}
\qquad 926&=\frac{11111+1}{11+1}=\frac{22222+2}{22+2}=\frac{33333+3}{33+3}=
\frac{44444+4}{44+4}=\frac{55555+5}{55+5}&\\\\
\qquad &=\frac{66666+6}{66+6}=\frac{77777+7}{77+7}=\frac{88888+8}{88+8}=\frac{99999+9}{99+9};&
\end{flalign*}

\textellipsis

\bigskip
In some cases, one has different symmetries:

\begin{flalign*}
\qquad 1&=\frac{2\times 2}{2+2}, \quad 2=\frac{4\times 4}{4+4}, \quad 3=\frac{6\times 6}{6+6},
\quad 4=\frac{8\times 8}{8+8},&
\end{flalign*}

\begin{flalign*}
\qquad   0&=\frac{1\times 1-1}{1+1}, \quad 1=\frac{3\times 3-3}{3+3}, \quad 2=\frac{5\times 5-5}{5+5}, \quad 3=\frac{7\times 7-7}{7+7}, \quad 4=\frac{9\times 9-9}{9+9}.&
 \end{flalign*}

\begin{flalign*}
\qquad 1&=\frac{1\times 1+1}{1+1},
\quad
2=\frac{3\times 3+3}{3+3},
\quad
3=\frac{5\times 5+5}{5+5},
\quad
4=\frac{7\times 7+7}{7+7},
\quad
5=\frac{9\times 9+9}{9+9}.&
 \end{flalign*}

\begin{flalign*}
\qquad 222&=111+\frac{111}{1},
\quad
333=222+\frac{222}{2},
\quad
444=333+\frac{333}{3},
\quad
555=444+\frac{444}{4}.&
 \end{flalign*}

\begin{flalign*}
\qquad 666&=555+\frac{555}{5},
\quad
777=666+\frac{666}{6},
\quad
888=777+\frac{777}{7},
\quad
999=888+\frac{888}{8}.&
 \end{flalign*}

\begin{flalign*}
\qquad 212&=111+\frac{1111}{11},
\quad
323=222+\frac{2222}{22},
\quad
434=333+\frac{3333}{33},
\quad
545=444+\frac{4444}{44}.&
 \end{flalign*}

\begin{flalign*}
\qquad 656&=555+\frac{5555}{55},
\quad
767=666+\frac{6666}{66},
\quad
878=777+\frac{7777}{77},
\quad
989=888+\frac{8888}{88}.&
 \end{flalign*}

\bigskip
The above symmetric expressions are beautiful to see, but the number of digits used are not necessarily minimum number, for example,

\[
5=4+\frac{4}{4}=\frac{44-4}{4+4}, \mbox{and}\,\, 6=4+\frac{4+4}{4}=\frac{44-4}{4+4}.
\]

\bigskip
In case of number five,  4 is used three times, for number six,  4 is used four times, while symmetric representations uses five times 4 in each case.

\bigskip
The following section deals with representation of natural numbers from 0 to 100 in with minimum possible digits in each case.

\newpage
\section{\textbf{Representations of Natural Numbers Using Different Digits}}

{\footnotesize
\begin{multicols}{4}
\begin{flalign*}
0&= 1-1&\\
&=2-2&\\
&=3-3&\\
&=4-4&\\
&=5-5&\\
&=6-6&\\
&=7-7&\\
&=8-8&\\
&=9-9.&
\end{flalign*}

\begin{flalign*}
1&=1&\\
&=2/2&\\
&=3/3&\\
&=4/4&\\
&=5/5&\\
&=6/6&\\
&=7/7&\\
&=8/8&\\
&=9/9.&
\end{flalign*}

\begin{flalign*}
2&=1+1&\\
&=2&\\
&=3-3/3&\\
&=(4+4)/4&\\
&=(5+5)/5&\\
&=(6+6)/6&\\
&=(7+7)/7&\\
&=(8+8)/8&\\
&=(9+9)/9.&
\end{flalign*}

\begin{flalign*}
3&=1+1+1&\\
&=2+2/2&\\
&=3&\\
&=4-4/4&\\
&=5-(5+5)/5&\\
&=6\times 6/(6+6)&\\
&=(7+7+7)/7&\\
&=88/8-8&\\
&=(9+9+9)/9.&
\end{flalign*}

\begin{flalign*}
4&=1+1+1+1&\\
&=2+2&\\
&=3+3/3&\\
&=4&\\
&=5-5/5&\\
&=6-(6+6)/6&\\
&=77/7-7&\\
&=8\times 8/(8+8)&\\
&=(9\times9-9)/(9+9).&
\end{flalign*}

\begin{flalign*}
5&=(11-1)/(1+1)&\\
&=2+2+2/2&\\
&=3+3-3/3&\\
&=4+4/4&\\
&=5&\\
&=6-6/6&\\
&=7-(7+7)/7&\\
&=(88-8)/(8+8)&\\
&=(99-9)/(9+9).&
\end{flalign*}

\begin{flalign*}
6&=(1+1)\times (1+1+1)&\\
&=2+2+2&\\
&=3+3&\\
&=4+(4+4)/4&\\
&=5+5/5&\\
&=6&\\
&=7-7/7&\\
&=8-(8+8)/8&\\
&=(99+9)/(9+9).&
\end{flalign*}

\begin{flalign*}
7&=(1+1)\times (1+1+1)+1&\\
&=2+2+2+2/2&\\
&=3+3+3/3&\\
&=4+4-4/4&\\
&=5+(5+5)/5&\\
&=6+6/6&\\
&=7&\\
&=8-8/8&\\
&=9-(9+9)/9.&
\end{flalign*}

\begin{flalign*}
8&=(1+1)^{(1+1+1)}&\\
&=2\times (2+2)&\\
&=3\times 3-3/3&\\
&=4+4&\\
&=5+5-(5+5)/5&\\
&=6+(6+6)/6&\\
&=7+7/7&\\
&=8&\\
&=9-9/9.&
\end{flalign*}

\begin{flalign*}
9&=11-1-1&\\
&=(2+2/2)^2&\\
&=3\times 3&\\
&=4+4+4/4&\\
&=5+5-5/5&\\
&=6+6\times 6/(6+6)&\\
&=7+(7+7)/7&\\
&=8+8/8&\\
&=9.&
\end{flalign*}

\begin{flalign*}
10&=11-1&\\
&=2\times 2 \times 2+2&\\
&=3\times 3+3/3&\\
&=(44-4)/4&\\
&=5+5&\\
&=(66-6)/6&\\
&=(77-7)/7&\\
&=(88-8)/8&\\
&=9+9/9.&
\end{flalign*}

\begin{flalign*}
11&=11&\\
&=22/2&\\
&=33/3&\\
&=44/4&\\
&=55/5&\\
&=66/6&\\
&=77/7&\\
&=88/8&\\
&=99/9.&
\end{flalign*}

\begin{flalign*}
12&=11+1&\\
&=2\times (2+2+2)&\\
&=3+3\times 3&\\
&=4+4+4&\\
&=6+6&\\
&=(55+5)/5&\\
&=(77+7)/7&\\
&=(88+8)/8&\\
&=(99+9)/9.&
\end{flalign*}

\begin{flalign*}
13&=11+1+1&\\
&=2+22/2&\\
&=3+3\times 3+3/3&\\
&=4+4+4+4/4&\\
&=(55+5+5)/5&\\
&=6+6+6/6&\\
&=7+7-7/7&\\
&=(88+8+8)/8&\\
&=(9+99+9)/9.&
\end{flalign*}

\begin{flalign*}
14&=11+1+1+1&\\
&=2^{(2+2)}-2&\\
&=3+33/3&\\
&=4+(44-4)/4&\\
&=5+5+5-5/5&\\
&=6+6+(6+6)/6&\\
&=7+7&\\
&=8+8-(8+8)/8&\\
&=9+(99-9)/(9+9).&
\end{flalign*}

\begin{flalign*}
15&=11+1+1+1+1&\\
&=2+2+22/2&\\
&=3+3+3\times 3&\\
&=4+44/4&\\
&=5+5+5&\\
&=6+6+6\times 6/(6+6)&\\
&=7+7+7/7&\\
&=8+8-8/8&\\
&=9+(99+9)/(9+9).&
\end{flalign*}

\begin{flalign*}
16&=(1+1)^{(1+1+1+1)}&\\
&=2^{(2+2)}&\\
&=3^3-33/3&\\
&=4\times 4&\\
&=5+55/5&\\
&=6+(66-6)/6&\\
&=7+7+(7+7)/7&\\
&=8+8&\\
&=9+9-(9+9)/9.&
\end{flalign*}

\begin{flalign*}
17&=(1+1)^{(1+1+1+1)}+1&\\
&=2^{(2+2)}+2/2&\\
&=3+3+33/3&\\
&=4\times 4+4/4&\\
&=5+(55+5)/5&\\
&=6+66/6&\\
&=7+(77-7)/7&\\
&=8+8+8/8&\\
&=9+9-9/9.&
\end{flalign*}

\begin{flalign*}
18&=(1+1)\times (11-1-1)&\\
&=2^{(2+2)}+2&\\
&=3\times (3+3)&\\
&=4\times 4+(4+4)/4&\\
&=5+(55+5+5)/5&\\
&=6+6+6&\\
&=7+77/7&\\
&=8+(88-8)/8&\\
&=9+9.&
\end{flalign*}

\begin{flalign*}
19&=(1+1)\times (11-1)-1&\\
&=22-2-2/2&\\
&=3\times (3+3)+3/3&\\
&=4+4+44/4&\\
&=5\times 5-5-5/5&\\
&=6+6+6+6/6&\\
&=7+(77+7)/7&\\
&=8+88/8&\\
&=9+9+9/9.&
\end{flalign*}
\end{multicols}
}

{\footnotesize
\begin{multicols}{3}
\begin{flalign*}
20&=(1+1)\times (11-1)&\\
&=22-2&\\
&=3\times 3+33/3&\\
&=4+4\times 4&\\
&=5\times 5-5&\\
&=6+6+6+(6+6)/6&\\
&=7+7+7-7/7&\\
&=8+(88+8)/8&\\
&=9+99/9.&
\end{flalign*}

\begin{flalign*}
21&=11+11-1&\\
&=22-2/2&\\
&=3\times (3+3)+3&\\
&=4+4\times 4+4/4&\\
&=5+5+55/5&\\
&=6\times (6\times 6+6)/(6+6)&\\
&=7+7+7&\\
&=(88+88-8)/8&\\
&=9+(99+9)/9.&
\end{flalign*}

\begin{flalign*}
22&=11+11&\\
&=22&\\
&=(33+33)/3&\\
&=(44+44)/4&\\
&=(55+55)/5&\\
&=(66+66)/6&\\
&=(77+77)/7&\\
&=(88+88)/8&\\
&=(99+99)/9.&
\end{flalign*}

\begin{flalign*}
23&=11+11+1&\\
&=22+2/2&\\
&=3^3-3-3/3&\\
&=4+4+4+44/4&\\
&=5\times 5-(5+5)/5&\\
&=6+6+66/6&\\
&=(77+77+7)/7&\\
&=8+8+8-8/8&\\
&=(99+99+9)/9.&
\end{flalign*}

\begin{flalign*}
24&=(1+1)\times (11+1)&\\
&=22+2&\\
&=3^3-3&\\
&=4+4+4\times 4&\\
&=5\times 5-5/5&\\
&=6+6+6+6&\\
&=7+7+(77-7)/7&\\
&=8+8+8&\\
&=(99+99+9+9)/9.&
\end{flalign*}

\begin{flalign*}
25&=(1+1)\times (11+1)+1&\\
&=22+2+2/2&\\
&=3^3-3+3/3&\\
&=4+4+4\times 4+4/4&\\
&=5\times 5&\\
&=6\times 6-66/6&\\
&=7+7+77/7&\\
&=8+8+8+8/8&\\
&=(9+9-(9+9)/9))+9.&
\end{flalign*}

\begin{flalign*}
26&=(1+1)\times (11+1+1)&\\
&=22+2+2&\\
&=3^3-3/3&\\
&=4+44\times 4/(4+4)&\\
&=5\times 5+5/5&\\
&=6\times 6-(66-6)/6&\\
&=7+7+(77+7)/7&\\
&=8+8+(88-8)/8&\\
&=9+9+9-9/9.&
\end{flalign*}

\begin{flalign*}
27&=(1+1+1)^{(1+1+1)}&\\
&=22+2+2+2/2&\\
&=3^3&\\
&=4\times 4+44/4&\\
&=5\times 5+(5+5)/5&\\
&=6\times 66/(6+6)-6&\\
&=77-7\times 7-7/7&\\
&=8+8+88/8&\\
&=9+9+9.&
\end{flalign*}

\begin{flalign*}
28&=(1+1+1)^{(1+1+1)}+1&\\
&=22+2+2+2&\\
&=3^3+3/3&\\
&=44-4\times 4&\\
&=5\times 5+5-(5+5)/5&\\
&=6+(66+66)/6&\\
&=7\times (77/7-7)&\\
&=8+8+(88+8)/8&\\
&=9+9+9+9/9.&
\end{flalign*}

\begin{flalign*}
29&=(1+1+1)\times (11-1)-1&\\
&=22+2+2+2+2/2&\\
&=3+3^3-3/3&\\
&=44-4\times 4+4/4&\\
&=5\times 5+5-5/5&\\
&=6\times 6-6-6/6&\\
&=77-7\times 7+7/7&\\
&=8+(88+88-8)/8&\\
&=9+9+99/9.&
\end{flalign*}

\begin{flalign*}
30&=(1+1+1)\times (11-1)&\\
&=22+2\times (2+2)&\\
&=3+3^3&\\
&=4\times (4+4)-(4+4)/4&\\
&=5\times 5+5&\\
&=6\times 6-6&\\
&=77-7\times 7+(7+7)/7&\\
&=8+(88+88)/8&\\
&=999/9-9\times9.&
\end{flalign*}

\begin{flalign*}
31&=(1+1+1)\times (11-1)+1&\\
&=22+(2+2/2)^2&\\
&=3+3^3+3/3&\\
&=4\times (4+4)-4/4&\\
&=5\times 5+5+5/5&\\
&=6\times 6-6+6/6&\\
&=7\times 7-7-77/7&\\
&=8+8+8+8-8/8&\\
&=9+(99+99)/9.&
\end{flalign*}

\begin{flalign*}
32&=11\times (1+1+1)-1&\\
&=2\times 2^{(2+2)}&\\
&=33-3/3&\\
&=4\times (4+4)&\\
&=((5+5)/5)^5&\\
&=6\times 6-6+(6+6)/6&\\
&=7+7+7+77/7&\\
&=8+8+8+8&\\
&=9+(99+99+9)/9.&
\end{flalign*}

\begin{flalign*}
33&=11\times (1+1+1)&\\
&=22+22/2&\\
&=33&\\
&=4\times (4+4)+4/4&\\
&=((5+5)/5)^5+5/5&\\
&=6\times 66/(6+6)&\\
&=(77+77+77)/7&\\
&=8+8+8+8+8/8&\\
&=99\times(9+9+9)/(9\times9).&
\end{flalign*}

\begin{flalign*}
34&=11\times (1+1+1)+1&\\
&=2+2\times 2^{(2+2)}&\\
&=33+3/3&\\
&=44-(44-4)/4&\\
&=5\times 5+5-5/5+5&\\
&=6\times 6-(6+6)/6&\\
&=777/7-77&\\
&=8+8+8+(88-8)/8&\\
&=((9+9)\times(9+9)+9)/9.&
\end{flalign*}

\begin{flalign*}
35&=11\times (1+1+1)+1+1&\\
&=22+2+22/2&\\
&=3+33-3/3&\\
&=4+4\times (4+4)-4/4&\\
&=5\times 5+5+5&\\
&=6\times 6-6/6&\\
&=7\times 7-7-7&\\
&=8+8+8+88/8&\\
&=9+9+9+9-9/9.&
\end{flalign*}

\begin{flalign*}
36&=(1+1+1)\times (11+1)&\\
&=(2+2+2)^2&\\
&=3+33&\\
&=4+4\times (4+4)&\\
&=5\times 5+55/5&\\
&=6\times 6&\\
&=7\times 7-7-7+7/7&\\
&=88\times 8/(8+8)-8&\\
&=9+9+9+9.&
\end{flalign*}

\begin{flalign*}
37&=111/(1+1+1)&\\
&=(2+2+2)^2+2/2&\\
&=3+33+3/3&\\
&=4+4\times (4+4)+4/4&\\
&=5+((5+5)/5)^5&\\
&=6\times 6+6/6&\\
&=777/(7+7+7)&\\
&=888/(8+8+8)&\\
&=999/(9+9+9).&
\end{flalign*}

\begin{flalign*}
38&=111/(1+1+1)+1&\\
&=(2+2+2)^2+2&\\
&=3^3+33/3&\\
&=44-4-(4+4)/4&\\
&=5+((5+5)/5)^5+5/5&\\
&=6\times 6+(6+6)/6&\\
&=7\times 7-77/7&\\
&=8+8+(88+88)/8&\\
&=9+9+9+99/9.&
\end{flalign*}

\begin{flalign*}
39&=(1+1+1)\times (11+1+1)&\\
&=2\times (22-2)-2/2&\\
&=3+3+33&\\
&=44-4-4/4&\\
&=55-5-55/5&\\
&=6+6\times 66/(6+6)&\\
&=7\times 7-(77-7)/7&\\
&=8\times 8-8-8-8-8/8&\\
&=9+9+9+(99+9)/9.&
\end{flalign*}

\begin{flalign*}
40&=(1+1)\times (1+1)\times (11-1)&\\
&=2\times (22-2)&\\
&=3+3+33+3/3&\\
&=44-4&\\
&=5\times 5+5+5+5&\\
&=6\times 6+6-(6+6)/6&\\
&=7\times 7-7-(7+7)/7&\\
&=8\times (8+8)-88&\\
&=(9\times9\times9-9)/(9+9).&
\end{flalign*}

\begin{flalign*}
41&=(1+1+1)\times (11-1)+11&\\
&=2\times (22-2)+2/2&\\
&=3+3^3+33/3&\\
&=44-4+4/4&\\
&=5\times 5+5+55/5&\\
&=6\times 6+6-6/6&\\
&=7\times 7-7-7/7&\\
&=8\times 8-8-8-8+8/8&\\
&=(9\times9\times9+9)/(9+9).&
\end{flalign*}

\begin{flalign*}
42&=(1+1)\times (11+11-1)&\\
&=2\times 22-2&\\
&=3\times 3+33&\\
&=44-(4+4)/4&\\
&=5+5+((5+5)/5)^5&\\
&=6\times 6+6&\\
&=7\times 7-7&\\
&=8\times 8-(88+88)/8&\\
&=9+9\times99/(9+9+9).&
\end{flalign*}

\begin{flalign*}
43&=(1+1)\times (11+11)-1&\\
&=2\times 22-2/2&\\
&=3\times 3+33+3/3&\\
&=44-4/4&\\
&=55-(55+5)/5&\\
&=6\times 6+6+6/6&\\
&=7\times 7-7+7/7&\\
&=8+8+8+8+88/8&\\
&=9\times9-9-9-9-99/9.&
\end{flalign*}

\begin{flalign*}
44&=(1+1)\times (11+11)&\\
&=2\times 22&\\
&=33+33/3&\\
&=44&\\
&=55-55/5&\\
&=6\times 6+6+(6+6)/6&\\
&=7\times 7-7+(7+7)/7&\\
&=88\times 8/(8+8)&\\
&=99\times(9-9/9)/(9+9).&
\end{flalign*}

\begin{flalign*}
45&=(1+1)\times (11+11)+1&\\
&=2\times 22+2/2&\\
&=3+3\times 3+33&\\
&=44+4/4&\\
&=55-5-5&\\
&=666/6-66&\\
&=7\times 7+7-77/7&\\
&=8\times 8-8-88/8&\\
&=9+9+9+9+9.&
\end{flalign*}

\begin{flalign*}
46&=(1+1)\times (11+11+1)&\\
&=2\times 22+2&\\
&=3+3\times 3+33+3/3&\\
&=44+(4+4)/4&\\
&=55-5-5+5/5&\\
&=6\times 6+(66-6)/6&\\
&=7\times 7-(7+7+7)/7&\\
&=8\times 8-8-8-(8+8)/8&\\
&=9+9+9+9+9+9/9.&
\end{flalign*}

\begin{flalign*}
47&=(1+1)\times (11+11+1)+1&\\
&=2\times 22+2+2/2&\\
&=3+33+33/3&\\
&=4+44-4/4&\\
&=5+5+5+((5+5)/5)^5&\\
&=66+6\times 6/6&\\
&=7\times 7-((7+7)/7)&\\
&=888/8-8\times 8&\\
&=9+9+9+9+99/9.&
\end{flalign*}

\begin{flalign*}
48&=(1+1)\times (1+1)\times (11+1)&\\
&=2\times (22+2)&\\
&=3\times 3^3-33&\\
&=4+44&\\
&=55-5-(5+5)/5&\\
&=6\times 6+6+6&\\
&=7\times 7-7/7&\\
&=8\times 8-8-8&\\
&=9+9-9\times9+999/9.&
\end{flalign*}

\begin{flalign*}
49&=(11-1)^{(1+1)}/(1+1)-1&\\
&=2\times (22+2)+2/2&\\
&=3^3+33-33/3&\\
&=4+44+4/4&\\
&=55-5-5/5&\\
&=6\times 6+6+6+6/6&\\
&=7\times 7&\\
&=8\times 8-8-8+8/8&\\
&=(9\times99-9)/(9+9).&
\end{flalign*}

\begin{flalign*}
50&=(11-1)^{(1+1)}/(1+1)&\\
&=2\times (22+2)+2&\\
&=3+3+33+33/3&\\
&=4+44+(4+4)/4&\\
&=5\times (5+5)&\\
&=6\times 6+6+6+(6+6)/6&\\
&=7\times 7+7/7&\\
&=8\times 8-8-8+(8+8)/8&\\
&=(9\times99+9)/(9+9).&
\end{flalign*}

\begin{flalign*}
51&=(11-1)^{(1+1)}/(1+1)+1&\\
&=2\times (22+2)+2+2/2&\\
&=3^3+3^3-3&\\
&=4+4+44-4/4&\\
&=55-5+5/5&\\
&=6-66+666/6&\\
&=7\times 7+(7+7)/7&\\
&=8\times 8-8-8-8+88/8&\\
&=(999-9\times9)/(9+9).&
\end{flalign*}

\begin{flalign*}
52&=(1+1)\times (1+1)\times (11+1+1)&\\
&=2\times (22+2+2)&\\
&=3^3+3^3-3+3/3&\\
&=4+4+44&\\
&=55-5+(5+5)/5&\\
&=((6+6)/6)^6-6-6&\\
&=7\times 7+(7+7+7)/7&\\
&=8+88\times 8/(8+8)&\\
&=9\times9-9-9-99/9.&
\end{flalign*}

\begin{flalign*}
53&=(111-1)/(1+1)-1-1&\\
&=2\times (22+2+2)+2/2&\\
&=3^3+3^3-3/3&\\
&=(4^4-44)/4&\\
&=55-(5+5)/5&\\
&=6\times 6+6+66/6&\\
&=7\times 7+77/7-7&\\
&=8\times 8-88/8&\\
&=9\times9-9-9-9-9/9.&
\end{flalign*}

\begin{flalign*}
54&=(111-1)/(1+1)-1&\\
&=2\times (22+2+2)+2&\\
&=3\times 3\times (3+3)&\\
&=44+(44-4)/4&\\
&=55-5/5&\\
&=66-6-6&\\
&=7\times 7+7-(7+7)/7&\\
&=8\times 8-(88-8)/8&\\
&=9\times9-9-9-9.&
\end{flalign*}

\begin{flalign*}
55&=(111-1)/(1+1)&\\
&=2\times 22+22/2&\\
&=3^3+3^3+3/3&\\
&=44+44/4&\\
&=55&\\
&=66-66/6&\\
&=7\times 7+7-7/7&\\
&=8\times 8-8-8/8&\\
&=(999-9)/(9+9).&
\end{flalign*}

\begin{flalign*}
56&=(111+1)/(1+1)&\\
&=(222+2)/(2+2)&\\
&=(333+3)/(3+3)&\\
&=44+4+4+4&\\
&=55+5/5&\\
&=(666+6)/(6+6).&\\
&=7\times 7+7&\\
&=8\times 8-8&\\
&=(999+9)/(9+9).&
\end{flalign*}

\begin{flalign*}
57&=(111+1)/(1+1)+1&\\
&=2\times 22+2+22/2&\\
&=3^3+3^3+3&\\
&=4+(4^4-44)/4&\\
&=55+(5+5)/5&\\
&=66+(6+6-66)/6&\\
&=7\times 7+7+7/7&\\
&=8\times 8-8+8/8&\\
&=(999+9)/(9+9)+9/9.&
\end{flalign*}

\begin{flalign*}
58&=(111+1)/(1+1)+1+1&\\
&=(2+2+2)^2+22&\\
&=(3+3/3)^3-3-3&\\
&=(4^4-4-4)/4-4&\\
&=55+5-(5+5)/5&\\
&=((6+6)/6)^6-6&\\
&=7\times 7+7+(7+7)/7&\\
&=8\times 8-8+(8+8)/8&\\
&=9+(9\times99-9)/(9+9).&
\end{flalign*}

\begin{flalign*}
59&=(11^{(1+1)}-1)/(1+1)-1&\\
&=(4^4-4)/4-4&\\
&=2\times (22+2)+22/2&\\
&=3^3+33-3/3&\\
&=55+5-5/5&\\
&=66+6-6/6&\\
&=77-7-77/7&\\
&=8\times 8-8-8+88/8&\\
&=9+(9\times99+9)/(9+9).&
\end{flalign*}

\begin{flalign*}
60&=(11^{(1+1)}-1)/(1+1)&\\
&=2\times (2\times (2+2)+22)&\\
&=3^3+33&\\
&=4\times 4+44&\\
&=55+5&\\
&=66-6&\\
&=7\times 7+77/7&\\
&=8+8+88\times 8/(8+8)&\\
&=9\times9-9-(99+9)/9.&
\end{flalign*}

\begin{flalign*}
61&=(11^{(1+1)}+1)/(1+1)&\\
&=(3+3/3)^3-3&\\
&=(4^4+4)/4-4&\\
&=2^{(2+2+2)}-2-2/2&\\
&=55+5+5/5&\\
&=66-6+6/6&\\
&=7\times 7+(77+7)/7&\\
&=8\times 8+8-88/8&\\
&=9\times9-99/9-9.&
\end{flalign*}

\begin{flalign*}
62&=(11^{(1+1)}+1)/(1+1)+1&\\
&=2^{(2+2+2)}-2&\\
&=3+3^3+33-3/3&\\
&=(4^4-4-4)/4&\\
&=55+5+(5+5)/5&\\
&=66-6+(6+6)/6&\\
&=777/7-7\times 7&\\
&=8\times 8-(8+8)/8&\\
&=9\times9-9-9-9/9.&
\end{flalign*}

\begin{flalign*}
63&=(1+1+1)\times (11+11-1)&\\
&=2^{(2+2+2)}-2/2&\\
&=3+3^3+33&\\
&=(4^4-4)/4&\\
&=(5+5^5/5)/(5+5)&\\
&=66-6\times 6/(6+6)&\\
&=77-7-7&\\
&=8\times 8-8/8&\\
&=9\times9-9-9.&\\
\end{flalign*}

\begin{flalign*}
64&=(1+1)^{((1+1)\times (1+1+1))}&\\
&=2^{(2+2+2)}&\\
&=4\times 4\times 4&\\
&=(3+3/3)^3&\\
&=55+5+5-5/5&\\
&=((6+6)/6)^6&\\
&=77-7-7+7/7&\\
&=8\times 8&\\
&=(9-9/9)\times(9-9/9).&
\end{flalign*}
\end{multicols}
}

{\footnotesize
\begin{multicols}{3}
\begin{flalign*}
65&=(1+1)^{((1+1)\times (1+1+1))}+1&\\
&=2^{(2+2+2)}+2/2&\\
&=(3+3/3)^3+3/3&\\
&=(4^4+4)/4&\\
&=55+5+5&\\
&=66-6/6&\\
&=77-(77+7)/7&\\
&=8\times 8+8/8&\\
&=9\times9-9-9+(9+9)/9.&
\end{flalign*}

\begin{flalign*}
66&=11\times (1+1)\times (1+1+1)&\\
&=2^{(2+2+2)}+2&\\
&=33+33&\\
&=(4^4+4+4)/4&\\
&=55+55/5&\\
&=66&\\
&=77-77/7&\\
&=8\times 8+(8+8)/8&\\
&=99\times(99+9)/(9\times(9+9)).&
\end{flalign*}

\begin{flalign*}
67&=11\times (1+1)\times (1+1+1)+1&\\
&=2^{(2+2+2)}+2+2/2&\\
&=(3+3/3)^3+3&\\
&=4+(4^4-4)/4&\\
&=55+(55+5)/5&\\
&=66+6/6&\\
&=77-(77-7)/7&\\
&=8\times 8-8+88/8&\\
&=9\times9-(99+9+9+9)/9.&
\end{flalign*}

\begin{flalign*}
68&=(1+1)\times (11\times (1+1+1)+1)&\\
&=2^{(2+2+2)}+2+2&\\
&=((3+3)^3-3)/3-3&\\
&=4+4\times 4\times 4&\\
&=5+(5^5/5+5)/(5+5)&\\
&=66+(6+6)/6&\\
&=77-7-(7+7)/7&\\
&=88-8-(88+8)/8&\\
&=9\times9-(99+9+9)/9.&
\end{flalign*}

\begin{flalign*}
69&=(1+1+1)\times (11+11+1)&\\
&=(22+2/2)\times (2+2/2)&\\
&=33+33+3&\\
&=4+(4^4+4)/4&\\
&=55+5+5+5-5/5&\\
&=66+6\times 6/(6+6)&\\
&=77-7-7/7&\\
&=88-8-88/8&\\
&=9\times9-(99+9)/9.&
\end{flalign*}

\begin{flalign*}
70&=(11-1-1)^{(1+1)}-11&\\
&=2\times (22+2)+22&\\
&=(3/3+3)^3+3+3&\\
&=4+(4^4+4+4)/4&\\
&=55+5+5+5&\\
&=6+((6+6)/6)^6&\\
&=77-7&\\
&=8\times 8+8-(8+8)/8&\\
&=9\times9-99/9.&
\end{flalign*}

\begin{flalign*}
71&=(11+1)^{(1+1)}/(1+1)-1&\\
&=2\times (2+2+2)^2-2/2&\\
&=((3+3)^3-3)/3&\\
&=4+4+(4^4-4)/4&\\
&=55+5+55/5&\\
&=66+6-6/6&\\
&=77-7+7/7&\\
&=8\times 8+8-8/8&\\
&=9\times9-9-9/9.&
\end{flalign*}

\begin{flalign*}
72&=(11+1)^{(1+1)}/(1+1)&\\
&=2\times (2+2+2)^2&\\
&=3\times (3^3-3)&\\
&=4+4+4\times 4\times 4&\\
&=55+5+(55+5)/5&\\
&=66+6&\\
&=77-7+(7+7)/7&\\
&=8+8\times 8&\\
&=9\times9-9.&
\end{flalign*}

\begin{flalign*}
73&=(11+1)^{(1+1)}/(1+1)+1&\\
&=2\times (2+2+2)^2+2/2&\\
&=((3+3)^3+3)/3&\\
&=4+4+(4^4+4)/4&\\
&=5\times (5+5+5)-(5+5)/5&\\
&=66+6+6/6&\\
&=77+7-77/7&\\
&=8\times 8+8+8/8&\\
&=9\times9-9+9/9.&
\end{flalign*}

\begin{flalign*}
74&=(1+1)\times 111/(1+1+1)&\\
&=2\times (2+2+2)^2+2&\\
&=3+((3+3)^3-3)/3&\\
&=(4^4+44-4)/4&\\
&=5\times (5+5+5)-5/5&\\
&=66+6+(6+6)/6&\\
&=77-(7+7+7)/7&\\
&=8\times 8+8+(8+8)/8&\\
&=9\times9-(9-(9+9)/9).&
\end{flalign*}

\begin{flalign*}
75&=(1+1)\times 111/(1+1+1)+1&\\
&=2^{(2+2+2)}+22/2&\\
&=3+3\times (3^3-3)&\\
&=(44+4^4)/4&\\
&=5\times (5+5+5)&\\
&=666/6-6\times 6&\\
&=77-(7+7)/7&\\
&=8\times 8+88/8&\\
&=9\times9-(99+9)/(9+9).&
\end{flalign*}

\begin{flalign*}
76&=(1+1)\times (111/(1+1+1)+1)&\\
&=2\times ((2+2+2)^2+2)&\\
&=3+(3+(3+3)^3)/3&\\
&=44+4\times (4+4)&\\
&=5\times (5+5+5)+5/5&\\
&=6+6+((6+6)/6)^6&\\
&=77-7/7&\\
&=88-(88+8)/8&\\
&=9\times9-(99-9)/(9+9).&
\end{flalign*}

\begin{flalign*}
77&=11\times ((1+1)\times (1+1+1)+1)&\\
&=2\times 2\times 22-22/2&\\
&=3\times 3^3-3-3/3&\\
&=(4-4/4)^4-4&\\
&=55+(55+55)/5&\\
&=66+66/6&\\
&=77&\\
&=88-88/8&\\
&=9\times9-(9+9+9+9)/9.&
\end{flalign*}

\begin{flalign*}
78&=111-11\times (1+1+1)&\\
&=2\times 2\times (22-2)-2&\\
&=3\times 3^3-3&\\
&=4+(4^4-4+44)/4&\\
&=5\times 5+55-(5+5)/5&\\
&=66+6+6&\\
&=77+7/7&\\
&=88-(88-8)/8&\\
&=9\times9-(9+9+9)/9.&
\end{flalign*}

\begin{flalign*}
79&=(11-1-1)^{(1+1)}-1-1&\\
&=(2+2/2)^{(2+2)}-2&\\
&=3\times 3^3-3+3/3&\\
&=4+(44+4^4)/4&\\
&=5\times 5+55-5/5&\\
&=66+6+6+6/6&\\
&=77+(7+7)/7&\\
&=88-8-8/8&\\
&=9\times9-(9+9)/9.&
\end{flalign*}

\begin{flalign*}
80&=(11-1-1)^{(1+1)}-1&\\
&=2\times 2\times (22-2)&\\
&=3\times 3^3-3/3&\\
&=4\times (4\times 4+4)&\\
&=5\times 5+55&\\
&=66+6+6+(6+6)/6&\\
&=77+(7+7+7)/7&\\
&=88-8&\\
&=9\times9-9/9.&
\end{flalign*}

\begin{flalign*}
81&=(11-1-1)^{(1+1)}&\\
&=(2+2/2)^{(2+2)}&\\
&=3\times 3^3&\\
&=(4-4/4)^4&\\
&=5\times 5+55+5/5&\\
&=6-66+6\times 66/6&\\
&=77/7-7+77&\\
&=88-8+8/8&\\
&=9\times9.&
\end{flalign*}

\begin{flalign*}
82&=(11-1-1)^{(1+1)}+1&\\
&=2\times 2\times (22-2)+2&\\
&=3\times 3^3+3/3&\\
&=(4-4/4)^4+4/4&\\
&=5\times 5+55+(5+5)/5&\\
&=6+6+6+((6+6)/6)^6&\\
&=77+7-(7+7)/7&\\
&=88-8+(8+8)/8&\\
&=9\times9+9/9.&
\end{flalign*}

\begin{flalign*}
83&=(11-1-1)^{(1+1)}+1+1&\\
&=(2+2/2)^{(2+2)}+2&\\
&=3+3\times 3^3-3/3&\\
&=4+4+(44+4^4)/4&\\
&=5\times 5+55+5-(5+5)/5&\\
&=66+6+66/6&\\
&=77+7-7/7&\\
&=88/8+8+8\times 8&\\
&=9\times9+(9+9)/9.&
\end{flalign*}

\begin{flalign*}
84&=(11-1-1)^{(1+1)}+1+1+1&\\
&=2\times (2\times 22-2)&\\
&=3+3\times 3^3&\\
&=4+4\times (4\times 4+4)&\\
&=5\times 5+55+5-5/5&\\
&=66+6+6+6&\\
&=77+7&\\
&=88-8\times 8/(8+8)&\\
&=9\times9+(9+9+9)/9.&
\end{flalign*}

\begin{flalign*}
85&=111-(1+1)\times (11+1+1)&\\
&=2\times 2\times 22-2-2/2&\\
&=3+3\times 3^3+3/3&\\
&=4+(4-4/4)^4&\\
&=5\times 5+55+5&\\
&=66+6+6+6+6/6&\\
&=77+7+7/7&\\
&=88+8-88/8&\\
&=9\times9+(9+9+9+9)/9.&
\end{flalign*}

\begin{flalign*}
86&=11\times (1+1)^{(1+1+1)}-1-1&\\
&=2\times 2\times 22-2&\\
&=3+3+3\times 3^3-3/3&\\
&=44+44-(4+4)/4&\\
&=555/5-5\times 5&\\
&=66+6+6+6+(6+6)/6&\\
&=77+7+(7+7)/7&\\
&=88-(8+8)/8&\\
&=99-(99+9+9)/9.&
\end{flalign*}

\begin{flalign*}
87&=111-(1+1)\times (11+1)&\\
&=2\times 2\times 22-2/2&\\
&=3+3+3\times 3^3&\\
&=44+44-4/4&\\
&=55+((5+5)/5)^5&\\
&=6+6-66+6\times 66/6&\\
&=77+(77-7)/7&\\
&=88-8/8&\\
&=99-(99+9)/9.&
\end{flalign*}

\begin{flalign*}
88&=11\times (1+1)^{(1+1+1)}&\\
&=2\times 2\times 22&\\
&=3\times 33-33/3&\\
&=44+44&\\
&=5\times 5+(5^5/5+5)/(5+5)&\\
&=66+(66+66)/6&\\
&=77+77/7&\\
&=88&\\
&=99-99/9.&
\end{flalign*}

\begin{flalign*}
89&=111-11-11&\\
&=2\times 2\times 22+2/2&\\
&=3\times (3^3+3)-3/3&\\
&=44+44+4/4&\\
&=5\times (5\times 5-5)-55/5&\\
&=66+6+6+66/6&\\
&=77+(77+7)/7&\\
&=88+8/8&\\
&=99-9-9/9.&
\end{flalign*}

\begin{flalign*}
90&=(11-1)\times (11-1-1)&\\
&=2\times 2\times 22+2&\\
&=3\times (3^3+3)&\\
&=44+44+(4+4)/4&\\
&=5\times 5+55+5+5&\\
&=66+6+6+6+6&\\
&=7+77+7-7/7&\\
&=88+(8+8)/8&\\
&=99-9.&
\end{flalign*}

\begin{flalign*}
91&=(11-1)\times (11-1-1)+1&\\
&=2\times 2\times 22+2+2/2&\\
&=3^3+(3+3/3)^3&\\
&=4+44+44-4/4&\\
&=5-5\times 5+555/5&\\
&=66+6\times 6-66/6&\\
&=77+7+7&\\
&=88-8+88/8&\\
&=99-9+9/9.&
\end{flalign*}

\begin{flalign*}
92&=(11-1-1)^{(1+1)}+11&\\
&=2\times (2\times 22+2)&\\
&=3\times 3^3+33/3&\\
&=4+44+44&\\
&=55+5+((5+5)/5)^5&\\
&=66+6\times 6-(66-6)/6&\\
&=77+7+7+7/7&\\
&=88+8\times 8/(8+8)&\\
&=9\times9+99/9.&
\end{flalign*}

\begin{flalign*}
93&=((1+1)^{(11-1)}-1)/11&\\
&=2\times (2\times 22+2)+2/2&\\
&=3+3\times (3^3+3)&\\
&=((4+4)^4-4)/44&\\
&=5\times 5\times 5-((5+5)/5)^5&\\
&=666/6-6-6-6&\\
&=(777-77)/7-7&\\
&=8888/88-8&\\
&=999/9-9-9.&
\end{flalign*}

\begin{flalign*}
94&=((1+1)^{(11-1)}-1)/11+1&\\
&=2\times (2\times 22+2)+2&\\
&=3+3\times (3^3+3)+3/3&\\
&=(444-4)/4-4\times 4&\\
&=5\times (5\times 5-5)-5-5/5&\\
&=6\times 6-6+((6+6)/6)^6&\\
&=7777/77-7&\\
&=88+8-(8+8)/8&\\
&=(999+9)/9-9-9.&
\end{flalign*}

\begin{flalign*}
95&=111-(1+1)^{(1+1+1+1)}&\\
&=2\times 2\times (22+2)-2/2&\\
&=3\times 33-3-3/3&\\
&=444/4-4\times 4&\\
&=5\times (5\times 5-5)-5&\\
&=66+6\times 6-6-6/6&\\
&=77+7+77/7&\\
&=88+8-8/8&\\
&=99-(9\times9-9)/(9+9).&
\end{flalign*}

\begin{flalign*}
96&=(11+1)\times (1+1)^{(1+1+1)}&\\
&=2\times 2\times (22+2)&\\
&=3\times 33-3&\\
&=4\times (4\times 4+4+4)&\\
&=5\times (5\times 5-5)-5+5/5&\\
&=66+6\times 6-6&\\
&=7\times (7+7)-(7+7)/7&\\
&=88+8&\\
&=99-(9+9+9)/9.&
\end{flalign*}

\begin{flalign*}
97&=111-11-1-1-1&\\
&=2\times 2\times (22+2)+2/2&\\
&=3\times 33-3+3/3&\\
&=4\times 4+(4-4/4)^4&\\
&=55+5+5+((5+5)/5)^5&\\
&=66+6\times 6-6+6/6&\\
&=7\times (7+7)-7/7&\\
&=88+8+8/8&\\
&=99-(9+9)/9.&
\end{flalign*}

\begin{flalign*}
98&=111-11-1-1&\\
&=2\times 2\times (22+2)+2&\\
&=3\times 33-3/3&\\
&=4-4\times 4+(444-4)/4&\\
&=5\times (5\times 5-5)-(5+5)/5&\\
&=(666-6)/6-6-6&\\
&=7\times (7+7)&\\
&=88+(88-8)/8&\\
&=99-9/9.&
\end{flalign*}

\begin{flalign*}
99&=11\times (11-1-1)&\\
&=(22/2)^2-22&\\
&=3\times 33&\\
&=4-4\times 4+444/4&\\
&=5\times (5\times 5-5)-5/5&\\
&=666/6-6-6&\\
&=7\times (7+7)+7/7&\\
&=88+88/8&\\
&=99.&
\end{flalign*}

\begin{flalign*}
100&=(11-1)^{(1+1)}&\\
&=(2\times (2+2)+2)^2&\\
&=3\times 33+3/3&\\
&=(444-44)/4&\\
&=5\times (5\times 5-5)&\\
&=(666-66)/6&\\
&=(777-77)/7&\\
&=(888-88)/8&\\
&=99+9/9.&
\end{flalign*}
\end{multicols}
}

\bigskip
Instead of writing jointly, as above, here below in the following sections, the natural numbers  from 101 to 1000 are represented separately in each case.

\bigskip
\section{\textbf{Representations Using Number 1}}

{\footnotesize
\begin{multicols}{2}
\begin{itemize}
\item [] $101=(11-1)^{(1+1)}+1.$
\item [] $102=(11-1)^{(1+1)}+1+1.$
\item [] $103=(11-1)^{(1+1)}+1+1+1.$
\item [] $104=(11-1)^{(1+1)}+1+1+1+1.$
\item [] $105=111-(1+1)\times (1+1+1).$
\item [] $106=111-(11-1)/(1+1).$
\item [] $107=111-(1+1+1+1).$
\item [] $108=111-(1+1+1).$
\item [] $109=111-(1+1).$
\item [] $110=111-1.$
\item [] $111=111.$
\item [] $112=111+1.$
\item [] $113=111+1+1.$
\item [] $114=111+1+1+1.$
\item [] $115=111+1+1+1+1.$
\item [] $116=111+1+1+1+1+1.$
\item [] $117=111+(1+1)\times (1+1+1).$
\item [] $118=11^{(1+1)}-1-1-1.$
\item [] $119=11^{(1+1)}-1-1.$
\item [] $120=11^{(1+1)}-1.$
\item [] $121=11^{(1+1)}.$
\item [] $122=11^{(1+1)}+1.$
\item [] $123=11^{(1+1)}+1+1.$
\item [] $124=11^{(1+1)}+1+1+1.$
\item [] $125=11^{(1+1)}+1+1+1+1.$
\item [] $126=11^{(1+1)}+1+1+1+1+1.$
\item [] $127=(1+1)^{(1+1+1+1)}+111.$
\item [] $128=(1+1)^{((1+1)\times (1+1+1)+1)}.$
\item [] $129=11\times (11+1)-1-1-1.$
\item [] $130=(11-1)\times (11+1+1).$
\item [] $131=11\times (11+1)-1.$
\item [] $132=11\times (11+1).$
\item [] $133=11\times (11+1)+1.$
\item [] $134=11\times (11+1)+1+1.$
\item [] $135=11\times (11+1)+1+1+1.$
\item [] $136=11\times (11+1)+1+1+1+1.$
\item [] $137=(1+1)\times (11+1+1)+111.$
\item [] $138=(1+1+1)^{(1+1+1)}+111.$
\item [] $139=(11-1)\times (11+1+1+1)-1.$
\item [] $140=(11-1)\times (11+1+1+1).$
\item [] $141=(11+1)^{(1+1)}-1-1-1.$
\item [] $142=(11+1)^{(1+1)}-1-1.$
\item [] $143=11\times (11+1+1).$
\item [] $144=(11+1)^{(1+1)}.$
\item [] $145=(11+1)^{(1+1)}+1.$
\item [] $146=(11+1)^{(1+1)}+1+1.$
\item [] $147=(11+1)^{(1+1)}+1+1+1.$
\item [] $148=(11+1)^{(1+1)}+1+1+1+1.$
\item [] $149=(11+1)^{(1+1)}+1+1+1+1+1.$
\item [] $150=(11-1)\times (11+1+1+1+1).$
\item [] $151=(11-1)\times (11+1+1+1+1)+1.$
\item [] $152=11\times (11+1+1+1)-1-1.$
\item [] $153=11\times (11+1+1+1)-1.$
\item [] $154=11\times (11+1+1+1).$
\item [] $155=(11+1)^{(1+1)}+11.$
\item [] $156=(11+1)\times (11+1+1).$
\item [] $157=(11+1)\times (11+1+1)+1.$
\item [] $158=(11+1+1)^{(1+1)}-11.$
\item [] $159=(11+1+1)^{(1+1)}-11+1.$
\item [] $160=(1+1)\times ((11-1-1)^{(1+1)}-1).$
\item [] $161=(1+1)\times (11-1-1)^{(1+1)}-1.$
\item [] $162=(1+1)\times (11-1-1)^{(1+1)}.$
\item [] $163=(1+1)\times (11-1-1)^{(1+1)}+1.$
\item [] $164=(1+1)\times ((11-1-1)^{(1+1)}+1).$
\item [] $165=11\times (11+1+1+1+1).$
\item [] $166=11\times (11+1+1+1+1)+1.$
\item [] $167=(11+1+1)^{(1+1)}-1-1.$
\item [] $168=(11+1+1)^{(1+1)}-1.$
\item [] $169=(11+1+1)^{(1+1)}.$
\item [] $170=(11+1+1)^{(1+1)}+1.$
\item [] $171=(11+1+1)^{(1+1)}+1+1.$
\item [] $172=(11+1+1)^{(1+1)}+1+1+1.$
\item [] $173=(11+1+1)^{(1+1)}+1+1+1+1.$
\item [] $174=(1+1)\times (111-(1+1)\times (11+1)).$
\item [] $175=11\times (1+1)^{(1+1+1+1)}-1.$
\item [] $176=11\times (1+1)^{(1+1+1+1)}.$
\item [] $177=11\times (1+1)^{(1+1+1+1)}+1.$
\item [] $178=(1+1)\times (111-11-11).$
\item [] $179=(11+1+1)^{(1+1)}+11-1.$
\item [] $180=(11+1+1)^{(1+1)}+11.$
\item [] $181=(11+1+1)^{(1+1)}+11+1.$
\item [] $182=(11+1+1)\times (11+1+1+1).$
\item [] $183=(11+1+1)\times (11+1+1+1)+1.$
\item [] $184=(1+1)\times ((11-1-1)^{(1+1)}+11).$
\item [] $185=(11+1+1+1)^{(1+1)}-11.$
\item [] $186=((1+1)^{11}-1-1)/11.$
\item [] $187=((1+1)^{11}-1-1)/11+1.$
\item [] $188=((1+1)^{11}-1-1)/11+1+1.$
\item [] $189=(11-1-1)\times (11+11-1).$
\item [] $190=(11-1)\times ((1+1)\times (11-1)-1).$
\item [] $191=(11-1)\times ((1+1)\times (11-1)-1)+1.$
\item [] $192=(11-1-1)^{(1+1)}+111.$
\item [] $193=(11-1-1)^{(1+1)}+111+1.$
\item [] $194=(11+1+1+1)^{(1+1)}-1-1.$
\item [] $195=(11+1+1+1)^{(1+1)}-1.$
\item [] $196=(11+1+1+1)^{(1+1)}.$
\item [] $197=(11+1+1+1)^{(1+1)}+1.$
\item [] $198=(1+1)\times 11\times (11-1-1).$
\item [] $199=(1+1)\times (11-1)^{(1+1)}-1.$
\item [] $200=(1+1)\times (11-1)^{(1+1)}.$
\item [] $201=(1+1)\times (11-1)^{(1+1)}+1.$
\item [] $202=(1+1)\times ((11-1)^{(1+1)}+1).$
\item [] $203=(1+1)\times ((11-1)^{(1+1)}+1)+1.$
\item [] $204=(1+1)\times ((11-1)^{(1+1)}+1+1).$
\item [] $205=((1+1)^{11}+1+1)/(11-1).$
\item [] $206=((1+1)^{11}+1+1)/(11-1)+1.$
\item [] $207=(11+1+1+1)^{(1+1)}+11.$
\item [] $208=(1+1)\times (111-1)-11-1.$
\item [] $209=11\times ((1+1)\times (11-1)-1).$
\item [] $210=(11-1)\times (11+11-1).$
\item [] $211=(1+1)\times 111-11.$
\item [] $212=(1+1)\times 111-11+1.$
\item [] $213=(1+1)\times (111+1)-11.$
\item [] $214=(1+1)\times (111+1)-11+1.$
\item [] $215=(1+1)\times (111+1+1)-11.$
\item [] $216=(1+1)\times (111-1-1-1).$
\item [] $217=(1+1)\times (111-1-1)-1.$
\item [] $218=(1+1)\times (111-1-1).$
\item [] $219=(1+1)\times (111-1)-1.$
\item [] $220=(1+1)\times (111-1).$
\item [] $221=(1+1)\times 111-1.$
\item [] $222=(1+1)\times 111.$
\item [] $223=(1+1)\times 111+1.$
\item [] $224=(1+1)\times (111+1).$
\item [] $225=(1+1)\times (111+1)+1.$
\item [] $226=(1+1)\times (111+1+1).$
\item [] $227=(1+1)\times (111+1+1)+1.$
\item [] $228=(1+1)\times (111+1+1+1).$
\item [] $229=(1+1)\times (111+1+1+1)+1.$
\item [] $230=(11-1)\times (11+11+1).$
\item [] $231=11\times (11+11-1).$
\item [] $232=11^{(1+1)}+111.$
\item [] $233=(1+1)\times 111+11.$
\item [] $234=(1+1)\times 111+11+1.$
\item [] $235=(1+1)\times (111+1)+11.$
\item [] $236=(1+1)\times (111+1)+11+1.$
\item [] $237=(1+1)\times (111+1+1)+11.$
\item [] $238=(1+1)\times (11^{(1+1)}-1-1).$
\item [] $239=(1+1)\times (11^{(1+1)}-1)-1.$
\item [] $240=(1+1)\times (11^{(1+1)}-1).$
\item [] $241=(1+1)\times 11^{(1+1)}-1.$
\item [] $242=(1+1)\times 11^{(1+1)}.$
\item [] $243=(1+1)\times 11^{(1+1)}+1.$
\item [] $244=(1+1)\times (11^{(1+1)}+1).$
\item [] $245=(1+1)\times (11^{(1+1)}+1)+1.$
\item [] $246=(1+1)\times (11^{(1+1)}+1+1).$
\item [] $247=(1+1)\times (11^{(1+1)}+1+1)+1.$
\item [] $248=(1+1)\times (11^{(1+1)}+1+1+1).$
\item [] $249=(1+1)\times (11^{(1+1)}+1+1+1)+1.$
\item [] $250=(11-1)\times ((1+1)\times (11+1)+1).$
\item [] $251=(1+1)\times (11^{(1+1)}-1)+11.$
\item [] $252=(11+1)\times (11+11-1).$
\item [] $253=11\times (11+11+1).$
\item [] $254=11\times (11+11+1)+1.$
\item [] $255=(11+1)^{(1+1)}+111.$
\item [] $256=(1+1)^{((1+1)^{(1+1+1)})}.$
\item [] $257=(1+1)^{((1+1)^{(1+1+1)})}+1.$
\item [] $258=(1+1)^{((1+1)^{(1+1+1)})}+1+1.$
\item [] $259=(1+1)^{((1+1)^{(1+1+1)})}+1+1+1.$
\item [] $260=(1+1)\times (11-1)\times (11+1+1).$
\item [] $261=(1+1)\times (11\times (11+1)-1)-1.$
\item [] $262=(1+1)\times (11\times (11+1)-1).$
\item [] $263=(1+1)\times 11\times (11+1)-1.$
\item [] $264=(1+1)\times 11\times (11+1).$
\item [] $265=(1+1)\times 11\times (11+1)+1.$
\item [] $266=(1+1)\times (11\times (11+1)+1).$
\item [] $267=(1+1)\times (11\times (11+1)+1)+1.$
\item [] $268=(1+1)\times (11\times (11+1)+1+1).$
\item [] $269=(1+1)\times (11\times (11+1)+1+1)+1.$
\item [] $270=(11-1)\times (1+1+1)^{(1+1+1)}.$
\item [] $271=(11-1)\times (1+1+1)^{(1+1+1)}+1.$
\item [] $272=(11+1+1)\times (11+11-1)-1.$
\item [] $273=(11+1+1)\times (11+11-1).$
\item [] $274=11\times ((1+1)\times (11+1)+1)-1.$
\item [] $275=11\times ((1+1)\times (11+1)+1).$
\item [] $276=(11+1)\times (11+11+1).$
\item [] $277=(11+1)\times (11+11+1)+1.$
\item [] $278=(1111+1)/(1+1+1+1).$
\item [] $279=(1111+1)/(1+1+1+1)+1.$
\item [] $280=(11+1+1)^{(1+1)}+111.$
\item [] $281=(11+1+1)^{(1+1)}+111+1.$
\item [] $282=(1+1)\times ((11+1)^{(1+1)}-1-1-1).$
\item [] $283=(1+1)\times ((11+1)^{(1+1)}-1-1)-1.$
\item [] $284=(1+1)\times ((11+1)^{(1+1)}-1-1).$
\item [] $285=(1+1)\times 11\times (11+1+1)-1.$
\item [] $286=(1+1)\times 11\times (11+1+1).$
\item [] $287=(1+1)\times (11+1)^{(1+1)}-1.$
\item [] $288=(1+1)\times (11+1)^{(1+1)}.$
\item [] $289=(1+1)\times (11+1)^{(1+1)}+1.$
\item [] $290=(1+1)\times ((11+1)^{(1+1)}+1).$
\item [] $291=(1+1)\times ((11+1)^{(1+1)}+1)+1.$
\item [] $292=(1+1)\times ((11+1)^{(1+1)}+1)+1+1.$
\item [] $293=(1+1)\times ((11+1)^{(1+1)}+1)+1+1+1.$
\item [] $294=(1+1)\times ((11+1)^{(1+1)}+1)+1+1+1+1.$
\item [] $295=11\times (1+1+1)^{(1+1+1)}-1-1.$
\item [] $296=11\times (1+1+1)^{(1+1+1)}-1.$
\item [] $297=11\times (1+1+1)^{(1+1+1)}.$
\item [] $298=11\times (1+1+1)^{(1+1+1)}+1.$
\item [] $299=(1+1)\times (11+1)^{(1+1)}+11.$
\item [] $300=(1+1+1)\times (11-1)^{(1+1)}.$
\item [] $301=(1+1+1)\times (11-1)^{(1+1)}+1.$
\item [] $302=(1+1+1)\times (11-1)^{(1+1)}+1+1.$
\item [] $303=(1+1+1)\times ((11-1)^{(1+1)}+1).$
\item [] $304=(1+1+1)\times ((11-1)^{(1+1)}+1)+1.$
\item [] $305=(11\times 111-1)/(1+1+1+1).$
\item [] $306=(1+1)\times (11\times (11+1+1+1)-1).$
\item [] $307=111+(11+1+1+1)^{(1+1)}.$
\item [] $308=11\times (1+1)\times (11+1+1+1).$
\item [] $309=11\times (1+1)\times (11+1+1+1)+1.$
\item [] $310=(1+1)\times ((11+1)^{(1+1)}+11).$
\item [] $311=(1+1)\times ((11+1)^{(1+1)}+11)+1.$
\item [] $312=(1+1)\times (11+1)\times (11+1+1).$
\item [] $313=(1+1)\times (11+1)\times (11+1+1)+1.$
\item [] $314=(1+1)\times ((11+1)\times (11+1+1)+1).$
\item [] $315=((1+1)^{(11+1)}-1)/(11+1+1).$
\item [] $316=(1+1)\times ((11+1+1)^{(1+1)}-11).$
\item [] $317=(1+1)\times ((11+1+1)^{(1+1)}-11)+1.$
\item [] $318=(1+1+1)\times (111-1)-11-1.$
\item [] $319=11\times ((1+1+1)\times (11-1)-1).$
\item [] $320=(11-1)\times (11\times (1+1+1)-1).$
\item [] $321=111\times (1+1+1)-11-1.$
\item [] $322=111\times (1+1+1)-11.$
\item [] $323=111\times (1+1+1)-11+1.$
\item [] $324=((1+1)\times (11-1-1))^{(1+1)}.$
\item [] $325=((1+1)\times (11-1-1))^{(1+1)}+1.$
\item [] $326=(1+1+1)\times (111-1-1)-1.$
\item [] $327=(1+1+1)\times (111-1-1).$
\item [] $328=(1+1+1)\times (111-1-1)+1.$
\item [] $329=(1+1+1)\times (111-1)-1.$
\item [] $330=(1+1+1)\times (111-1).$
\item [] $331=(1+1+1)\times (111-1)+1.$
\item [] $332=(1+1+1)\times 111-1.$
\item [] $333=(1+1+1)\times 111.$
\item [] $334=(1+1+1)\times 111+1.$
\item [] $335=(1+1+1)\times 111+1+1.$
\item [] $336=(1+1+1)\times (111+1).$
\item [] $337=(1+1+1)\times (111+1)+1.$
\item [] $338=(1+1)\times (11+1+1)^{(1+1)}.$
\item [] $339=(1+1+1)\times (111+1+1).$
\item [] $340=(1+1+1)\times (111+1+1)+1.$
\item [] $341=(1+1+1)\times (111-1)+11.$
\item [] $342=(1+1+1)\times (111+1+1+1).$
\item [] $343=(1+1+1)\times 111+11-1.$
\item [] $344=(1+1+1)\times 111+11.$
\item [] $345=(1+1+1)\times 111+11+1.$
\item [] $346=(1+1+1)\times 111+11+1+1.$
\item [] $347=(1+1+1)\times (111+1)+11.$
\item [] $348=(1+1+1)\times (111+1)+11+1.$
\item [] $349=(1+1)\times (11+1+1)^{(1+1)}+11.$
\item [] $350=(1+1+1)\times (111+1+1)+11.$
\item [] $351=11\times (11\times (1+1+1)-1)-1.$
\item [] $352=11\times (11\times (1+1+1)-1).$
\item [] $353=11\times (11\times (1+1+1)-1)+1.$
\item [] $354=11\times (11\times (1+1+1)-1)+1+1.$
\item [] $355=(1+1+1)\times 111+11+11.$
\item [] $356=(1+1+1)\times (11^{(1+1)}-1-1)-1.$
\item [] $357=(1+1+1)\times (11^{(1+1)}-1-1).$
\item [] $358=(1+1+1)\times (11^{(1+1)}-1-1)+1.$
\item [] $359=(1+1+1)\times (11^{(1+1)}-1)-1.$
\item [] $360=(1+1+1)\times (11^{(1+1)}-1).$
\item [] $361=((1+1)\times (11-1)-1)^{(1+1)}.$
\item [] $362=(1+1+1)\times 11\times 11-1.$
\item [] $363=(1+1+1)\times 11\times 11.$
\item [] $364=(1+1+1)\times 11\times 11+1.$
\item [] $365=(1+1+1)\times 11\times 11+1+1.$
\item [] $366=(1+1+1)\times (11^{(1+1)}+1).$
\item [] $367=(1+1+1)\times (11^{(1+1)}+1)+1.$
\item [] $368=(1+1+1)\times (11^{(1+1)}+1)+1+1.$
\item [] $369=(1+1+1)\times (11^{(1+1)}+1+1).$
\item [] $370=(1111-1)/(1+1+1).$
\item [] $371=(1111-1)/(1+1+1)+1.$
\item [] $372=(1111-1)/(1+1+1)+1+1.$
\item [] $373=(11+11)^{(1+1)}-111.$
\item [] $374=11\times (11\times (1+1+1)+1).$
\item [] $375=11\times (11\times (1+1+1)+1)+1.$
\item [] $376=11\times (11\times (1+1+1)+1)+1+1.$
\item [] $377=(1+1+1)\times (11^{(1+1)}+1)+11.$
\item [] $378=(1+1+1)\times (11^{(1+1)}+1)+11+1.$
\item [] $379=(11-1)\times (111/(1+1+1)+1)-1.$
\item [] $380=(11-1)\times (111/(1+1+1)+1).$
\item [] $381=(1111-1)/(1+1+1)+11.$
\item [] $382=(1111-1)/(1+1+1)+11+1.$
\item [] $383=(11+1)\times (11\times (1+1+1)-1)-1.$
\item [] $384=(11+1)\times (11\times (1+1+1)-1).$
\item [] $385=11\times (11\times (1+1+1)+1+1).$
\item [] $386=11\times (11\times (1+1+1)+1+1)+1.$
\item [] $387=11\times (11\times (1+1+1)+1+1)+1+1.$
\item [] $388=((1+1)\times (11-1))^{(1+1)}-11-1.$
\item [] $389=((1+1)\times (11-1))^{(1+1)}-11.$
\item [] $390=((1+1)\times (11-1))^{(1+1)}-11+1.$
\item [] $391=(1+1)^{(11-1-1)}-11^{(1+1)}.$
\item [] $392=(1+1)\times (11+1+1+1)^{(1+1)}.$
\item [] $393=(1+1+1)\times (11\times (11+1)-1).$
\item [] $394=(1+1+1)\times (11\times (11+1)-1)+1.$
\item [] $395=(1+1+1)\times (11+1)\times 11-1.$
\item [] $396=(1+1+1)\times (11+1)\times 11.$
\item [] $397=(1+1+1)\times (11+1)\times 11+1.$
\item [] $398=((1+1)\times (11-1))^{(1+1)}-1-1.$
\item [] $399=((1+1)\times (11-1))^{(1+1)}-1.$
\item [] $400=((1+1)\times (11-1))^{(1+1)}.$
\item [] $401=((1+1)\times (11-1))^{(1+1)}+1.$
\item [] $402=((1+1)\times (11-1))^{(1+1)}+1+1.$
\item [] $403=((1+1)\times (11-1))^{(1+1)}+1+1+1.$
\item [] $404=(1+1)\times (1+1)\times ((11-1)^{(1+1)}+1).$
\item [] $405=11\times 111/(1+1+1)-1-1.$
\item [] $406=11\times 111/(1+1+1)-1.$
\item [] $407=11\times 111/(1+1+1).$
\item [] $408=11\times 111/(1+1+1)+1.$
\item [] $409=11\times 111/(1+1+1)+1+1.$
\item [] $410=((1+1)\times (11-1))^{(1+1)}+11-1.$
\item [] $411=((1+1)\times (11-1))^{(1+1)}+11.$
\item [] $412=((1+1)\times (11-1))^{(1+1)}+11+1.$
\item [] $413=((1+1)\times (11-1))^{(1+1)}+11+1+1.$
\item [] $414=(1+1)\times ((11+1+1+1)^{(1+1)}+11).$
\item [] $415=(11+1+1)\times (11\times (1+1+1)-1)-1.$
\item [] $416=(11+1+1)\times (11\times (1+1+1)-1).$
\item [] $417=11\times (111/(1+1+1)+1)-1.$
\item [] $418=11\times (111/(1+1+1)+1).$
\item [] $419=11\times (111/(1+1+1)+1)+1.$
\item [] $420=(1+1)\times (11-1)\times (11+11-1).$
\item [] $421=(1+1+1)\times (11+1)^{(1+1)}-11.$
\item [] $422=(1+1)\times ((1+1)\times 111-11).$
\item [] $423=(1+1)\times ((1+1)\times 111-11)+1.$
\item [] $424=(1+1)\times ((1+1)\times 111-11+1).$
\item [] $425=(1+1)\times ((1+1)\times 111-11+1)+1.$
\item [] $426=(1+1)\times ((1+1)\times (111+1)-11).$
\item [] $427=(1+1)\times ((1+1)\times (111+1)-11)+1.$
\item [] $428=(1+1+1)\times (11+1+1)\times 11-1.$
\item [] $429=(1+1+1)\times (11+1+1)\times 11.$
\item [] $430=(11+11-1)^{(1+1)}-11.$
\item [] $431=(1+1+1)\times (11+1)^{(1+1)}-1.$
\item [] $432=(1+1+1)\times (11+1)^{(1+1)}.$
\item [] $433=(1+1+1)\times (11+1)^{(1+1)}+1.$
\item [] $434=(1+1+1)\times (11+1)^{(1+1)}+1+1.$
\item [] $435=(1+1+1)\times ((11+1)^{(1+1)}+1).$
\item [] $436=(1+1)\times (1+1)\times (111-1-1).$
\item [] $437=(1+1)\times (1+1)\times (111-1-1)+1.$
\item [] $438=(1+1)\times ((1+1)\times (111-1)-1).$
\item [] $439=(11+11-1)^{(1+1)}-1-1.$
\item [] $440=(1+1)\times (1+1)\times (111-1).$
\item [] $441=(11+11-1)^{(1+1)}.$
\item [] $442=(11+11-1)^{(1+1)}+1.$
\item [] $443=(1+1)\times (1+1)\times 111-1.$
\item [] $444=(1+1)\times (1+1)\times 111.$
\item [] $445=(1+1)\times (1+1)\times 111+1.$
\item [] $446=(1+1)\times ((1+1)\times 111+1).$
\item [] $447=(1+1)\times ((1+1)\times 111+1)+1.$
\item [] $448=(1+1)\times (1+1)\times (111+1).$
\item [] $449=(1+1)\times (1+1)\times (111+1)+1.$
\item [] $450=(1+1)\times ((1+1)\times (111+1)+1).$
\item [] $451=(1+1)\times ((1+1)\times (111+1)+1)+1.$
\item [] $452=(11+11-1)^{(1+1)}+11.$
\item [] $453=(11+11-1)^{(1+1)}+11+1.$
\item [] $454=(1+1)\times (1+1)\times 111+11-1.$
\item [] $455=(1+1)\times (1+1)\times 111+11.$
\item [] $456=(1+1)\times (1+1)\times 111+11+1.$
\item [] $457=(1+1)\times ((1+1)\times 111+1)+11.$
\item [] $458=(1+1)\times ((1+1)\times 111+1)+11+1.$
\item [] $459=(1+1)\times (1+1)\times (111+1)+11.$
\item [] $460=(1+1)\times (11-1)\times (11+11+1).$
\item [] $461=(1+1)\times 11\times (11+11-1)-1.$
\item [] $462=(1+1)\times 11\times (11+11-1).$
\item [] $463=(1+1)\times 11\times (11+11-1)+1.$
\item [] $464=(1+1)\times (11^{(1+1)}+111).$
\item [] $465=(1+1)\times (11^{(1+1)}+111)+1.$
\item [] $466=(1+1)\times ((1+1)\times 111+11).$
\item [] $467=(1+1)\times ((1+1)\times 111+11)+1.$
\item [] $468=(1+1)\times ((1+1)\times 111+11+1).$
\item [] $469=(1+1)\times ((1+1)\times 111+11+1)+1.$
\item [] $470=(1+1)\times ((1+1)\times (111+1)+11).$
\item [] $471=(11+11)^{(1+1)}-11-1-1.$
\item [] $472=(11+11)^{(1+1)}-11-1.$
\item [] $473=(11+11)^{(1+1)}-11.$
\item [] $474=(11+11)^{(1+1)}-11+1.$
\item [] $475=(11+11)^{(1+1)}-11+1+1.$
\item [] $476=(1+1)\times ((1+1)\times (11^{(1+1)}-1-1)).$
\item [] $477=(1+1)\times ((1+1)\times (11^{(1+1)}-1-1))+1.$
\item [] $478=(1+1)\times (((1+1)\times (11^{(1+1)}-1-1))+1).$
\item [] $479=(1+1)\times ((1+1)\times (11^{(1+1)}-1))-1.$
\item [] $480=(1+1)\times (1+1)\times (11^{(1+1)}-1).$
\item [] $481=(11+11)^{(1+1)}-1-1-1.$
\item [] $482=(11+11)^{(1+1)}-1-1.$
\item [] $483=(11+11)^{(1+1)}-1.$
\item [] $484=(11+11)^{(1+1)}.$
\item [] $485=(11+11)^{(1+1)}+1.$
\item [] $486=(11+11)^{(1+1)}+1+1.$
\item [] $487=(11+11)^{(1+1)}+1+1+1.$
\item [] $488=(1+1)\times (1+1)\times (11^{(1+1)}+1).$
\item [] $489=(1+1)\times (1+1)\times (11^{(1+1)}+1)+1.$
\item [] $490=(1+1)\times ((1+1)\times (11^{(1+1)}+1)+1).$
\item [] $491=(1+1)\times ((1+1)\times (11^{(1+1)}+1)+1)+1.$
\item [] $492=(1+1)\times (1+1)\times (11^{(1+1)}+1+1).$
\item [] $493=(11+11)^{(1+1)}+11-1-1.$
\item [] $494=(11+11)^{(1+1)}+11-1.$
\item [] $495=(11+11)^{(1+1)}+11.$
\item [] $496=(11+11)^{(1+1)}+11+1.$
\item [] $497=(11+11)^{(1+1)}+11+1+1.$
\item [] $498=(11-1)^{(1+1+1)}/(1+1)-1-1.$
\item [] $499=(11-1)^{(1+1+1)}/(1+1)-1.$
\item [] $500=(11-1)^{(1+1+1)}/(1+1).$
\item [] $501=(1+1)^{(11-1-1)}-11.$
\item [] $502=(1+1)^{(11-1-1)}-11+1.$
\item [] $503=(1+1)^{(11-1-1)}-11+1+1.$
\item [] $504=(1+1)\times (11+1)\times (11+11-1).$
\item [] $505=(11111-1)/(11+11).$
\item [] $506=(1+1)\times 11\times (11+11+1).$
\item [] $507=(1+1+1)\times (11+1+1)^{(1+1)}.$
\item [] $508=(1+1+1)\times (11+1+1)^{(1+1)}+1.$
\item [] $509=(1+1)^{(11-1-1)}-1-1-1.$
\item [] $510=(1+1)^{(11-1-1)}-1-1.$
\item [] $511=(1+1)^{(11-1-1)}-1.$
\item [] $512=(1+1)^{(11-1-1)}.$
\item [] $513=(1+1)^{(11-1-1)}+1.$
\item [] $514=(1+1)^{(11-1-1)}+1+1.$
\item [] $515=(1+1)^{(11-1-1)}+1+1+1.$
\item [] $516=(1+1)^{(11-1-1)}+1+1+1+1.$
\item [] $517=11\times ((1+1)\times (11+11+1)+1).$
\item [] $518=(11+11+1)^{(1+1)}-11.$
\item [] $519=(11+11+1)^{(1+1)}-11+1.$
\item [] $520=(11+11+1)^{(1+1)}-11+1+1.$
\item [] $521=(1+1)^{(11-1-1)}+11-1-1.$
\item [] $522=(1+1)^{(11-1-1)}+11-1.$
\item [] $523=(1+1)^{(11-1-1)}+11.$
\item [] $524=(1+1)^{(11-1-1)}+11+1.$
\item [] $525=(1+1)^{(11-1-1)}+11+1+1.$
\item [] $526=(1+1)\times ((1+1)\times 11\times (11+1)-1).$
\item [] $527=(11+11+1)^{(1+1)}-1-1.$
\item [] $528=(11+11+1)^{(1+1)}-1.$
\item [] $529=(11+11+1)^{(1+1)}.$
\item [] $530=(11+11+1)^{(1+1)}+1.$
\item [] $531=(11+11+1)^{(1+1)}+1+1.$
\item [] $532=(11+11+1)^{(1+1)}+1+1+1.$
\item [] $533=(11+11+1)^{(1+1)}+1+1+1+1.$
\item [] $534=(1+1)^{(11-1-1)}+11+11.$
\item [] $535=(1+1)^{(11-1-1)}+11+11+1.$
\item [] $536=(1+1)\times (1+1)\times (11\times (11+1)+1+1).$
\item [] $537=(1+1)^{(11-1-1)}+(1+1)\times (11+1)+1.$
\item [] $538=(11+11+1)^{(1+1)}+11-1-1.$
\item [] $539=(11+11+1)^{(1+1)}+11-1.$
\item [] $540=(11+11+1)^{(1+1)}+11.$
\item [] $541=(11+11+1)^{(1+1)}+11+1.$
\item [] $542=(11+11+1)^{(1+1)}+11+1+1.$
\item [] $543=(1111-1)/(1+1)-11-1.$
\item [] $544=(1111-1)/(1+1)-11.$
\item [] $545=(1111+1)/(1+1)-11.$
\item [] $546=(1111+1)/(1+1)-11+1.$
\item [] $547=(1111+1)/(1+1)-11+1+1.$
\item [] $548=(1111-11)/(1+1)-1-1.$
\item [] $549=(1111-11)/(1+1)-1.$
\item [] $550=(1111-11)/(1+1).$
\item [] $551=(1111-11)/(1+1)+1.$
\item [] $552=(1111-11)/(1+1)+1+1.$
\item [] $553=(1111-1)/(1+1)-1-1.$
\item [] $554=(1111-1)/(1+1)-1.$
\item [] $555=(1111-1)/(1+1).$
\item [] $556=(1111+1)/(1+1).$
\item [] $557=(1111+1)/(1+1)+1.$
\item [] $558=(1111+1)/(1+1)+1+1.$
\item [] $559=(1111+1)/(1+1)+1+1+1.$
\item [] $560=(11-1)\times (111+1)/(1+1).$
\item [] $561=(1111+11)/(1+1).$
\item [] $562=(1111+11)/(1+1)+1.$
\item [] $563=(1111+11)/(1+1)+1+1.$
\item [] $564=((1+1)\times (11+1))^{(1+1)}-11-1.$
\item [] $565=((1+1)\times (11+1))^{(1+1)}-11.$
\item [] $566=(1111-1)/(1+1)+11.$
\item [] $567=(1111+1)/(1+1)+11.$
\item [] $568=(1111+1)/(1+1)+11+1.$
\item [] $569=(1111+1)/(1+1)+11+1+1.$
\item [] $570=(11-1)\times ((111+1)/(1+1)+1).$
\item [] $571=(11-1)\times ((111+1)/(1+1)+1)+1.$
\item [] $572=(1111+11)/(1+1)+11.$
\item [] $573=((1+1)\times (11+1))^{(1+1)}-1-1-1.$
\item [] $574=((1+1)\times (11+1))^{(1+1)}-1-1.$
\item [] $575=((1+1)\times (11+1))^{(1+1)}-1.$
\item [] $576=((1+1)\times (11+1))^{(1+1)}.$
\item [] $577=((1+1)\times (11+1))^{(1+1)}+1.$
\item [] $578=((1+1)\times (11+1))^{(1+1)}+1+1.$
\item [] $579=((1+1)\times (11+1))^{(1+1)}+1+1+1.$
\item [] $580=(1+1)\times (1+1)\times ((11+1)^{(1+1)}+1).$
\item [] $581=(1+1)\times (1+1)\times ((11+1)^{(1+1)}+1)+1.$
\item [] $582=11\times ((111-1)/(1+1)-1-1)-1.$
\item [] $583=11\times ((111-1)/(1+1)-1-1).$
\item [] $584=11\times ((111-1)/(1+1)-1-1)+1.$
\item [] $585=((1+1)\times (11+1))^{(1+1)}+11-1-1.$
\item [] $586=((1+1)\times (11+1))^{(1+1)}+11-1.$
\item [] $587=((1+1)\times (11+1))^{(1+1)}+11.$
\item [] $588=((1+1)\times (11+1))^{(1+1)}+11+1.$
\item [] $589=((1+1)\times (11+1))^{(1+1)}+11+1+1.$
\item [] $590=(11-1)\times ((11^{(1+1)}-1)/(1+1)-1).$
\item [] $591=(1+1+1)\times ((11+1+1+1)^{(1+1)}+1).$
\item [] $592=(1+1)\times (11\times (1+1+1)^{(1+1+1)}-1).$
\item [] $593=11\times ((111-1)/(1+1)-1)-1.$
\item [] $594=11\times ((111-1)/(1+1)-1).$
\item [] $595=(11+11)^{(1+1)}+111.$
\item [] $596=(11+11)^{(1+1)}+111+1.$
\item [] $597=(11+11)^{(1+1)}+111+1+1.$
\item [] $598=((1+1)\times (11+1))^{(1+1)}+11+11.$
\item [] $599=1111-(1+1)^{(11-1-1)}.$
\item [] $600=(1+1)\times (1+1+1)\times (11-1)^{(1+1)}.$
\item [] $601=(1+1)\times (1+1+1)\times (11-1)^{(1+1)}+1.$
\item [] $602=(1+1)\times ((1+1+1)\times (11-1)^{(1+1)}+1).$
\item [] $603=11\times (111-1)/(1+1)-1-1.$
\item [] $604=11\times (111-1)/(1+1)-1.$
\item [] $605=11\times (111-1)/(1+1).$
\item [] $606=11\times (111-1)/(1+1)+1.$
\item [] $607=11\times (111-1)/(1+1)+1+1.$
\item [] $608=(11\times 111-1)/(1+1)-1-1.$
\item [] $609=(11\times 111-1)/(1+1)-1.$
\item [] $610=(11\times 111-1)/(1+1).$
\item [] $611=(11\times 111+1)/(1+1).$
\item [] $612=(11\times 111+1)/(1+1)+1.$
\item [] $613=(11\times 111+1)/(1+1)+1+1.$
\item [] $614=11\times (111+1)/(1+1)-1-1.$
\item [] $615=11\times (111+1)/(1+1)-1.$
\item [] $616=11\times (111+1)/(1+1).$
\item [] $617=11\times (111+1)/(1+1)+1.$
\item [] $618=11\times (111+1)/(1+1)+1+1.$
\item [] $619=11\times (111+1)/(1+1)+1+1+1.$
\item [] $620=(11\times 111-1)/(1+1)+11-1.$
\item [] $621=(11\times 111-1)/(1+1)+11.$
\item [] $622=(11\times 111+1)/(1+1)+11.$
\item [] $623=(1+1)^{(11-1-1)}+111.$
\item [] $624=((1+1)\times (11+1)+1)^{(1+1)}-1.$
\item [] $625=((1+1)\times (11+1)+1)^{(1+1)}.$
\item [] $626=((1+1)\times (11+1)+1)^{(1+1)}+1.$
\item [] $627=11\times ((111+1)/(1+1)+1).$
\item [] $628=11\times ((111+1)/(1+1)+1)+1.$
\item [] $629=11\times ((111+1)/(1+1)+1)+1+1.$
\item [] $630=(1+1+1)\times (11-1)\times (11+11-1).$
\item [] $631=(1+1+1)\times (11-1)\times (11+11-1)+1.$
\item [] $632=11^{(1+1)}+(1+1)^{(11-1-1)}-1.$
\item [] $633=(1+1+1)\times ((1+1)\times 111-11).$
\item [] $634=11^{(1+1)}+(1+1)^{(11-1-1)}+1.$
\item [] $635=((1+1)\times (11+1)+1)^{(1+1)}+11-1.$
\item [] $636=((1+1)\times (11+1)+1)^{(1+1)}+11.$
\item [] $637=((1+1)\times (11+1)+1)^{(1+1)}+11+1.$
\item [] $638=11\times ((111+1)/(1+1)+1+1).$
\item [] $639=11\times ((111+1)/(1+1)+1+1)+1.$
\item [] $640=(11+11+1)^{(1+1)}+111.$
\item [] $641=(11+11+1)^{(1+1)}+111+1.$
\item [] $642=(1+1)\times ((1+1+1)\times 111-11-1).$
\item [] $643=(1+1)\times ((1+1+1)\times 111-11)-1.$
\item [] $644=(1+1)\times ((1+1+1)\times 111-11).$
\item [] $645=(1+1)\times ((1+1+1)\times 111-11)+1.$
\item [] $646=(1+1)\times ((1+1+1)\times 111-11+1).$
\item [] $647=(1+1)\times ((1+1)\times (11-1-1))^{(1+1)}-1.$
\item [] $648=(1+1)\times ((1+1)\times (11-1-1))^{(1+1)}.$
\item [] $649=11\times ((11^{(1+1)}-1)/(1+1)-1).$
\item [] $650=11\times ((11^{(1+1)}-1)/(1+1)-1)+1.$
\item [] $651=(1+1+1)\times ((1+1)\times (111-1-1)-1).$
\item [] $652=(1+1)\times ((1+1+1)\times (111-1-1)-1).$
\item [] $653=(11^{(1+1+1)}-1)/(1+1)-11-1.$
\item [] $654=(1+1)\times (1+1+1)\times (111-1-1).$
\item [] $655=(1+1)\times (1+1+1)\times 111-11.$
\item [] $656=(1+1)\times (1+1+1)\times 111-11+1.$
\item [] $657=(1+1+1)\times ((1+1)\times (111-1)-1).$
\item [] $658=(1+1)\times ((1+1+1)\times (111-1)-1).$
\item [] $659=(1+1)\times (1+1+1)\times (111-1)-1.$
\item [] $660=(1+1)\times (1+1+1)\times (111-1).$
\item [] $661=(1+1)\times (1+1+1)\times (111-1)+1.$
\item [] $662=(1+1)\times ((1+1+1)\times (111-1)+1).$
\item [] $663=(1+1+1)\times ((1+1)\times 111-1).$
\item [] $664=(1+1)\times ((1+1+1)\times 111-1).$
\item [] $665=(11^{(1+1+1)}-1)/(1+1).$
\item [] $666=(1+1)\times (1+1+1)\times 111.$
\item [] $667=(1+1)\times (1+1+1)\times 111+1.$
\item [] $668=(1+1)\times ((1+1+1)\times 111+1).$
\item [] $669=(1+1+1)\times ((1+1)\times 111+1).$
\item [] $670=(1+1+1)\times ((1+1)\times 111+1)+1.$
\item [] $671=11\times (11^{(1+1)}+1)/(1+1).$
\item [] $672=(1+1)\times (1+1+1)\times (111+1).$
\item [] $673=(1+1)\times (1+1+1)\times (111+1)+1.$
\item [] $674=(1+1)\times ((1+1+1)\times (111+1)+1).$
\item [] $675=((1+1)\times (11+1+1))^{(1+1)}-1.$
\item [] $676=((1+1)\times (11+1+1))^{(1+1)}.$
\item [] $677=((1+1)\times (11+1+1))^{(1+1)}+1.$
\item [] $678=((1+1)\times (11+1+1))^{(1+1)}+1+1.$
\item [] $679=((1+1)^{11}-11)/(1+1+1).$
\item [] $680=((1+1)^{11}-11)/(1+1+1)+1.$
\item [] $681=((1+1)^{11}+1)/(1+1+1)-1-1.$
\item [] $682=((1+1)^{11}-1-1)/(1+1+1).$
\item [] $683=((1+1)^{11}+1)/(1+1+1).$
\item [] $684=((1+1)^{11}+1)/(1+1+1)+1.$
\item [] $685=((1+1)^{11}+1)/(1+1+1)+1+1.$
\item [] $686=((1+1)^{11}+11-1)/(1+1+1).$
\item [] $687=((1+1)\times (11+1+1))^{(1+1)}+11.$
\item [] $688=(1+1)\times ((1+1+1)\times 111+11).$
\item [] $689=(1+1)\times ((1+1+1)\times 111+11)+1.$
\item [] $690=((1+1)^{11}-11)/(1+1+1)+11.$
\item [] $691=(1+1)^{(11-1)}-(1+1+1)\times 111.$
\item [] $692=11\times (1+1+1)\times (11+11-1)-1.$
\item [] $693=11\times (1+1+1)\times (11+11-1).$
\item [] $694=((1+1)^{11}+1)/(1+1+1)+11.$
\item [] $695=((1+1)^{11}+1)/(1+1+1)+11+1.$
\item [] $696=(1+1+1)\times (11^{(1+1)}+111).$
\item [] $697=(1+1+1)\times (11^{(1+1)}+111)+1.$
\item [] $698=(1+1+1)\times ((1+1)\times 111+11)-1.$
\item [] $699=(1+1+1)\times ((1+1)\times 111+11).$
\item [] $700=(1+1+1)\times ((1+1)\times 111+11)+1.$
\item [] $701=(1+1+1)\times ((1+1)\times 111+11)+1+1.$
\item [] $702=(1+1)\times (11\times (11\times (1+1+1)-1)-1).$
\item [] $703=11\times (1+1)^{((1+1)\times (1+1+1))}-1.$
\item [] $704=11\times (1+1)^{((1+1)\times (1+1+1))}.$
\item [] $705=11\times (1+1)^{((1+1)\times (1+1+1))}+1.$
\item [] $706=11\times (1+1)^{((1+1)\times (1+1+1))}+1+1.$
\item [] $707=(11-1-1)^{(1+1+1)}-11-11.$
\item [] $708=(11+1)\times ((11^{(1+1)}-1)/(1+1)-1).$
\item [] $709=11\times 111-(1+1)^{(11-1-1)}.$
\item [] $710=(11-1)\times ((11+1)^{(1+1)}/(1+1)-1).$
\item [] $711=1111-((1+1)\times (11-1))^{(1+1)}.$
\item [] $712=1111-((1+1)\times (11-1))^{(1+1)}+1.$
\item [] $713=(11+11+1)\times ((1+1+1)\times (11-1)+1).$
\item [] $714=(1+1)\times (1+1+1)\times (11^{(1+1)}-1-1).$
\item [] $715=11\times ((1+1)^{((1+1)\times (1+1+1))}+1).$
\item [] $716=(1+1)^{11}-11^{(1+1+1)}-1.$
\item [] $717=(1+1)^{11}-11^{(1+1+1)}.$
\item [] $718=(11-1-1)^{(1+1+1)}-11.$
\item [] $719=(11-1-1)^{(1+1+1)}-11+1.$
\item [] $720=(1+1)\times (1+1+1)\times (11^{(1+1)}-1).$
\item [] $721=(11^{(1+1+1)}+111)/(1+1).$
\item [] $722=(1+1)\times ((1+1)\times (11-1)-1)^{(1+1)}.$
\item [] $723=(1+1+1)\times ((1+1)\times 11^{(1+1)}-1).$
\item [] $724=(1+1)\times (11\times 11\times (1+1+1)-1).$
\item [] $725=(1+1)\times 11\times 11\times (1+1+1)-1.$
\item [] $726=(1+1)\times 11\times 11\times (1+1+1).$
\item [] $727=(11-1-1)^{(1+1+1)}-1-1.$
\item [] $728=(11-1-1)^{(1+1+1)}-1.$
\item [] $729=(11-1-1)^{(1+1+1)}.$
\item [] $730=(11-1-1)^{(1+1+1)}+1.$
\item [] $731=(11-1-1)^{(1+1+1)}+1+1.$
\item [] $732=(11-1-1)^{(1+1+1)}+1+1+1.$
\item [] $733=(11-1-1)^{(1+1+1)}+1+1+1+1.$
\item [] $734=(1+1)\times ((1+1+1)\times (11^{(1+1)}+1)+1).$
\item [] $735=(1+1+1)\times ((1+1)\times (11^{(1+1)}+1)+1).$
\item [] $736=((1+1)\times (11+1)+1)^{(1+1)}+111.$
\item [] $737=11\times (11\times (1+1)\times (1+1+1)+1).$
\item [] $738=11\times (11\times (1+1)\times (1+1+1)+1)+1.$
\item [] $739=(11-1-1)^{(1+1+1)}+11-1.$
\item [] $740=(11-1-1)^{(1+1+1)}+11.$
\item [] $741=(11-1-1)^{(1+1+1)}+11+1.$
\item [] $742=(11-1-1)^{(1+1+1)}+11+1+1.$
\item [] $743=(11-1-1)^{(1+1+1)}+11+1+1+1.$
\item [] $744=(11+1)\times ((11^{(1+1)}+1)/(1+1)+1).$
\item [] $745=((1+1)^{(11+1+1)}+1+1+1)/11.$
\item [] $746=(1+1)\times ((11+11)^{(1+1)}-111).$
\item [] $747=(1+1)\times 11\times (11\times (1+1+1)+1)-1.$
\item [] $748=(1+1)\times 11\times (11\times (1+1+1)+1).$
\item [] $749=(1+1)\times 11\times (11\times (1+1+1)+1)+1.$
\item [] $750=(1+1)\times (11\times (11\times (1+1+1)+1)+1).$
\item [] $751=(11-1-1)^{(1+1+1)}+11+11.$
\item [] $752=(11-1-1)^{(1+1+1)}+11+11+1.$
\item [] $753=(11+1)^{(1+1+1)}/(1+1)-111.$
\item [] $754=(11+1)^{(1+1+1)}/(1+1)-111+1.$
\item [] $755=11^{(1+1+1)}-((1+1)\times (11+1))^{(1+1)}.$
\item [] $756=(1+1+1)\times (11+1)\times (11+11-1).$
\item [] $757=(1+1+1)\times (11+1)\times (11+11-1)+1.$
\item [] $758=11\times (1+1+1)\times (11+11+1)-1.$
\item [] $759=11\times (1+1+1)\times (11+11+1).$
\item [] $760=11\times (1+1+1)\times (11+11+1)+1.$
\item [] $761=11\times (1+1+1)\times (11+11+1)+1+1.$
\item [] $762=(1+1+1)\times (11\times (11+11+1)+1).$
\item [] $763=(111-1-1)\times (11-1-1-1-1).$
\item [] $764=(111-1-1)\times (11-1-1-1-1)+1.$
\item [] $765=(1+1+1)\times ((11+1)^{(1+1)}+111).$
\item [] $766=111\times ((1+1)\times (1+1+1)+1)-11.$
\item [] $767=(1+1+1)\times (1+1)^{((1+1)^{(1+1+1)})}-1.$
\item [] $768=(1+1+1)\times (1+1)^{((1+1)^{(1+1+1)})}.$
\item [] $769=(1+1+1)\times (1+1)^{((1+1)^{(1+1+1)})}+1.$
\item [] $770=11\times ((11-1-1)^{(1+1)}-11).$
\item [] $771=11\times ((11-1-1)^{(1+1)}-11)+1.$
\item [] $772=11\times ((11-1-1)^{(1+1)}-11)+1+1.$
\item [] $773=((1+1+1)^{(1+1+1)}+1)^{(1+1)}-11.$
\item [] $774=((1+1+1)^{(1+1+1)}+1)^{(1+1)}-11+1.$
\item [] $775=1111-(1+1+1)\times (111+1).$
\item [] $776=111\times ((1+1)\times (1+1+1)+1)-1.$
\item [] $777=111\times ((1+1)\times (1+1+1)+1).$
\item [] $778=111\times ((1+1)\times (1+1+1)+1)+1.$
\item [] $779=111\times ((1+1)\times (1+1+1)+1)+1+1.$
\item [] $780=(11-1)\times (111-11\times (1+1+1)).$
\item [] $781=11\times ((11+1)^{(1+1)}/(1+1)-1).$
\item [] $782=11\times ((11+1)^{(1+1)}/(1+1)-1)+1.$
\item [] $783=((1+1+1)^{(1+1+1)}+1)^{(1+1)}-1.$
\item [] $784=((1+1+1)^{(1+1+1)}+1)^{(1+1)}.$
\item [] $785=((1+1+1)^{(1+1+1)}+1)^{(1+1)}+1.$
\item [] $786=((1+1+1)^{(1+1+1)}+1)^{(1+1)}+1+1.$
\item [] $787=((1+1)\times (11+1+1))^{(1+1)}+111.$
\item [] $788=111\times ((1+1)\times (1+1+1)+1)+11.$
\item [] $789=(1+1+1)\times ((1+1)\times 11\times (11+1)-1).$
\item [] $790=(1+1)\times (11\times (1+1+1)\times (11+1)-1).$
\item [] $791=11\times (11+1)^{(1+1)}/(1+1)-1.$
\item [] $792=11\times (11+1)^{(1+1)}/(1+1).$
\item [] $793=11\times (11+1)^{(1+1)}/(1+1)+1.$
\item [] $794=11\times (11+1)^{(1+1)}/(1+1)+1+1.$
\item [] $795=11+((1+1+1)^{(1+1+1)}+1)^{(1+1)}.$
\item [] $796=(1+1)\times (((1+1)\times (11-1))^{(1+1)}-1-1).$
\item [] $797=(11\times ((11+1)^{(1+1)}+1)-1)/(1+1).$
\item [] $798=(1+1)\times (((1+1)\times (11-1))^{(1+1)}-1).$
\item [] $799=(1+1)\times ((1+1)\times (11-1))^{(1+1)}-1.$
\item [] $800=(1+1)\times ((1+1)\times (11-1))^{(1+1)}.$
\item [] $801=(1+1)\times ((1+1)\times (11-1))^{(1+1)}+1.$
\item [] $802=(1+1)\times (((1+1)\times (11-1))^{(1+1)}+1).$
\item [] $803=11\times ((11+1)^{(1+1)}/(1+1)+1).$
\item [] $804=11\times ((11+1)^{(1+1)}/(1+1)+1)+1.$
\item [] $805=(1+1)^{11}-11\times (111+1+1).$
\item [] $806=(1+1)^{11}-11\times (111+1+1)+1.$
\item [] $807=(11-1)\times (11-1-1)^{(1+1)}-1-1-1.$
\item [] $808=(11-1)\times (11-1-1)^{(1+1)}-1-1.$
\item [] $809=(11-1)\times (11-1-1)^{(1+1)}-1.$
\item [] $810=(11-1)\times (11-1-1)^{(1+1)}.$
\item [] $811=(11-1)\times (11-1-1)^{(1+1)}+1.$
\item [] $812=(11-1)\times (11-1-1)^{(1+1)}+1+1.$
\item [] $813=(1+1)\times 11\times 111/(1+1+1)-1.$
\item [] $814=(1+1)\times 11\times 111/(1+1+1).$
\item [] $815=(1+1)\times 11\times 111/(1+1+1)+1.$
\item [] $816=(1+1)^{11}-11\times (111+1).$
\item [] $817=(1+1)^{11}-11\times (111+1)+1.$
\item [] $818=(1+1)^{11}-11\times (111+1)+1+1.$
\item [] $819=11^{(1+1+1)}-(1+1)^{(11-1-1)}.$
\item [] $820=(11-1)\times ((11-1-1)^{(1+1)}+1).$
\item [] $821=(11-1)\times ((11-1-1)^{(1+1)}+1)+1.$
\item [] $822=(1+1)\times (((1+1)\times (11-1))^{(1+1)}+11).$
\item [] $823=1111-(1+1)\times (11+1)^{(1+1)}.$
\item [] $824=1111-(1+1)\times (11+1)^{(1+1)}+1.$
\item [] $825=11\times ((1+1)\times 111/(1+1+1)+1).$
\item [] $826=(1+1)^{11}-11\times 111-1.$
\item [] $827=(1+1)^{11}-11\times 111.$
\item [] $828=(1+1)^{11}-11\times 111+1.$
\item [] $829=(1+1)^{11}-11\times 111+1+1.$
\item [] $830=(11-1)\times ((11-1-1)^{(1+1)}+1+1).$
\item [] $831=(11-1)\times ((11-1-1)^{(1+1)}+1+1)+1.$
\item [] $832=(11+1+1)\times (1+1)^{((1+1)\times (1+1+1))}.$
\item [] $833=(11+1+1)\times (1+1)^{((1+1)\times (1+1+1))}+1.$
\item [] $834=(1+1)\times (11\times (111/(1+1+1)+1)-1).$
\item [] $835=1111-(11+1)\times (11+11+1).$
\item [] $836=(1+1)\times 11\times (111/(1+1+1)+1).$
\item [] $837=(1+1)^{11}-11\times (111-1)-1.$
\item [] $838=(1+1)^{11}-11\times (111-1).$
\item [] $839=(1+1)^{11}-11\times (111-1)+1.$
\item [] $840=(11-1-1)^{(1+1+1)}+111.$
\item [] $841=((1+1+1)\times (11-1)-1)^{(1+1)}.$
\item [] $842=((1+1+1)\times (11-1)-1)^{(1+1)}+1.$
\item [] $843=((1+1+1)\times (11-1)-1)^{(1+1)}+1+1.$
\item [] $844=(1+1)\times (1+1)\times ((1+1)\times 111-11).$
\item [] $845=(1+1+1+1+1)\times (11+1+1)^{(1+1)}.$
\item [] $846=11\times 11\times ((1+1)\times (1+1+1)+1)-1.$
\item [] $847=11\times 11\times ((1+1)\times (1+1+1)+1).$
\item [] $848=11\times 11\times ((1+1)\times (1+1+1)+1)+1.$
\item [] $849=((1+1)^{11})+(11\times ((1-111)+1)).$
\item [] $850=11^{(1+1)}+(11-1-1)^{(1+1+1)}.$
\item [] $851=(11+11+1)\times 111/(1+1+1).$
\item [] $852=((1+1+1)\times (11-1)-1)^{(1+1)}+11.$
\item [] $853=(11+1)^{(1+1+1)}/(1+1)-11.$
\item [] $854=(11+1)^{(1+1+1)}/(1+1)-11+1.$
\item [] $855=1111-(1+1)^{((1+1)^{(1+1+1)})}.$
\item [] $856=(11-1)^{(1+1+1)}-(11+1)^{(1+1)}.$
\item [] $857=(1+1)\times (11+11)^{(1+1)}-111.$
\item [] $858=11\times (111-11\times (1+1+1)).$
\item [] $859=11\times (111-11\times (1+1+1))+1.$
\item [] $860=(1+1)\times ((11+11-1)^{(1+1)}-11).$
\item [] $861=(11+1)^{(1+1+1)}/(1+1)-1-1-1.$
\item [] $862=(11+1)^{(1+1+1)}/(1+1)-1-1.$
\item [] $863=(11+1)^{(1+1+1)}/(1+1)-1.$
\item [] $864=(11+1)^{(1+1+1)}/(1+1).$
\item [] $865=(11+1)^{(1+1+1)}/(1+1)+1.$
\item [] $866=(11+1)^{(1+1+1)}/(1+1)+1+1.$
\item [] $867=(11+1)^{(1+1+1)}/(1+1)+1+1+1.$
\item [] $868=1111-(1+1)\times 11^{(1+1)}-1.$
\item [] $869=11\times ((11-1-1)^{(1+1)}-1-1).$
\item [] $870=11\times ((11-1-1)^{(1+1)}-1-1)+1.$
\item [] $871=1111-(1+1)\times (11^{(1+1)}-1).$
\item [] $872=(111-1-1)\times (11-1-1-1).$
\item [] $873=(111-1-1)\times (11-1-1-1)+1.$
\item [] $874=11+((11+1)^{(1+1+1)}/(1+1)-1).$
\item [] $875=11+(11+1)^{(1+1+1)}/(1+1).$
\item [] $876=11+(11+1)^{(1+1+1)}/(1+1)+1.$
\item [] $877=(1+1)^{(1+1+1)}\times 111-11.$
\item [] $878=(1+1)^{(1+1+1)}\times 111-11+1.$
\item [] $879=(11-1)^{(1+1+1)}-11^{(1+1)}.$
\item [] $880=11\times ((11-1-1)^{(1+1)}-1).$
\item [] $881=11\times ((11-1-1)^{(1+1)}-1)+1.$
\item [] $882=(1+1)\times (11+11-1)^{(1+1)}.$
\item [] $883=(1+1)\times (11+11-1)^{(1+1)}+1.$
\item [] $884=(1+1)\times ((11+11-1)^{(1+1)}+1).$
\item [] $885=(1+1)\times ((11+11-1)^{(1+1)}+1)+1.$
\item [] $886=(1+1)\times ((1+1)\times (1+1)\times 111-1).$
\item [] $887=(1+1)^{(1+1+1)}\times 111-1.$
\item [] $888=(1+1)^{(1+1+1)}\times 111.$
\item [] $889=(1+1)^{(1+1+1)}\times 111+1.$
\item [] $890=11\times (11-1-1)^{(1+1)}-1.$
\item [] $891=11\times (11-1-1)^{(1+1)}.$
\item [] $892=11\times (11-1-1)^{(1+1)}+1.$
\item [] $893=11\times (11-1-1)^{(1+1)}+1+1.$
\item [] $894=11\times (11-1-1)^{(1+1)}+1+1+1.$
\item [] $895=(111+1)\times (1+1)^{(1+1+1)}-1.$
\item [] $896=(111+1)\times (1+1)^{(1+1+1)}.$
\item [] $897=(111+1)\times (1+1)^{(1+1+1)}+1.$
\item [] $898=((1+1+1)\times (11-1))^{(1+1)}-1-1.$
\item [] $899=((1+1+1)\times (11-1))^{(1+1)}-1.$
\item [] $900=((1+1+1)\times (11-1))^{(1+1)}.$
\item [] $901=((1+1+1)\times (11-1))^{(1+1)}+1.$
\item [] $902=11\times ((11-1-1)^{(1+1)}+1).$
\item [] $903=(1+1)^{(11-1)}-11^{(1+1)}.$
\item [] $904=(1+1)^{(11-1)}-11^{(1+1)}+1.$
\item [] $905=(1+1)^{(11-1)}-11^{(1+1)}+1+1.$
\item [] $906=(1+1)^{(11-1)}-11^{(1+1)}+1+1+1.$
\item [] $907=(111+1)\times (1+1)^{(1+1+1)}+11.$
\item [] $908=(11-1-1)\times ((11-1)^{(1+1)}+1)-1.$
\item [] $909=(11-1-1)\times ((11-1)^{(1+1)}+1).$
\item [] $910=(11-1-1)\times ((11-1)^{(1+1)}+1)+1.$
\item [] $911=((1+1+1)\times (11-1))^{(1+1)}+11.$
\item [] $912=(1+1)^{(11-1)}-111-1.$
\item [] $913=(1+1)^{(11-1)}-111.$
\item [] $914=(1+1)^{(11-1)}-111+1.$
\item [] $915=(1+1)^{(11-1)}-111+1+1.$
\item [] $916=(1+1)^{(11-1)}-111+1+1+1.$
\item [] $917=(1+1)^{(11-1)}-111+1+1+1+1.$
\item [] $918=(11-1-1)\times ((11-1)^{(1+1)}+1+1).$
\item [] $919=(11-1-1)\times ((11-1)^{(1+1)}+1+1)+1.$
\item [] $920=(11-1)\times ((11-1-1)^{(1+1)}+11).$
\item [] $921=(11-1)\times ((11-1-1)^{(1+1)}+11)+1.$
\item [] $922=((1+1+1)\times (11-1))^{(1+1)}+11+11.$
\item [] $923=(1+1)^{(11-1)}+11-111-1.$
\item [] $924=(1+1)^{(11-1)}+11-111.$
\item [] $925=(11111-11)/(11+1).$
\item [] $926=(11111+1)/(11+1).$
\item [] $927=(11111+1)/(11+1)+1.$
\item [] $928=(11111+1)/(11+1)+1+1.$
\item [] $929=(11111+1)/(11+1)+1+1+1.$
\item [] $930=(11-1)\times ((1+1)^{(11-1)}-1)/11.$
\item [] $931=((11-1)\times (1+1)^{(11-1)}+1)/11.$
\item [] $932=(1+1)\times (1+1)\times ((1+1)\times 111+11).$
\item [] $933=(1+1)\times (1+1)\times ((1+1)\times 111+11)+1.$
\item [] $934=(1+1)^{11}-1111-1-1-1.$
\item [] $935=(1+1)^{11}-1111-1-1.$
\item [] $936=(1+1)^{11}-1111-1.$
\item [] $937=(1+1)^{11}-1111.$
\item [] $938=(1+1)^{11}-1111+1.$
\item [] $939=(1+1)^{11}-1111+1+1.$
\item [] $940=(1+1)^{11}-1111+1+1+1.$
\item [] $941=1111-(11+1+1)^{(1+1)}-1.$
\item [] $942=1111-(11+1+1)^{(1+1)}.$
\item [] $943=1111-(11+1+1)^{(1+1)}+1.$
\item [] $944=(1+1)\times ((11+11)^{(1+1)}-11-1).$
\item [] $945=(1+1)\times ((11+11)^{(1+1)}-11)-1.$
\item [] $946=(1+1)\times ((11+11)^{(1+1)}-11).$
\item [] $947=(1+1)\times ((11+11)^{(1+1)}-11)+1.$
\item [] $948=(1+1)^{11}-1111+11.$
\item [] $949=(1+1)^{11}-1111+11+1.$
\item [] $950=((1+1+1)\times (11-1)+1)^{(1+1)}-11.$
\item [] $951=((1+1+1)\times (11-1)+1)^{(1+1)}-11+1.$
\item [] $952=(1+1)^{(1+1+1)}\times (11^{(1+1)}-1-1).$
\item [] $953=(1+1)^{(1+1+1)}\times (11^{(1+1)}-1-1)+1.$
\item [] $954=(11-1-1)\times (111-(11-1)/(1+1)).$
\item [] $955=1111-(11+1)\times (11+1+1).$
\item [] $956=1111-(11+1)^{(1+1)}-11.$
\item [] $957=11\times (111-(1+1)\times (11+1)).$
\item [] $958=11\times (111-(1+1)\times (11+1))+1.$
\item [] $959=(1+1)^{11}-(11\times (1+1+1))^{(1+1)}.$
\item [] $960=(11+1)\times ((11-1-1)^{(1+1)}-1).$
\item [] $961=((1+1+1)\times (11-1)+1)^{(1+1)}.$
\item [] $962=((1+1+1)\times (11-1)+1)^{(1+1)}+1.$
\item [] $963=((1+1+1)\times (11-1)+1)^{(1+1)}+1+1.$
\item [] $964=(1+1)\times ((11+11)^{(1+1)}-1-1).$
\item [] $965=(1+1)\times ((11+11)^{(1+1)}-1)-1.$
\item [] $966=(1+1)\times ((11+11)^{(1+1)}-1).$
\item [] $967=1111-(11+1)^{(1+1)}.$
\item [] $968=(1+1)\times (11+11)^{(1+1)}.$
\item [] $969=(1+1)\times (11+11)^{(1+1)}+1.$
\item [] $970=(1+1)\times ((11+11)^{(1+1)}+1).$
\item [] $971=(1+1)\times ((11+11)^{(1+1)}+1)+1.$
\item [] $972=(11+1)\times (11-1-1)^{(1+1)}.$
\item [] $973=(11+1)\times (11-1-1)^{(1+1)}+1.$
\item [] $974=(11+1)\times (11-1-1)^{(1+1)}+1+1.$
\item [] $975=(11+1)^{(1+1+1)}/(1+1)+111.$
\item [] $976=(1+1)^{(1+1+1)}\times (11^{(1+1)}+1).$
\item [] $977=(1+1)^{(1+1+1)}\times (11^{(1+1)}+1)+1.$
\item [] $978=11\times (111-11-11)-1.$
\item [] $979=11\times (111-11-11).$
\item [] $980=11\times (111-11-11)+1.$
\item [] $981=(11-1-1)\times (111-1-1).$
\item [] $982=(11-1-1)\times (111-1-1)+1.$
\item [] $983=(11-1-1)\times (111-1-1)+1+1.$
\item [] $984=(11+1)\times ((11-1-1)^{(1+1)}+1).$
\item [] $985=(11+1)\times ((11-1-1)^{(1+1)}+1)+1.$
\item [] $986=(11-1)^{(1+1+1)}-11-1-1-1.$
\item [] $987=(11-1)^{(1+1+1)}-11-1-1.$
\item [] $988=(11-1)^{(1+1+1)}-11-1.$
\item [] $989=(11-1)^{(1+1+1)}-11.$
\item [] $990=(11-1-1)\times (111-1).$
\item [] $991=(11-1-1)\times (111-1)+1.$
\item [] $992=(11-1-1)\times (111-1)+1+1.$
\item [] $993=(11-1-1)\times (111-1)+1+1+1.$
\item [] $994=(1+1)^{(11-1)}-(1+1+1)\times (11-1).$
\item [] $995=(11-1)^{(1+1+1)}-(11-1)/(1+1).$
\item [] $996=(1+1+1)\times ((1+1+1)\times 111-1).$
\item [] $997=(11-1)^{(1+1+1)}-1-1-1.$
\item [] $998=(11-1)^{(1+1+1)}-1-1.$
\item [] $999=111\times (11-1-1).$
\item [] $1000=(11-1)^{(1+1+1)}.$
\end{itemize}
\end{multicols}
}

\section{\textbf{Representations Using Number 2}}

{\footnotesize
\begin{multicols}{3}
\begin{itemize}
\item [] $101=2222/22.$
\item [] $102=2+(2\times (2+2)+2)^2.$
\item [] $103=2+2222/22.$
\item [] $104=2\times 2\times (22+2+2).$
\item [] $105=222/2-2-2-2.$
\item [] $106=2+2\times 2\times (22+2+2).$
\item [] $107=222/2-2-2.$
\item [] $108=(222-2)/2-2.$
\item [] $109=222/2-2.$
\item [] $110=(222-2)/2.$
\item [] $111=222/2.$
\item [] $112=(222+2)/2.$
\item [] $113=2+222/2.$
\item [] $114=2+(222+2)/2.$
\item [] $115=2+2+222/2.$
\item [] $116=2+2+(222+2)/2.$
\item [] $117=(22/2)^2-2-2.$
\item [] $118=2\times 2\times (22+2)+22.$
\item [] $119=(22/2)^2-2.$
\item [] $120=(2+2+2)\times (22-2).$
\item [] $121=(22/2)^2.$
\item [] $122=(22/2)^2+2/2.$
\item [] $123=(22/2)^2+2.$
\item [] $124=2\times (2^{(2+2+2)}-2).$
\item [] $125=(22/2)^2+2+2.$
\item [] $126=2\times 2^{(2+2+2)}-2.$
\item [] $127=(2^{2\times (2+2)}-2)/2.$
\item [] $128=2\times 2^{(2+2+2)}.$
\item [] $129=(2^{2\times (2+2)}+2)/2.$
\item [] $130=2\times 2^{(2+2+2)}+2.$
\item [] $131=22\times (2+2+2)-2/2.$
\item [] $132=22\times (2+2+2).$
\item [] $133=22+222/2.$
\item [] $134=22\times (2+2+2)+2.$
\item [] $135=22+2+222/2.$
\item [] $136=22\times (2+2+2)+2+2.$
\item [] $137=2^{(2+2)}+(22/2)^2.$
\item [] $138=(2+2/2)\times (2\times 22+2).$
\item [] $139=(2^{2\times (2+2)}+22)/2.$
\item [] $140=2\times (2\times (22+2)+22).$
\item [] $141=22-2+(22/2)^2.$
\item [] $142=(2\times (2+2+2))^2-2.$
\item [] $143=22+(22/2)^2.$
\item [] $144=(2\times (2+2+2))^2.$
\item [] $145=(22/2)^2+22+2.$
\item [] $146=(2\times (2+2+2))^2+2.$
\item [] $147=(2+22/2)^2-22.$
\item [] $148=(2\times (2+2+2))^2+2+2.$
\item [] $149=(2+22/2)^2-22+2.$
\item [] $150=2\times 2^{(2+2+2)}+22.$
\item [] $151=(2^{2\times (2+2)}+2)/2+22.$
\item [] $152=2\times 2\times ((2+2+2)^2+2).$
\item [] $153=2\times 22-2+222/2.$
\item [] $154=22+22\times (2+2+2).$
\item [] $155=2\times 22+222/2.$
\item [] $156=2\times (2\times 2\times (22-2)-2).$
\item [] $157=2\times 22+2+222/2.$
\item [] $158=2\times ((2+2/2)^{(2+2)}-2).$
\item [] $159=2\times (22+2)+222/2.$
\item [] $160=2\times 2\times 2\times (22-2).$
\item [] $161=((2^{(2+2)}+2)^2-2)/2.$
\item [] $162=2\times (2+2/2)^{(2+2)}.$
\item [] $163=((2^{(2+2)}+2)^2+2)/2.$
\item [] $164=2\times (2+2/2)^{(2+2)}+2.$
\item [] $165=2\times 22+(22/2)^2.$
\item [] $166=2\times ((2+2/2)^{(2+2)}+2).$
\item [] $167=(2+22/2)^2-2.$
\item [] $168=2\times 2\times (2\times 22-2).$
\item [] $169=(2+22/2)^2.$
\item [] $170=2\times 2\times (2\times 22-2)+2.$
\item [] $171=(2+22/2)^2+2.$
\item [] $172=2\times (2\times 2\times 22-2).$
\item [] $173=(2+22/2)^2+2+2.$
\item [] $174=2\times 2\times 2\times 22-2.$
\item [] $175=2\times 2\times 2\times 22-2/2.$
\item [] $176=2\times 2\times 2\times 22.$
\item [] $177=2\times 2\times 2\times 22+2/2.$
\item [] $178=2\times 2\times 2\times 22+2.$
\item [] $179=2\times 2\times 2\times 22+2+2/2.$
\item [] $180=2\times (2\times 2\times 22+2).$
\item [] $181=2\times (2\times 2\times 22+2)+2/2.$
\item [] $182=2\times (2\times 2\times 22+2)+2.$
\item [] $183=2\times 2\times (2\times 22+2)-2/2.$
\item [] $184=2\times 2\times (2\times 22+2).$
\item [] $185=2\times 2\times (2\times 22+2)+2/2.$
\item [] $186=2\times 2\times (2\times 22+2)+2.$
\item [] $187=2\times 2\times 2\times 22+22/2.$
\item [] $188=2\times (2\times (2\times 22+2)+2).$
\item [] $189=((22-2)^2-22)/2.$
\item [] $190=222-2\times 2^{(2+2)}.$
\item [] $191=22+(2+22/2)^2.$
\item [] $192=2\times 2\times 2\times (22+2).$
\item [] $193=(2+22/2)^2+22+2.$
\item [] $194=(2^{(2+2)}-2)^2-2.$
\item [] $195=(2^{(2+2)}-2)^2-2/2.$
\item [] $196=(2^{(2+2)}-2)^2.$
\item [] $197=(2^{(2+2)}-2)^2+2/2.$
\item [] $198=(2^{(2+2)}-2)^2+2.$
\item [] $199=((22-2)^2-2)/2.$
\item [] $200=222-22.$
\item [] $201=((22-2)^2+2)/2.$
\item [] $202=222-22+2.$
\item [] $203=((22-2)^2+2)/2+2.$
\item [] $204=222-22+2+2.$
\item [] $205=((22-2)^2+2)/2+2+2.$
\item [] $206=222-2^{(2+2)}.$
\item [] $207=222-2^{(2+2)}+2/2.$
\item [] $208=222-2^{(2+2)}+2.$
\item [] $209=222-2-22/2.$
\item [] $210=222-2\times (2+2+2).$
\item [] $211=222-22/2.$
\item [] $212=222-(22-2)/2.$
\item [] $213=222+2-22/2.$
\item [] $214=222-2\times (2+2).$
\item [] $215=222+2+2-22/2.$
\item [] $216=(2+2+2)^{(2+2/2)}.$
\item [] $217=222-2-2-2/2.$
\item [] $218=222-2-2.$
\item [] $219=222-2-2/2.$
\item [] $220=222-2.$
\item [] $221=222-2/2.$
\item [] $222=222.$
\item [] $223=222+2/2.$
\item [] $224=222+2.$
\item [] $225=222+2+2/2.$
\item [] $226=222+2+2.$
\item [] $227=222+2+2+2/2.$
\item [] $228=222+2+2+2.$
\item [] $229=(22^2-22)/2-2.$
\item [] $230=222+2\times (2+2).$
\item [] $231=(22^2-22)/2.$
\item [] $232=222+2\times (2+2)+2.$
\item [] $233=222+22/2.$
\item [] $234=2^{2\times (2+2)}-22.$
\item [] $235=222+2+22/2.$
\item [] $236=2^{2\times (2+2)}-22+2.$
\item [] $237=(22^2-2)/2-2-2.$
\item [] $238=2\times ((22/2)^2-2).$
\item [] $239=(22^2-2)/2-2.$
\item [] $240=22^2/2-2.$
\item [] $241=(22^2-2)/2.$
\item [] $242=22^2/2.$
\item [] $243=(22^2+2)/2.$
\item [] $244=22^2/2+2.$
\item [] $245=(22^2+2)/2+2.$
\item [] $246=22^2/2+2+2.$
\item [] $247=(22^2+2)/2+2+2.$
\item [] $248=22^2/2+2+2+2.$
\item [] $249=(22^2+2)/2+2+2+2.$
\item [] $250=2\times ((22/2)^2+2+2).$
\item [] $251=(22^2+22)/2-2.$
\item [] $252=2^{2\times (2+2)}-2-2.$
\item [] $253=(22^2+22)/2.$
\item [] $254=2^{2\times (2+2)}-2.$
\item [] $255=2^{2\times (2+2)}-2/2.$
\item [] $256=2^{2\times (2+2)}.$
\item [] $257=2^{2\times (2+2)}+2/2.$
\item [] $258=2^{2\times (2+2)}+2.$
\item [] $259=2^{2\times (2+2)}+2+2/2.$
\item [] $260=2^{2\times (2+2)}+2+2.$
\item [] $261=2^{2\times (2+2)}+2+2+2/2.$
\item [] $262=22^2-222.$
\item [] $263=22+(22^2-2)/2.$
\item [] $264=2\times 22\times (2+2+2).$
\item [] $265=22+(22^2+2)/2.$
\item [] $266=222+2\times 22.$
\item [] $267=2^{2\times (2+2)}+22/2.$
\item [] $268=222+2\times 22+2.$
\item [] $269=2+2^{2\times (2+2)}+22/2.$
\item [] $270=222+2\times (22+2).$
\item [] $271=2^{2\times (2+2)}+2+2+22/2.$
\item [] $272=2^{2\times (2+2)}+2^{(2+2)}.$
\item [] $273=(2+22/2)\times (22-2/2).$
\item [] $274=2^{2\times (2+2)}+2^{(2+2)}+2.$
\item [] $275=22+(22^2+22)/2.$
\item [] $276=22+2^{2\times (2+2)}-2.$
\item [] $277=((22+2)^2-22)/2.$
\item [] $278=22+2^{2\times (2+2)}.$
\item [] $279=((22+2)^2-22)/2+2.$
\item [] $280=2^{2\times (2+2)}+22+2.$
\item [] $281=((22+2)^2-22)/2+2+2.$
\item [] $282=2^{2\times (2+2)}+22+2+2.$
\item [] $283=((22+2)^2-2)/2-2-2.$
\item [] $284=22\times (2+22/2)-2.$
\item [] $285=((22+2)^2-2)/2-2.$
\item [] $286=22\times (2+22/2).$
\item [] $287=((22+2)^2-2)/2.$
\item [] $288=(22+2)^2/2.$
\item [] $289=(2^{(2+2)}+2/2)^2.$
\item [] $290=(22+2)^2/2+2.$
\item [] $291=(2^{(2+2)}+2/2)^2+2.$
\item [] $292=(22+2)^2/2+2+2.$
\item [] $293=(2^{(2+2)}+2/2)^2+2+2.$
\item [] $294=(22+2)^2/2+2+2+2.$
\item [] $295=(2^{(2+2)}+2/2)^2+2+2+2.$
\item [] $296=2\times ((2\times (2+2+2))^2+2+2).$
\item [] $297=((22+2)^2+22)/2-2.$
\item [] $298=2\times 22+2^{2\times (2+2)}-2.$
\item [] $299=((22+2)^2+22)/2.$
\item [] $300=2\times 22+2^{2\times (2+2)}.$
\item [] $301=((22+2)^2+22)/2+2.$
\item [] $302=(2^{(2+2)}+2)^2-22.$
\item [] $303=(2+2/2)^{(2+2)}+222.$
\item [] $304=(2^{(2+2)}+2)^2-22+2.$
\item [] $305=(2+2/2)^{(2+2)}+222+2.$
\item [] $306=22\times (2^{(2+2)}-2)-2.$
\item [] $307=(22+2/2)^2-222.$
\item [] $308=22\times (2^{(2+2)}-2).$
\item [] $309=((22+2)^2-2)/2+22.$
\item [] $310=22\times (2^{(2+2)}-2)+2.$
\item [] $311=(2^{(2+2)}+2/2)^2+22.$
\item [] $312=(22+2)\times (2+22/2).$
\item [] $313=(2^{(2+2)}+2)^2-22/2.$
\item [] $314=(22+2)\times (2+22/2)+2.$
\item [] $315=22^2-(2+22/2)^2.$
\item [] $316=2\times 2\times ((2+2/2)^{(2+2)}-2).$
\item [] $317=22^2-(2+22/2)^2+2.$
\item [] $318=2^{(2+2)}\times (22-2)-2.$
\item [] $319=2^{(2+2)}\times (22-2)-2/2.$
\item [] $320=(22-2)\times 2^{(2+2)}.$
\item [] $321=(22-2)\times 2^{(2+2)}+2/2.$
\item [] $322=(2^{(2+2)}+2)^2-2.$
\item [] $323=(2^{(2+2)}+2)^2-2/2.$
\item [] $324=(2^{(2+2)}+2)^2-2.$
\item [] $325=(2^{(2+2)}+2)^2+2/2.$
\item [] $326=(2^{(2+2)}+2)^2+2.$
\item [] $327=(2^{(2+2)}+2)^2+2+2/2.$
\item [] $328=(2^{(2+2)}+2)^2+2+2.$
\item [] $329=(2^{(2+2)}+2)^2+2+2+2/2.$
\item [] $330=22\times (2+2+22/2).$
\item [] $331=222-2+222/2.$
\item [] $332=22\times (2+2+22/2)+2.$
\item [] $333=222+222/2.$
\item [] $334=2\times ((2+22/2)^2-2).$
\item [] $335=222+2+222/2.$
\item [] $336=2\times 2\times 2\times (2\times 22-2).$
\item [] $337=((22+2+2)^2-2)/2.$
\item [] $338=2\times (2+22/2)^2.$
\item [] $339=((22+2+2)^2+2)/2.$
\item [] $340=2\times (2+22/2)^2+2.$
\item [] $341=((22+2+2)^2+2)/2+2.$
\item [] $342=2\times ((2+22/2)^2+2).$
\item [] $343=(22/2)^2+222.$
\item [] $344=2\times 2\times (2\times 2\times 22-2).$
\item [] $345=(22/2)^2+222+2.$
\item [] $346=(2^{(2+2)}+2)^2+22.$
\item [] $347=(22/2)^2+222+2+2.$
\item [] $348=2\times (2\times 2\times 2\times 22-2).$
\item [] $349=22\times 2^{(2+2)}-2-2/2.$
\item [] $350=22\times 2^{(2+2)}-2.$
\item [] $351=22\times 2^{(2+2)}-2/2.$
\item [] $352=22\times 2^{(2+2)}.$
\item [] $353=22\times 2^{(2+2)}+2/2.$
\item [] $354=22\times 2^{(2+2)}+2.$
\item [] $355=22\times 2^{(2+2)}+2+2/2.$
\item [] $356=22\times 2^{(2+2)}+2+2.$
\item [] $357=(2+2/2)\times ((22/2)^2-2).$
\item [] $358=22\times 2^{(2+2)}+2+2+2.$
\item [] $359=(22-2-2/2)^2-2.$
\item [] $360=2\times 2\times (2\times 2\times 22+2).$
\item [] $361=(22-2-2/2)^2.$
\item [] $362=2\times 2\times (2\times 2\times 22+2)+2.$
\item [] $363=(22-2-2/2)^2+2.$
\item [] $364=2\times (2\times (2\times 2\times 22+2)+2).$
\item [] $365=(22-2-2/2)^2+2+2.$
\item [] $366=(2\times (2+2+2))^2+222.$
\item [] $367=2^{2\times (2+2)}+222/2.$
\item [] $368=2\times 2\times 2\times (2\times 22+2).$
\item [] $369=(2+2/2)\times (2+(22/2)^2).$
\item [] $370=2\times 2\times 2\times (2\times 22+2)+2.$
\item [] $371=22^2-2-222/2.$
\item [] $372=2\times (2\times 2\times (2\times 22+2)+2).$
\item [] $373=22^2-222/2.$
\item [] $374=22\times 2^{(2+2)}+22.$
\item [] $375=22^2+2-222/2.$
\item [] $376=(22-2)^2-22-2.$
\item [] $377=(22/2)^2+2^{2\times (2+2)}.$
\item [] $378=(22-2)^2-22.$
\item [] $379=(22-2)^2-22+2/2.$
\item [] $380=(22-2)^2-22+2.$
\item [] $381=(22-2)^2-22+2+2/2.$
\item [] $382=2^{(2+2)}\times (22+2)-2.$
\item [] $383=(22-2-2/2)^2+22.$
\item [] $384=2^{(2+2)}\times (22+2).$
\item [] $385=2^{(2+2)}\times (22+2)+2/2.$
\item [] $386=2^{(2+2)}\times (22+2)+2.$
\item [] $387=(22-2)^2-2-22/2.$
\item [] $388=2\times ((2^{(2+2)}-2)^2-2).$
\item [] $389=(22-2)^2-22/2.$
\item [] $390=2\times (2^{(2+2)}-2)^2-2.$
\item [] $391=(22-2)^2+2-22/2.$
\item [] $392=2\times (2^{(2+2)}-2)^2.$
\item [] $393=2\times (2^{(2+2)}-2)^2+2/2.$
\item [] $394=2\times (2^{(2+2)}-2)^2+2.$
\item [] $395=(22-2)^2-2-2-2/2.$
\item [] $396=22\times (2^{(2+2)}+2).$
\item [] $397=(22-2)^2-2-2/2.$
\item [] $398=(22-2)^2-2.$
\item [] $399=(22-2)^2-2/2.$
\item [] $400=(22-2)^2.$
\item [] $401=(22-2)^2+2/2.$
\item [] $402=(22-2)^2+2.$
\item [] $403=(22-2)^2+2+2/2.$
\item [] $404=(22-2)^2+2+2.$
\item [] $405=(22-2)^2+2+2+2/2.$
\item [] $406=(22-2)^2+2+2+2.$
\item [] $407=(22-2)^2+2+2+2+2/2.$
\item [] $408=2\times (2+2)+(22-2)^2.$
\item [] $409=(22-2)^2-2+22/2.$
\item [] $410=2\times (2+2)+(22-2)^2+2.$
\item [] $411=(22-2)^2+22/2.$
\item [] $412=2\times (222-2^{(2+2)}).$
\item [] $413=(22-2)^2+2+22/2.$
\item [] $414=2\times (222-2^{(2+2)})+2.$
\item [] $415=(22-2)^2+2+2+22/2.$
\item [] $416=2^{(2+2)}+(22-2)^2.$
\item [] $417=(22-2/2)^2-22-2.$
\item [] $418=22\times (22-2-2/2).$
\item [] $419=(22-2/2)^2-22.$
\item [] $420=(22-2)^2+22-2.$
\item [] $421=(22-2/2)^2-22+2.$
\item [] $422=(22-2)^2+22.$
\item [] $423=(22-2)^2+22+2/2.$
\item [] $424=(22-2)^2+22+2.$
\item [] $425=(22-2/2)^2-2^{(2+2)}.$
\item [] $426=2\times (222+2)-22.$
\item [] $427=(22-2/2)^2-2^{(2+2)}+2.$
\item [] $428=2\times (222-2\times (2+2)).$
\item [] $429=2\times (222-2)-22/2.$
\item [] $430=2\times (222-2\times (2+2))+2.$
\item [] $431=2\times 222-2-22/2.$
\item [] $432=2\times (2+2+2)^{(2+2/2)}.$
\item [] $433=2\times 222-22/2.$
\item [] $434=2\times (222-2-2)-2.$
\item [] $435=2\times 222+2-22/2.$
\item [] $436=2\times (222-2-2).$
\item [] $437=(22-2/2)^2-2-2.$
\item [] $438=2\times (222-2)-2.$
\item [] $439=(22-2/2)^2-2.$
\item [] $440=2\times (222-2).$
\item [] $441=(22-2/2)^2.$
\item [] $442=2\times 222-2.$
\item [] $443=2\times 222-2/2.$
\item [] $444=2\times 222.$
\item [] $445=2\times 222+2/2.$
\item [] $446=2\times 222+2.$
\item [] $447=2\times 222+2+2/2.$
\item [] $448=2\times (222+2).$
\item [] $449=2\times (222+2)+2/2.$
\item [] $450=2\times (222+2)+2.$
\item [] $451=2\times (222+2)+2+2/2.$
\item [] $452=2\times (222+2+2).$
\item [] $453=2\times (222+2+2)+2/2.$
\item [] $454=2\times (222+2+2)+2.$
\item [] $455=2\times 222+22/2.$
\item [] $456=2\times (222+2+2+2).$
\item [] $457=2\times 222+2+22/2.$
\item [] $458=22^2-22-2-2.$
\item [] $459=2\times (222+2)+22/2.$
\item [] $460=22^2-22-2.$
\item [] $461=22^2-22-2/2.$
\item [] $462=22^2-22.$
\item [] $463=(22-2/2)^2+22.$
\item [] $464=22^2-22+2.$
\item [] $465=22^2-22+2+2/2.$
\item [] $466=2\times 222+22.$
\item [] $467=2\times 222+22+2/2.$
\item [] $468=22^2-2^{(2+2)}.$
\item [] $469=22^2-2^{(2+2)}+2/2.$
\item [] $470=22^2-2^{(2+2)}+2.$
\item [] $471=22^2-2-22/2.$
\item [] $472=22^2-2\times (2+2+2).$
\item [] $473=22^2-22/2.$
\item [] $474=22^2-(22-2)/2.$
\item [] $475=22^2+2-22/2.$
\item [] $476=22^2-2\times (2+2).$
\item [] $477=22^2+2+2-22/2.$
\item [] $478=22^2-2-2-2.$
\item [] $479=22^2-2-2-2/2.$
\item [] $480=22^2-2-2.$
\item [] $481=22^2-2-2/2.$
\item [] $482=22^2-2.$
\item [] $483=22^2-2/2.$
\item [] $484=22^2.$
\item [] $485=22^2+2/2.$
\item [] $486=22^2+2.$
\item [] $487=22^2+2+2/2.$
\item [] $488=22^2+2+2.$
\item [] $489=22^2+2+2+2/2.$
\item [] $490=22^2+2+2+2.$
\item [] $491=22^2+2+2+2+2/2.$
\item [] $492=22^2+2\times (2+2).$
\item [] $493=22^2-2+22/2.$
\item [] $494=22^2+2\times (2+2)+2.$
\item [] $495=22^2+22/2.$
\item [] $496=22^2+2\times (2+2+2).$
\item [] $497=22^2+2+22/2.$
\item [] $498=22^2+2^{(2+2)}-2.$
\item [] $499=22^2+2+2+22/2.$
\item [] $500=22^2+2^{(2+2)}.$
\item [] $501=22^2+2^{(2+2)}+2/2.$
\item [] $502=22^2+2^{(2+2)}+2.$
\item [] $503=22^2+22-2-2/2.$
\item [] $504=22^2+22-2.$
\item [] $505=22^2+22-2/2.$
\item [] $506=22^2+22.$
\item [] $507=22^2+22+2/2.$
\item [] $508=22^2+22+2.$
\item [] $509=22^2+22+2+2/2.$
\item [] $510=2^{((2+2/2)^{2})}-2.$
\item [] $511=2^{((2+2/2)^{2})}-2/2.$
\item [] $512=2^{((2+2/2)^{2})}.$
\item [] $513=2^{((2+2/2)^{2})}+2/2.$
\item [] $514=2^{((2+2/2)^{2})}+2.$
\item [] $515=2^{((2+2/2)^{2})}+2+2/2.$
\item [] $516=2\times (2^{2\times (2+2)}+2).$
\item [] $517=22^2+22+22/2.$
\item [] $518=2\times (2^{2\times (2+2)}+2)+2.$
\item [] $519=2\times (2^{2\times (2+2)}+2)+2+2/2.$
\item [] $520=2\times (2^{2\times (2+2)}+2+2).$
\item [] $521=(22-2)^2+(22/2)^2.$
\item [] $522=+22^2+(2+2+2)^2+2.$
\item [] $523=2^{((2+2/2)^{2})}+22/2.$
\item [] $524=2\times (22^2-222).$
\item [] $525=(22+2/2)^2-2-2.$
\item [] $526=22\times (22+2)-2.$
\item [] $527=(22+2/2)^2-2.$
\item [] $528=22\times (22+2).$
\item [] $529=(22+2/2)^2.$
\item [] $530=22\times (22+2)+2.$
\item [] $531=(22+2/2)^2+2.$
\item [] $532=22\times (22+2)+2+2.$
\item [] $533=(22+2/2)^2+2+2.$
\item [] $534=2^{((2+2/2)^{2})}+22.$
\item [] $535=(22+2/2)^2+2+2+2.$
\item [] $536=2^{((2+2/2)^{2})}+22+2.$
\item [] $537=2\times (2+2)+(22+2/2)^2.$
\item [] $538=2\times (2^{2\times (2+2)}+2)+22.$
\item [] $539=22\times (22+2)+22/2.$
\item [] $540=2\times (2\times (22+2)+222).$
\item [] $541=22\times (22+2)+2+22/2.$
\item [] $542=(22+2/2)^2+2+22/2.$
\item [] $543=(22+2/2)^2+2^{(2+2)}-2.$
\item [] $544=2\times (2^{2\times (2+2)}+2^{(2+2)}).$
\item [] $545=(22+2/2)^2+2^{(2+2)}.$
\item [] $546=(2^{(2+2)}+2)^2+222.$
\item [] $547=(22+2/2)^2+2^{(2+2)}+2.$
\item [] $548=22^2+2^{(2+2+2)}.$
\item [] $549=(22+2/2)^2+22-2.$
\item [] $550=22\times (22+2)+22.$
\item [] $551=(22+2/2)^2+22.$
\item [] $552=(22+2)\times (22+2/2).$
\item [] $553=(22+2/2)^2+22+2.$
\item [] $554=(22+2)^2-22.$
\item [] $555=(2222-2)/(2+2).$
\item [] $556=(2222+2)/(2+2).$
\item [] $557=(2222-2)/(2+2)+2.$
\item [] $558=(2222+2)/(2+2)+2.$
\item [] $559=(2222-2)/(2+2)+2+2.$
\item [] $560=(22+2)^2-2^{(2+2)}.$
\item [] $561=(2222+22)/(2+2).$
\item [] $562=(22+2)^2-2^{(2+2)}+2.$
\item [] $563=(22+2)^2-22/2-2.$
\item [] $564=22^2+2\times 2\times (22-2).$
\item [] $565=(22+2)^2-22/2.$
\item [] $566=(22+2)^2-(22-2)/2.$
\item [] $567=(22+2)^2+2-22/2.$
\item [] $568=(22+2)^2-2\times (2+2).$
\item [] $569=(22+2)^2+2+2-22/2.$
\item [] $570=(22+2)^2-2-2-2.$
\item [] $571=(22+2)^2-2-2-2/2.$
\item [] $572=22\times (22+2+2).$
\item [] $573=(22+2)^2-2-2/2.$
\item [] $574=(22+2)^2-2.$
\item [] $575=(22+2)^2-2/2.$
\item [] $576=(22+2)^2.$
\item [] $577=(22+2)^2+2/2.$
\item [] $578=(22+2)^2+2.$
\item [] $579=(22+2)^2+2+2/2.$
\item [] $580=(22+2)^2+2+2.$
\item [] $581=(22+2)^2+2+2+2/2.$
\item [] $582=(22+2)^2+2+2+2.$
\item [] $583=(22+2)^2+2+2+2+2/2.$
\item [] $584=(22+2)^2+2\times (2+2).$
\item [] $585=(22+2)^2-2+22/2.$
\item [] $586=(22+2)^2+2\times (2+2)+2.$
\item [] $587=(22+2)^2+22/2.$
\item [] $588=(2+2/2)\times (2^{(2+2)}-2)^2.$
\item [] $589=(22+2)^2+2+22/2.$
\item [] $590=(22+2)^2+2^{(2+2)}-2.$
\item [] $591=(22+2)^2+2+2+22/2.$
\item [] $592=2^{(2+2)}+(22+2)^2.$
\item [] $593=22^2-2+222/2.$
\item [] $594=2^{(2+2)}+(22+2)^2+2.$
\item [] $595=22^2+222/2.$
\item [] $596=(22+2)^2+22-2.$
\item [] $597=22^2+2+222/2.$
\item [] $598=(22+2)^2+22.$
\item [] $599=(22+2)^2+22+2/2.$
\item [] $600=(22+2)^2+22+2.$
\item [] $601=(22+2)^2+22+2+2/2.$
\item [] $602=(22+2)^2+22+2+2.$
\item [] $603=22^2+(22/2)^2-2.$
\item [] $604=2\times ((2^{(2+2)}+2)^2-22).$
\item [] $605=22^2+(22/2)^2.$
\item [] $606=2\times ((2^{(2+2)}+2)^2-22)+2.$
\item [] $607=22^2+(22/2)^2+2.$
\item [] $608=2^{(2+2)}\times ((2+2+2)^2+2).$
\item [] $609=22^2+(22/2)^2+2+2.$
\item [] $610=(22+2)^2+2\times 2^{(2+2)}+2.$
\item [] $611=22^2+(2^{2\times (2+2)}-2)/2.$
\item [] $612=2\times (22\times (2^{(2+2)}-2)-2).$
\item [] $613=22^2+(2^{2\times (2+2)}+2)/2.$
\item [] $614=2\times 22\times (2^{(2+2)}-2)-2.$
\item [] $615=(2/2+2+2)\times ((22/2)^2+2).$
\item [] $616=2\times 22\times (2^{(2+2)}-2).$
\item [] $617=22^2+22+222/2.$
\item [] $618=2\times 22\times (2^{(2+2)}-2)+2.$
\item [] $619=2\times 22+(22+2)^2-2/2.$
\item [] $620=2\times 22+(22+2)^2.$
\item [] $621=(2+2+2/2)^{(2+2)}-2-2.$
\item [] $622=(22-2)^2+222.$
\item [] $623=(2+2+2/2)^{(2+2)}-2.$
\item [] $624=(22+2)\times (22+2+2).$
\item [] $625=(2+2+2/2)^{(2+2)}.$
\item [] $626=(22+2)\times (22+2+2)+2.$
\item [] $627=(2+2+2/2)^{(2+2)}+2.$
\item [] $628=22^2+(2\times (2+2+2))^2.$
\item [] $629=(2+2+2/2)^{(2+2)}+2+2.$
\item [] $630=22^2+(2\times (2+2+2))^2+2.$
\item [] $631=(2+2+2/2)^{(2+2)}+2+2+2.$
\item [] $632=(22+2+2)^2-2\times 22.$
\item [] $633=2\times (2+2)+(2+2+2/2)^{(2+2)}.$
\item [] $634=(22+2+2)^2-2\times 22+2.$
\item [] $635=((2+2+2)^{(2+2)}-22)/2-2.$
\item [] $636=2\times (2^{(2+2)}\times (22-2)-2).$
\item [] $637=((2+2+2)^{(2+2)}-22)/2.$
\item [] $638=2\times 2^{(2+2)}\times (22-2)-2.$
\item [] $639=((2+2+2)^{(2+2)}-22)/2+2.$
\item [] $640=2\times 2^{(2+2)}\times (22-2).$
\item [] $641=2\times 2^{(2+2)}\times (22-2)+2/2.$
\item [] $642=2\times 2^{(2+2)}\times (22-2)+2.$
\item [] $643=(22-2)^2+(22^2+2)/2.$
\item [] $644=2\times ((2^{(2+2)}+2)^2-2).$
\item [] $645=((2+2+2)^{(2+2)}-2)/2-2.$
\item [] $646=2\times (2^{(2+2)}+2)^2-2.$
\item [] $647=((2+2+2)^{(2+2)}-2)/2.$
\item [] $648=2\times (2^{(2+2)}+2)^2.$
\item [] $649=((2+2+2)^{(2+2)}+2)/2.$
\item [] $650=2\times (2^{(2+2)}+2)^2+2.$
\item [] $651=((2+2+2)^{(2+2)}+2)/2+2.$
\item [] $652=2\times ((2^{(2+2)}+2)^2+2).$
\item [] $653=22^2+(2+22/2)^2.$
\item [] $654=(22+2+2)^2-22.$
\item [] $655=22^2+(2+22/2)^2+2.$
\item [] $656=(22+2+2)^2-22+2.$
\item [] $657=(2+2/2)\times (222-2-2/2).$
\item [] $658=22\times (2\times (2+2)+22)-2.$
\item [] $659=((2+2+2)^{(2+2)}+22)/2.$
\item [] $660=22\times (2\times (2+2)+22).$
\item [] $661=((2+2+2)^{(2+2)}+22)/2+2.$
\item [] $662=2+22\times (2\times (2+2)+22).$
\item [] $663=(22-2/2)^2+222.$
\item [] $664=(2+2/2)\times 222-2.$
\item [] $665=(22+2+2)^2-22/2.$
\item [] $666=(2+2/2)\times 222.$
\item [] $667=(2+2/2)\times 222+2/2.$
\item [] $668=(2+2/2)\times 222+2.$
\item [] $669=(2+2/2)\times (222+2/2).$
\item [] $670=(2+2/2)\times 222+2+2.$
\item [] $671=(2+2/2)\times (222+2/2)+2.$
\item [] $672=(2+2/2)\times (222+2).$
\item [] $673=(22+2+2)^2-2-2/2.$
\item [] $674=(22+2+2)^2-2.$
\item [] $675=(22+2+2)^2-2/2.$
\item [] $676=(22+2+2)^2.$
\item [] $677=(22+2+2)^2+2/2.$
\item [] $678=(22+2+2)^2+2.$
\item [] $679=(22+2+2)^2+2+2/2.$
\item [] $680=(22+2+2)^2+2+2.$
\item [] $681=(22+2+2)^2+2+2+2/2.$
\item [] $682=(22+2+2)^2+2+2+2.$
\item [] $683=22^2+((22-2)^2-2)/2.$
\item [] $684=2\times 2\times ((2+22/2)^2+2).$
\item [] $685=(22+2+2)^2-2+22/2.$
\item [] $686=2\times ((22/2)^2+222).$
\item [] $687=(22+2+2)^2+22/2.$
\item [] $688=2\times 2\times 2\times (2\times 2\times 22-2).$
\item [] $689=(22+2+2)^2+2+22/2.$
\item [] $690=((2+2/2)\times 222+22)+2.$
\item [] $691=(22+2+2)^2+2+2+22/2.$
\item [] $692=2\times ((2^{(2+2)}+2)^2+22).$
\item [] $693=(2+2/2)\times (22^2-22)/2.$
\item [] $694=(22+2+2)^2+2^{(2+2)}+2.$
\item [] $695=22^2+222-22/2.$
\item [] $696=2\times 2\times (2\times 2\times 2\times 22-2).$
\item [] $697=(22+2)^2+(22/2)^2.$
\item [] $698=(22+2+2)^2+22.$
\item [] $699=(22+2)^2+2+(22/2)^2.$
\item [] $700=2\times (22\times 2^{(2+2)}-2).$
\item [] $701=2\times (22\times 2^{(2+2)}-2)+2/2.$
\item [] $702=2\times 22\times 2^{(2+2)}-2.$
\item [] $703=2\times 22\times 2^{(2+2)}-2/2.$
\item [] $704=2\times 22\times 2^{(2+2)}.$
\item [] $705=2\times 22\times 2^{(2+2)}+2/2.$
\item [] $706=22^2+222.$
\item [] $707=22^2+222+2/2.$
\item [] $708=22^2+222+2.$
\item [] $709=22^2+222+2+2/2.$
\item [] $710=22^2+222+2+2.$
\item [] $711=2222/2-(22-2)^2.$
\item [] $712=2\times (22\times 2^{(2+2)}+2+2).$
\item [] $713=(2+2/2)^{(2+2+2)}-2^{(2+2)}.$
\item [] $714=(2+2+2)\times ((22/2)^2-2).$
\item [] $715=2\times 22\times 2^{(2+2)}+22/2.$
\item [] $716=(2+2+2)\times ((22/2)^2-2)+2.$
\item [] $717=22^2+222+22/2.$
\item [] $718=2\times ((22-2-2/2)^2-2).$
\item [] $719=(22-2)\times (2+2+2)^2-2/2.$
\item [] $720=(22-2)\times (2+2+2)^2.$
\item [] $721=(((2+2+2)^2+2)^2-2)/2.$
\item [] $722=2\times (22-2-2/2)^2.$
\item [] $723=(2+2/2)\times (22^2-2)/2.$
\item [] $724=2\times (22-2/2-2)^2+2.$
\item [] $725=22^2+(22^2-2)/2.$
\item [] $726=22\times (22+22/2).$
\item [] $727=(2+2/2)^{(2+2+2)}-2.$
\item [] $728=22^2/2+22^2+2.$
\item [] $729=(2+2/2)^{(2+2+2)}.$
\item [] $730=(2+2/2)^{(2+2+2)}+2/2.$
\item [] $731=(2+2/2)^{(2+2+2)}+2.$
\item [] $732=(2+2/2)\times (2+22^2/2).$
\item [] $733=(2+2/2)^{(2+2+2)}+2+2.$
\item [] $734=2^{((2+2/2)^{2})}+222.$
\item [] $735=2^{(2+2)}\times (2\times 22+2)-2/2.$
\item [] $736=2^{(2+2)}\times (2\times 22+2).$
\item [] $737=2^{(2+2)}\times (2\times 22+2)+2/2.$
\item [] $738=2^{(2+2)}\times (2\times 22+2)+2.$
\item [] $739=22^2+2^{2\times (2+2)}-2/2.$
\item [] $740=22^2+2^{2\times (2+2)}.$
\item [] $741=22^2+2^{2\times (2+2)}+2/2.$
\item [] $742=22^2+2^{2\times (2+2)}+2.$
\item [] $743=22^2+2^{2\times (2+2)}+2+2/2.$
\item [] $744=2\times 22^2-222-2.$
\item [] $745=2^{(2+2)}+(2+2/2)^{(2+2+2)}.$
\item [] $746=2\times 22^2-222.$
\item [] $747=2\times 22^2-222+2/2.$
\item [] $748=22\times (2\times 2^{(2+2)}+2).$
\item [] $749=22+(2+2/2)^{(2+2+2)}-2.$
\item [] $750=2+22\times (2\times 2^{(2+2)}+2).$
\item [] $751=(2+2/2)^{(2+2+2)}+22.$
\item [] $752=2\times ((22-2)^2-22-2).$
\item [] $753=(2+2/2)^{(2+2+2)}+22+2.$
\item [] $754=2\times ((22-2)^2-22)-2.$
\item [] $755=2\times ((22-2)^2-22)-2/2.$
\item [] $756=2\times ((22-2)^2-22).$
\item [] $757=2\times ((22-2)^2-22)+2/2.$
\item [] $758=2\times ((22-2)^2-22)+2.$
\item [] $759=(2+2/2)\times (22^2+22)/2.$
\item [] $760=2\times ((22-2)^2-22+2).$
\item [] $761=(2+2/2)\times (22^2+22)/2+2.$
\item [] $762=(2+2/2)\times (2^{2\times (2+2)}-2).$
\item [] $763=(2+2/2)\times (2^{2\times (2+2)}-2)+2/2.$
\item [] $764=2\times (2^{(2+2)}\times (22+2)-2).$
\item [] $765=(2+2/2)\times (2^{2\times (2+2)}-2/2).$
\item [] $766=2\times 2^{(2+2)}\times (22+2)-2.$
\item [] $767=2\times 2\times 222-(22/2)^2.$
\item [] $768=2\times 2^{(2+2)}\times (22+2).$
\item [] $769=2\times 2^{(2+2)}\times (22+2)+2/2.$
\item [] $770=2\times 2^{(2+2)}\times (22+2)+2.$
\item [] $771=22^2+((22+2)^2-2)/2.$
\item [] $772=2\times 2^{(2+2)}\times (22+2)+2+2.$
\item [] $773=2\times 22+(2+2/2)^{(2+2+2)}.$
\item [] $774=(2+2/2)\times (2^{2\times (2+2)}+2).$
\item [] $775=22^2+(2^{(2+2)}+2/2)^2+2.$
\item [] $776=2\times 2\times ((2^{(2+2)}-2)^2-2).$
\item [] $777=222\times (2+2+2+2/2)/2.$
\item [] $778=2\times (22-2)^2-22.$
\item [] $779=2\times (22-2)^2-22+2/2.$
\item [] $780=2\times (22-2)^2-22+2.$
\item [] $781=(22+2+2+2)^2-2-2/2.$
\item [] $782=(22+2+2+2)^2-2.$
\item [] $783=(22+2+2+2)^2-2/2.$
\item [] $784=(22+2+2+2)^2.$
\item [] $785=(22+2+2+2)^2+2/2.$
\item [] $786=(22+2+2+2)^2+2.$
\item [] $787=(22+2+2+2)^2+2+2/2.$
\item [] $788=(22+2+2+2)^2+2+2.$
\item [] $789=2\times (22-2)^2-22/2.$
\item [] $790=22\times (2+2+2)^2-2.$
\item [] $791=22\times (2+2+2)^2-2/2.$
\item [] $792=22\times (2+2+2)^2.$
\item [] $793=22\times (2+2+2)^2+2/2.$
\item [] $794=22\times (2+2+2)^2+2.$
\item [] $795=2\times ((22-2)^2-2)-2/2.$
\item [] $796=2\times ((22-2)^2-2).$
\item [] $797=2\times ((22-2)^2-2)+2/2.$
\item [] $798=2\times (22-2)^2-2.$
\item [] $799=2\times (22-2)^2-2/2.$
\item [] $800=2\times (22-2)^2.$
\item [] $801=2\times (22-2)^2+2/2.$
\item [] $802=2\times (22-2)^2+2.$
\item [] $803=2\times (22-2)^2+2+2/2.$
\item [] $804=2\times ((22-2)^2+2).$
\item [] $805=2\times ((22-2)^2+2)+2/2.$
\item [] $806=2\times ((22-2)^2+2)+2.$
\item [] $807=2\times ((22-2)^2+2)+2+2/2.$
\item [] $808=2\times ((22-2)^2+2+2).$
\item [] $809=2\times ((22-2)^2+2+2)+2/2.$
\item [] $810=2\times ((22-2)^2+2+2)+2.$
\item [] $811=2\times (22-2)^2+22/2.$
\item [] $812=2\times ((22-2)^2+2+2+2).$
\item [] $813=2\times (22-2)^2+22/2+2.$
\item [] $814=22\times (2+2+2)^2+22.$
\item [] $815=2\times ((22-2)^2+2)+22/2.$
\item [] $816=2\times ((22-2)^2+2\times (2+2)).$
\item [] $817=(22+2)^2+(22^2-2)/2.$
\item [] $818=2\times ((22-2)^2-2)+22.$
\item [] $819=(22+2)^2+(22^2+2)/2.$
\item [] $820=2\times (22-2)^2+22-2.$
\item [] $821=2\times (22-2)^2+22-2/2.$
\item [] $822=2\times (22-2)^2+22.$
\item [] $823=2\times (22-2)^2+22+2/2.$
\item [] $824=2\times (22-2)^2+22+2.$
\item [] $825=(2\times 22)^2-2222/2.$
\item [] $826=2\times ((22-2)^2+2)+22.$
\item [] $827=(2\times 22)^2+2-2222/2.$
\item [] $828=(2^{(2+2)}+2)\times (2\times 22+2).$
\item [] $829=(2\times 22+2)\times (2^{(2+2)}+2)+2/2.$
\item [] $830=(2^{(2+2)}+2)\times (2\times 22+2)+2.$
\item [] $831=(2+2/2)\times ((22+2)^2-22)/2.$
\item [] $832=2\times ((22-2)^2+2^{(2+2)}).$
\item [] $833=2\times (22-2)^2+22+22/2.$
\item [] $834=22\times ((2+2+2)^2+2)-2.$
\item [] $835=22\times ((2+2+2)^2+2)-2/2.$
\item [] $836=22\times ((2+2+2)^2+2).$
\item [] $837=22\times ((2+2+2)^2+2)+2-2/2.$
\item [] $838=22\times ((2+2+2)^2+2)+2.$
\item [] $839=(22-2)\times (2\times 22-2)-2/2.$
\item [] $840=(22-2)\times (2\times 22-2).$
\item [] $841=(22+2+2+2+2/2)^2.$
\item [] $842=(22-2)\times (2\times 22-2)+2.$
\item [] $843=(22+2+2+2+2/2)^2+2.$
\item [] $844=2\times ((22-2)^2+22).$
\item [] $845=2\times ((22-2)^2+22)+2/2.$
\item [] $846=2\times ((22-2)^2+22)+2.$
\item [] $847=2\times 22^2-(22/2)^2.$
\item [] $848=2\times ((22-2)^2+22+2).$
\item [] $849=2+2\times 22^2-(22/2)^2.$
\item [] $850=2\times ((22-2)^2+22+2)+2.$
\item [] $851=2\times (22^2+2)-(22/2)^2.$
\item [] $852=2\times (2\times (222+2)-22).$
\item [] $853=2\times (22^2-2)-222/2.$
\item [] $854=2\times (2\times (222+2)-22)+2.$
\item [] $855=2\times 22^2-2-222/2.$
\item [] $856=2\times 2\times (222-2\times (2+2)).$
\item [] $857=2\times 22^2-222/2.$
\item [] $858=22\times (2\times (22-2)-2/2).$
\item [] $859=2\times 22^2+2-222/2.$
\item [] $860=(22-2)\times (2\times 22-2/2).$
\item [] $861=(2+2/2)\times ((22+2)^2-2)/2.$
\item [] $862=2\times (2\times 222-2)-22.$
\item [] $863=(22+2)\times (2+2+2)^2-2/2.$
\item [] $864=(22+2)\times (2+2+2)^2.$
\item [] $865=(22+2)\times (2+2+2)^2+2/2.$
\item [] $866=2\times 2\times 222-22.$
\item [] $867=2\times 2\times 222-22+2/2.$
\item [] $868=2\times 2\times 222-22+2.$
\item [] $869=(2222-22^2)/2.$
\item [] $870=2\times (2\times 222+2)-22.$
\item [] $871=((2\times 22-2)^2-22)/2.$
\item [] $872=2\times 2\times (222-2-2).$
\item [] $873=((2\times 22-2)^2-22)/2+2.$
\item [] $874=2\times 2\times (222-2-2)+2.$
\item [] $875=2\times 2\times 222-2-22/2.$
\item [] $876=2\times (2\times (222-2)-2).$
\item [] $877=2\times 2\times 222-22/2.$
\item [] $878=2\times 2\times (222-2)-2.$
\item [] $879=2\times 2\times (222-2)-2/2.$
\item [] $880=2\times 2\times (222-2).$
\item [] $881=((2\times 22-2)^2-2)/2.$
\item [] $882=2\times (22-2/2)^2.$
\item [] $883=((2\times 22-2)^2+2)/2.$
\item [] $884=2\times (2\times 222-2).$
\item [] $885=2\times 2\times 222-2-2/2.$
\item [] $886=2\times 2\times 222-2.$
\item [] $887=2\times 2\times 222-2/2.$
\item [] $888=2\times 2\times 222.$
\item [] $889=2\times 2\times 222+2/2.$
\item [] $890=2\times 2\times 222+2.$
\item [] $891=2\times 2\times 222+2+2/2.$
\item [] $892=2\times (2\times 222+2).$
\item [] $893=2\times (2\times 222+2)+2/2.$
\item [] $894=2\times (2\times 222+2)+2.$
\item [] $895=2\times 2\times (222+2)-2/2.$
\item [] $896=2\times 2\times (222+2).$
\item [] $897=2\times 2\times (222+2)+2/2.$
\item [] $898=2\times 2\times (222+2)+2.$
\item [] $899=2\times 2\times 222+22/2.$
\item [] $900=(2\times (2+2)+22)^2.$
\item [] $901=(2\times (2+2)+22)^2+2/2.$
\item [] $902=(2\times (2+2)+22)^2+2.$
\item [] $903=(2\times (2+2)+22)^2+2+2/2.$
\item [] $904=2\times 2\times (222+2+2).$
\item [] $905=((2\times 22-2)^2+2)/2+22.$
\item [] $906=2\times 2\times (222+2+2)+2.$
\item [] $907=2\times 2\times (222+2)+22/2.$
\item [] $908=2\times (2\times (222+2+2)+2).$
\item [] $909=2\times 2\times 222+22-2/2.$
\item [] $910=2\times 2\times 222+22.$
\item [] $911=2\times 2\times 222+22+2/2.$
\item [] $912=2\times 2\times 222+22+2.$
\item [] $913=(2^(22/2)-222)/2.$
\item [] $914=2\times (2\times 222+2)+22.$
\item [] $915=(2^(22/2)-222)/2+2.$
\item [] $916=2\times (22^2-22-2-2).$
\item [] $917=(2^(22/2)-222)/2+2+2.$
\item [] $918=2\times 2\times (222+2)+22.$
\item [] $919=2\times (22^2-22-2)-2/2.$
\item [] $920=2\times (22^2-22-2).$
\item [] $921=2\times (22^2-22-2)+2/2.$
\item [] $922=2\times (22^2-22)-2.$
\item [] $923=2\times (22^2-22)-2/2.$
\item [] $924=2\times (22^2-22).$
\item [] $925=2\times (22^2-22)+2/2.$
\item [] $926=2\times (22^2-22)+2.$
\item [] $927=2\times (22^2-22)+2+2/2.$
\item [] $928=2\times (22^2-22+2).$
\item [] $929=2\times 222+22^2+2/2.$
\item [] $930=2\times 222+22^2+2.$
\item [] $931=2\times (2\times 222+22)-2/2.$
\item [] $932=2\times (2\times 222+22).$
\item [] $933=2\times (2\times 222+22)+2/2.$
\item [] $934=2\times (2\times 222+22)+2.$
\item [] $935=2\times (22^2-22)+22/2.$
\item [] $936=2\times (22^2-2^{(2+2)}).$
\item [] $937=2\times (22^2-2^{(2+2)})+2/2.$
\item [] $938=2\times (22^2-2^{(2+2)})+2.$
\item [] $939=((2+2/2)^2+22)^2-22.$
\item [] $940=2\times (22^2-2^{(2+2)}+2).$
\item [] $941=2\times (22^2-2)-22-2/2.$
\item [] $942=2\times (22^2-2)-22.$
\item [] $943=2\times (22^2-2)-22+2/2.$
\item [] $944=2\times 22^2-22-2.$
\item [] $945=2\times 22^2-22-2/2.$
\item [] $946=2\times 22^2-22.$
\item [] $947=2\times 22^2-22+2/2.$
\item [] $948=2\times 22^2-22+2.$
\item [] $949=2\times 22^2-22+2+2/2.$
\item [] $950=2\times (22^2+2)-22.$
\item [] $951=2\times (22^2+2)-22+2/2.$
\item [] $952=2\times (22^2-2\times (2+2)).$
\item [] $953=2\times (22^2-2)-22/2.$
\item [] $954=2\times (22^2-2\times (2+2))+2.$
\item [] $955=2\times 22^2-2-22/2.$
\item [] $956=2\times (22^2-2-2-2).$
\item [] $957=2\times 22^2-22/2.$
\item [] $958=2\times (22^2-2-2-2)+2.$
\item [] $959=2\times 22^2+2-22/2.$
\item [] $960=2\times (22^2-2-2).$
\item [] $961=(22+(2+2/2)^2)^2.$
\item [] $962=2\times (22^2-2)-2.$
\item [] $963=2\times (22^2-2)-2/2.$
\item [] $964=2\times (22^2-2).$
\item [] $965=2\times 22^2-2-2/2.$
\item [] $966=2\times 22^2-2.$
\item [] $967=2\times 22^2-2/2.$
\item [] $968=2\times 22^2.$
\item [] $969=2\times 22^2+2/2.$
\item [] $970=2\times 22^2+2.$
\item [] $971=2\times 22^2+2+2/2.$
\item [] $972=2\times (22^2+2).$
\item [] $973=2\times (22^2+2)+2/2.$
\item [] $974=2\times (22^2+2)+2.$
\item [] $975=2\times (22^2+2)+2+2/2.$
\item [] $976=2\times (22^2+2+2).$
\item [] $977=2\times (22^2+2+2)+2/2.$
\item [] $978=2\times (22^2+2+2)+2.$
\item [] $979=2\times 22^2+22/2.$
\item [] $980=2\times (22^2+2+2+2).$
\item [] $981=2\times 22^2+2+22/2.$
\item [] $982=2\times (22^2+2+2+2)+2.$
\item [] $983=2\times (22^2+2)+22/2.$
\item [] $984=2\times (2\times (2+2)+22^2).$
\item [] $985=2\times (22^2+2)+2+22/2.$
\item [] $986=2\times (22^2-2)+22.$
\item [] $987=2\times (22^2-2)+22+2/2.$
\item [] $988=2\times 22^2+22-2.$
\item [] $989=2\times 22^2+22-2/2.$
\item [] $990=2\times 22^2+22.$
\item [] $991=2\times 22^2+22+2/2.$
\item [] $992=2\times 22^2+22+2.$
\item [] $993=2\times 22^2+22+2+2/2.$
\item [] $994=22+2\times (22^2+2).$
\item [] $995=22+2\times (22^2+2)+2/2.$
\item [] $996=22+2\times (22^2+2)+2.$
\item [] $997=(2+2/2)^2\times 222/2-2.$
\item [] $998=2\times (22^2+2+2)+22.$
\item [] $999=(2+2/2)^2\times 222/2.$
\item [] $1000=2\times (22^2+2^{(2+2)}).$
\end{itemize}
\end{multicols}
}

\section{\textbf{Representations Using Number 3}}

{\footnotesize
\begin{multicols}{3}
\begin{itemize}
\item [] $101=3+3\times 33-3/3.$
\item [] $102=3+3\times 33.$
\item [] $103=3+3\times 33+3/3.$
\item [] $104=3+3+3\times 33-3/3.$
\item [] $105=3+(3\times 33+3).$
\item [] $106=3+3+3\times 33+3/3.$
\item [] $107=3\times (33+3)-3/3.$
\item [] $108=3\times (33+3).$
\item [] $109=3\times (33+3)+3/3.$
\item [] $110=(333-3)/3.$
\item [] $111=333/3.$
\item [] $112=(333+3)/3.$
\item [] $113=3+(333-3)/3.$
\item [] $114=3+333/3.$
\item [] $115=3+(333+3)/3.$
\item [] $116=3+3+(333-3)/3.$
\item [] $117=3\times (33+3+3).$
\item [] $118=3+3+(333+3)/3.$
\item [] $119=3\times 3+(333-3)/3.$
\item [] $120=3+3\times (33+3+3).$
\item [] $121=(33/3)^{(3-3/3)}.$
\item [] $122=(3^{(3+3)}+3)/(3+3).$
\item [] $123=3^3+3\times 33-3.$
\item [] $124=3+(33/3)^{(3-3/3)}.$
\item [] $125=(3+3-3/3)^3.$
\item [] $126=3\times (3\times 3+33).$
\item [] $127=3^3+3\times 33+3/3.$
\item [] $128=3+(3+3-3/3)^3.$
\item [] $129=3+3\times 33+3^3.$
\item [] $130=3+3\times 33+3^3+3/3.$
\item [] $131=3+3+(3+3-3/3)^3.$
\item [] $132=33+3\times 33.$
\item [] $133=3\times 33+33+3/3.$
\item [] $134=3\times 3+(3+3-3/3)^3.$
\item [] $135=3+3\times 33+33.$
\item [] $136=(3+3/3)\times (33+3/3).$
\item [] $137=3^3+(333-3)/3.$
\item [] $138=3^3+333/3.$
\item [] $139=3^3+(333+3)/3.$
\item [] $140=3+3^3+(333-3)/3.$
\item [] $141=33+3\times (33+3).$
\item [] $142=3+3^3+(333+3)/3.$
\item [] $143=33+(333-3)/3.$
\item [] $144=(3+3)\times (3^3-3).$
\item [] $145=(3+3)\times (3^3-3)+3/3.$
\item [] $146=3+33+(333-3)/3.$
\item [] $147=3+(3+3)\times (3^3-3).$
\item [] $148=(33\times 3^3-3)/(3+3).$
\item [] $149=(33\times 3^3+3)/(3+3).$
\item [] $150=3+3+(3+3)\times (3^3-3).$
\item [] $151=(3+3)\times 3^3-33/3.$
\item [] $152=3^3+(3+3-3/3)^3.$
\item [] $153=3\times (3^3+3^3-3).$
\item [] $154=33\times (3+33/3)/3.$
\item [] $155=3^3+3+(3+3-3/3)^3.$
\item [] $156=(3+3)\times (3^3-3/3).$
\item [] $157=3+33\times (3+33/3)/3.$
\item [] $158=33+(3+3-3/3)^3.$
\item [] $159=(3+3)\times 3^3-3.$
\item [] $160=(3+3)\times 3^3-3+3/3.$
\item [] $161=(3+3)\times 3^3-3/3.$
\item [] $162=(3+3)\times 3^3.$
\item [] $163=(3+3)\times 3^3+3/3.$
\item [] $164=(3+3)\times 3^3+3-3/3.$
\item [] $165=(3+3)\times 3^3+3.$
\item [] $166=(3+3)\times 3^3+3+3/3.$
\item [] $167=(3\times 333+3)/(3+3).$
\item [] $168=(3+3)\times 3^3+3+3.$
\item [] $169=(3+3)\times 3^3+3+3+3/3.$
\item [] $170=3+(3\times 333+3)/(3+3).$
\item [] $171=3\times (3^3+3^3+3).$
\item [] $172=3\times 3+(3+3)\times 3^3+3/3.$
\item [] $173=(3+3)\times 3^3+33/3.$
\item [] $174=(3+3)\times 3^3+3\times 3+3.$
\item [] $175=(3+3/3)^3+333/3.$
\item [] $176=(3+3)\times 3^3+3+33/3.$
\item [] $177=3\times (3^3+33)-3.$
\item [] $178=3\times (3^3+33)-3+3/3.$
\item [] $179=3\times (3^3+33)-3/3.$
\item [] $180=3\times (3^3+33).$
\item [] $181=3\times (3^3+33)+3/3.$
\item [] $182=(3+3)^3-33-3/3.$
\item [] $183=(3+3)^3-33.$
\item [] $184=(3+3)^3-33+3/3.$
\item [] $185=(3+3)^3+3-33-3/3.$
\item [] $186=(3+3)^3+3-33.$
\item [] $187=33\times (3+3)-33/3.$
\item [] $188=(3+3)^3-3^3-3/3.$
\item [] $189=(3+3)^3-3^3.$
\item [] $190=(3+3)^3-3^3+3/3.$
\item [] $191=3\times (3+3/3)^3-3/3.$
\item [] $192=3\times (3+3/3)^3.$
\item [] $193=3\times (3+3/3)^3+3/3.$
\item [] $194=33\times (3+3)-3-3/3.$
\item [] $195=33\times (3+3)-3.$
\item [] $196=33\times (3+3)-3+3/3.$
\item [] $197=33\times (3+3)-3/3.$
\item [] $198=33\times (3+3).$
\item [] $199=33\times (3+3)+3/3.$
\item [] $200=33\times (3+3)+3-3/3.$
\item [] $201=33\times (3+3)+3.$
\item [] $202=33\times (3+3)+3+3/3.$
\item [] $203=33\times (3+3)+3+3-3/3.$
\item [] $204=33\times (3+3)+3+3.$
\item [] $205=(3+3)^3-33/3.$
\item [] $206=(3+3)^3-(33-3)/3.$
\item [] $207=(3+3)^3-3\times 3.$
\item [] $208=(3+3)^3-3\times 3+3/3.$
\item [] $209=(3+3)^3-3-3-3/3.$
\item [] $210=(3+3)^3-3-3.$
\item [] $211=(3+3)^3-3-3+3/3.$
\item [] $212=(3+3)^3-3-3/3.$
\item [] $213=(3+3)^3-3.$
\item [] $214=(3+3)^3-3+3/3.$
\item [] $215=(3+3)^3-3/3.$
\item [] $216=(3+3)^3.$
\item [] $217=(3+3)^3+3/3.$
\item [] $218=(3+3)^3+3-3/3.$
\item [] $219=(3+3)^3+3.$
\item [] $220=(3+3)^3+3+3/3.$
\item [] $221=(3+3)^3+3+3-3/3.$
\item [] $222=(3+3)^3+3+3.$
\item [] $223=(3+3)^3+3+3+3/3.$
\item [] $224=(3+3)^3+3\times 3-3/3.$
\item [] $225=(3+3)^3+3\times 3.$
\item [] $226=(3+3)^3+3\times 3+3/3.$
\item [] $227=(3+3)^3+33/3.$
\item [] $228=(3+3)^3+3\times 3+3.$
\item [] $229=(3^{(3+3)}-33)/3-3.$
\item [] $230=(3+3)^3+3+33/3.$
\item [] $231=33\times (3+3)+33.$
\item [] $232=(3^{(3+3)}-33)/3.$
\item [] $233=3\times (3\times 3^3-3)-3/3.$
\item [] $234=3\times (3\times 3^3-3).$
\item [] $235=3+(3^{(3+3)}-33)/3.$
\item [] $236=(3^{(3+3)}-3)/3-3-3.$
\item [] $237=3\times (3\times 3^3-3)+3.$
\item [] $238=(3^{(3+3)}+3)/3-3-3.$
\item [] $239=(3^{(3+3)}-3)/3-3.$
\item [] $240=3\times 3\times 3^3-3.$
\item [] $241=(3^{(3+3)}+3)/3-3.$
\item [] $242=(3^{(3+3)}-3)/3.$
\item [] $243=3\times 3\times 3^3.$
\item [] $244=(3^{(3+3)}+3)/3.$
\item [] $245=3+(3^{(3+3)}-3)/3.$
\item [] $246=3\times 3\times 3^3+3.$
\item [] $247=3+(3^{(3+3)}+3)/3.$
\item [] $248=3+3+(3^{(3+3)}-3)/3.$
\item [] $249=33+(3+3)^3.$
\item [] $250=3+3+(3^{(3+3)}+3)/3.$
\item [] $251=3\times 3+(3^{(3+3)}-3)/3.$
\item [] $252=3\times (3\times 3^3+3).$
\item [] $253=3\times 3+(3^{(3+3)}+3)/3.$
\item [] $254=(3^{(3+3)}+33)/3.$
\item [] $255=3\times (3\times 3^3+3)+3.$
\item [] $256=(3+3/3)^{(3+3/3)}.$
\item [] $257=3+(3^{(3+3)}+33)/3.$
\item [] $258=3+3+3\times (3\times 3^3+3).$
\item [] $259=3+(3+3/3)^{(3+3/3)}.$
\item [] $260=3+3+(3^{(3+3)}+33)/3.$
\item [] $261=3\times (3\times 3^3+3+3).$
\item [] $262=3+3+(3+3/3)^{(3+3/3)}.$
\item [] $263=(33\times (3^3-3)-3)/3.$
\item [] $264=33\times (3\times 3-3/3).$
\item [] $265=(33\times (3^3-3)+3)/3.$
\item [] $266=3+(33\times (3^3-3)-3)/3.$
\item [] $267=3\times 3\times (3^3+3)-3.$
\item [] $268=3+(33\times (3^3-3)+3)/3.$
\item [] $269=3^3+(3^{(3+3)}-3)/3.$
\item [] $270=3\times 3\times (3^3+3).$
\item [] $271=3^3+(3^{(3+3)}+3)/3.$
\item [] $272=3^3+3+(3^{(3+3)}-3)/3.$
\item [] $273=3\times 3\times (3^3+3)+3.$
\item [] $274=3\times 3\times (3^3+3)+3+3/3.$
\item [] $275=33+(3^{(3+3)}-3)/3.$
\item [] $276=33+3\times 3\times 3^3.$
\item [] $277=33+(3^{(3+3)}+3)/3.$
\item [] $278=33+3+(3^{(3+3)}-3)/3.$
\item [] $279=3\times (3\times (3^3+3)+3).$
\item [] $280=(3+3)^3+(3+3/3)^3.$
\item [] $281=3^3+(3^{(3+3)}+33)/3.$
\item [] $282=3+3\times (3\times (3^3+3)+3).$
\item [] $283=(3+3/3)^3+(3+3)^3+3.$
\item [] $284=(3+3/3)\times ((3+3)^3-3)/3.$
\item [] $285=3\times (3\times 33-3)-3.$
\item [] $286=33\times (3^3-3/3)/3.$
\item [] $287=3\times (3\times 33-3)-3/3.$
\item [] $288=3\times (3\times 33-3).$
\item [] $289=3\times (3\times 33-3)+3/3.$
\item [] $290=3\times (3\times 33-3)+3-3/3.$
\item [] $291=3\times (3\times 33-3)+3.$
\item [] $292=3\times (3\times 33-3)+3+3/3.$
\item [] $293=3\times 3\times 33-3-3/3.$
\item [] $294=3\times 3\times 33-3.$
\item [] $295=3\times 3\times 33-3+3/3.$
\item [] $296=3\times 3\times 33-3/3.$
\item [] $297=3\times 3\times 33.$
\item [] $298=3\times 3\times 33+3/3.$
\item [] $299=3\times 3\times 33+3-3/3.$
\item [] $300=3\times 3\times 33+3.$
\item [] $301=3\times 3\times 33+3+3/3.$
\item [] $302=3\times 3\times 33+3+3-3/3.$
\item [] $303=3\times 3\times 33+3+3.$
\item [] $304=3\times 3\times 33+3+3+3/3.$
\item [] $305=333-3^3-3/3.$
\item [] $306=3\times (3\times 33+3).$
\item [] $307=3\times (3\times 33+3)+3/3.$
\item [] $308=3\times 3\times 33+33/3.$
\item [] $309=3\times (3\times 33+3)+3.$
\item [] $310=(3+3+3/3)^3-33.$
\item [] $311=3\times 3\times 33+3+33/3.$
\item [] $312=3\times (3\times 33+3)+3+3.$
\item [] $313=3-33+(3+3+3/3)^3.$
\item [] $314=3\times 33+(3+3)^3-3/3.$
\item [] $315=3\times (3\times 33+3+3).$
\item [] $316=(3+3+3/3)^3-3^3.$
\item [] $317=333-3^3+33/3.$
\item [] $318=3\times 33+(3+3)^3+3.$
\item [] $319=333-3-33/3.$
\item [] $320=333-3-(33-3)/3.$
\item [] $321=333-3\times 3-3.$
\item [] $322=333-33/3.$
\item [] $323=333-(33-3)/3.$
\item [] $324=3\times 3\times (33+3).$
\item [] $325=3+333-33/3.$
\item [] $326=333-3-3-3/3.$
\item [] $327=333-3-3.$
\item [] $328=333-3-3+3/3.$
\item [] $329=333-3-3/3.$
\item [] $330=333-3.$
\item [] $331=333-3+3/3.$
\item [] $332=333-3/3.$
\item [] $333=333.$
\item [] $334=333+3/3.$
\item [] $335=333+3-3/3.$
\item [] $336=333+3.$
\item [] $337=333+3+3/3.$
\item [] $338=333+3+3-3/3.$
\item [] $339=333+3+3.$
\item [] $340=(3+3+3/3)^3-3.$
\item [] $341=333+3\times 3-3/3.$
\item [] $342=333+3\times 3.$
\item [] $343=(3+3+3/3)^3.$
\item [] $344=333+33/3.$
\item [] $345=333+3\times 3+3.$
\item [] $346=3+(3+3+3/3)^3.$
\item [] $347=3+333+33/3.$
\item [] $348=3+3+333+3\times 3.$
\item [] $349=3+3+(3+3+3/3)^3.$
\item [] $350=3+3+333+33/3.$
\item [] $351=3\times 3\times (33+3+3).$
\item [] $352=3\times 3+(3+3+3/3)^3.$
\item [] $353=3\times 3+333+33/3.$
\item [] $354=3+333+3\times (3+3).$
\item [] $355=3+3\times 3+(3+3+3/3)^3.$
\item [] $356=3^3+333-3-3/3.$
\item [] $357=3^3+333-3.$
\item [] $358=3^3+333-3+3/3.$
\item [] $359=3^3+333-3/3.$
\item [] $360=3^3+333.$
\item [] $361=3^3+333+3/3.$
\item [] $362=(33\times 33-3)/3.$
\item [] $363=33\times 33/3.$
\item [] $364=(33\times 33+3)/3.$
\item [] $365=3+(33\times 33-3)/3.$
\item [] $366=33+333.$
\item [] $367=3+(33\times 33+3)/3.$
\item [] $368=3+3+(33\times 33-3)/3.$
\item [] $369=3+333+33.$
\item [] $370=3^3+(3+3+3/3)^3.$
\item [] $371=3^3+333+33/3.$
\item [] $372=3+3+333+33.$
\item [] $373=3+3^3+(3+3+3/3)^3.$
\item [] $374=33\times (33+3/3)/3.$
\item [] $375=3\times (3+3-3/3)^3.$
\item [] $376=33+(3+3+3/3)^3.$
\item [] $377=3+33\times (33+3/3)/3.$
\item [] $378=3\times (3\times 33+3^3).$
\item [] $379=3+33+(3+3+3/3)^3.$
\item [] $380=3+3+33\times (33+3/3)/3.$
\item [] $381=3+3\times (3\times 33+3^3).$
\item [] $382=3333/3-3^{(3+3)}.$
\item [] $383=(3+3)\times (3+3/3)^3-3/3.$
\item [] $384=(3+3)\times (3+3/3)^3.$
\item [] $385=(3+3)\times (3+3/3)^3+3/3.$
\item [] $386=3^{(3+3)}-(3+3+3/3)^3.$
\item [] $387=3+(3+3)\times (3+3/3)^3.$
\item [] $388=3+(3+3)\times (3+3/3)^3+3/3.$
\item [] $389=3^3+(33\times 33-3)/3.$
\item [] $390=3^3+33\times 33/3.$
\item [] $391=3^3+(33\times 33+3)/3.$
\item [] $392=33\times (3\times 3+3)-3-3/3.$
\item [] $393=33\times (3\times 3+3)-3.$
\item [] $394=33\times (3\times 3+3)-3+3/3.$
\item [] $395=33\times (3\times 3+3)-3/3.$
\item [] $396=33\times (3\times 3+3).$
\item [] $397=33\times (3\times 3+3)+3/3.$
\item [] $398=33\times (3\times 3+3)+3-3/3.$
\item [] $399=33\times (3\times 3+3)+3.$
\item [] $400=33\times (3\times 3+3)+3+3/3.$
\item [] $401=33\times (3\times 3+3)+3+3-3/3.$
\item [] $402=33\times (3\times 3+3)+3+3.$
\item [] $403=33\times (3\times 3+3)+3+3+3/3.$
\item [] $404=333+((3+3)^3-3)/3.$
\item [] $405=3^3\times (3\times 3+3+3).$
\item [] $406=3^3\times (3\times 3+3+3)+3/3.$
\item [] $407=33\times 333/(3^3).$
\item [] $408=3^3\times (3\times 3+3+3)+3.$
\item [] $409=3^3\times (3\times 3+3+3)+3+3/3.$
\item [] $410=33\times 333/(3^3)+3.$
\item [] $411=3\times 3^3+333-3.$
\item [] $412=3\times 3^3+333-3+3/3.$
\item [] $413=3\times 3^3+333-3/3.$
\item [] $414=3\times 3^3+333.$
\item [] $415=3\times 3^3+333+3/3.$
\item [] $416=3\times 3^3+333+3-3/3.$
\item [] $417=3\times 3^3+333+3.$
\item [] $418=33\times (3^3+33/3)/3.$
\item [] $419=(3^3+3)\times (3+33/3)-3/3.$
\item [] $420=(3^3+3)\times (3+33/3).$
\item [] $421=((3+3)^{(3+3/3)}-33)/3.$
\item [] $422=3^3+33\times (3\times 3+3)-3/3.$
\item [] $423=3\times (3\times (33+3)+33).$
\item [] $424=3\times 3^3+(3+3+3/3)^3.$
\item [] $425=((3+3)\times ((3+3)^3-3)-3)/3.$
\item [] $426=(3+3)\times ((3+3)^3-3)/3.$
\item [] $427=((3+3)\times ((3+3)^3-3)+3)/3.$
\item [] $428=((3+3)^{(3+3/3)}-3)/3-3.$
\item [] $429=33+33\times (3\times 3+3).$
\item [] $430=(3-3/3)\times ((3+3)^3-3/3).$
\item [] $431=((3+3)^{(3+3/3)}-3)/3.$
\item [] $432=3\times (3+3)\times (3^3-3).$
\item [] $433=((3+3)^{(3+3/3)}+3)/3.$
\item [] $434=3+((3+3)^{(3+3/3)}-3)/3.$
\item [] $435=3+3\times 33+333.$
\item [] $436=3+((3+3)^{(3+3/3)}+3)/3.$
\item [] $437=((3+3)\times ((3+3)^3+3)-3)/3.$
\item [] $438=(3+3)\times ((3+3)^3+3)/3.$
\item [] $439=((3+3)\times ((3+3)^3+3)+3)/3.$
\item [] $440=(3+3/3)\times (333-3)/3.$
\item [] $441=3\times ((3+3)\times (3^3-3)+3).$
\item [] $442=3\times 33+(3+3+3/3)^3.$
\item [] $443=333+(333-3)/3.$
\item [] $444=333+333/3.$
\item [] $445=333+(333+3)/3.$
\item [] $446=(3\times 33\times 3^3+3)/(3+3).$
\item [] $447=3+333+333/3.$
\item [] $448=(3+3/3)\times (333+3)/3.$
\item [] $449=(3\times 33\times 3^3+3)/(3+3)+3.$
\item [] $450=(3+3)\times (3\times (3^3-3)+3).$
\item [] $451=(3+3/3)\times (333+3)/3+3.$
\item [] $452=(3+3/3)\times (3+(333-3)/3).$
\item [] $453=3\times (3+3)\times 3^3-33.$
\item [] $454=3\times (3+3)\times 3^3-33+3/3.$
\item [] $455=(3+3)^3+(3^{(3+3)}-3)/3-3.$
\item [] $456=3\times (3+3)\times 3^3-33+3.$
\item [] $457=(3+3)^3+(3^{(3+3)}+3)/3-3.$
\item [] $458=(3+3)^3+(3^{(3+3)}-3)/3.$
\item [] $459=3\times 3\times (3^3-3+3^3).$
\item [] $460=(3+3)^3+(3^{(3+3)}+3)/3.$
\item [] $461=33\times (3+33/3)-3/3.$
\item [] $462=33\times (3+33/3).$
\item [] $463=33\times (3+33/3)+3/3.$
\item [] $464=3+33\times (3+33/3)-3/3.$
\item [] $465=3+33\times (3+33/3).$
\item [] $466=3+33\times (3+33/3)+3/3.$
\item [] $467=(3+3)\times (3\times 3^3-3)-3/3.$
\item [] $468=(3+3)\times (3\times 3^3-3).$
\item [] $469=(3+3)\times (3\times 3^3-3)+3/3.$
\item [] $470=(3+3)\times (3\times 3^3-3)+3-3/3.$
\item [] $471=(3+3)\times (3\times 3^3-3)+3.$
\item [] $472=(3+3)\times (3\times 3^3-3)+3+3/3.$
\item [] $473=33\times (3+33/3)+33/3.$
\item [] $474=3\times ((3+3)\times 3^3-3)-3.$
\item [] $475=3\times (3+3)\times 3^3-33/3.$
\item [] $476=(3+33/3)\times (33+3/3).$
\item [] $477=3\times ((3+3)\times 3^3-3).$
\item [] $478=3\times ((3+3)\times 3^3-3)+3/3.$
\item [] $479=(3-3/3)^{(3\times 3)}-33.$
\item [] $480=3\times ((3+3)\times 3^3-3)+3.$
\item [] $481=3\times ((3+3)\times 3^3-3)+3+3/3.$
\item [] $482=(3-3/3)^{(3\times 3)}+3-33.$
\item [] $483=3\times (3+3)\times 3^3-3.$
\item [] $484=3\times (3+3)\times 3^3-3+3/3.$
\item [] $485=(3-3/3)^{(3\times 3)}-3^3.$
\item [] $486=3\times (3+3)\times 3^3.$
\item [] $487=3\times (3+3)\times 3^3+3/3.$
\item [] $488=(3-3/3)^{(3\times 3)}+3-3^3.$
\item [] $489=3\times (3+3)\times 3^3+3.$
\item [] $490=3\times (3+3)\times 3^3+3+3/3.$
\item [] $491=(3-3/3)^{(3\times 3)}+3+3-3^3.$
\item [] $492=3\times (3+3)\times 3^3+3+3.$
\item [] $493=3\times (3+3)\times 3^3+3+3+3/3.$
\item [] $494=(3-3/3)^{(3\times 3)}-3\times (3+3).$
\item [] $495=3\times ((3+3)\times 3^3+3).$
\item [] $496=3\times ((3+3)\times 3^3+3)+3/3.$
\item [] $497=3\times (3+3)\times 3^3+33/3.$
\item [] $498=3\times ((3+3)\times 3^3+3)+3.$
\item [] $499=(3\times 3\times 333-3)/(3+3).$
\item [] $500=(3+3/3)\times (3+3-3/3)^3.$
\item [] $501=(3+3)\times (3\times 3^3+3)-3.$
\item [] $502=(3\times 3\times 333-3)/(3+3)+3.$
\item [] $503=(3-3/3)^{(3\times 3)}-3\times 3.$
\item [] $504=(3+3)\times (3\times 3^3+3).$
\item [] $505=(3+3)\times (3\times 3^3+3)+3/3.$
\item [] $506=(3-3/3)^{(3\times 3)}-3-3.$
\item [] $507=(3+3)\times (3\times 3^3+3)+3.$
\item [] $508=(3-3/3)^{(3\times 3)}-3-3/3.$
\item [] $509=(3-3/3)^{(3\times 3)}-3.$
\item [] $510=3^{(3+3)}-(3+3)^3-3.$
\item [] $511=(3-3/3)^{(3\times 3)}-3/3.$
\item [] $512=(3-3/3)^{(3\times 3)}.$
\item [] $513=(3-3/3)^{(3\times 3)}+3/3.$
\item [] $514=(3-3/3)^{(3\times 3)}+3-3/3.$
\item [] $515=(3-3/3)^{(3\times 3)}+3.$
\item [] $516=3^{(3+3)}-(3+3)^3+3.$
\item [] $517=(3-3/3)^{(3\times 3)}+3+3-3/3.$
\item [] $518=(3-3/3)^{(3\times 3)}+3+3.$
\item [] $519=3\times (3+3)\times 3^3+33.$
\item [] $520=3\times 3+(3-3/3)^{(3\times 3)}-3/3.$
\item [] $521=3\times 3+(3-3/3)^{(3\times 3)}.$
\item [] $522=(3+3)\times (3\times 3^3+3+3).$
\item [] $523=(3-3/3)^{(3\times 3)}+33/3.$
\item [] $524=3\times 3+(3-3/3)^{(3\times 3)}+3.$
\item [] $525=(3+3)\times (3\times 3^3+3+3)+3.$
\item [] $526=(3-3/3)^{(3\times 3)}+3+33/3.$
\item [] $527=3\times 3+(3-3/3)^{(3\times 3)}+3+3.$
\item [] $528=33\times (3^3-33/3).$
\item [] $529=(3^3-3-3/3)^{(3-3/3)}.$
\item [] $530=3\times (3+3)+(3-3/3)^{(3\times 3)}.$
\item [] $531=3\times (3\times (3^3+33)-3).$
\item [] $532=(3^3-3-3/3)^{(3-3/3)}+3.$
\item [] $533=3\times (3+3)+(3-3/3)^{(3\times 3)}+3.$
\item [] $534=33\times (3+3)+333+3.$
\item [] $535=(3333-(3\times 3+3)^3)/3.$
\item [] $536=3^3+(3-3/3)^{(3\times 3)}-3.$
\item [] $537=3\times 3\times (3^3+33)-3.$
\item [] $538=3^3+(3-3/3)^{(3\times 3)}-3/3.$
\item [] $539=3^3+(3-3/3)^{(3\times 3)}.$
\item [] $540=3\times 3\times (3^3+33).$
\item [] $541=3\times 3\times (3^3+33)+3/3.$
\item [] $542=3^3+(3-3/3)^{(3\times 3)}+3.$
\item [] $543=3\times 3\times (3^3+33)+3.$
\item [] $544=((3\times 3+3)^3+3)/3-33.$
\item [] $545=(3-3/3)^{(3\times 3)}+33.$
\item [] $546=(3+3)^3+333-3.$
\item [] $547=((3\times 3+3)^3+3)/3-33+3.$
\item [] $548=(3-3/3)^{(3\times 3)}+33+3.$
\item [] $549=3\times ((3+3)^3-33).$
\item [] $550=3\times ((3+3)^3-33)+3/3.$
\item [] $551=(3333-3^3)/(3+3).$
\item [] $552=(3+3)^3+333+3.$
\item [] $553=(3333+3)/(3+3)-3.$
\item [] $554=(3333-3\times 3)/(3+3).$
\item [] $555=(3333-3)/(3+3).$
\item [] $556=(3333+3)/(3+3).$
\item [] $557=(3333+3\times 3)/(3+3).$
\item [] $558=3\times ((3+3)^3+3-33).$
\item [] $559=(3333+3)/(3+3)+3.$
\item [] $560=(3333+3^3)/(3+3).$
\item [] $561=33\times (3+3+33/3).$
\item [] $562=((3\times 3+3+3)^3-3)/(3+3).$
\item [] $563=((3\times 3+3+3)^3+3)/(3+3).$
\item [] $564=3\times ((3+3)^3-3^3)-3.$
\item [] $565=((3\times 3+3)^3-33)/3.$
\item [] $566=3\times ((3+3)^3-3^3)-3/3.$
\item [] $567=3\times ((3+3)^3-3^3).$
\item [] $568=((3\times 3+3)^3-33)/3+3.$
\item [] $569=((3\times 3+3)^3-3)/3-3-3.$
\item [] $570=3\times ((3+3)^3-3^3)+3.$
\item [] $571=((3\times 3+3)^3+3)/3-3-3.$
\item [] $572=((3\times 3+3)^3-3)/3-3.$
\item [] $573=((3\times 3+3)^3)/3-3.$
\item [] $574=((3\times 3+3)^3+3)/3-3.$
\item [] $575=((3\times 3+3)^3-3)/3.$
\item [] $576=(3\times 3+3)^3/3.$
\item [] $577=((3\times 3+3)^3+3)/3.$
\item [] $578=((3\times 3+3)^3-3)/3+3.$
\item [] $579=(3\times 3+3)^3/3+3.$
\item [] $580=((3\times 3+3)^3+3)/3+3.$
\item [] $581=((3\times 3+3)^3-3)/3+3+3.$
\item [] $582=(3\times 3+3)^3/3+3+3.$
\item [] $583=((3\times 3+3)^3+3)/3+3+3.$
\item [] $584=3\times 3+((3\times 3+3)^3-3)/3.$
\item [] $585=3\times (33\times (3+3)-3).$
\item [] $586=3\times 3+((3\times 3+3)^3+3)/3.$
\item [] $587=((3\times 3+3)^3+33)/3.$
\item [] $588=3\times (33\times (3+3)-3)+3.$
\item [] $589=3\times (33\times (3+3)-3)+3+3/3.$
\item [] $590=((3\times 3+3)^3+33)/3+3.$
\item [] $591=3\times 33\times (3+3)-3.$
\item [] $592=3\times 33\times (3+3)-3+3/3.$
\item [] $593=3\times 33\times (3+3)-3/3.$
\item [] $594=3\times 33\times (3+3).$
\item [] $595=3\times 33\times (3+3)+3/3.$
\item [] $596=3\times 33\times (3+3)+3-3/3.$
\item [] $597=3\times 33\times (3+3)+3.$
\item [] $598=3\times 33\times (3+3)+3+3/3.$
\item [] $599=(33/3)^3-3^{(3+3)}-3.$
\item [] $600=3\times 33\times (3+3)+3+3.$
\item [] $601=3\times 33\times (3+3)+3+3+3/3.$
\item [] $602=(33/3)^3-3^{(3+3)}.$
\item [] $603=3\times (33\times (3+3)+3).$
\item [] $604=3^3+((3\times 3+3)^3+3)/3.$
\item [] $605=3+(33/3)^3-3^{(3+3)}.$
\item [] $606=3\times (33\times (3+3)+3)+3.$
\item [] $607=3\times (33\times (3+3)+3)+3+3/3.$
\item [] $608=33+((3\times 3+3)^3-3)/3.$
\item [] $609=33+(3\times 3+3)^3/3.$
\item [] $610=33+((3\times 3+3)^3+3)/3.$
\item [] $611=3\times 33+(3-3/3)^{(3\times 3)}.$
\item [] $612=(3+3)\times (3\times 33+3).$
\item [] $613=(3+3)\times (3\times 33+3)+3/3.$
\item [] $614=3\times (3+3)^3-33-3/3.$
\item [] $615=3\times (3+3)^3-33.$
\item [] $616=3\times (3+3)^3-33+3/3.$
\item [] $617=3^{(3+3)}-(333+3)/3.$
\item [] $618=+3\times (3+3)^3+3-33.$
\item [] $619=3^{(3+3)}-(333-3)/3.$
\item [] $620=3\times (3+3)^3-3^3-3/3.$
\item [] $621=3\times ((3+3)^3-3\times 3).$
\item [] $622=3\times ((3+3)^3-3\times 3)+3/3.$
\item [] $623=3\times (3+3)^3-3^3+3-3/3.$
\item [] $624=3\times ((3+3)^3-3\times 3)+3.$
\item [] $625=(3+3-3/3)^{(3+3/3)}.$
\item [] $626=(3+3-3/3)^{(3+3/3)}+3/3.$
\item [] $627=3\times 33\times (3+3)+33.$
\item [] $628=(3+3-3/3)^{(3+3/3)}+3.$
\item [] $629=3^{(3+3)}-3\times 33-3/3.$
\item [] $630=3^{(3+3)}-3\times 33.$
\item [] $631=3^{(3+3)}-3\times 33+3/3.$
\item [] $632=3^{(3+3)}-3\times 33+3-3/3.$
\item [] $633=3^{(3+3)}-3\times 33+3.$
\item [] $634=3\times (3+3)^3-3-33/3.$
\item [] $635=3\times ((3+3)^3-3)-3-3/3.$
\item [] $636=3\times ((3+3)^3-3)-3.$
\item [] $637=3\times (3+3)^3-33/3.$
\item [] $638=3\times ((3+3)^3-3)-3/3.$
\item [] $639=3\times ((3+3)^3-3).$
\item [] $640=3\times ((3+3)^3-3)+3/3.$
\item [] $641=3\times ((3+3)^3-3)+3-3/3.$
\item [] $642=3\times ((3+3)^3-3)+3.$
\item [] $643=3\times ((3+3)^3-3)+3+3/3.$
\item [] $644=3\times (3+3)^3-3-3/3.$
\item [] $645=3\times (3+3)^3-3.$
\item [] $646=3\times (3+3)^3-3+3/3.$
\item [] $647=3\times (3+3)^3-3/3.$
\item [] $648=3\times (3+3)^3.$
\item [] $649=3\times (3+3)^3+3/3.$
\item [] $650=3\times (3+3)^3+3-3/3.$
\item [] $651=3\times (3+3)^3+3.$
\item [] $652=3\times (3+3)^3+3+3/3.$
\item [] $653=3\times (3+3)^3+3+3-3/3.$
\item [] $654=3\times (3+3)^3+3+3.$
\item [] $655=3\times ((3+3)^3+3)-3+3/3.$
\item [] $656=3\times ((3+3)^3+3)-3/3.$
\item [] $657=3\times ((3+3)^3+3).$
\item [] $658=3\times ((3+3)^3+3)+3/3.$
\item [] $659=3\times (3+3)^3+33/3.$
\item [] $660=3\times ((3+3)^3+3)+3.$
\item [] $661=3\times ((3+3)^3+3)+3+3/3.$
\item [] $662=3\times (3+3)^3+3+33/3.$
\item [] $663=3\times ((3+3)^3+3)+3+3.$
\item [] $664=(3-3/3)\times (333-3/3).$
\item [] $665=3^{(3+3)}-(3+3/3)^3.$
\item [] $666=3\times ((3+3)^3+3+3).$
\item [] $667=((3+3)\times 333+3)/3.$
\item [] $668=3^{(3+3)}-(3+3/3)^3+3.$
\item [] $669=333+333+3.$
\item [] $670=3+((3+3)\times 333+3)/3.$
\item [] $671=33\times ((3/3+3)^3-3)/3.$
\item [] $672=3^3+3\times (3+3)^3-3.$
\item [] $673=(3^3-3/3)^{(3-3/3)}-3.$
\item [] $674=3^3+3\times (3+3)^3-3/3.$
\item [] $675=3^3+3\times (3+3)^3.$
\item [] $676=(3^3-3/3)^{(3-3/3)}.$
\item [] $677=((3+3)\times 333+33)/3.$
\item [] $678=3\times (3+3)^3+3^3+3.$
\item [] $679=3+(3^3-3/3)^{(3-3/3)}.$
\item [] $680=33+3\times (3+3)^3-3/3.$
\item [] $681=33+3\times (3+3)^3.$
\item [] $682=33+3\times (3+3)^3+3/3.$
\item [] $683=(33/3)^3-3\times (3+3)^3.$
\item [] $684=3\times (3+3)^3+33+3.$
\item [] $685=3^{(3+3)}-33-33/3.$
\item [] $686=(3-3/3)\times (3+3+3/3)^3.$
\item [] $687=3^{(3+3)}-33-3\times 3.$
\item [] $688=3^{(3+3)}+3-33-33/3.$
\item [] $689=(3-3/3)\times (3+3+3/3)^3+3.$
\item [] $690=3\times ((3+3)^3+3)+33.$
\item [] $691=3^{(3+3)}-3^3-33/3.$
\item [] $692=3^{(3+3)}-3-33-3/3.$
\item [] $693=33\times (3\times (3+3)+3).$
\item [] $694=33\times (3\times (3+3)+3)+3/3.$
\item [] $695=3^{(3+3)}-33-3/3.$
\item [] $696=3^{(3+3)}-33.$
\item [] $697=3^{(3+3)}-33+3/3.$
\item [] $698=3^{(3+3)}-33+3-3/3.$
\item [] $699=3^{(3+3)}-33+3.$
\item [] $700=3^{(3+3)}-33+3+3/3.$
\item [] $701=3^{(3+3)}-3^3-3/3.$
\item [] $702=3^{(3+3)}-3^3.$
\item [] $703=3^{(3+3)}-3^3+3/3.$
\item [] $704=33\times (3+3/3)^3/3.$
\item [] $705=3^{(3+3)}-3^3+3.$
\item [] $706=3^{(3+3)}-3^3+3+3/3.$
\item [] $707=33\times (3/3+3)^3/3+3.$
\item [] $708=3^{(3+3)}-3^3+3+3.$
\item [] $709=3^{(3+3)}-3\times 3-33/3.$
\item [] $710=3^{(3+3)}-3\times (3+3)-3/3.$
\item [] $711=3^{(3+3)}-3\times (3+3).$
\item [] $712=3^{(3+3)}-3\times (3+3)+3/3.$
\item [] $713=3^{(3+3)}-3^3+33/3.$
\item [] $714=3^{(3+3)}+3-3\times (3+3).$
\item [] $715=3^{(3+3)}-3-33/3.$
\item [] $716=3^{(3+3)}-3-(33-3)/3.$
\item [] $717=3^{(3+3)}-3-3\times 3.$
\item [] $718=3^{(3+3)}-33/3.$
\item [] $719=3^{(3+3)}-(33-3)/3.$
\item [] $720=3^{(3+3)}-3\times 3.$
\item [] $721=3^{(3+3)}+3-33/3.$
\item [] $722=3^{(3+3)}-3-3-3/3.$
\item [] $723=3^{(3+3)}-3-3.$
\item [] $724=3^{(3+3)}-3-3+3/3.$
\item [] $725=3^{(3+3)}-3-3/3.$
\item [] $726=3^{(3+3)}-3.$
\item [] $727=3^{(3+3)}-3+3/3.$
\item [] $728=3^{(3+3)}-3/3.$
\item [] $729=3^{(3+3)}.$
\item [] $730=3^{(3+3)}+3/3.$
\item [] $731=3^{(3+3)}+3-3/3.$
\item [] $732=3^{(3+3)}+3.$
\item [] $733=3^{(3+3)}+3+3/3.$
\item [] $734=3^{(3+3)}+3+3-3/3.$
\item [] $735=3^{(3+3)}+3+3.$
\item [] $736=3^{(3+3)}+3+3+3/3.$
\item [] $737=3^{(3+3)}+3\times 3-3/3.$
\item [] $738=3^{(3+3)}+3\times 3.$
\item [] $739=3^{(3+3)}+3\times 3+3/3.$
\item [] $740=3^{(3+3)}+33/3.$
\item [] $741=3^{(3+3)}+3\times 3+3.$
\item [] $742=3^{(3+3)}+3\times 3+3+3/3.$
\item [] $743=3^{(3+3)}+3+33/3.$
\item [] $744=3^{(3+3)}+3\times 3+3+3.$
\item [] $745=3^{(3+3)}+3^3-33/3.$
\item [] $746=3^{(3+3)}+3+3+33/3.$
\item [] $747=3\times ((3+3)^3+33).$
\item [] $748=3\times ((3+3)^3+33)+3/3.$
\item [] $749=3^{(3+3)}+3\times 3+33/3.$
\item [] $750=3^{(3+3)}+3\times (3+3)+3.$
\item [] $751=3^{(3+3)}+33-33/3.$
\item [] $752=3^{(3+3)}+3^3-3-3/3.$
\item [] $753=3^{(3+3)}+3^3-3.$
\item [] $754=3^{(3+3)}+3^3-3+3/3.$
\item [] $755=3^{(3+3)}+3^3-3/3.$
\item [] $756=3^{(3+3)}+3^3.$
\item [] $757=3^{(3+3)}+3^3+3/3.$
\item [] $758=3^{(3+3)}+3^3+3-3/3.$
\item [] $759=3^{(3+3)}+3^3+3.$
\item [] $760=3^{(3+3)}+3^3+3+3/3.$
\item [] $761=3^{(3+3)}+33-3/3.$
\item [] $762=3^{(3+3)}+33.$
\item [] $763=3^{(3+3)}+33+3/3.$
\item [] $764=3^{(3+3)}+33+3-3/3.$
\item [] $765=3^{(3+3)}+33+3.$
\item [] $766=3^{(3+3)}+33+3+3/3.$
\item [] $767=3^{(3+3)}+3^3+33/3.$
\item [] $768=3^{(3+3)}+33+3+3.$
\item [] $769=3^{(3+3)}+33+3+3+3/3.$
\item [] $770=3^{(3+3)}+3^3+3+33/3.$
\item [] $771=3^{(3+3)}+3\times 3+33.$
\item [] $772=3^{(3+3)}+3\times 3+33+3/3.$
\item [] $773=3^{(3+3)}+33+33/3.$
\item [] $774=3^{(3+3)}+3\times 3+33+3.$
\item [] $775=3333/3-333-3.$
\item [] $776=3^{(3+3)}+33+3+33/3.$
\item [] $777=3^3\times (3^3+3)-33.$
\item [] $778=3333/3-333.$
\item [] $779=33\times 3^3-(333+3)/3.$
\item [] $780=(3^3+3)\times (3^3-3/3).$
\item [] $781=33\times ((3+3)^3-3)/(3\times 3).$
\item [] $782=3^{(3+3)}+3^3+3^3-3/3.$
\item [] $783=3\times 3\times (3\times 3^3+3+3).$
\item [] $784=(3^3+3/3)^{(3-3/3)}.$
\item [] $785=33\times (3^3-3)-3-3-3/3.$
\item [] $786=33\times (3^3-3)-3-3.$
\item [] $787=(3^3+3/3)^{(3-3/3)}+3.$
\item [] $788=33\times (3^3-3)-3-3/3.$
\item [] $789=33\times (3^3-3)-3.$
\item [] $790=33\times (3^3-3)-3+3/3.$
\item [] $791=33\times (3^3-3)-3/3.$
\item [] $792=33\times (3^3-3).$
\item [] $793=33\times (3^3-3)+3/3.$
\item [] $794=33\times (3^3-3)+3-3/3.$
\item [] $795=33\times (3^3-3)+3.$
\item [] $796=3+33\times (3^3-3)+3/3.$
\item [] $797=33\times (3^3-3)+3+3-3/3.$
\item [] $798=33\times (3^3-3)+3+3.$
\item [] $799=3^3\times (3^3+3)-33/3.$
\item [] $800=3^{(3+3)}+((3+3)^3-3)/3.$
\item [] $801=3\times (3\times 3\times (3^3+3)-3).$
\item [] $802=3^{(3+3)}+((3+3)^3+3)/3.$
\item [] $803=33\times (3^3-3)+33/3.$
\item [] $804=3^3\times (3^3+3)-3-3.$
\item [] $805=3^{(3+3)}+3+((3+3)^3+3)/3.$
\item [] $806=3^3\times (3^3+3)-3-3/3.$
\item [] $807=3^3\times (3^3+3)-3.$
\item [] $808=3^3\times (3^3+3)-3+3/3.$
\item [] $809=3^3\times (3^3+3)-3/3.$
\item [] $810=3^3\times (3^3+3).$
\item [] $811=3^3\times (3^3+3)+3/3.$
\item [] $812=3^3\times (3^3+3)+3-3/3.$
\item [] $813=3^3\times (3^3+3)+3.$
\item [] $814=3^3\times (3^3+3)+3+3/3.$
\item [] $815=3^3\times (3^3+3)+3+3-3/3.$
\item [] $816=3^3\times (3^3+3)+3+3.$
\item [] $817=3^3\times (3^3+3)+3+3+3/3.$
\item [] $818=3^3+33\times (3^3-3)-3/3.$
\item [] $819=3^3+33\times (3^3-3).$
\item [] $820=3^3+33\times (3^3-3)+3/3.$
\item [] $821=3^3\times (3^3+3)+33/3.$
\item [] $822=3^3+33\times (3^3-3)+3.$
\item [] $823=3^3+33\times (3^3-3)+3+3/3.$
\item [] $824=3^3\times (3^3+3)+3+33/3.$
\item [] $825=33+33\times (3^3-3).$
\item [] $826=33+33\times (3^3-3)+3/3.$
\item [] $827=3\times 33+3^{(3+3)}-3/3.$
\item [] $828=3\times 33+3^{(3+3)}.$
\item [] $829=3\times 33+3^{(3+3)}+3/3.$
\item [] $830=3\times 33+3^{(3+3)}+3-3/3.$
\item [] $831=3\times 33+3^{(3+3)}+3.$
\item [] $832=(3^3-3/3)\times (33-3/3).$
\item [] $833=(3333-3/3)/(3+3/3).$
\item [] $834=3\times 33+3^{(3+3)}+3+3.$
\item [] $835=(3^3-3/3)\times (33-3/3)+3.$
\item [] $836=3^3+3^3\times (3^3+3)-3/3.$
\item [] $837=3^3+3^3\times (3^3+3).$
\item [] $838=3^3+3^3\times (3^3+3)+3/3.$
\item [] $839=3^{(3+3)}+(333-3)/3.$
\item [] $840=3^{(3+3)}+333/3.$
\item [] $841=3^{(3+3)}+(333+3)/3.$
\item [] $842=3^{(3+3)}+3+(333-3)/3.$
\item [] $843=3^3\times (3^3+3)+33.$
\item [] $844=3^{(3+3)}+3+(333+3)/3.$
\item [] $845=(3-3/3)^{(3\times 3)}+333.$
\item [] $846=3^3\times (3^3+3)+33+3.$
\item [] $847=33\times (3\times 3^3-3-3/3)/3.$
\item [] $848=(3-3/3)^{(3\times 3)}+3+333.$
\item [] $849=33\times 3^3-3\times 3-33.$
\item [] $850=3^{(3+3)}+(33/3)^{(3-3/3)}.$
\item [] $851=3^{(3+3)}+(3^{(3+3)}+3)/(3+3).$
\item [] $852=(3+3/3)\times ((3+3)^3-3).$
\item [] $853=(3+3/3)\times ((3+3)^3-3)+3/3.$
\item [] $854=3^{(3+3)}+(3+3-3/3)^3.$
\item [] $855=3\times (3\times (3\times 33-3)-3).$
\item [] $856=3\times 3^{(3+3)}-(33/3)^3.$
\item [] $857=33\times 3^3-33-3/3.$
\item [] $858=33\times (3^3-3/3).$
\item [] $859=33\times (3^3-3/3)+3/3.$
\item [] $860=(3+3/3)\times ((3+3)^3-3/3).$
\item [] $861=33\times (3^3-3/3)+3.$
\item [] $862=33\times (3^3-3/3)+3+3/3.$
\item [] $863=33\times 3^3-3^3-3/3.$
\item [] $864=3\times 3\times (3\times 33-3).$
\item [] $865=3\times 3\times (3\times 33-3)+3/3.$
\item [] $866=33\times 3^3-3^3+3-3/3.$
\item [] $867=33\times 3^3-3^3+3.$
\item [] $868=(3+3/3)\times ((3+3)^3+3/3).$
\item [] $869=33\times (3^3-3/3)+33/3.$
\item [] $870=3\times 3\times (3\times 33-3)+3+3.$
\item [] $871=(3+3/3)\times ((3+3)^3+3/3)+3.$
\item [] $872=(3+3/3)\times ((3+3)^3+3-3/3).$
\item [] $873=3\times (3\times (3\times 33-3)+3).$
\item [] $874=3\times (3\times (3\times 33-3)+3)+3/3.$
\item [] $875=3\times 3\times (3\times 33-3)+33/3.$
\item [] $876=(3+3/3)\times ((3+3)^3+3).$
\item [] $877=33\times 3^3-3-33/3.$
\item [] $878=33\times 3^3-3-(33-3)/3.$
\item [] $879=33\times 3^3-3\times 3-3.$
\item [] $880=33\times 3^3-33/3.$
\item [] $881=33\times 3^3-(33-3)/3.$
\item [] $882=3\times (3\times 3\times 33-3).$
\item [] $883=33\times 3^3+3-33/3.$
\item [] $884=(3^3-3/3)\times (33+3/3).$
\item [] $885=33\times 3^3-3-3.$
\item [] $886=33\times 3^3-3-3+3/3.$
\item [] $887=33\times 3^3-3-3/3.$
\item [] $888=33\times 3^3-3.$
\item [] $889=33\times 3^3-3+3/3.$
\item [] $890=33\times 3^3-3/3.$
\item [] $891=33\times 3^3.$
\item [] $892=33\times 3^3+3/3.$
\item [] $893=33\times 3^3+3-3/3.$
\item [] $894=33\times 3^3+3.$
\item [] $895=33\times 3^3+3+3/3.$
\item [] $896=33\times 3^3+3+3-3/3.$
\item [] $897=33\times 3^3+3+3.$
\item [] $898=33\times 3^3+3+3+3/3.$
\item [] $899=3\times 3+33\times 3^3-3/3.$
\item [] $900=3\times 3+33\times 3^3.$
\item [] $901=33\times 3^3+3\times 3+3/3.$
\item [] $902=33\times 3^3+33/3.$
\item [] $903=33\times 3^3+3\times 3+3.$
\item [] $904=33\times 3^3+3\times 3+3+3/3.$
\item [] $905=33\times 3^3+3+33/3.$
\item [] $906=33\times 3^3+3\times 3+3+3.$
\item [] $907=33\times 3^3+3^3-33/3.$
\item [] $908=33\times 3^3+3+3+33/3.$
\item [] $909=3\times (3\times 3\times 33+3+3).$
\item [] $910=33\times 3^3+3\times (3+3)+3/3.$
\item [] $911=33\times 3^3+3\times 3+33/3.$
\item [] $912=3\times (3\times 3\times 33+3+3)+3.$
\item [] $913=33\times 3^3+33-33/3.$
\item [] $914=33\times 3^3+3^3-3-3/3.$
\item [] $915=33\times 3^3+3^3-3.$
\item [] $916=33\times 3^3+3^3-3+3/3.$
\item [] $917=33\times 3^3+3^3-3/3.$
\item [] $918=3\times 3\times (3\times 33+3).$
\item [] $919=33\times 3^3+3^3+3/3.$
\item [] $920=33\times 3^3+3^3+3-3/3.$
\item [] $921=33\times 3^3+3^3+3.$
\item [] $922=33\times 3^3+3^3+3+3/3.$
\item [] $923=33+33\times 3^3-3/3.$
\item [] $924=33+33\times 3^3.$
\item [] $925=33+33\times 3^3+3/3.$
\item [] $926=33+33\times 3^3+3-3/3.$
\item [] $927=33+33\times 3^3+3.$
\item [] $928=33+33\times 3^3+3+3/3.$
\item [] $929=33\times 3^3+3^3+33/3.$
\item [] $930=33+33\times 3^3+3+3.$
\item [] $931=33+33\times 3^3+3+3+3/3.$
\item [] $932=3\times 333-3-(3+3/3)^3.$
\item [] $933=33\times 3^3+33+3\times 3.$
\item [] $934=3^{(3+3)}+(3+3)^3-33/3.$
\item [] $935=3\times 333-(3+3/3)^3.$
\item [] $936=(3\times 3+3)\times (3\times 3^3-3).$
\item [] $937=(3\times 3+3)\times (3\times 3^3-3)+3/3.$
\item [] $938=3\times 333+3-(3+3/3)^3.$
\item [] $939=3^3\times (33+3)-33.$
\item [] $940=3^3\times (33+3)-33+3/3.$
\item [] $941=3^{(3+3)}+(3+3)^3-3-3/3.$
\item [] $942=3^{(3+3)}+(3+3)^3-3.$
\item [] $943=3^{(3+3)}+(3+3)^3-3+3/3.$
\item [] $944=3^{(3+3)}+(3+3)^3-3/3.$
\item [] $945=3^{(3+3)}+(3+3)^3.$
\item [] $946=3^{(3+3)}+(3+3)^3+3/3.$
\item [] $947=3^{(3+3)}+(3+3)^3+3-3/3.$
\item [] $948=3^{(3+3)}+(3+3)^3+3.$
\item [] $949=3^{(3+3)}+(3+3)^3+3+3/3.$
\item [] $950=3^{(3+3)}+(3+3)^3+3+3-3/3.$
\item [] $951=3^{(3+3)}+(3+3)^3+3+3.$
\item [] $952=(3^3+3/3)\times (33+3/3).$
\item [] $953=(3+3)\times ((3+3)\times 3^3-3)-3/3.$
\item [] $954=(3+3)\times ((3+3)\times 3^3-3).$
\item [] $955=33\times 3^3+(3+3/3)^3.$
\item [] $956=3^{(3+3)}+(3+3)^3+33/3.$
\item [] $957=33\times (3^3+3-3/3).$
\item [] $958=3+33\times 3^3+(3+3/3)^3.$
\item [] $959=(3+3^3)\times (33-3/3)-3/3.$
\item [] $960=(3+3^3)\times (33-3/3).$
\item [] $961=(3^3+3+3/3)^{(3-3/3)}.$
\item [] $962=33\times 3^3+((3+3)^3-3)/3.$
\item [] $963=3\times (333-3\times 3-3).$
\item [] $964=(3\times 3+3/3)^3-33-3.$
\item [] $965=3\times 333-33-3/3.$
\item [] $966=3\times 333-33.$
\item [] $967=(3\times 3+3/3)^3-33.$
\item [] $968=(3+3/3)\times (3^{(3+3)}-3)/3.$
\item [] $969=3^3\times (33+3)-3.$
\item [] $970=3^3\times (33+3)-3+3/3.$
\item [] $971=3^3\times (33+3)-3/3.$
\item [] $972=3^3\times (33+3).$
\item [] $973=3^3\times (33+3)+3/3.$
\item [] $974=3^3\times (33+3)+3-3/3.$
\item [] $975=3^3\times (33+3)+3.$
\item [] $976=3^3\times (33+3)+3+3/3.$
\item [] $977=3^3\times (33+3)+3+3-3/3.$
\item [] $978=3^3\times (33+3)+3+3.$
\item [] $979=3\times (333-3)-33/3.$
\item [] $980=3\times (333-3-3)-3/3.$
\item [] $981=3\times (333-3-3).$
\item [] $982=3\times (333-3-3)+3/3.$
\item [] $983=3^3\times (33+3)+33/3.$
\item [] $984=3\times (333-3-3)+3.$
\item [] $985=3\times 333-3-33/3.$
\item [] $986=3\times (333-3)-3-3/3.$
\item [] $987=3\times (333-3)-3.$
\item [] $988=3\times 333-33/3.$
\item [] $989=3\times (333-3)-3/3.$
\item [] $990=3\times (333-3).$
\item [] $991=3\times (333-3)+3/3.$
\item [] $992=3\times (333-3)+3-3/3.$
\item [] $993=3\times (333-3)+3.$
\item [] $994=(3\times 3+3/3)^3-3-3.$
\item [] $995=3\times 333-3-3/3.$
\item [] $996=3\times 333-3.$
\item [] $997=(3\times 3+3/3)^3-3.$
\item [] $998=3\times 333-3/3.$
\item [] $999=3\times 333.$
\item [] $1000=(3\times 3+3/3)^3.$
\end{itemize}
\end{multicols}
}

\section{\textbf{Representations Using Number 4}}

{\footnotesize
\begin{multicols}{3}
\begin{itemize}
\item [] $101=4444/44.$
\item [] $102=(444-4)/4-4-4.$
\item [] $103=444/4-4-4.$
\item [] $104=4\times 4+44+44.$
\item [] $105=4+4444/44.$
\item [] $106=(444-4)/4-4.$
\item [] $107=444/4-4.$
\item [] $108=44+4\times 4\times 4.$
\item [] $109=44+(4^4+4)/4.$
\item [] $110=(444-4)/4.$
\item [] $111=444/4.$
\item [] $112=4\times (44-4\times 4).$
\item [] $113=(444+4+4)/4.$
\item [] $114=4+(444-4)/4.$
\item [] $115=4+444/4.$
\item [] $116=4+4\times (44-4\times 4).$
\item [] $117=4+(444+4+4)/4.$
\item [] $118=4+4+(444-4)/4.$
\item [] $119=4+4+444/4.$
\item [] $120=(4+4)\times (44/4+4).$
\item [] $121=(44/4)^{((4+4)/4)}.$
\item [] $122=(444+44)/4.$
\item [] $123=4+4+4+444/4.$
\item [] $124=4\times 4\times (4+4)-4.$
\item [] $125=44+(4-4/4)^4.$
\item [] $126=4\times (4^4-4)/(4+4).$
\item [] $127=4\times 4+444/4.$
\item [] $128=4\times 4\times (4+4).$
\item [] $129=4\times 4\times (4+4)+4/4.$
\item [] $130=4\times (4^4+4)/(4+4).$
\item [] $131=4\times 4+4+444/4.$
\item [] $132=4\times 4\times (4+4)+4.$
\item [] $133=4\times 4\times (4+4)+4+4/4.$
\item [] $134=4+4\times (4^4+4)/(4+4).$
\item [] $135=(4-4/4)\times (44+4/4).$
\item [] $136=4\times 4\times (4+4)+4+4.$
\item [] $137=4\times 4\times (4+4)+4+4+4/4.$
\item [] $138=4\times 4\times (4+4)+(44-4)/4.$
\item [] $139=4\times 4\times (4+4)+44/4.$
\item [] $140=4\times (4\times (4+4)+4)-4.$
\item [] $141=4^4-4-444/4.$
\item [] $142=4\times (4+4)+(444-4)/4.$
\item [] $143=4\times (4+4)+444/4.$
\item [] $144=4\times (4\times (4+4)+4).$
\item [] $145=4^4-444/4.$
\item [] $146=4^4+(4-444)/4.$
\item [] $147=4\times (4+4)+4+444/4.$
\item [] $148=4\times (4\times (4+4)+4)+4.$
\item [] $149=4^4+4-444/4.$
\item [] $150=4\times (44+4^4)/(4+4).$
\item [] $151=44-4+444/4.$
\item [] $152=4\times (44-4)-4-4.$
\item [] $153=4^4+4+4-444/4.$
\item [] $154=44+(444-4)/4.$
\item [] $155=44+444/4.$
\item [] $156=4\times (44-4)-4.$
\item [] $157=4\times (44-4)-4+4/4.$
\item [] $158=4\times (44-4)-(4+4)/4.$
\item [] $159=4\times (44-4)-4/4.$
\item [] $160=4\times (44-4).$
\item [] $161=4\times (44-4)+4/4.$
\item [] $162=4\times (44-4)+(4+4)/4.$
\item [] $163=4\times (44-4)+4-4/4.$
\item [] $164=4\times (44-4)+4.$
\item [] $165=4\times 44-44/4.$
\item [] $166=4\times 44-(44-4)/4.$
\item [] $167=4\times 44-4-4-4/4.$
\item [] $168=4\times 44-4-4.$
\item [] $169=4\times 44+4-44/4.$
\item [] $170=4\times 44-4-(4+4)/4.$
\item [] $171=4\times 44-4-4/4.$
\item [] $172=4\times 44-4.$
\item [] $173=4\times 44-4+4/4.$
\item [] $174=4\times 44-(4+4)/4.$
\item [] $175=4\times 44-4/4.$
\item [] $176=4\times 44.$
\item [] $177=4\times 44+4/4.$
\item [] $178=4\times 44+(4+4)/4.$
\item [] $179=4\times 44+4-4/4.$
\item [] $180=4\times 44+4.$
\item [] $181=4\times 44+4+4/4.$
\item [] $182=4\times 44+4+(4+4)/4.$
\item [] $183=4\times 44+4+4-4/4.$
\item [] $184=4\times 44+4+4.$
\item [] $185=4\times 44+4+4+4/4.$
\item [] $186=4\times 44+(44-4)/4.$
\item [] $187=4\times 44+44/4.$
\item [] $188=444-4^4.$
\item [] $189=444-4^4+4/4.$
\item [] $190=4^4-(4^4+4+4)/4.$
\item [] $191=4^4-(4^4+4)/4.$
\item [] $192=4\times (44+4).$
\item [] $193=4\times (44+4)+4/4.$
\item [] $194=4^4-(4^4-4-4)/4.$
\item [] $195=4^4+4-(4^4+4)/4.$
\item [] $196=4\times (44+4)+4.$
\item [] $197=4^4+4-(4^4-4)/4.$
\item [] $198=4^4+4-(4^4-4-4)/4.$
\item [] $199=4^4+4+4-(4^4+4)/4.$
\item [] $200=4\times (44+4)+4+4.$
\item [] $201=4^4-44-44/4.$
\item [] $202=4^4-44+(44-4)/4.$
\item [] $203=4^4-(4^4-44)/4.$
\item [] $204=44+4\times (44-4).$
\item [] $205=((4+4)^4+4)/(4\times 4+4).$
\item [] $206=4^4-4-44-(4+4)/4.$
\item [] $207=4^4-4-44-4/4.$
\item [] $208=4\times (44+4+4).$
\item [] $209=4\times (44+4+4)+4/4.$
\item [] $210=4^4-44-(4+4)/4.$
\item [] $211=4^4-44-4/4.$
\item [] $212=4^4-44.$
\item [] $213=4^4-44+4/4.$
\item [] $214=4^4-44+(4+4)/4.$
\item [] $215=4^4+4-44-4/4.$
\item [] $216=4^4+4-44.$
\item [] $217=4^4+4-44+4/4.$
\item [] $218=4\times 444/(4+4)-4.$
\item [] $219=44+4\times 44-4/4.$
\item [] $220=44+4\times 44.$
\item [] $221=44+4\times 44+4/4.$
\item [] $222=4\times 444/(4+4).$
\item [] $223=4^4-44+44/4.$
\item [] $224=4^4-4\times (4+4).$
\item [] $225=4^4-4\times (4+4)+4/4.$
\item [] $226=4+4\times 444/(4+4).$
\item [] $227=4^4+4-44+44/4.$
\item [] $228=4^4+4-4\times (4+4).$
\item [] $229=4^4-4\times 4-44/4.$
\item [] $230=4+4+4\times 444/(4+4).$
\item [] $231=44+4\times 44+44/4.$
\item [] $232=4^4-4\times 4-4-4.$
\item [] $233=4+4^4-4\times 4-44/4.$
\item [] $234=4^4-4\times 44/(4+4).$
\item [] $235=4^4-4-4\times 4-4/4.$
\item [] $236=4^4-4\times 4-4.$
\item [] $237=4^4-4\times 4-4+4/4.$
\item [] $238=4^4-4\times 4-(4+4)/4.$
\item [] $239=4^4-4\times 4-4/4.$
\item [] $240=4^4-4\times 4.$
\item [] $241=4^4-4\times 4+4/4.$
\item [] $242=44\times 44/(4+4).$
\item [] $243=(4-4/4)^{(4+4/4)}.$
\item [] $244=4^4-4\times 4+4.$
\item [] $245=4^4-44/4.$
\item [] $246=4^4-(44-4)/4.$
\item [] $247=4^4-4-4-4/4.$
\item [] $248=4^4-4-4.$
\item [] $249=4^4+4-44/4.$
\item [] $250=4^4-4-(4+4)/4.$
\item [] $251=4^4-4-4/4.$
\item [] $252=4^4-4.$
\item [] $253=4^4-4+4/4.$
\item [] $254=4^4-(4+4)/4.$
\item [] $255=4^4-4/4.$
\item [] $256=4^4.$
\item [] $257=4^4+4/4.$
\item [] $258=4^4+(4+4)/4.$
\item [] $259=4^4+4-4/4.$
\item [] $260=4^4+4.$
\item [] $261=4^4+4+4/4.$
\item [] $262=4^4+4+(4+4)/4.$
\item [] $263=4^4+4+4-4/4.$
\item [] $264=4^4+4+4.$
\item [] $265=4^4+4+4+4/4.$
\item [] $266=4^4+(44-4)/4.$
\item [] $267=4^4+44/4.$
\item [] $268=4^4+4+4+4.$
\item [] $269=4^4+4+4+4+4/4.$
\item [] $270=4^4+4+(44-4)/4.$
\item [] $271=4^4+4+44/4.$
\item [] $272=4^4+4\times 4.$
\item [] $273=4^4+4\times 4+4/4.$
\item [] $274=4^4+4\times 4+(4+4)/4.$
\item [] $275=4^4+4+4+44/4.$
\item [] $276=4^4+4\times 4+4.$
\item [] $277=4^4+4\times 4+4+4/4.$
\item [] $278=4^4+4\times 44/(4+4).$
\item [] $279=4^4+4+4+4+44/4.$
\item [] $280=4^4+4\times 4+4+4.$
\item [] $281=4^4+4\times 4+4+4+4/4.$
\item [] $282=4^4+4+4\times 44/(4+4).$
\item [] $283=4\times 4+4^4+44/4.$
\item [] $284=44+4^4-4\times 4.$
\item [] $285=44+4^4-4\times 4-4/4.$
\item [] $286=44+44\times 44/(4+4).$
\item [] $287=4^4+4\times (4+4)-4/4.$
\item [] $288=4^4+4\times (4+4).$
\item [] $289=44+4^4-44/4.$
\item [] $290=44+4^4-(44-4)/4.$
\item [] $291=4\times (4+4)+4^4+4-4/4.$
\item [] $292=4\times (4+4)+4^4+4.$
\item [] $293=4\times (4+4)+4^4+4+4/4.$
\item [] $294=44+4^4-4-(4+4)/4.$
\item [] $295=44+4^4-4-4/4.$
\item [] $296=44+4^4-4.$
\item [] $297=44+4^4-4+4/4.$
\item [] $298=44+4^4-(4+4)/4.$
\item [] $299=44+4^4-4/4.$
\item [] $300=44+4^4.$
\item [] $301=4^4+44+4/4.$
\item [] $302=4^4+44+(4+4)/4.$
\item [] $303=4^4+44+4-4/4.$
\item [] $304=4^4+44+4.$
\item [] $305=4^4+44+4+4/4.$
\item [] $306=4^4+44+4+(4+4)/4.$
\item [] $307=4^4+44+4+4-4/4.$
\item [] $308=4^4+44+4+4.$
\item [] $309=4^4+(4^4-44)/4.$
\item [] $310=4^4+44+(44-4)/4.$
\item [] $311=4^4+44+44/4.$
\item [] $312=4^4+44+4+4+4.$
\item [] $313=4+4^4+(4^4-44)/4.$
\item [] $314=4^4-4+(4^4-4-4)/4.$
\item [] $315=4^4-4+(4^4-4)/4.$
\item [] $316=4^4+44+4\times 4.$
\item [] $317=4^4-4+(4^4+4)/4.$
\item [] $318=4^4+(4^4-4-4)/4.$
\item [] $319=4^4+(4^4-4)/4.$
\item [] $320=4\times 4\times (4\times 4+4).$
\item [] $321=4^4+(4^4+4)/4.$
\item [] $322=4^4+(4^4+4+4)/4.$
\item [] $323=4^4+4+(4^4-4)/4.$
\item [] $324=4\times (4-4/4)^4.$
\item [] $325=4^4+4+(4^4+4)/4.$
\item [] $326=4^4+4+(4^4+4+4)/4.$
\item [] $327=4^4+4+4+(4^4-4)/4.$
\item [] $328=4+4\times (4-4/4)^4.$
\item [] $329=4^4+4+4+(4^4+4)/4.$
\item [] $330=4^4+(4^4+44-4)/4.$
\item [] $331=4^4+(44+4^4)/4.$
\item [] $332=4+4+4\times (4-4/4)^4.$
\item [] $333=4^4-4+(4-4/4)^4.$
\item [] $334=444-(444-4)/4.$
\item [] $335=4^4+4+(44+4^4)/4.$
\item [] $336=4\times (4\times (4\times 4+4)+4).$
\item [] $337=4^4+(4-4/4)^4.$
\item [] $338=4^4+(4-4/4)^4+4/4.$
\item [] $339=4\times ((4-4/4)^4+4)-4/4.$
\item [] $340=4\times ((4-4/4)^4+4).$
\item [] $341=4^4+4+(4-4/4)^4.$
\item [] $342=44\times (4+4)-(44-4)/4.$
\item [] $343=(4+4-4/4)^{(4-4/4)}.$
\item [] $344=4^4+44+44.$
\item [] $345=4^4+44+44+4/4.$
\item [] $346=44\times (4+4)-4-(4+4)/4.$
\item [] $347=44\times (4+4)-4-4/4.$
\item [] $348=44\times (4+4)-4.$
\item [] $349=44\times (4+4)-4+4/4.$
\item [] $350=44\times (4+4)-(4+4)/4.$
\item [] $351=44\times (4+4)-4/4.$
\item [] $352=44\times (4+4).$
\item [] $353=44\times (4+4)+4/4.$
\item [] $354=44\times (4+4)+(4+4)/4.$
\item [] $355=44\times (4+4)+4-4/4.$
\item [] $356=44\times (4+4)+4.$
\item [] $357=44\times (4+4)+4+4/4.$
\item [] $358=44\times (4+4)+4+(4+4)/4.$
\item [] $359=44\times (4+4)+4+4-4/4.$
\item [] $360=44\times (4+4)+4+4.$
\item [] $361=44\times (4+4)+4+4+4/4.$
\item [] $362=4^4-4+(444-4)/4.$
\item [] $363=4^4-4+444/4.$
\item [] $364=44\times (4+4)+4+4+4.$
\item [] $365=(4+4/4)^4-4^4-4.$
\item [] $366=4^4+(444-4)/4.$
\item [] $367=4^4+444/4.$
\item [] $368=4\times (44+44+4).$
\item [] $369=(4+4/4)^4-4^4.$
\item [] $370=4^4+4+(444-4)/4.$
\item [] $371=4^4+4+444/4.$
\item [] $372=44\times (4+4)+4\times 4+4.$
\item [] $373=(4+4/4)^4+4-4^4.$
\item [] $374=44\times (4\times 4\times 4+4)/(4+4).$
\item [] $375=4^4+4+4+444/4.$
\item [] $376=(4+4)\times (44+4-4/4).$
\item [] $377=(4+4/4)^4-4^4+4+4.$
\item [] $378=4^4+(444+44)/4.$
\item [] $379=444-(4^4+4)/4.$
\item [] $380=444-4\times 4\times 4.$
\item [] $381=444+(4-4^4)/4.$
\item [] $382=4^4+4\times (4^4-4)/(4+4).$
\item [] $383=(4+4)\times (44+4)-4/4.$
\item [] $384=(4+4)\times (44+4).$
\item [] $385=(4+4)\times (44+4)+4/4.$
\item [] $386=4^4+4\times (4^4+4)/(4+4).$
\item [] $387=(4+4)\times (44+4)+4-4/4.$
\item [] $388=(4+4)\times (44+4)+4.$
\item [] $389=(4+4)\times (44+4)+4+4/4.$
\item [] $390=(4+(4+4)/4)\times (4^4+4)/4.$
\item [] $391=444-(4^4-44)/4.$
\item [] $392=(4+4)\times (44+4)+4+4.$
\item [] $393=(4+4)\times (44+4)+4+4+4/4.$
\item [] $394=(4+(4+4)/4)\times (4^4+4)/4+4.$
\item [] $395=44\times (4+4)+44-4/4.$
\item [] $396=44\times (4+4)+44.$
\item [] $397=44\times (4+4)+44+4/4.$
\item [] $398=444+44-(4+4)/4.$
\item [] $399=4\times 4^4-(4+4/4)^4.$
\item [] $400=444-44.$
\item [] $401=444-44+4/4.$
\item [] $402=444-44+(4+4)/4.$
\item [] $403=4\times 4^4+4-(4+4/4)^4.$
\item [] $404=444-44+4.$
\item [] $405=(4+4/4)\times (4-4/4)^4.$
\item [] $406=444-44+4+(4+4)/4.$
\item [] $407=4\times 4^4+4+4-(4+4/4)^4.$
\item [] $408=444-44+4+4.$
\item [] $409=(4+4/4)\times (4-4/4)^4+4.$
\item [] $410=4\times ((4+4)^4+4)/(44-4).$
\item [] $411=4^4+44+444/4.$
\item [] $412=444-4\times (4+4).$
\item [] $413=44-4^4+(4+4/4)^4.$
\item [] $414=4\times ((4+4)^4+4)/(44-4)+4.$
\item [] $415=4^4+4\times (44-4)-4/4.$
\item [] $416=4^4+4\times (44-4).$
\item [] $417=4^4+4\times (44-4)+4/4.$
\item [] $418=4^4+4\times (44-4)+(4+4)/4.$
\item [] $419=(44\times 44-4^4-4)/4.$
\item [] $420=4^4+4+4\times (44-4).$
\item [] $421=4^4+4\times 44-44/4.$
\item [] $422=444-4\times 44/(4+4).$
\item [] $423=444-4\times 4-4-4/4.$
\item [] $424=444-4\times 4-4.$
\item [] $425=444-4\times 4-4+4/4.$
\item [] $426=444-4\times 4-(4+4)/4.$
\item [] $427=444-4\times 4-4/4.$
\item [] $428=444-4\times 4.$
\item [] $429=444-4\times 4+4/4.$
\item [] $430=4^4+4\times 44-(4+4)/4.$
\item [] $431=4^4+4\times 44-4/4.$
\item [] $432=4^4+4\times 44.$
\item [] $433=444-44/4.$
\item [] $434=444-(44-4)/4.$
\item [] $435=444-4-4-4/4.$
\item [] $436=444-4-4.$
\item [] $437=444+4-44/4.$
\item [] $438=444-4-(4+4)/4.$
\item [] $439=444-4-4/4.$
\item [] $440=444-4.$
\item [] $441=444-4+4/4.$
\item [] $442=444-(4+4)/4.$
\item [] $443=444-4/4.$
\item [] $444=444.$
\item [] $445=444+4/4.$
\item [] $446=444+(4+4)/4.$
\item [] $447=444+4-4/4.$
\item [] $448=444+4.$
\item [] $449=444+4+4/4.$
\item [] $450=444+4+(4+4)/4.$
\item [] $451=444+4+4-4/4.$
\item [] $452=444+4+4.$
\item [] $453=444+4+4+4/4.$
\item [] $454=444+(44-4)/4.$
\item [] $455=444+44/4.$
\item [] $456=444+4+4+4.$
\item [] $457=444+4+4+4+4/4.$
\item [] $458=444+4+(44-4)/4.$
\item [] $459=444+4+44/4.$
\item [] $460=444+4\times 4.$
\item [] $461=444+4\times 4+4/4.$
\item [] $462=444+4\times 4+(4+4)/4.$
\item [] $463=444+4+4+44/4.$
\item [] $464=444+4\times 4+4.$
\item [] $465=444+4\times 4+4+4/4.$
\item [] $466=444+4\times 44/(4+4).$
\item [] $467=4^4+4^4-44-4/4.$
\item [] $468=4^4+4^4-44.$
\item [] $469=4^4+4^4-44+4/4.$
\item [] $470=444+4+4\times 44/(4+4).$
\item [] $471=444+4\times 4+44/4.$
\item [] $472=4^4+4^4+4-44.$
\item [] $473=44\times (44-4/4)/4.$
\item [] $474=(44\times (44-4/4)+4)/4.$
\item [] $475=444+4\times (4+4)-4/4.$
\item [] $476=444+4\times (4+4).$
\item [] $477=4+44\times (44-4/4)/4.$
\item [] $478=4^4+4\times 444/(4+4).$
\item [] $479=(44\times 44-4)/4-4.$
\item [] $480=(4+4)\times (4\times 4+44).$
\item [] $481=(44\times 44+4)/4-4.$
\item [] $482=(44\times 44-4-4)/4.$
\item [] $483=(44\times 44-4)/4.$
\item [] $484=44\times 44/4.$
\item [] $485=(44\times 44+4)/4.$
\item [] $486=(44\times 44+4+4)/4.$
\item [] $487=4+(44\times 44-4)/4.$
\item [] $488=44+444.$
\item [] $489=4+(44\times 44+4)/4.$
\item [] $490=4+(44\times 44+4+4)/4.$
\item [] $491=4+4+(44\times 44-4)/4.$
\item [] $492=444+44+4.$
\item [] $493=4+4+(44\times 44+4)/4.$
\item [] $494=(44\times 44+44-4)/4.$
\item [] $495=44\times (44+4/4)/4.$
\item [] $496=4^4+4^4-4\times 4.$
\item [] $497=4^4+4^4-4\times 4+4/4.$
\item [] $498=4^4+44\times 44/(4+4).$
\item [] $499=4+44\times (44+4/4)/4.$
\item [] $500=4^4+4^4-4\times 4+4.$
\item [] $501=4^4+4^4-44/4.$
\item [] $502=4^4+4^4-(44-4)/4.$
\item [] $503=((4+4)\times (4^4-4)-4)/4.$
\item [] $504=4^4+4^4-4-4.$
\item [] $505=((4+4)\times (4^4-4)+4)/4.$
\item [] $506=4^4+4^4-4-(4+4)/4.$
\item [] $507=4^4+4^4-4-4/4.$
\item [] $508=4^4+4^4-4.$
\item [] $509=4^4+4^4-4+4/4.$
\item [] $510=4^4+4^4-(4+4)/4.$
\item [] $511=4^4+4^4-4/4.$
\item [] $512=4^4+4^4.$
\item [] $513=4^4+4^4+4/4.$
\item [] $514=4^4+4^4+(4+4)/4.$
\item [] $515=4^4+4^4+4-4/4.$
\item [] $516=4^4+4^4+4.$
\item [] $517=4^4+4^4+4+4/4.$
\item [] $518=4^4+4^4+4+(4+4)/4.$
\item [] $519=((4+4)\times (4^4+4)-4)/4.$
\item [] $520=4^4+4^4+4+4.$
\item [] $521=((4+4)\times (4^4+4)+4)/4.$
\item [] $522=4^4+4^4+(44-4)/4.$
\item [] $523=4^4+4^4+44/4.$
\item [] $524=44\times (4+4+4)-4.$
\item [] $525=444+(4-4/4)^4.$
\item [] $526=44\times (4+4+4)-(4+4)/4.$
\item [] $527=44\times (4+4+4)-4/4.$
\item [] $528=44\times (4+4+4).$
\item [] $529=44\times (4+4+4)+4/4.$
\item [] $530=44\times (4+4+4)+(4+4)/4.$
\item [] $531=44\times (4+4+4)+4-4/4.$
\item [] $532=44\times (4+4+4)+4.$
\item [] $533=44\times (4+4+4)+4+4/4.$
\item [] $534=(4+4)\times (4^4+44/4)/4.$
\item [] $535=(4+4/4)\times (444/4-4).$
\item [] $536=44\times (4+4+4)+4+4.$
\item [] $537=(4+4/4)^4-44-44.$
\item [] $538=(4+4)\times (4^4+44/4)/4+4.$
\item [] $539=44\times (4+4+4)+44/4.$
\item [] $540=(4-4/4)\times (4\times 44+4).$
\item [] $541=(4+4/4)^4+4-44-44.$
\item [] $542=(4+4)\times (4^4+4+44/4)/4.$
\item [] $543=(4+4)\times (4\times (4\times 4)+4)-4/4.$
\item [] $544=(4+4)\times (4\times (4\times 4)+4).$
\item [] $545=(4+4)\times (4\times (4\times 4)+4)+4/4.$
\item [] $546=(44\times 44+4^4-4-4)/4.$
\item [] $547=(44\times 44+4^4-4)/4.$
\item [] $548=(44\times 44+4^4)/4.$
\item [] $549=(44\times 44+4^4+4)/4.$
\item [] $550=(4+4/4)\times (444-4)/4.$
\item [] $551=444-4+444/4.$
\item [] $552=4^4+4^4+44-4.$
\item [] $553=444+(444-4-4)/4.$
\item [] $554=444+(444-4)/4.$
\item [] $555=444+444/4.$
\item [] $556=4^4+4^4+44.$
\item [] $557=4^4+4^4+44+4/4.$
\item [] $558=444+4+(444-4)/4.$
\item [] $559=444+4+444/4.$
\item [] $560=4^4+4^4+44+4.$
\item [] $561=(4+4/4)^4-4\times 4\times 4.$
\item [] $562=(4+4/4)^4-(4^4-4)/4.$
\item [] $563=(4^4\times 44-4)/(4\times 4+4).$
\item [] $564=4\times (4^4-4)-444.$
\item [] $565=(4+4/4)^4+4-4\times 4\times 4.$
\item [] $566=(4+4/4)^4+4-(4^4-4)/4.$
\item [] $567=(4+4+4/4)\times (4^4-4)/4.$
\item [] $568=4\times (4^4-4)+4-444.$
\item [] $569=(4+4/4)^4+4+4-4\times 4\times 4.$
\item [] $570=444+4\times (4^4-4)/(4+4).$
\item [] $571=(44\times (44+4+4)-4)/4.$
\item [] $572=44\times (4+4+4)+44.$
\item [] $573=(4+4/4)^4-44-4-4.$
\item [] $574=4^4+4^4+(4^4-4-4)/4.$
\item [] $575=4^4+4^4+(4^4-4)/4.$
\item [] $576=4\times 4\times (4\times (4+4)+4).$
\item [] $577=(4+4/4)^4-44-4.$
\item [] $578=(4+4/4)^4-44-4+4/4.$
\item [] $579=4\times 4^4-444-4/4.$
\item [] $580=4\times 4^4-444.$
\item [] $581=(4+4/4)^4-44.$
\item [] $582=(4+4/4)^4-44+4/4.$
\item [] $583=44\times (4^4-44)/(4\times 4).$
\item [] $584=4\times 4^4-444+4.$
\item [] $585=(4+4/4)^4-44+4.$
\item [] $586=(4+4/4)^4-44+4+4/4.$
\item [] $587=44\times (4^4-44)/(4\times 4)+4.$
\item [] $588=4\times 4^4-444+4+4.$
\item [] $589=(4+4/4)^4-44+4+4.$
\item [] $590=(44-4)\times ((4^4-4)/4-4)/4.$
\item [] $591=4\times 4^4-444+44/4.$
\item [] $592=4\times (4\times (4\times (4+4)+4)+4).$
\item [] $593=(4+4/4)^4-4\times (4+4).$
\item [] $594=(4+4/4)^4-4\times (4+4)+4/4.$
\item [] $595=(44\times 44+444)/4.$
\item [] $596=4\times (4^4+4)-444.$
\item [] $597=(4+4/4)^4+4-4\times (4+4).$
\item [] $598=(4+4)\times (44+4^4-4/4)/4.$
\item [] $599=((4+4)\times (44+4^4)-4)/4.$
\item [] $600=(4+4)\times (44+4^4)/4.$
\item [] $601=(4+4/4)^4-4\times 4-4-4.$
\item [] $602=(4+4)\times (44+4^4+4/4)/4.$
\item [] $603=4+((4+4)\times (44+4^4)-4)/4.$
\item [] $604=4^4+(4+4)\times 44-4.$
\item [] $605=(4+4/4)^4-4\times 4-4.$
\item [] $606=(4+4/4)^4-4\times 4-4+4/4.$
\item [] $607=4^4+(4+4)\times 44-4/4.$
\item [] $608=4^4+(4+4)\times 44.$
\item [] $609=(4+4/4)^4-4\times 4.$
\item [] $610=(4+4/4)^4-4\times 4+4/4.$
\item [] $611=(4+4)\times 44+4^4+4-4/4.$
\item [] $612=(4+4)\times 44+4^4+4.$
\item [] $613=(4+4/4)^4+4-4\times 4.$
\item [] $614=(4+4/4)^4-44/4.$
\item [] $615=(4+4/4)^4-(44-4)/4.$
\item [] $616=(4+4)\times ((4-4/4)^4-4).$
\item [] $617=(4+4/4)^4-4-4.$
\item [] $618=(4+4/4)^4+4-44/4.$
\item [] $619=444+4\times 44-4/4.$
\item [] $620=444+4\times 44.$
\item [] $621=(4+4/4)^4-4.$
\item [] $622=(4+4/4)^4-4+4/4.$
\item [] $623=(4+4/4)^4-(4+4)/4.$
\item [] $624=4\times (4\times (44-4)-4).$
\item [] $625=(4+4/4)^4.$
\item [] $626=(4+4/4)^4+4/4.$
\item [] $627=(4+4/4)^4+(4+4)/4.$
\item [] $628=4\times (4\times (44-4)-4)+4.$
\item [] $629=(4+4/4)^4+4.$
\item [] $630=(4+4/4)^4+4+4/4.$
\item [] $631=(4+4/4)^4+4+(4+4)/4.$
\item [] $632=4\times 4\times (44-4)-4-4.$
\item [] $633=(4+4/4)^4+4+4.$
\item [] $634=(4+4/4)^4+4+4+4/4.$
\item [] $635=(4+4/4)^4+(44-4)/4.$
\item [] $636=4\times 4\times (44-4)-4.$
\item [] $637=(4+4/4)^4+4+4+4.$
\item [] $638=4\times 4\times (44-4)-(4+4)/4.$
\item [] $639=4\times 4\times (44-4)-4/4.$
\item [] $640=4\times 4\times (44-4).$
\item [] $641=4\times 4+(4+4/4)^4.$
\item [] $642=4\times 4+(4+4/4)^4+4/4.$
\item [] $643=4\times 4\times (44-4)+4-4/4.$
\item [] $644=4\times 4\times (44-4)+4.$
\item [] $645=4\times 4+(4+4/4)^4+4.$
\item [] $646=4\times 4+(4+4/4)^4+4+4/4.$
\item [] $647=(4+4)\times (4-4/4)^4-4/4.$
\item [] $648=(4+4)\times (4-4/4)^4.$
\item [] $649=4\times 4+(4+4/4)^4+4+4.$
\item [] $650=(44-4)\times (4^4+4)/(4\times 4).$
\item [] $651=4\times 4\times (44-4)+44/4.$
\item [] $652=(4+4)\times (4-4/4)^4+4.$
\item [] $653=(4+4/4)^4+44-4\times 4.$
\item [] $654=(44-4)\times (4^4+4)/(4\times 4)+4.$
\item [] $655=4\times (4\times (44-4)+4)-4/4.$
\item [] $656=4\times (4\times (44-4)+4).$
\item [] $657=4\times (4+4)+(4+4/4)^4.$
\item [] $658=(4+4/4)^4+44-44/4.$
\item [] $659=4\times 4\times 44-44-4/4.$
\item [] $660=4\times 4\times 44-44.$
\item [] $661=4\times (4+4)+4+(4+4/4)^4.$
\item [] $662=44\times (4+44/4)+(4+4)/4.$
\item [] $663=4\times 4\times 44+4-44-4/4.$
\item [] $664=44\times (4+44/4)+4.$
\item [] $665=(4+4/4)^4+44-4.$
\item [] $666=444\times (4+(4+4)/4)/4.$
\item [] $667=4444/4-444.$
\item [] $668=4\times (4\times 44-4-4)-4.$
\item [] $669=44+(4+4/4)^4.$
\item [] $670=44+(4+4/4)^4+4/4.$
\item [] $671=4\times (4\times 44-4-4)-4/4.$
\item [] $672=4\times (4\times 44-4-4).$
\item [] $673=44+4+(4+4/4)^4.$
\item [] $674=44+4+(4+4/4)^4+4/4.$
\item [] $675=(4+44/4)\times (44+4/4).$
\item [] $676=4+4\times (4\times 44-4-4).$
\item [] $677=4+4\times (4\times 44-4-4)+4/4.$
\item [] $678=4\times (4\times 44-4)-(44-4)/4.$
\item [] $679=(4+44/4)\times (44+4/4)+4.$
\item [] $680=(4+4)\times ((4-4/4)^4+4).$
\item [] $681=(4+4/4)^4+44+4+4+4.$
\item [] $682=44\times (4^4-4-4)/(4\times 4).$
\item [] $683=4\times (4\times 44-4)-4-4/4.$
\item [] $684=4\times (4\times 44-4)-4.$
\item [] $685=4\times 4+44+(4+4/4)^4.$
\item [] $686=4\times (4\times 44-4)-(4+4)/4.$
\item [] $687=4\times (4\times 44-4)-4/4.$
\item [] $688=4\times (4\times 44-4).$
\item [] $689=4\times (4\times 44-4)+4/4.$
\item [] $690=4\times (4\times 44-4)+(4+4)/4.$
\item [] $691=4\times (4\times 44-4)+4-4/4.$
\item [] $692=4\times (4\times 44-4)+4.$
\item [] $693=44\times (4^4-4)/(4\times 4).$
\item [] $694=4\times 4\times 44-(44-4)/4.$
\item [] $695=4^4+444-4-4/4.$
\item [] $696=4^4+444-4.$
\item [] $697=4^4+444-4+4/4.$
\item [] $698=4^4+444-(4+4)/4.$
\item [] $699=4^4+444-4/4.$
\item [] $700=4^4+444.$
\item [] $701=444+4^4+4/4.$
\item [] $702=4\times 4\times 44-(4+4)/4.$
\item [] $703=4\times 4\times 44-4/4.$
\item [] $704=4\times 4\times 44.$
\item [] $705=4\times 4\times 44+4/4.$
\item [] $706=4\times 4\times 44+(4+4)/4.$
\item [] $707=4\times 4\times 44+4-4/4.$
\item [] $708=4\times 4\times 44+4.$
\item [] $709=4\times 4\times 44+4+4/4.$
\item [] $710=4\times 4\times 44+4+(4+4)/4.$
\item [] $711=4\times 4\times 44+4+4-4/4.$
\item [] $712=4\times 4\times 44+4+4.$
\item [] $713=4\times 4\times 44+4+4+4/4.$
\item [] $714=4\times 4\times 44+(44-4)/4.$
\item [] $715=4\times 4\times 44+44/4.$
\item [] $716=4\times (4\times 44+4)-4.$
\item [] $717=4\times (4\times 44+4)-4+4/4.$
\item [] $718=4\times (4\times 44+4)-(4+4)/4.$
\item [] $719=4\times (4\times 44+4)-4/4.$
\item [] $720=4\times (4\times 44+4).$
\item [] $721=4\times (4\times 44+4)+4/4.$
\item [] $722=4\times (4\times 44+4)+(4+4)/4.$
\item [] $723=4\times (4\times 44+4)+4-4/4.$
\item [] $724=4\times (4\times 44+4)+4.$
\item [] $725=4\times (4\times 44+4)+4+4/4.$
\item [] $726=44\times (4^4+4+4)/(4\times 4).$
\item [] $727=4\times (4\times 44+4)+4+4-4/4.$
\item [] $728=4\times (4\times 44+4)+4+4.$
\item [] $729=(4-4/4)^{(4+(4+4)/4)}.$
\item [] $730=44\times (4^4+4+4)/(4\times 4)+4.$
\item [] $731=4\times (4\times 44+4)+44/4.$
\item [] $732=44+4\times (4\times 44-4).$
\item [] $733=(4-4/4)^{(4+(4+4)/4)}+4.$
\item [] $734=(4+4)\times (4^4+444/4)/4.$
\item [] $735=(4-4/4)\times (4^4-44/4).$
\item [] $736=4\times (4\times 44+4+4).$
\item [] $737=4\times (4\times 44+4+4)+4/4.$
\item [] $738=4^4+(44\times 44-4-4)/4.$
\item [] $739=4^4+(44\times 44-4)/4.$
\item [] $740=4\times (4\times 44+4+4)+4.$
\item [] $741=4^4+(44\times 44+4)/4.$
\item [] $742=4^4+(44\times 44+4+4)/4.$
\item [] $743=4^4+4+(44\times 44-4)/4.$
\item [] $744=44+444+4^4.$
\item [] $745=4^4+4+(44\times 44+4)/4.$
\item [] $746=44+4\times 4\times 44-(4+4)/4.$
\item [] $747=44+4\times 4\times 44-4/4.$
\item [] $748=44+4\times 4\times 44.$
\item [] $749=44+4\times 4\times 44+4/4.$
\item [] $750=44+4\times 4\times 44+(4+4)/4.$
\item [] $751=4\times (444-4^4)-4/4.$
\item [] $752=4\times (444-4^4).$
\item [] $753=4\times (444-4^4)+4/4.$
\item [] $754=4\times (444-4^4)+(4+4)/4.$
\item [] $755=(4-4/4)\times (4^4-4)-4/4.$
\item [] $756=(4-4/4)\times (4^4-4).$
\item [] $757=(4-4/4)\times (4^4-4)+4/4.$
\item [] $758=(4-4/4)\times (4^4-4)+(4+4)/4.$
\item [] $759=(4-4/4)\times (4^4-4+4/4).$
\item [] $760=(4-4/4)\times (4^4-4)+4.$
\item [] $761=(4-4/4)\times (4^4-4/4)-4.$
\item [] $762=(4-4/4)\times (4^4-(4+4)/4).$
\item [] $763=4\times 4^4-4^4-4-4/4.$
\item [] $764=4\times 4^4-4^4-4.$
\item [] $765=(4-4/4)\times (4^4-4/4).$
\item [] $766=4\times 4^4-4^4-(4+4)/4.$
\item [] $767=4\times 4^4-4^4-4/4.$
\item [] $768=4\times 4\times (44+4).$
\item [] $769=4\times 4\times (44+4)+4/4.$
\item [] $770=4\times 4\times (44+4)+(4+4)/4.$
\item [] $771=(4-4/4)\times (4^4+4/4).$
\item [] $772=4\times 4\times (44+4)+4.$
\item [] $773=4\times 4\times (44+4)+4+4/4.$
\item [] $774=(4-4/4)\times (4^4+(4+4)/4).$
\item [] $775=4+((4-4/4)\times (4/4+(4^4))).$
\item [] $776=4\times 4\times (44+4)+4+4.$
\item [] $777=(4-4/4)\times (((4^4)-4/4)+4).$
\item [] $778=(4-4/4)\times (4^4+(4+4)/4)+4.$
\item [] $779=44/4+4\times 4\times (44+4).$
\item [] $780=(4-4/4)\times (4^4+4).$
\item [] $781=(4-4/4)\times (4^4+4)+4/4.$
\item [] $782=(4-4/4)\times (4^4+4)+(4+4)/4.$
\item [] $783=(4-4/4)\times (4^4+4+4/4).$
\item [] $784=4\times (4\times (44+4)+4).$
\item [] $785=4\times (4\times (44+4)+4)+4/4.$
\item [] $786=4\times (4\times (44+4)+4)+(4+4)/4.$
\item [] $787=(4-4/4)\times (4^4+4+4/4)+4.$
\item [] $788=4\times (4\times (44+4)+4)+4.$
\item [] $789=4\times (4\times (44+4)+4)+4+4/4.$
\item [] $790=(4-4/4)\times (4^4+4)+(44-4)/4.$
\item [] $791=(4-4/4)\times (4^4+4)+44/4.$
\item [] $792=44\times (4\times 4+(4+4)/4).$
\item [] $793=44\times (4\times 4+(4+4)/4)+4/4.$
\item [] $794=444+(4+4)\times 44-(4+4)/4.$
\item [] $795=444+(4+4)\times 44-4/4.$
\item [] $796=444+(4+4)\times 44.$
\item [] $797=4\times 44-4+(4+4/4)^4.$
\item [] $798=(4-4/4)\times (4^4+(44-4)/4).$
\item [] $799=(4\times 4+4)\times (44-4)-4/4.$
\item [] $800=(4\times 4+4)\times (44-4).$
\item [] $801=4\times 44+(4+4/4)^4.$
\item [] $802=4\times 44+(4+4/4)^4+4/4.$
\item [] $803=(4\times 4+4)\times (44-4)+4-4/4.$
\item [] $804=(4\times 4+4)\times (44-4)+4.$
\item [] $805=4\times 44+4+(4+4/4)^4.$
\item [] $806=(44-4)\times (4-4/4)^4/4-4.$
\item [] $807=(4\times 4+4)\times (44-4)+4+4-4/4.$
\item [] $808=(4\times 4+4)\times (44-4)+4+4.$
\item [] $809=4\times 44+4+4+(4+4/4)^4.$
\item [] $810=(44-4)\times (4-4/4)^4/4.$
\item [] $811=4\times 4^4+44-4^4-4/4.$
\item [] $812=4\times 4\times (44+4)+44.$
\item [] $813=4\times 4\times (44+4)+44+4/4.$
\item [] $814=(44-4)\times (4-4/4)^4/4+4.$
\item [] $815=4\times 4\times 44+444/4.$
\item [] $816=4\times (4\times (44-4)+44).$
\item [] $817=4\times (44+4)+(4+4/4)^4.$
\item [] $818=(44-4)\times (4-4/4)^4/4+4+4.$
\item [] $819=((4+4)^4-4/4)/(4+4/4).$
\item [] $820=((4+4)^4+4)/(4+4/4).$
\item [] $821=4\times (44+4)+(4+4/4)^4+4.$
\item [] $822=(44\times (44+4^4)/4+4)/4-4.$
\item [] $823=((4+4)^4-4/4)/(4+4/4)+4.$
\item [] $824=((4+4)^4+4)/(4+4/4)+4.$
\item [] $825=44\times (44+4^4)/(4\times 4).$
\item [] $826=(44\times (44+4^4)/4+4)/4.$
\item [] $827=4\times 4\times (44+4+4)-4-4/4.$
\item [] $828=4\times 4\times (44+4+4)-4.$
\item [] $829=44\times (44+4^4)/(4\times 4)+4.$
\item [] $830=4\times 4\times (44+4+4)-(4+4)/4.$
\item [] $831=4\times 4\times (44+4+4)-4/4.$
\item [] $832=4\times 4\times (44+4+4).$
\item [] $833=4\times 4\times (44+4+4)+4/4.$
\item [] $834=4\times 4\times (44+4+4)+(4+4)/4.$
\item [] $835=4\times 4\times (44+4+4)+4-4/4.$
\item [] $836=4\times 4\times (44+4+4)+4.$
\item [] $837=4^4-44+(4+4/4)^4.$
\item [] $838=4\times (4^4-44)-(44-4)/4.$
\item [] $839=4444/4-4\times 4-4^4.$
\item [] $840=4\times (4^4-44)-4-4.$
\item [] $841=4^4-44+4+(4/4+4)^4.$
\item [] $842=4\times (4^4-44)-4-(4+4)/4.$
\item [] $843=4\times (4^4-44)-4-4/4.$
\item [] $844=4\times (4^4-44)-4.$
\item [] $845=4\times (4^4-44)-4+4/4.$
\item [] $846=4\times (4^4-44)-(4+4)/4.$
\item [] $847=4\times (4^4-44)-4/4.$
\item [] $848=4\times (4^4-44).$
\item [] $849=4\times (4^4-44)+4/4.$
\item [] $850=4\times (4^4-44)+(4+4)/4.$
\item [] $851=4\times (4^4-44)+4-4/4.$
\item [] $852=4\times (4^4-44)+4.$
\item [] $853=4\times (4^4-44)+4+4/4.$
\item [] $854=(4444-4)/4-4^4.$
\item [] $855=4444/4-4^4.$
\item [] $856=4\times (4^4-44)+4+4.$
\item [] $857=4\times (4^4-44)+4+4+4/4.$
\item [] $858=(4444-4)/4-4^4+4.$
\item [] $859=4444/4-4^4+4.$
\item [] $860=(4+4/4)\times (4\times 44-4).$
\item [] $861=(4+4/4)\times (4\times 44-4)+4/4.$
\item [] $862=4\times (4^4-44+4)-(4+4)/4.$
\item [] $863=4\times (4^4-44+4)-4/4.$
\item [] $864=4\times (4^4-44+4).$
\item [] $865=4\times (4^4-44+4)+4/4.$
\item [] $866=4\times (4^4-44+4)+(4+4)/4.$
\item [] $867=4\times (4^4-44+4)+4-4/4.$
\item [] $868=4\times (4^4-44+4)+4.$
\item [] $869=44\times (4\times 4+4)-44/4.$
\item [] $870=4^4+(4+4/4)^4-44/4.$
\item [] $871=4\times 4-4^4+4444/4.$
\item [] $872=44\times (4\times 4+4)-4-4.$
\item [] $873=4^4+(4+4/4)^4-4-4.$
\item [] $874=44\times (4\times 4+4)-4-(4+4)/4.$
\item [] $875=4\times 44\times (4+4/4)-4/4.$
\item [] $876=44\times (4\times 4+4)-4.$
\item [] $877=4^4+(4+4/4)^4-4.$
\item [] $878=44\times (4\times 4+4)-(4+4)/4.$
\item [] $879=44\times (4\times 4+4)-4/4.$
\item [] $880=44\times (4\times 4+4).$
\item [] $881=4^4+(4+4/4)^4.$
\item [] $882=4^4+(4+4/4)^4+4/4.$
\item [] $883=44\times (4\times 4+4)+4-4/4.$
\item [] $884=44\times (4\times 4+4)+4.$
\item [] $885=(4+4/4)^4+4^4+4.$
\item [] $886=(4+4)/4\times (444-4/4).$
\item [] $887=((4+4)\times 444-4)/4.$
\item [] $888=(4+4)\times 444/4.$
\item [] $889=((4+4)\times 444+4)/4.$
\item [] $890=(4+4)\times (444+4/4)/4.$
\item [] $891=44\times (4-4/4)^4/4.$
\item [] $892=444+444+4.$
\item [] $893=((4+4)\times 444+4)/4+4.$
\item [] $894=(4+4)\times (444+4/4)/4+4.$
\item [] $895=44\times (4-4/4)^4/4+4.$
\item [] $896=4\times (4^4-4\times (4+4)).$
\item [] $897=4\times 4+4^4+(4+4/4)^4.$
\item [] $898=4\times (4^4-4\times (4+4))+(4+4)/4.$
\item [] $899=((4+4)\times 444+44)/4.$
\item [] $900=4\times (4^4-4\times (4+4))+4.$
\item [] $901=4\times (4^4-4\times (4+4))+4+4/4.$
\item [] $902=4\times 4^4-(444+44)/4.$
\item [] $903=((4+4)\times 444+44)/4+4.$
\item [] $904=4\times (4^4-4\times (4+4))+4+4.$
\item [] $905=4\times 4+((4+4)\times 444+4)/4.$
\item [] $906=4+4\times 4^4-(444+44)/4.$
\item [] $907=4\times 4+44\times (4-4/4)^4/4.$
\item [] $908=44+4\times (4^4-44+4).$
\item [] $909=4\times 4^4-4-444/4.$
\item [] $910=4\times 4^4-4-(444-4)/4.$
\item [] $911=4\times 4^4-(444+4+4)/4.$
\item [] $912=4\times (4^4-4\times (4+4)+4).$
\item [] $913=4\times 4^4-444/4.$
\item [] $914=4\times 4^4-(444-4)/4.$
\item [] $915=4\times 4^4-(444-4-4)/4.$
\item [] $916=4\times (4^4-4\times 4)-44.$
\item [] $917=4\times 4^4-444/4+4.$
\item [] $918=4\times 4^4-(444-4)/4+4.$
\item [] $919=4444/4-4\times (44+4).$
\item [] $920=(4+4)\times (4+444/4).$
\item [] $921=4\times 4^4-444/4+4+4.$
\item [] $922=4\times 4^4-(444-4)/4+4+4.$
\item [] $923=44+44\times (4\times 4+4)-4/4.$
\item [] $924=44+(44\times ((4\times 4)+4)).$
\item [] $925=44+4^4+(4+4/4)^4.$
\item [] $926=44+4^4+(4+4/4)^4+4/4.$
\item [] $927=4\times (4^4-4)-(4-4/4)^4.$
\item [] $928=4\times (4^4-4\times 4-4-4).$
\item [] $929=4\times (4^4+4)-444/4.$
\item [] $930=4\times (4^4+4)-(444-4)/4.$
\item [] $931=4\times (4^4-4)-(4-4/4)^4+4.$
\item [] $932=4\times (4^4-4\times 4-4-4)+4.$
\item [] $933=4\times (4^4+4)+4-444/4.$
\item [] $934=(4444-4)/4-4\times 44.$
\item [] $935=44\times ((4-4/4)^4+4)/4.$
\item [] $936=4\times 4^4-44-44.$
\item [] $937=4\times 4^4-44-44+4/4.$
\item [] $938=(4^4\times 44-4-4)/(4+4+4).$
\item [] $939=4\times 4^4-(4-4/4)^4-4.$
\item [] $940=4\times (4^4-4\times 4-4)-4.$
\item [] $941=4\times (4^4-4\times 4-4)-4+4/4.$
\item [] $942=4\times 4^4-(4-4/4)^4-4/4.$
\item [] $943=4\times 4^4-(4-4/4)^4.$
\item [] $944=4\times (4^4-4\times 4-4).$
\item [] $945=4\times (4^4-4\times 4-4)+4/4.$
\item [] $946=44\times (4\times 44-4)/(4+4).$
\item [] $947=4\times 4^4-(4-4/4)^4+4.$
\item [] $948=4\times (4^4-4\times 4-4)+4.$
\item [] $949=4\times 4^4-(44+4^4)/4.$
\item [] $950=44\times (4\times 44-4)/(4+4)+4.$
\item [] $951=4\times 4^4+4+4-(4-4/4)^4.$
\item [] $952=4\times (4^4-4\times 4)-4-4.$
\item [] $953=4\times 4^4+4-(44+4^4)/4.$
\item [] $954=4\times 4^4-(4^4+4+4)/4-4.$
\item [] $955=4\times (4^4-4\times 4)-4-4/4.$
\item [] $956=4\times (4^4-4\times 4)-4.$
\item [] $957=4\times (4^4-4\times 4)-4+4/4.$
\item [] $958=4\times 4^4-(4^4+4+4)/4.$
\item [] $959=4\times (4^4-4\times 4)-4/4.$
\item [] $960=4\times (4^4-4\times 4).$
\item [] $961=4\times (4^4-4\times 4)+4/4.$
\item [] $962=4\times (4^4-4\times 4)+(4+4)/4.$
\item [] $963=4\times 4^4+4-(4^4+4)/4.$
\item [] $964=4\times (4^4-4\times 4)+4.$
\item [] $965=4\times (4^4-4\times 4)+4+4/4.$
\item [] $966=4\times (44\times 44-4)/(4+4).$
\item [] $967=(44\times (44+44)-4)/4.$
\item [] $968=4\times (4^4-4\times 4)+4+4.$
\item [] $969=4\times 4^4-44-44/4.$
\item [] $970=4\times (44\times 44+4)/(4+4).$
\item [] $971=4\times (4^4-4\times 4)+44/4.$
\item [] $972=4\times (4-4/4)^{(4+4/4)}.$
\item [] $973=4\times 4^4-44+4-44/4.$
\item [] $974=4\times (44\times 44+4)/(4+4)+4.$
\item [] $975=4\times 4^4-44-4-4/4.$
\item [] $976=4\times (4^4-4\times 4+4).$
\item [] $977=4\times (4^4-4\times 4+4)+4/4.$
\item [] $978=4\times 4^4-44-(4+4)/4.$
\item [] $979=4\times 4^4-44-4/4.$
\item [] $980=4\times 4^4-44.$
\item [] $981=4\times 4^4-44+4/4.$
\item [] $982=4\times 4^4-44+(4+4)/4.$
\item [] $983=4\times 4^4-44+4-4/4.$
\item [] $984=4\times 4^4-44+4.$
\item [] $985=4\times 4^4-44+4+4/4.$
\item [] $986=4\times 4^4-44+4+(4+4)/4.$
\item [] $987=4\times 4^4-44+4+4-4/4.$
\item [] $988=4\times (4^4-4-4)-4.$
\item [] $989=4\times (4^4-4-4)-4+4/4.$
\item [] $990=4\times (4^4-4-4)-(4+4)/4.$
\item [] $991=4\times (4^4-4-4)-4/4.$
\item [] $992=4\times (4^4-4-4).$
\item [] $993=4\times (4^4-4-4)+4/4.$
\item [] $994=4\times (4^4-4-4)+(4+4)/4.$
\item [] $995=4\times (4^4-4-4)+4-4/4.$
\item [] $996=4\times (4^4-4-4)+4.$
\item [] $997=4\times (4^4-4)-44/4.$
\item [] $998=4\times (4^4-4)-(44-4)/4.$
\item [] $999=(4+4+4/4)\times 444/4.$
\item [] $1000=4\times (4^4-4)-4-4.$
\end{itemize}
\end{multicols}
}

\section{\textbf{Representations Using Number 5}}

{\footnotesize
\begin{multicols}{3}
\begin{itemize}
\item [] $101=5\times (5\times 5-5)+5/5.$
\item [] $102=5\times (5\times 5-5)+(5+5)/5.$
\item [] $103=5\times (5\times 5-5)+5-(5+5)/5.$
\item [] $104=(5^5-5)/(5\times 5+5).$
\item [] $105=5\times (5\times 5-5)+5.$
\item [] $106=555/5-5.$
\item [] $107=(555+5)/5-5.$
\item [] $108=55+55-(5+5)/5.$
\item [] $109=55+55-5/5.$
\item [] $110=55+55.$
\item [] $111=555/5.$
\item [] $112=(555+5)/5.$
\item [] $113=(555+5+5)/5.$
\item [] $114=5\times 5\times 5-55/5.$
\item [] $115=55+55+5.$
\item [] $116=5+555/5.$
\item [] $117=5+(555+5)/5.$
\item [] $118=5+(555+5+5)/5.$
\item [] $119=5\times 5\times 5-5-5/5.$
\item [] $120=5\times 5\times 5-5.$
\item [] $121=5+5+555/5.$
\item [] $122=(555+55)/5.$
\item [] $123=5\times 5\times 5-(5+5)/5.$
\item [] $124=5\times 5\times 5-5/5.$
\item [] $125=5\times 5\times 5.$
\item [] $126=5\times 5\times 5+5/5.$
\item [] $127=5\times 5\times 5+(5+5)/5.$
\item [] $128=5\times 5\times 5+5-(5+5)/5.$
\item [] $129=5\times 5\times 5+5-5/5.$
\item [] $130=5\times 5\times 5+5.$
\item [] $131=5\times 5\times 5+5+5/5.$
\item [] $132=5\times 5\times 5+5+(5+5)/5.$
\item [] $133=5\times 5\times 5+5+5-(5+5)/5.$
\item [] $134=5\times 5\times 5+5+5-5/5.$
\item [] $135=5\times 5\times 5+5+5.$
\item [] $136=5\times 5+555/5.$
\item [] $137=5\times 5+(555+5)/5.$
\item [] $138=5\times (5\times 55+5)/(5+5).$
\item [] $139=5\times (5\times 5+5)-55/5.$
\item [] $140=5\times 5\times 5+5+5+5.$
\item [] $141=5\times 5+5+555/5.$
\item [] $142=5\times 5+5+(555+5)/5.$
\item [] $143=5+(5\times 5\times 55+5)/(5+5).$
\item [] $144=(5+5/5)\times (5\times 5-5/5).$
\item [] $145=5\times (5\times 5+5)-5.$
\item [] $146=5\times (5\times 5+5)-5+5/5.$
\item [] $147=5\times 5+(555+55)/5.$
\item [] $148=5\times (5\times 5+5)-(5+5)/5.$
\item [] $149=5\times (5\times 5+5)-5/5.$
\item [] $150=5\times (5\times 5+5).$
\item [] $151=5\times (5\times 5+5)+5/5.$
\item [] $152=5\times (5\times 5+5)+(5+5)/5.$
\item [] $153=5\times (5\times 5+5)+5-(5+5)/5.$
\item [] $154=5\times (5\times 5+5)+5-5/5.$
\item [] $155=5\times (5\times 5+5)+5.$
\item [] $156=(5^5-5)/(5\times 5-5).$
\item [] $157=5\times 5\times 5+((5+5)/5)^5.$
\item [] $158=5\times ((5+5)/5)^5-(5+5)/5.$
\item [] $159=(55+5^5)/(5\times 5-5).$
\item [] $160=5\times ((5+5)/5)^5.$
\item [] $161=5+(5^5-5)/(5\times 5-5).$
\item [] $162=5\times 5\times 5+5+((5+5)/5)^5.$
\item [] $163=55\times 5-(555+5)/5.$
\item [] $164=55\times 5-555/5.$
\item [] $165=5+5\times ((5+5)/5)^5.$
\item [] $166=55+555/5.$
\item [] $167=55+(555+5)/5.$
\item [] $168=55+(555+5+5)/5.$
\item [] $169=55\times 5+5-555/5.$
\item [] $170=5\times (5\times 5+5+5)-5.$
\item [] $171=55+5+555/5.$
\item [] $172=55+5+(555+5)/5.$
\item [] $173=5\times (5\times 5+5+5)-(5+5)/5.$
\item [] $174=5\times (5\times 5+5+5)-5/5.$
\item [] $175=5\times (5\times 5+5+5).$
\item [] $176=5\times (5\times 5+5+5)+5/5.$
\item [] $177=55+(555+55)/5.$
\item [] $178=55+5\times 5\times 5-(5+5)/5.$
\item [] $179=55+5\times 5\times 5-5/5.$
\item [] $180=55+5\times 5\times 5.$
\item [] $181=55+5\times 5\times 5+5/5.$
\item [] $182=55+5\times 5\times 5+(5+5)/5.$
\item [] $183=(5-(5+5)/5)^5-55-5.$
\item [] $184=5\times 5\times 5+55+5-5/5.$
\item [] $185=5\times 5\times 5+55+5.$
\item [] $186=5\times 5\times 5+55+5+5/5.$
\item [] $187=55\times (5+(55+5)/5)/5.$
\item [] $188=(5-(5+5)/5)^5-55.$
\item [] $189=5\times (55+5)-555/5.$
\item [] $190=5\times 5\times 5+55+5+5.$
\item [] $191=5\times 5+55+555/5.$
\item [] $192=(5+5/5)\times ((5+5)/5)^5.$
\item [] $193=(5-(5+5)/5)^5-55+5.$
\item [] $194=5\times 5\times (5+5)-55-5/5.$
\item [] $195=5\times 5\times (5+5)-55.$
\item [] $196=5\times 5\times (5+5)-55+5/5.$
\item [] $197=5+(5+5/5)\times ((5+5)/5)^5.$
\item [] $198=(5+5)\times (5\times 5-5)-(5+5)/5.$
\item [] $199=(5+5)\times (5\times 5-5)-5/5.$
\item [] $200=(5+5)\times (5\times 5-5).$
\item [] $201=(5+5)\times (5\times 5-5)+5/5.$
\item [] $202=(5+5)\times (5\times 5-5)+(5+5)/5.$
\item [] $203=(5^5-5)/(5+5+5)-5.$
\item [] $204=(5+5)\times (5\times 5-5)+5-5/5.$
\item [] $205=(5+5)\times (5\times 5-5)+5.$
\item [] $206=(5+5)\times (5\times 5-5)+5+5/5.$
\item [] $207=(5^5-5\times 5+5)/(5+5+5).$
\item [] $208=(5^5-5)/(5+5+5).$
\item [] $209=(5^5+5+5)/(5+5+5).$
\item [] $210=(5+5)\times (5\times 5-5)+5+5.$
\item [] $211=55+(5^5-5)/(5\times 5-5).$
\item [] $212=(55+5^5)/(5+5+5).$
\item [] $213=(5^5-5)/(5+5+5)+5.$
\item [] $214=(5^5+5+5)/(5+5+5)+5.$
\item [] $215=55\times 5-55-5.$
\item [] $216=(5-5/5)\times (55-5/5).$
\item [] $217=(55+5^5)/(5+5+5)+5.$
\item [] $218=(5-(5+5)/5)^5-5\times 5.$
\item [] $219=55\times 5-55-5/5.$
\item [] $220=55\times (5-5/5).$
\item [] $221=55\times (5-5/5)+5/5.$
\item [] $222=555\times (5+5)/(5\times 5).$
\item [] $223=(5-(5+5)/5)^5-5\times 5+5.$
\item [] $224=(5-5/5)\times (55+5/5).$
\item [] $225=5\times (55-5-5).$
\item [] $226=5\times (55-5-5)+5/5.$
\item [] $227=(5+5)\times 555/(5\times 5)+5.$
\item [] $228=(5-5/5)\times (55+(5+5)/5).$
\item [] $229=(5-5/5)\times (55+5/5)+5.$
\item [] $230=5\times (55-5-5)+5.$
\item [] $231=5\times (55-5-5)+5+5/5.$
\item [] $232=(5+5)\times (5+555/5)/5.$
\item [] $233=(5-(5+5)/5)^5-5-5.$
\item [] $234=(5+5-5/5)\times (5\times 5+5/5).$
\item [] $235=5\times (55-5-5)+5+5.$
\item [] $236=5\times 5\times 5+555/5.$
\item [] $237=5\times 5\times 5+(555+5)/5.$
\item [] $238=(5-(5+5)/5)^5-5.$
\item [] $239=5\times 5\times (5+5)-55/5.$
\item [] $240=(5+5)\times (5\times 5-5/5).$
\item [] $241=5\times 5\times 5+5+555/5.$
\item [] $242=(5-(5+5)/5)^5-5/5.$
\item [] $243=(5-(5+5)/5)^5.$
\item [] $244=(5-(5+5)/5)^5+5/5.$
\item [] $245=5\times 5\times (5+5)-5.$
\item [] $246=5\times 5\times (5+5)-5+5/5.$
\item [] $247=5+(5-(5+5)/5)^5-5/5.$
\item [] $248=5+(5-(5+5)/5)^5.$
\item [] $249=5\times 5\times (5+5)-5/5.$
\item [] $250=5\times 5\times (5+5).$
\item [] $251=5\times 5\times (5+5)+5/5.$
\item [] $252=5\times 5\times (5+5)+(5+5)/5.$
\item [] $253=(5-(5+5)/5)^5+5+5.$
\item [] $254=5\times 5\times (5+5)+5-5/5.$
\item [] $255=5\times 5\times (5+5)+5.$
\item [] $256=(5-5/5)^{(5-5/5)}.$
\item [] $257=(5^5-5)/(5+5)-55.$
\item [] $258=(5^5+5)/(5+5)-55.$
\item [] $259=5\times 55-5-55/5.$
\item [] $260=5\times 5\times (5+5)+5+5.$
\item [] $261=(5-5/5)^{(5-5/5)}+5.$
\item [] $262=(5^5-5)/(5+5)+5-55.$
\item [] $263=5\times 55-(55+5)/5.$
\item [] $264=5\times 55-55/5.$
\item [] $265=5\times 55-5-5.$
\item [] $266=5\times 55-5-5+5/5.$
\item [] $267=5\times 55-5-5+(5+5)/5.$
\item [] $268=5\times 5+(5-(5+5)/5)^5.$
\item [] $269=5\times 55-5-5/5.$
\item [] $270=5\times 55-5.$
\item [] $271=5\times 55-5+5/5.$
\item [] $272=5\times 55-5+(5+5)/5.$
\item [] $273=5\times 55-((5+5)/5).$
\item [] $274=5\times 55-5/5.$
\item [] $275=5\times 55.$
\item [] $276=5\times 55+5/5.$
\item [] $277=5\times 55+(5+5)/5.$
\item [] $278=5\times 55+5-(5+5)/5.$
\item [] $279=5\times 55+5-5/5.$
\item [] $280=5\times 55+5.$
\item [] $281=5\times 55+5+5/5.$
\item [] $282=5\times 55+5+(5+5)/5.$
\item [] $283=5\times 55+5+5-(5+5)/5.$
\item [] $284=5\times 55+5+5-5/5.$
\item [] $285=5\times 55+5+5.$
\item [] $286=5\times 55+55/5.$
\item [] $287=5\times 55+(55+5)/5.$
\item [] $288=(5^5+5)/(5+5)-5\times 5.$
\item [] $289=5\times (55+5)-55/5.$
\item [] $290=5\times 55+5+5+5.$
\item [] $291=5\times 55+5+55/5.$
\item [] $292=5\times 55+5+(55+5)/5.$
\item [] $293=(5^5+5)/(5+5)+5-5\times 5.$
\item [] $294=5\times (55+5)-5-5/5.$
\item [] $295=5\times (55+5)-5.$
\item [] $296=5\times (55+5)-5+5/5.$
\item [] $297=5\times 55+(55+55)/5.$
\item [] $298=55+(5-(5+5)/5)^5.$
\item [] $299=5\times (55+5)-5/5.$
\item [] $300=5\times (55+5).$
\item [] $301=5\times (55+5)+5/5.$
\item [] $302=5\times (55+5)+(5+5)/5.$
\item [] $303=(5^5+5)/(5+5)-5-5.$
\item [] $304=5\times (55+5)+5-5/5.$
\item [] $305=5\times (55+5)+5.$
\item [] $306=5\times (55+5)+5+5/5.$
\item [] $307=(5^5-55)/(5+5).$
\item [] $308=(5^5+5)/(5+5)-5.$
\item [] $309=5\times (55+5)+5+5-5/5.$
\item [] $310=5\times (55+5)+5+5.$
\item [] $311=5\times (55+5)+55/5.$
\item [] $312=(5^5-5)/(5+5).$
\item [] $313=(5^5+5)/(5+5).$
\item [] $314=(5^5+5)/(5+5)+5/5.$
\item [] $315=(5\times 5+5^5)/(5+5).$
\item [] $316=5\times (55+5)+5+55/5.$
\item [] $317=(5^5-5)/(5+5)+5.$
\item [] $318=(5^5+5)/(5+5)+5.$
\item [] $319=(5^5+5)/(5+5)+5+5/5.$
\item [] $320=(5+5)\times ((5+5)/5)^5.$
\item [] $321=(5+5)\times ((5+5)/5)^5+5/5.$
\item [] $322=(5^5-5)/(5+5)+5+5.$
\item [] $323=(5^5+5)/(5+5)+5+5.$
\item [] $324=(5+5/5)\times (55-5/5).$
\item [] $325=5\times (55+5+5).$
\item [] $326=5\times (55+5+5)+5/5.$
\item [] $327=5+5+5+(5^5-5)/(5+5).$
\item [] $328=5+5+5+(5^5+5)/(5+5).$
\item [] $329=55+5\times 55-5/5.$
\item [] $330=55+5\times 55.$
\item [] $331=55+5\times 55+5/5.$
\item [] $332=5\times 5+(5^5-55)/(5+5).$
\item [] $333=5\times 5-5+(5^5+5)/(5+5).$
\item [] $334=5\times 55+55+5-5/5.$
\item [] $335=5\times 55+55+5.$
\item [] $336=(5+5/5)\times (55+5/5).$
\item [] $337=5\times 5+(5^5-5)/(5+5).$
\item [] $338=5\times 5+(5^5+5)/(5+5).$
\item [] $339=(5^5-55)/5-5\times 55.$
\item [] $340=5\times 55+55+5+5.$
\item [] $341=5+(5+5/5)\times (55+5/5).$
\item [] $342=5\times 5+5+(5^5-5)/(5+5).$
\item [] $343=5\times 5+5+(5^5+5)/(5+5).$
\item [] $344=(5^5-5)/5-5\times 55-5.$
\item [] $345=5^5-5\times 555-5.$
\item [] $346=(5^5+5)/5-5\times 55-5.$
\item [] $347=55\times ((5+5)/5)^5/5-5.$
\item [] $348=(5^5-5-5)/5-5\times 55.$
\item [] $349=(5^5-5)/5-5\times 55.$
\item [] $350=5\times (55+5+5+5).$
\item [] $351=(5^5+5)/5-5\times 55.$
\item [] $352=55\times ((5+5)/5)^5/5.$
\item [] $353=(55\times ((5+5)/5)^5+5)/5.$
\item [] $354=(5^5-5)/5+5-5\times 55.$
\item [] $355=55+5\times (55+5).$
\item [] $356=(5^5+5)/5+5-5\times 55.$
\item [] $357=55\times ((5+5)/5)^5/5+5.$
\item [] $358=(55\times ((5+5)/5)^5+5)/5+5.$
\item [] $359=(5+5/5)\times (55+5)-5/5.$
\item [] $360=(5+5/5)\times (55+5).$
\item [] $361=(5+5/5)\times (55+5)+5/5.$
\item [] $362=55+(5^5-55)/(5+5).$
\item [] $363=55-5+(5^5+5)/(5+5).$
\item [] $364=5\times 5\times (5+5+5)-55/5.$
\item [] $365=(5+5/5)\times (55+5)+5.$
\item [] $366=(5+5/5)\times (55+5+5/5).$
\item [] $367=55+(5^5-5)/(5+5).$
\item [] $368=55+(5^5+5)/(5+5).$
\item [] $369=5\times 5\times (5+5+5)-5-5/5.$
\item [] $370=5\times 5\times (5+5+5)-5.$
\item [] $371=5\times 5\times (5+5+5)-5+5/5.$
\item [] $372=55+5+(5^5-5)/(5+5).$
\item [] $373=55+5+(5^5+5)/(5+5).$
\item [] $374=5\times 5\times (5+5+5)-5/5.$
\item [] $375=5\times 5\times (5+5+5).$
\item [] $376=5\times 5\times (5+5+5)+5/5.$
\item [] $377=5\times 5\times (5+5+5)+(5+5)/5.$
\item [] $378=(5+(5+5)/5)\times (55-5/5).$
\item [] $379=5\times 5\times (5+5+5)+5-5/5.$
\item [] $380=5\times 5\times (5+5+5)+5.$
\item [] $381=5\times 5\times (5+5+5)+5+5/5.$
\item [] $382=5^5/5-(5-(5+5)/5)^5.$
\item [] $383=(5+(5+5)/5)\times (55-5/5)+5.$
\item [] $384=(55+5)\times ((5+5)/5)^5/5.$
\item [] $385=55\times (5+(5+5)/5).$
\item [] $386=5\times 55+555/5.$
\item [] $387=5\times 55+(555+5)/5.$
\item [] $388=5\times 55+(555+5+5)/5.$
\item [] $389=5\times (5\times 5+55)-55/5.$
\item [] $390=55\times (5+(5+5)/5)+5.$
\item [] $391=5\times 55+5+555/5.$
\item [] $392=(5+(5+5)/5)\times (55+5/5).$
\item [] $393=5\times 5+55+(5^5+5)/(5+5).$
\item [] $394=(5-5/5)^5-5^5/5-5.$
\item [] $395=5\times (5\times 5+55)-5.$
\item [] $396=5\times (5\times 5+55)-5+5/5.$
\item [] $397=(5+(5+5)/5)\times (55+5/5)+5.$
\item [] $398=(5-5/5)^5-(5^5+5)/5.$
\item [] $399=(5-5/5)^5-5^5/5.$
\item [] $400=5\times (5\times 5+55).$
\item [] $401=5\times (5\times 5+55)+5/5.$
\item [] $402=5\times (5\times 5+55)+(5+5)/5.$
\item [] $403=5+(5-5/5)^5-(5^5+5)/5.$
\item [] $404=5+(5-5/5)^5-5^5/5.$
\item [] $405=5\times (5\times 5+55)+5.$
\item [] $406=5\times (5\times 5+55)+5+5/5.$
\item [] $407=55\times (((5+5)/5)^5+5)/5.$
\item [] $408=(5\times 5-5/5)\times (5+(55+5)/5).$
\item [] $409=5+5+(5-5/5)^5-5^5/5.$
\item [] $410=5\times (5\times 5+55)+5+5.$
\item [] $411=5\times (5\times 5+55)+55/5.$
\item [] $412=55\times (5+((5+5)/5)^5)/5+5.$
\item [] $413=5\times (5\times 5-5)+(5^5+5)/(5+5).$
\item [] $414=(5-5/5)^5-555-55.$
\item [] $415=5\times (5\times 5+55)+5+5+5.$
\item [] $416=(5+55/5)\times (5\times 5+5/5).$
\item [] $417=(5^5+5^5+5)/(5+5+5).$
\item [] $418=55+55-5+(5^5+5)/(5+5).$
\item [] $419=555-5\times 5-555/5.$
\item [] $420=(55+5)\times (5+(5+5)/5).$
\item [] $421=5+(5+55/5)\times (5\times 5+5/5).$
\item [] $422=5+(5^5+5^5+5)/(5+5+5).$
\item [] $423=55+55+(5^5+5)/(5+5).$
\item [] $424=(5-5/5)\times (555/5-5).$
\item [] $425=5\times (5\times 5+55+5).$
\item [] $426=5\times (5\times 5+55+5)+5/5.$
\item [] $427=5\times (5\times 5+55+5)+(5+5)/5.$
\item [] $428=(5-5/5)\times ((555+5)/5-5).$
\item [] $429=555-5\times 5\times 5-5/5.$
\item [] $430=555-5\times 5\times 5.$
\item [] $431=555-5\times 5\times 5+5/5.$
\item [] $432=5\times 5\times 5+(5^5-55)/(5+5).$
\item [] $433=5\times 5\times 5-5+(5^5+5)/(5+5).$
\item [] $434=555-5\times 5\times 5+5-5/5.$
\item [] $435=555-5\times 5\times 5+5.$
\item [] $436=555-5\times 5\times 5+5+5/5.$
\item [] $437=5\times 5\times 5+(5^5-5)/(5+5).$
\item [] $438=5\times 5\times 5+(5^5+5)/(5+5).$
\item [] $439=555-5-555/5.$
\item [] $440=(5+5)\times (55-55/5).$
\item [] $441=(5+5+55/5)^{((5+5)/5)}.$
\item [] $442=5\times 5\times 5+5+(5^5-5)/(5+5).$
\item [] $443=555-(555+5)/5.$
\item [] $444=(5-5/5)\times 555/5.$
\item [] $445=555-55-55.$
\item [] $446=5+(55/5+5+5)^{((5+5)/5)}.$
\item [] $447=(5^5+5-5/5)/(5+(5+5)/5).$
\item [] $448=(5-5/5)\times (555+5)/5.$
\item [] $449=5+(5-5/5)\times 555/5.$
\item [] $450=(5+5)\times (55-5-5).$
\item [] $451=(5+5)\times (55-5-5)+5/5.$
\item [] $452=(5+5)\times (55-5-5)+(5+5)/5.$
\item [] $453=5+(5-5/5)\times (555+5)/5.$
\item [] $454=5+5+(5-5/5)\times 555/5.$
\item [] $455=(5+5)\times (55-5-5)+5.$
\item [] $456=(5+5)\times (55-5-5)+5+5/5.$
\item [] $457=((5+5)/5)^{(5+5-5/5)}-55.$
\item [] $458=(5\times 5\times 5\times 55-5)/(5+5+5).$
\item [] $459=(5^5-555)/5-55.$
\item [] $460=(5+5)\times (55-5-5)+5+5.$
\item [] $461=(5+5)\times (55-5-5)+55/5.$
\item [] $462=55\times (((5+5)/5)^5+5+5)/5.$
\item [] $463=5\times (5\times 5+5)+(5^5+5)/(5+5).$
\item [] $464=(5-5/5)\times (5+555/5).$
\item [] $465=5\times (5\times 5\times 5-((5+5)/5)^5).$
\item [] $466=(5^5+5)/5-5\times ((5+5)/5)^5.$
\item [] $467=5\times (5\times 5+5)+5+(5^5-5)/(5+5).$
\item [] $468=(5-5/5)\times (5+(555+5)/5).$
\item [] $469=(5-5/5)^5-555.$
\item [] $470=5\times (5\times (5\times 5-5)-5)-5.$
\item [] $471=5\times (5\times (5\times 5-5)-5)-5+5/5.$
\item [] $472=55+(5^5+5^5+5)/(5+5+5).$
\item [] $473=55\times (55-(55+5)/5)/5.$
\item [] $474=5+(5-5/5)^5-555.$
\item [] $475=5\times (5\times (5\times 5-5)-5).$
\item [] $476=5\times (5\times (5\times 5-5)-5)+5/5.$
\item [] $477=5\times (5\times (5\times 5-5)-5)+(5+5)/5.$
\item [] $478=5+55\times (55-(55+5)/5)/5.$
\item [] $479=(5-5/5)^5+5+5-555.$
\item [] $480=5\times (5\times (5\times 5-5)-5)+5.$
\item [] $481=(5555-5^5)/5-5.$
\item [] $482=(5555-5^5+5)/5-5.$
\item [] $483=5\times 55+(5^5-5)/(5+5+5).$
\item [] $484=55\times (55-55/5)/5.$
\item [] $485=5\times (5\times (5\times 5-5)-5)+5+5.$
\item [] $486=(5555-5^5)/5.$
\item [] $487=(5555+5-5^5)/5.$
\item [] $488=(5-5/5)\times (555+55)/5.$
\item [] $489=555-55-55/5.$
\item [] $490=555-55-5-5.$
\item [] $491=5+(5555-5^5)/5.$
\item [] $492=5+(5555-5^5+5)/5.$
\item [] $493=555-55-5-(5+5)/5.$
\item [] $494=555-55-5-5/5.$
\item [] $495=555-55-5.$
\item [] $496=555-55-5+5/5.$
\item [] $497=555-55-5+(5+5)/5.$
\item [] $498=555-55-(5+5)/5.$
\item [] $499=555-55-5/5.$
\item [] $500=555-55.$
\item [] $501=5\times 5\times (5\times 5-5)+5/5.$
\item [] $502=555-55+(5+5)/5.$
\item [] $503=(5^5-555-55)/5.$
\item [] $504=555-55+5-5/5.$
\item [] $505=5\times 5\times (5\times 5-5)+5.$
\item [] $506=5\times 5\times (5\times 5-5)+5+5/5.$
\item [] $507=((5+5)/5)^{(5+5-5/5)}-5.$
\item [] $508=(5^5-555-5)/5-5.$
\item [] $509=(5^5-555)/5-5.$
\item [] $510=5\times 5\times (5\times 5-5)+5+5.$
\item [] $511=555-55+55/5.$
\item [] $512=((5+5)/5)^{(5+5-5/5)}.$
\item [] $513=(5^5-555-5)/5.$
\item [] $514=(5^5-555)/5.$
\item [] $515=5^5/5-55-55.$
\item [] $516=(5^5+5)/5-55-55.$
\item [] $517=5+((5+5)/5)^{(5+5-5/5)}.$
\item [] $518=5+(5^5-555-5)/5.$
\item [] $519=5+(5^5-555)/5.$
\item [] $520=(5^5-5)/(5+5/5).$
\item [] $521=(5^5+5/5)/(5+5/5).$
\item [] $522=5+5+((5+5)/5)^{(5+5-5/5)}.$
\item [] $523=555-((5+5)/5)^5.$
\item [] $524=5+5+(5^5-555)/5.$
\item [] $525=5\times (5\times (5\times 5-5)+5).$
\item [] $526=5+(5^5+5/5)/(5+5/5).$
\item [] $527=555-5\times 5-5+(5+5)/5.$
\item [] $528=555+5-((5+5)/5)^5.$
\item [] $529=555-5\times 5-5/5.$
\item [] $530=555-5\times 5.$
\item [] $531=555-5\times 5+5/5.$
\item [] $532=555-5\times 5+(5+5)/5.$
\item [] $533=555-(55+55)/5.$
\item [] $534=555-5\times 5+5-5/5.$
\item [] $535=555-5\times 5+5.$
\item [] $536=555-5\times 5+5+5/5.$
\item [] $537=555-5\times 5+5+(5+5)/5.$
\item [] $538=555-5-(55+5)/5.$
\item [] $539=555-5-55/5.$
\item [] $540=(5+5)\times (55-5/5).$
\item [] $541=(5+5)\times (55-5/5)+5/5.$
\item [] $542=555-((55+5+5)/5).$
\item [] $543=555-(55+5)/5.$
\item [] $544=555-55/5.$
\item [] $545=555-5-5.$
\item [] $546=555-5-5+5/5.$
\item [] $547=555-5-5+(5+5)/5.$
\item [] $548=555-5-(5+5)/5.$
\item [] $549=555-5-5/5.$
\item [] $550=555-5.$
\item [] $551=555-5+5/5.$
\item [] $552=555-5+(5+5)/5.$
\item [] $553=555-(5+5)/5.$
\item [] $554=555-5/5.$
\item [] $555=555.$
\item [] $556=555+5/5.$
\item [] $557=555+(5+5)/5.$
\item [] $558=555+5-(5+5)/5.$
\item [] $559=555+5-5/5.$
\item [] $560=555+5.$
\item [] $561=555+5+5/5.$
\item [] $562=555+5+(5+5)/5.$
\item [] $563=555+5+5-(5+5)/5.$
\item [] $564=555+5+5-5/5.$
\item [] $565=555+5+5.$
\item [] $566=555+55/5.$
\item [] $567=555+(55+5)/5.$
\item [] $568=(5^5-5-5)/5-55.$
\item [] $569=(5^5-5)/5-55.$
\item [] $570=5^5/5-55.$
\item [] $571=(5^5+5)/5-55.$
\item [] $572=(5^5+5+5)/5-55.$
\item [] $573=5-55+(5^5-5-5)/5.$
\item [] $574=5-55+(5^5-5)/5.$
\item [] $575=5-55+5^5/5.$
\item [] $576=5-55+(5^5+5)/5.$
\item [] $577=5-55+(5^5+5+5)/5.$
\item [] $578=5\times 5+555-(5+5)/5.$
\item [] $579=5\times 5+555-5/5.$
\item [] $580=5\times 5+555.$
\item [] $581=5\times 5+555+5/5.$
\item [] $582=5\times 5+555+(5+5)/5.$
\item [] $583=55\times (55-(5+5)/5)/5.$
\item [] $584=555+5\times 5+5-5/5.$
\item [] $585=555+5\times 5+5.$
\item [] $586=555+5\times 5+5+5/5.$
\item [] $587=555+((5+5)/5)^5.$
\item [] $588=5\times 55+(5^5+5)/(5+5).$
\item [] $589=(5^5-55)/5-5\times 5.$
\item [] $590=555+5\times 5+5+5.$
\item [] $591=555+5\times 5+55/5.$
\item [] $592=555+5+((5+5)/5)^5.$
\item [] $593=5^5/5-((5+5)/5)^5.$
\item [] $594=55\times (55-5/5)/5.$
\item [] $595=5^5/5-5\times 5-5.$
\item [] $596=(5^5+5)/5-5\times 5-5.$
\item [] $597=(5^5+5+5)/5-5\times 5-5.$
\item [] $598=5\times (5\times 5\times 5-5)-(5+5)/5.$
\item [] $599=5\times (5\times 5\times 5-5)-5/5.$
\item [] $600=5\times (5\times 5\times 5-5).$
\item [] $601=(5^5+5)/5-5\times 5.$
\item [] $602=(5^5+5+5)/5-5\times 5.$
\item [] $603=(5^5-55-55)/5.$
\item [] $604=(55\times 55-5)/5.$
\item [] $605=55\times 55/5.$
\item [] $606=(55\times 55+5)/5.$
\item [] $607=(55\times 55+5+5)/5.$
\item [] $608=(5^5-55-5)/5-5.$
\item [] $609=(5^5-55)/5-5.$
\item [] $610=555+55.$
\item [] $611=5+(55\times 55+5)/5.$
\item [] $612=(5^5-55-5-5)/5.$
\item [] $613=(5^5-55-5)/5.$
\item [] $614=(5^5-55)/5.$
\item [] $615=5^5/5-5-5.$
\item [] $616=(5^5+5)/5-5-5.$
\item [] $617=(5^5+5+5)/5-5-5.$
\item [] $618=(5^5-5-5)/5-5.$
\item [] $619=(5^5-5)/5-5.$
\item [] $620=5^5/5-5.$
\item [] $621=(5^5+5)/5-5.$
\item [] $622=(5^5+5+5)/5-5.$
\item [] $623=(5^5-5-5)/5.$
\item [] $624=(5^5-5)/5.$
\item [] $625=5^5/5.$
\item [] $626=(5^5+5)/5.$
\item [] $627=(5^5+5+5)/5.$
\item [] $628=5+(5^5-5-5)/5.$
\item [] $629=5+(5^5-5)/5.$
\item [] $630=5+5^5/5.$
\item [] $631=5+(5^5+5)/5.$
\item [] $632=5+(5^5+5+5)/5.$
\item [] $633=5+5+(5^5-5-5)/5.$
\item [] $634=5+5+(5^5-5)/5.$
\item [] $635=5+5+5^5/5.$
\item [] $636=(55+5^5)/5.$
\item [] $637=(55+5^5+5)/5.$
\item [] $638=(55+5^5+5+5)/5.$
\item [] $639=5+5+5+(5^5-5)/5.$
\item [] $640=5+5+5+5^5/5.$
\item [] $641=5+(55+5^5)/5.$
\item [] $642=5+(55+5^5+5)/5.$
\item [] $643=5+(55+5^5+5+5)/5.$
\item [] $644=5\times 5-5+(5^5-5)/5.$
\item [] $645=5\times 5-5+5^5/5.$
\item [] $646=5+5+(55+5^5)/5.$
\item [] $647=(55+55+5^5)/5.$
\item [] $648=5\times 5+(5^5-5-5)/5.$
\item [] $649=5\times 5+(5^5-5)/5.$
\item [] $650=5\times (5\times 5\times 5+5).$
\item [] $651=5\times 5+(5^5+5)/5.$
\item [] $652=5\times 5+(5^5+5+5)/5.$
\item [] $653=5+5\times 5+(5^5-5-5)/5.$
\item [] $654=5+5\times 5+(5^5-5)/5.$
\item [] $655=5+5\times 5+5^5/5.$
\item [] $656=5+5\times 5+(5^5+5)/5.$
\item [] $657=5^5/5+((5+5)/5)^5.$
\item [] $658=(5^5+5)/5+((5+5)/5)^5.$
\item [] $659=(55\times (55+5)-5)/5.$
\item [] $660=55\times (55+5)/5.$
\item [] $661=5\times 5+(55+5^5)/5.$
\item [] $662=5+((5+5)/5)^5+5^5/5.$
\item [] $663=5+(5^5+5)/5+((5+5)/5)^5.$
\item [] $664=5+(55\times (55+5)-5)/5.$
\item [] $665=5+55\times (55+5)/5.$
\item [] $666=555+555/5.$
\item [] $667=555+(555+5)/5.$
\item [] $668=55+(5^5-55-5)/5.$
\item [] $669=55+(5^5-55)/5.$
\item [] $670=55-5-5+5^5/5.$
\item [] $671=5+555+555/5.$
\item [] $672=(5+5/5)\times (555+5)/5.$
\item [] $673=55-5+(5^5-5-5)/5.$
\item [] $674=55-5+(5^5-5)/5.$
\item [] $675=5\times (5\times 5\times 5+5+5).$
\item [] $676=55-5+(5^5+5)/5.$
\item [] $677=55-5+(5^5+5+5)/5.$
\item [] $678=55+(5^5-5-5)/5.$
\item [] $679=55+(5^5-5)/5.$
\item [] $680=55+5^5/5.$
\item [] $681=55+(5^5+5)/5.$
\item [] $682=55+(5^5+5+5)/5.$
\item [] $683=55+5+(5^5-5-5)/5.$
\item [] $684=55+5+(5^5-5)/5.$
\item [] $685=55+5+5^5/5.$
\item [] $686=55+5+(5^5+5)/5.$
\item [] $687=55+5+(5^5+5+5)/5.$
\item [] $688=(55\times 5\times 5\times 5+5)/(5+5).$
\item [] $689=55+5+5+(5^5-5)/5.$
\item [] $690=55+5+5+5^5/5.$
\item [] $691=55+(55+5^5)/5.$
\item [] $692=55+(55+5^5+5)/5.$
\item [] $693=5+(5\times 5\times 5\times 55+5)/(5+5).$
\item [] $694=55+5+5+5+(5^5-5)/5.$
\item [] $695=5\times 5\times (5\times 5+5)-55.$
\item [] $696=55+5+(55+5^5)/5.$
\item [] $697=55+5+(55+5^5+5)/5.$
\item [] $698=5\times (5+5+5)+(5^5-5-5)/5.$
\item [] $699=5\times (5+5+5)+(5^5-5)/5.$
\item [] $700=5\times (5\times 5\times 5+5+5+5).$
\item [] $701=5\times (5+5+5)+(5^5+5)/5.$
\item [] $702=55+(55+55+5^5)/5.$
\item [] $703=5\times 5+55+(5^5-5-5)/5.$
\item [] $704=5\times 5+55+(5^5-5)/5.$
\item [] $705=5\times 5+55+5^5/5.$
\item [] $706=5\times 5+55+(5^5+5)/5.$
\item [] $707=5\times 5+55+(5^5+5+5)/5.$
\item [] $708=(55+5)\times (55+5-5/5)/5.$
\item [] $709=5\times 5+55+5+(5^5-5)/5.$
\item [] $710=5\times 5+55+5+5^5/5.$
\item [] $711=(555+5^5)/5-5\times 5.$
\item [] $712=55+((5+5)/5)^5+5^5/5.$
\item [] $713=(55\times (55+5+5)-5-5)/5.$
\item [] $714=(55\times (55+5+5)-5)/5.$
\item [] $715=55\times (55+5+5)/5.$
\item [] $716=(55\times (55+5+5)+5)/5.$
\item [] $717=(55\times (55+5+5)+5+5)/5.$
\item [] $718=5\times 5\times (5\times 5+5)-((5+5)/5)^5.$
\item [] $719=((55+5)^{((5+5)/5)}-5)/5.$
\item [] $720=(5\times 5\times 5-5)\times (5+5/5).$
\item [] $721=((55+5)^{((5+5)/5)}+5)/5.$
\item [] $722=((55+5)^{((5+5)/5)}+5+5)/5.$
\item [] $723=5\times (5\times 5-5)+(5^5-5-5)/5.$
\item [] $724=5\times (5\times 5-5)+(5^5-5)/5.$
\item [] $725=5\times (5\times (5\times 5+5)-5).$
\item [] $726=5\times (5\times (5\times 5+5)-5)+5/5.$
\item [] $727=(555+5^5+5)/5-5-5.$
\item [] $728=(55+5/5)\times (55+5+5)/5.$
\item [] $729=(5-(5+5)/5)^{(5+5/5)}.$
\item [] $730=5+5\times (5\times (5\times 5+5)-5).$
\item [] $731=(555+5^5)/5-5.$
\item [] $732=(555+5^5+5)/5-5.$
\item [] $733=55+55+(5^5-5-5)/5.$
\item [] $734=5+(5-(5+5)/5)^{(5+5/5)}.$
\item [] $735=55+55+5^5/5.$
\item [] $736=(555+5^5)/5.$
\item [] $737=(555+5^5+5)/5.$
\item [] $738=(555+5^5+5+5)/5.$
\item [] $739=5\times 5\times 5+(5^5-55)/5.$
\item [] $740=5\times 5\times (5\times 5+5)-5-5.$
\item [] $741=5+(555+5^5)/5.$
\item [] $742=5+(555+5^5+5)/5.$
\item [] $743=5+(555+5^5+5+5)/5.$
\item [] $744=5\times 5\times 5-5+(5^5-5)/5.$
\item [] $745=5\times 5\times (5\times 5+5)-5.$
\item [] $746=5+5+(555+5^5)/5.$
\item [] $747=(555+55+5^5)/5.$
\item [] $748=5\times 5\times 5+(5^5-5-5)/5.$
\item [] $749=5\times 5\times 5+(5^5-5)/5.$
\item [] $750=5\times 5\times (5\times 5+5).$
\item [] $751=5\times 5\times (5\times 5+5)+5/5.$
\item [] $752=5\times 5\times 5+(5^5+5+5)/5.$
\item [] $753=5\times 5\times 5+5+(5^5-5-5)/5.$
\item [] $754=5\times 5\times 5+5+(5^5-5)/5.$
\item [] $755=5+5\times 5\times (5\times 5+5).$
\item [] $756=5+5\times 5\times (5\times 5+5)+5/5.$
\item [] $757=5+5\times 5\times (5\times 5+5)+(5+5)/5.$
\item [] $758=5\times 5\times 5+5+5+(5^5-5-5)/5.$
\item [] $759=5\times 5\times 5+5+5+(5^5-5)/5.$
\item [] $760=5+5+5\times 5\times (5\times 5+5).$
\item [] $761=5\times 5+(555+5^5)/5.$
\item [] $762=5\times 5+(555+5^5+5)/5.$
\item [] $763=555+(5^5-5)/(5+5+5).$
\item [] $764=5\times (5\times 5+5)+(5^5-55)/5.$
\item [] $765=5+5+5+5\times 5\times (5\times 5+5).$
\item [] $766=5+5\times 5+(555+5^5)/5.$
\item [] $767=5+5\times 5+(555+5^5+5)/5.$
\item [] $768=(5\times 5-5/5)\times ((5+5)/5)^5.$
\item [] $769=(5-5/5)^5-5\times (55-5)-5.$
\item [] $770=55\times (5+5+5-5/5).$
\item [] $771=55\times (5+5+5-5/5)+5/5.$
\item [] $772=555\times (5+(5+5)/5)/5-5.$
\item [] $773=5+(5\times 5-5/5)\times ((5+5)/5)^5.$
\item [] $774=(5-5/5)^5-5\times (55-5).$
\item [] $775=5\times (5\times (5\times 5+5)+5).$
\item [] $776=5\times (5\times (5\times 5+5)+5)+5/5.$
\item [] $777=555\times (5+(5+5)/5)/5.$
\item [] $778=(5^5-5)/(5-5/5)-(5+5)/5.$
\item [] $779=(5^5-5)/(5-5/5)-5/5.$
\item [] $780=(5^5-5)/(5-5/5).$
\item [] $781=(5^5-5/5)/(5-5/5).$
\item [] $782=5+555\times (5+(5+5)/5)/5.$
\item [] $783=(5^5-5/5)/(5-5/5)+(5+5)/5.$
\item [] $784=(5^5+55/5)/(5-5/5).$
\item [] $785=5+(5^5-5)/(5-5/5).$
\item [] $786=5+(5^5-5/5)/(5-5/5).$
\item [] $787=5+5+555\times (5+(5+5)/5)/5.$
\item [] $788=5+(5^5+5+(5+5)/5)/(5-5/5).$
\item [] $789=5+(5^5+55/5)/(5-5/5).$
\item [] $790=5+5+(5^5-5)/(5-5/5).$
\item [] $791=55+(555+5^5)/5.$
\item [] $792=55+(555+5^5+5)/5.$
\item [] $793=555-5+(5-(5+5)/5)^5.$
\item [] $794=5\times 5\times ((5+5)/5)^5-5-5/5.$
\item [] $795=(55+5^5)/(5-5/5).$
\item [] $796=55+5+(555+5^5)/5.$
\item [] $797=55+5+(555+5^5+5)/5.$
\item [] $798=555+(5-(5+5)/5)^5.$
\item [] $799=5\times 5\times ((5+5)/5)^5-5/5.$
\item [] $800=5\times 5\times ((5+5)/5)^5.$
\item [] $801=5\times 5\times ((5+5)/5)^5+5/5.$
\item [] $802=5\times 5\times ((5+5)/5)^5+(5+5)/5.$
\item [] $803=5+555+(5-(5+5)/5)^5.$
\item [] $804=5+5\times 5\times ((5+5)/5)^5-5/5.$
\item [] $805=5+5\times 5\times ((5+5)/5)^5.$
\item [] $806=5+5\times 5\times ((5+5)/5)^5+5/5.$
\item [] $807=5+5\times 5\times ((5+5)/5)^5+(5+5)/5.$
\item [] $808=555+5+5+(5-(5+5)/5)^5.$
\item [] $809=55\times (5+5+5)-5-55/5.$
\item [] $810=(5+5+5)\times (55-5/5).$
\item [] $811=(5+5+5)\times (55-5/5)+5/5.$
\item [] $812=(5+(5+5)/5)\times (5+555/5).$
\item [] $813=55\times (5+5+5)-(55+5)/5.$
\item [] $814=55\times (5+5+5)-55/5.$
\item [] $815=55\times (5+5+5)-5-5.$
\item [] $816=55\times (5+5+5)-5-5+5/5.$
\item [] $817=55\times (5+5+5)-5-5+(5+5)/5.$
\item [] $818=55\times (5+5+5)-5-(5+5)/5.$
\item [] $819=55\times (5+5+5)-5-5/5.$
\item [] $820=55\times (5+5+5)-5.$
\item [] $821=55\times (5+5+5)-5+5/5.$
\item [] $822=55\times (5+5+5)-5+(5+5)/5.$
\item [] $823=55\times (5+5+5)-(5+5)/5.$
\item [] $824=55\times (5+5+5)-5/5.$
\item [] $825=55\times (5+5+5).$
\item [] $826=55\times (5+5+5)+5/5.$
\item [] $827=55\times (5+5+5)+(5+5)/5.$
\item [] $828=55\times (5+5+5)+5-(5+5)/5.$
\item [] $829=55\times (5+5+5)+5-5/5.$
\item [] $830=55\times (5+5+5)+5.$
\item [] $831=55\times (5+5+5)+5+5/5.$
\item [] $832=(5\times 5+5/5)\times ((5+5)/5)^5.$
\item [] $833=5^5/5+(5^5-5)/(5+5+5).$
\item [] $834=55\times (5+5+5)+5+5-5/5.$
\item [] $835=55\times (5+5+5)+5+5.$
\item [] $836=55\times (5+5+5)+55/5.$
\item [] $837=55\times (5+5+5)+(55+5)/5.$
\item [] $838=5+5^5/5+(5^5-5)/(5+5+5).$
\item [] $839=(5+5+5)\times (55+5/5)-5/5.$
\item [] $840=(5+5+5)\times (55+5/5).$
\item [] $841=(5\times 5+5-5/5)^{((5+5)/5)}.$
\item [] $842=5+5+(5\times 5+5/5)\times ((5+5)/5)^5.$
\item [] $843=5^5/5-5\times 5+(5-(5+5)/5)^5.$
\item [] $844=5\times 55-55+(5^5-5)/5.$
\item [] $845=(5+5+5)\times (55+5/5)+5.$
\item [] $846=(5\times 5+5-5/5)^{((5+5)/5)}+5.$
\item [] $847=(555+555+5^5)/5.$
\item [] $848=(5+55/5)\times (55-(5+5)/5).$
\item [] $849=5\times 5+55\times (5+5+5)-5/5.$
\item [] $850=5\times (5\times (5\times 5+5+5)-5).$
\item [] $851=5\times 5+55\times (5+5+5)+5/5.$
\item [] $852=5+(555+555+5^5)/5.$
\item [] $853=5+(55/5+5)\times (55-(5+5)/5).$
\item [] $854=555+5\times (55+5)-5/5.$
\item [] $855=555+5\times (55+5).$
\item [] $856=555+5\times (55+5)+5/5.$
\item [] $857=55\times (5+5+5)+((5+5)/5)^5.$
\item [] $858=5^5/5+(5-(5+5)/5)^5-5-5.$
\item [] $859=(5+55/5)\times (55-5/5)-5.$
\item [] $860=5+5\times (55+5)+555.$
\item [] $861=5\times 5\times 5+(555+5^5)/5.$
\item [] $862=555+(5^5-55)/(5+5).$
\item [] $863=5^5/5+(5-(5+5)/5)^5-5.$
\item [] $864=(5+55/5)\times (55-5/5).$
\item [] $865=5\times 5\times (5\times 5+5+5)-5-5.$
\item [] $866=5\times 5\times 5+5+(555+5^5)/5.$
\item [] $867=555+(5^5-5)/(5+5).$
\item [] $868=5^5/5+(5-(5+5)/5)^5.$
\item [] $869=5+(5+55/5)\times (55-5/5).$
\item [] $870=5\times 5\times (5\times 5+5+5)-5.$
\item [] $871=5\times 5\times (5\times 5+5+5)-5+5/5.$
\item [] $872=555+5+(5^5-5)/(5+5).$
\item [] $873=5^5/5+5+(5-(5+5)/5)^5.$
\item [] $874=5\times 5\times (5+5)+(5^5-5)/5.$
\item [] $875=5\times 5\times (5\times 5+5+5).$
\item [] $876=5\times 5\times (5\times 5+5+5)+5/5.$
\item [] $877=5\times 5\times (5\times 5+5+5)+(5+5)/5.$
\item [] $878=55\times (5+55/5)-(5+5)/5.$
\item [] $879=55\times (5+55/5)-5/5.$
\item [] $880=55\times (5+55/5).$
\item [] $881=55\times (5+55/5)+5/5.$
\item [] $882=55\times (5+55/5)+(5+5)/5.$
\item [] $883=55\times (5+55/5)+5-(5+5)/5.$
\item [] $884=55\times (5+55/5)+5-5/5.$
\item [] $885=55\times (5+55/5)+5.$
\item [] $886=55\times (5+55/5)+5+5/5.$
\item [] $887=55\times (5+55/5)+5+(5+5)/5.$
\item [] $888=5\times 55+(5^5-55-5)/5.$
\item [] $889=5\times 55+(5^5-55)/5.$
\item [] $890=5+5+55\times (5+55/5).$
\item [] $891=5\times 55-5-5+(5^5+5)/5.$
\item [] $892=5\times 5+555+(5^5-5)/(5+5).$
\item [] $893=5\times 55-5+(5^5-5-5)/5.$
\item [] $894=5\times 55-5+(5^5-5)/5.$
\item [] $895=5\times 55-5+5^5/5.$
\item [] $896=5\times 55-5+(5^5+5)/5.$
\item [] $897=5\times 55-5+(5^5+5+5)/5.$
\item [] $898=5\times 55+(5^5-5-5)/5.$
\item [] $899=5\times 55+(5^5-5)/5.$
\item [] $900=5\times (5\times 5\times 5+55).$
\item [] $901=5\times 55+(5^5+5)/5.$
\item [] $902=5\times 55+(5^5+5+5)/5.$
\item [] $903=5\times 55+5+(5^5-5-5)/5.$
\item [] $904=5\times 55+5+(5^5-5)/5.$
\item [] $905=5\times 55+5+5^5/5.$
\item [] $906=5\times 55+5+(5^5+5)/5.$
\item [] $907=5\times 55+5+(5^5+5+5)/5.$
\item [] $908=(5-5/5)^5-5-555/5.$
\item [] $909=5\times 55+5+5+(5^5-5)/5.$
\item [] $910=5\times 55+5+5+5^5/5.$
\item [] $911=5\times 55+(55+5^5)/5.$
\item [] $912=5\times 55+(55+5^5+5)/5.$
\item [] $913=(5-5/5)^5-555/5.$
\item [] $914=(5-5/5)^5-55-55.$
\item [] $915=(5+5+5)\times (55+5+5/5).$
\item [] $916=5\times 55+5+(55+5^5)/5.$
\item [] $917=5\times 55+5+(55+5^5+5)/5.$
\item [] $918=(5-5/5)^5-555/5+5.$
\item [] $919=(5-5/5)^5-55-55+5.$
\item [] $920=5\times (55+5)-5+5^5/5.$
\item [] $921=5\times (55+5)-5+(5^5+5)/5.$
\item [] $922=555+55+(5^5-5)/(5+5).$
\item [] $923=5\times (55+5)+(5^5-5-5)/5.$
\item [] $924=(5-5/5)^5-5\times (5\times 5-5).$
\item [] $925=5\times (5\times 5\times 5+55+5).$
\item [] $926=5\times (55+5)+(5^5+5)/5.$
\item [] $927=5\times (55+5)+(5^5+5+5)/5.$
\item [] $928=((5+5)/5)^5\times (5\times 5+5-5/5).$
\item [] $929=5-5\times (5\times 5-5)+(5-5/5)^5.$
\item [] $930=5+5\times (55+5)+5^5/5.$
\item [] $931=5+5\times (55+5)+(5^5+5)/5.$
\item [] $932=5^5/5+(5^5-55)/(5+5).$
\item [] $933=5^5/5+(5^5+5)/(5+5)-5.$
\item [] $934=55\times (5+(55+5)/5)-5/5.$
\item [] $935=55\times (5+(55+5)/5).$
\item [] $936=55\times (5+(55+5)/5)+5/5.$
\item [] $937=(5^5-5)/(5+5)+5^5/5.$
\item [] $938=(5^5+5)/(5+5)+5^5/5.$
\item [] $939=(5^5+5)/5+(5^5+5)/(5+5).$
\item [] $940=55\times (5+(55+5)/5)+5.$
\item [] $941=55\times (5+(55+5)/5)+5+5/5.$
\item [] $942=(5^5-5)/(5+5)+5^5/5+5.$
\item [] $943=(5^5+5)/(5+5)+5^5/5+5.$
\item [] $944=(5-5/5)^5-5\times 5-55.$
\item [] $945=5\times 5\times (55+5)-555.$
\item [] $946=55\times (555/5-5\times 5)/5.$
\item [] $947=(5^5-5)/(5+5)+5^5/5+5+5.$
\item [] $948=(5^5+5)/(5+5)+5^5/5+5+5.$
\item [] $949=(5-5/5)^5-5\times (5+5+5).$
\item [] $950=(5+5)\times (5\times (5\times 5-5)-5).$
\item [] $951=(5+5)\times (5\times (5\times 5-5)-5)+5/5.$
\item [] $952=(5-5/5)\times ((5-(5+5)/5)^5-5).$
\item [] $953=55+5\times 55+(5^5-5-5)/5.$
\item [] $954=(5-5/5)^5+5-5\times (5+5+5).$
\item [] $955=5+(5+5)\times (5\times (5\times 5-5)-5).$
\item [] $956=55+5\times 55+(5^5+5)/5.$
\item [] $957=55\times (55+((5+5)/5)^5)/5.$
\item [] $958=(5-5/5)^5-55-55/5.$
\item [] $959=(5-5/5)^5-55-5-5.$
\item [] $960=(5\times 5+5)\times ((5+5)/5)^5.$
\item [] $961=(5\times 5+5+5/5)^{((5+5)/5)}.$
\item [] $962=(5\times 5+5/5)\times (5+((5+5)/5)^5).$
\item [] $963=(5-5/5)^5-55-5-5/5.$
\item [] $964=(5-5/5)^5-55-5.$
\item [] $965=5+(5\times 5+5)\times ((5+5)/5)^5.$
\item [] $966=5+(5\times 5+5+5/5)^{((5+5)/5)}.$
\item [] $967=(5-5/5)^5-(55+(5+5)/5).$
\item [] $968=(5-5/5)^5-55-5/5.$
\item [] $969=(5-5/5)^5-55.$
\item [] $970=(5-5/5)^5-55+5/5.$
\item [] $971=(5-5/5)^{(5+5/5)}-5^5.$
\item [] $972=(5-5/5)\times (5-(5+5)/5)^5.$
\item [] $973=(5-5/5)^5+5-55-5/5.$
\item [] $974=(5-5/5)^5+5-55.$
\item [] $975=5\times (5\times 5\times (5+5)-55).$
\item [] $976=5-5^5+(5-5/5)^{(5+5/5)}.$
\item [] $977=5+(5-5/5)\times (5-(5+5)/5)^5.$
\item [] $978=(5-5/5)^5+5+5-55-5/5.$
\item [] $979=(5-5/5)^5+5+5-55.$
\item [] $980=5+5\times (5\times 5\times (5+5)-55).$
\item [] $981=5555/5-5\times 5\times 5-5.$
\item [] $982=(5-5/5)\times (5-(5+5)/5)^5+5+5.$
\item [] $983=(5-5/5)^5-55/5-5\times 5-5.$
\item [] $984=(5-5/5)^5-55+5+5+5.$
\item [] $985=5+5\times (5\times 5\times (5+5)-55)+5.$
\item [] $986=5555/5-5\times 5\times 5.$
\item [] $987=(5-5/5)^5-5-((5+5)/5)^5.$
\item [] $988=(5-5/5)^5-5\times 5-55/5.$
\item [] $989=(5-5/5)^5-5\times 5-5-5.$
\item [] $990=(5+5)\times (5\times (5\times 5-5)-5/5).$
\item [] $991=5+5555/5-5\times 5\times 5.$
\item [] $992=(5-5/5)^5-((5+5)/5)^5.$
\item [] $993=(5-5/5)^5-5\times 5-5-5/5.$
\item [] $994=(5-5/5)^5-5\times 5-5.$
\item [] $995=5\times (5+5)\times (5\times 5-5)-5.$
\item [] $996=5\times (5+5)\times (5\times 5-5)-5+5/5.$
\item [] $997=(5-5/5)^5-((5+5)/5)^5+5.$
\item [] $998=(5-5/5)^5-5\times 5-5/5.$
\item [] $999=(5-5/5)^5-5\times 5.$
\item [] $1000=5\times (5+5)\times (5\times 5-5).$
\end{itemize}
\end{multicols}
}
\section{\textbf{Representations Using Number 6}}

{\footnotesize
\begin{multicols}{3}
\begin{itemize}
\item [] $101=66+6\times 6-6/6.$
\item [] $102=66+6\times 6.$
\item [] $103=66+6\times 6+6/6.$
\item [] $104=(666-6)/6-6.$
\item [] $105=666/6-6.$
\item [] $106=(666+6)/6-6.$
\item [] $107=6\times (6+6+6)-6/6.$
\item [] $108=6\times (6+6+6).$
\item [] $109=(666-6-6)/6.$
\item [] $110=(666-6)/6.$
\item [] $111=666/6.$
\item [] $112=(666+6)/6.$
\item [] $113=(666+6+6)/6.$
\item [] $114=6\times (6+6+6)+6.$
\item [] $115=6\times (6+6+6)+6+6/6.$
\item [] $116=6+(666-6)/6.$
\item [] $117=6+666/6.$
\item [] $118=6+(666+6)/6.$
\item [] $119=6+(666+6+6)/6.$
\item [] $120=6\times (6+6+6)+6+6.$
\item [] $121=66\times 66/(6\times 6).$
\item [] $122=(666+66)/6.$
\item [] $123=6+6+666/6.$
\item [] $124=6+6+(666+6)/6.$
\item [] $125=66+66-6-6/6.$
\item [] $126=66+66-6.$
\item [] $127=6+66\times 66/(6\times 6).$
\item [] $128=((6+6)/6)^{(6+6/6)}.$
\item [] $129=6+6+6+666/6.$
\item [] $130=66+((6+6)/6)^6.$
\item [] $131=66+66-6/6.$
\item [] $132=66+66.$
\item [] $133=66+66+6/6.$
\item [] $134=66+66+(6+6)/6.$
\item [] $135=6+6+6+6+666/6.$
\item [] $136=66+6+((6+6)/6)^6.$
\item [] $137=66+6+66-6/6.$
\item [] $138=66+66+6.$
\item [] $139=66+66+6+6/6.$
\item [] $140=66+66+6+(6+6)/6.$
\item [] $141=6\times 6-6+666/6.$
\item [] $142=(6+6)\times (6+6)-(6+6)/6.$
\item [] $143=(6+6)\times (6+6)-6/6.$
\item [] $144=(6+6)\times (6+6).$
\item [] $145=(6+6)\times (6+6)+6/6.$
\item [] $146=6\times 6+(666-6)/6.$
\item [] $147=66+6\times 66/6.$
\item [] $148=6\times 6+(666+6)/6.$
\item [] $149=(6+6)\times (6+6)+6-6/6.$
\item [] $150=(6+6)\times (6+6)+6.$
\item [] $151=(6+6)\times (6+6)+6+6/6.$
\item [] $152=6\times 6\times 6-((6+6)/6)^6.$
\item [] $153=6\times 6+6+666/6.$
\item [] $154=6\times 6+6+(666+6)/6.$
\item [] $155=(6+6)\times (6+6)+66/6.$
\item [] $156=(6+6)\times (6+6)+6+6.$
\item [] $157=(6+6)\times (6+6)+6+6+6/6.$
\item [] $158=6\times 6\times 6+6-((6+6)/6)^6.$
\item [] $159=6\times 6+6+6+666/6.$
\item [] $160=6\times 6+6+6+(666+6)/6.$
\item [] $161=(6+6)\times (6+6)+6+66/6.$
\item [] $162=(6+6)\times (6+6)+6+6+6.$
\item [] $163=6\times (6\times 6-6)-6-66/6.$
\item [] $164=6\times 6+((6+6)/6)^{(6+6/6)}.$
\item [] $165=66\times (6\times 6-6)/(6+6).$
\item [] $166=66+6\times 6+((6+6)/6)^6.$
\item [] $167=66+6\times 6+66-6/6.$
\item [] $168=66+6\times 6+66.$
\item [] $169=(6+6+6/6)^{((6+6)/6)}.$
\item [] $170=66-6+(666-6)/6.$
\item [] $171=66-6+666/6.$
\item [] $172=66-6+(666+6)/6.$
\item [] $173=6\times (6\times 6-6)-6-6/6.$
\item [] $174=6\times (6\times 6-6)-6.$
\item [] $175=6\times (6\times 6-6)-6+6/6.$
\item [] $176=66+(666-6)/6.$
\item [] $177=66+666/6.$
\item [] $178=66+(666+6)/6.$
\item [] $179=6\times (6\times 6-6)-6/6.$
\item [] $180=6\times (6\times 6-6).$
\item [] $181=6\times (6\times 6-6)+6/6.$
\item [] $182=6\times (6\times 6-6)+(6+6)/6.$
\item [] $183=66+6+666/6.$
\item [] $184=66+6+(666+6)/6.$
\item [] $185=6\times (6\times 6-6)+6-6/6.$
\item [] $186=6\times (6\times 6-6)+6.$
\item [] $187=6\times (6\times 6-6)+6+6/6.$
\item [] $188=6\times (6\times 6-6)+6+(6+6)/6.$
\item [] $189=66+6+6+666/6.$
\item [] $190=(6-6/6)\times (6\times 6+(6+6)/6).$
\item [] $191=6\times (6\times 6-6)+66/6.$
\item [] $192=6\times (6\times 6-6)+6+6.$
\item [] $193=6\times (6\times 6-6)+6+6+6/6.$
\item [] $194=66+((6+6)/6)^{(6+6/6)}.$
\item [] $195=6\times (6\times 66-6)/(6+6).$
\item [] $196=66+66+((6+6)/6)^6.$
\item [] $197=(66\times (6+6+6)-6)/6.$
\item [] $198=6\times 6\times 66/(6+6).$
\item [] $199=6\times 6\times 6-6-66/6.$
\item [] $200=6\times 6\times 6-6-(66-6)/6.$
\item [] $201=6\times (6\times 66+6)/(6+6).$
\item [] $202=6\times 6\times 6-6-6-(6+6)/6.$
\item [] $203=6\times 6\times 6-6-6-6/6.$
\item [] $204=6\times 6\times 6-6-6.$
\item [] $205=6\times 6\times 6-66/6.$
\item [] $206=6\times 6\times 6-(66-6)/6.$
\item [] $207=6+6\times (6\times 66+6)/(6+6).$
\item [] $208=6\times 6\times 6-6-(6+6)/6.$
\item [] $209=6\times 6\times 6-6-6/6.$
\item [] $210=6\times 6\times 6-6.$
\item [] $211=6\times 6\times 6-6+6/6.$
\item [] $212=6\times 6\times 6-6+(6+6)/6.$
\item [] $213=6\times 6\times 6-6\times 6/(6+6).$
\item [] $214=6\times 6\times 6-(6+6)/6.$
\item [] $215=6\times 6\times 6-6/6.$
\item [] $216=6\times 6\times 6.$
\item [] $217=6\times 6\times 6+6/6.$
\item [] $218=6\times 6\times 6+(6+6)/6.$
\item [] $219=6\times 6\times 6+6\times 6/(6+6).$
\item [] $220=6+6\times 6\times 6-(6+6)/6.$
\item [] $221=6+6\times 6\times 6-6/6.$
\item [] $222=6+6\times 6\times 6.$
\item [] $223=6+6\times 6\times 6+6/6.$
\item [] $224=6+6\times 6\times 6+(6+6)/6.$
\item [] $225=6\times 6\times 6+6+6\times 6/(6+6).$
\item [] $226=6\times 6\times 6+(66-6)/6.$
\item [] $227=6\times 66+6\times 6/6.$
\item [] $228=6\times 6\times 6+6+6.$
\item [] $229=6\times 6\times 6+6+6+6/6.$
\item [] $230=6\times 6\times 6+6+6+(6+6)/6.$
\item [] $231=66\times (6\times 6+6)/(6+6).$
\item [] $232=6\times 6\times 6+6+(66-6)/6.$
\item [] $233=6\times 6\times 6+6+66/6.$
\item [] $234=6\times 6\times 6+6+6+6.$
\item [] $235=6\times 6\times 6+6+6+6+6/6.$
\item [] $236=66\times 66/(6+6+6)-6.$
\item [] $237=66\times (6\times 6+6)/(6+6)+6.$
\item [] $238=6\times 6\times 6+(66+66)/6.$
\item [] $239=6\times 6\times 6+6+6+66/6.$
\item [] $240=6\times (6\times 6+6)-6-6.$
\item [] $241=6\times (6\times 6+6)-66/6.$
\item [] $242=66\times 66/(6+6+6).$
\item [] $243=(6\times 6/(6+6))^(6-6/6).$
\item [] $244=6\times 6\times 6-6+((6+6)/6)^6.$
\item [] $245=(6+6/6)\times (6\times 6-6/6).$
\item [] $246=6\times (6\times 6+6)-6.$
\item [] $247=6\times (6\times 6+6)-6+6/6.$
\item [] $248=66\times 66/(6+6+6)+6.$
\item [] $249=(6\times 6/(6+6))^(6-6/6)+6.$
\item [] $250=6\times (6\times 6+6)-(6+6)/6.$
\item [] $251=6\times (6\times 6+6)-6/6.$
\item [] $252=6\times (6\times 6+6).$
\item [] $253=6\times (6\times 6+6)+6/6.$
\item [] $254=6\times (6\times 6+6)+(6+6)/6.$
\item [] $255=6\times (6\times 6+6)+6\times 6/(6+6).$
\item [] $256=((6+6)/6)^{(6+(6+6)/6)}.$
\item [] $257=6\times (6\times 6+6)+6-6/6.$
\item [] $258=6\times (6\times 6+6)+6.$
\item [] $259=6\times (6\times 6+6)+6+6/6.$
\item [] $260=6\times (6\times 6+6)+6+(6+6)/6.$
\item [] $261=6\times (6\times 6+6)+6+6\times 6/(6+6).$
\item [] $262=((6+6)/6)^{(6+(6+6)/6)}+6.$
\item [] $263=6\times (6\times 6+6)+66/6.$
\item [] $264=6\times (6\times 6+6)+6+6.$
\item [] $265=6\times (6\times 6+6)+6+6+6/6.$
\item [] $266=(6+6/6)\times (6\times 6+(6+6)/6).$
\item [] $267=6\times 666/(6+6)-66.$
\item [] $268=(66+6/6)\times (6-(6+6)/6).$
\item [] $269=6\times (6\times 6+6)+6+66/6.$
\item [] $270=666-6\times 66.$
\item [] $271=66+6\times 6\times 6-66/6.$
\item [] $272=(6+6/6)\times (6\times 6+(6+6)/6)+6.$
\item [] $273=6-66+6\times 666/(6+6).$
\item [] $274=6\times 6\times 6-6+((6+6)/6)^6.$
\item [] $275=66+6\times 6\times 6-6-6/6.$
\item [] $276=66+6\times 6\times 6-6.$
\item [] $277=66+6\times 6\times 6-6+6/6.$
\item [] $278=66+6\times 6\times 66/(6+6+6).$
\item [] $279=6\times 66-6-666/6.$
\item [] $280=6\times 6\times 6+((6+6)/6)^6.$
\item [] $281=66+6\times 6\times 6-6/6.$
\item [] $282=66+6\times 6\times 6.$
\item [] $283=66+6\times 6\times 6+6/6.$
\item [] $284=66+6\times 6\times 6+(6+6)/6.$
\item [] $285=6\times 66-666/6.$
\item [] $286=6\times 6\times 6+((6+6)/6)^6+6.$
\item [] $287=6\times (6\times 6+6+6)-6/6.$
\item [] $288=6\times (6\times 6+6+6).$
\item [] $289=6\times (6\times 6+6+6)+6/6.$
\item [] $290=6\times (6\times 6+6+6)+(6+6)/6.$
\item [] $291=6+6\times 66-666/6.$
\item [] $292=6\times 6\times 6+6+6+((6+6)/6)^6.$
\item [] $293=6\times (6\times 6+6+6)+6-6/6.$
\item [] $294=6\times (6\times 6+6+6)+6.$
\item [] $295=6\times (6\times 6+6+6)+6+6/6.$
\item [] $296=(6+(6+6)/6)\times (6\times 6+6/6).$
\item [] $297=6\times 6\times 66/(6+(6+6)/6).$
\item [] $298=6\times (6\times 6+6+6)+(66-6)/6.$
\item [] $299=6\times (6\times 6+6+6)+66/6.$
\item [] $300=(66-6)\times (6-6/6).$
\item [] $301=(66-6)\times (6-6/6)+6/6.$
\item [] $302=(6+(6+6)/6)\times (6\times 6+6/6)+6.$
\item [] $303=6\times 6\times 66/(6+(6+6)/6)+6.$
\item [] $304=(6+(6+6)/6)\times (6\times 6+(6+6)/6).$
\item [] $305=(6-6/6)\times (66-6+6/6).$
\item [] $306=(66-6)\times (6-6/6)+6.$
\item [] $307=(66-6)\times (6-6/6)+6+6/6.$
\item [] $308=66+66\times 66/(6+6+6).$
\item [] $309=66+(6\times 6/(6+6))^(6-6/6).$
\item [] $310=(6-6/6)\times (66-6+(6+6)/6).$
\item [] $311=(6-6/6)\times (66-6+6/6)+6.$
\item [] $312=6\times (((6+6)/6)^6-6-6).$
\item [] $313=6\times (66-6-6)-66/6.$
\item [] $314=(6-6/6)\times ((6+6)/6)^6-6.$
\item [] $315=6\times (666-6\times 6)/(6+6).$
\item [] $316=6\times (6\times 6+6)+((6+6)/6)^6.$
\item [] $317=66+6\times (6\times 6+6)-6/6.$
\item [] $318=66+6\times (6\times 6+6).$
\item [] $319=66+6\times (6\times 6+6)+6/6.$
\item [] $320=(6-6/6)\times ((6+6)/6)^6.$
\item [] $321=6\times 6\times 6-6+666/6.$
\item [] $322=6\times (66-6-6)-(6+6)/6.$
\item [] $323=6\times (66-6-6)-6/6.$
\item [] $324=6\times (66-6-6).$
\item [] $325=6\times (66-6-6)+6/6.$
\item [] $326=(6-6/6)\times ((6+6)/6)^6+6.$
\item [] $327=6\times 66+6\times 66/6.$
\item [] $328=6\times 6\times 6+(666+6)/6.$
\item [] $329=6\times 66-66-6/6.$
\item [] $330=6\times 66-66.$
\item [] $331=6\times 66-66+6/6.$
\item [] $332=6\times 66-((6+6)/6)^6.$
\item [] $333=6\times 666/(6+6).$
\item [] $334=6\times 666/(6+6)+6/6.$
\item [] $335=(6-6/6)\times (66+6/6).$
\item [] $336=66\times (6-6/6)+6.$
\item [] $337=66\times (6-6/6)+6+6/6.$
\item [] $338=6\times 66+6-((6+6)/6)^6.$
\item [] $339=6\times 666/(6+6)+6.$
\item [] $340=(6-6/6)\times (66+(6+6)/6).$
\item [] $341=(6-6/6)\times (66+6/6)+6.$
\item [] $342=66\times (6-6/6)+6+6.$
\item [] $343=(6+6/6)^(6\times 6/(6+6)).$
\item [] $344=6\times 66+6+6-((6+6)/6)^6.$
\item [] $345=6\times 666/(6+6)+6+6.$
\item [] $346=(6-6/6)\times (66+(6+6)/6)+6.$
\item [] $347=6\times (66-6)-6-6-6/6.$
\item [] $348=6\times (66-6)-6-6.$
\item [] $349=6\times (66-6)-66/6.$
\item [] $350=(6-6/6)\times (6+((6+6)/6)^6).$
\item [] $351=6\times (666+6\times 6)/(6+6).$
\item [] $352=6\times (66-6)-6-(6+6)/6.$
\item [] $353=6\times (66-6)-6-6/6.$
\item [] $354=6\times (66-6)-6.$
\item [] $355=6\times (66-6)-6+6/6.$
\item [] $356=6\times (66-6)-6+(6+6)/6.$
\item [] $357=66\times 66/(6+6)-6.$
\item [] $358=6\times (66-6)-(6+6)/6.$
\item [] $359=6\times (66-6)-6/6.$
\item [] $360=6\times (66-6).$
\item [] $361=6\times (66-6)+6/6.$
\item [] $362=6\times (66-6)+(6+6)/6.$
\item [] $363=66\times 66/(6+6).$
\item [] $364=6\times (66-6)+6-(6+6)/6.$
\item [] $365=6\times (66-6)+6-6/6.$
\item [] $366=6\times (66-6)+6.$
\item [] $367=6\times (66-6)+6+6/6.$
\item [] $368=6\times (66-6)+6+(6+6)/6.$
\item [] $369=66\times 66/(6+6)+6.$
\item [] $370=6\times (66-6)+(66-6)/6.$
\item [] $371=6\times (66-6)+66/6.$
\item [] $372=6\times (66-6)+6+6.$
\item [] $373=6\times (66-6)+6+6+6/6.$
\item [] $374=66\times (6\times 6-(6+6)/6)/6.$
\item [] $375=66\times 66/(6+6)+6+6.$
\item [] $376=6\times (66-6)+6+(66-6)/6.$
\item [] $377=6\times (66-6)+6+66/6.$
\item [] $378=6\times 66-6-6-6.$
\item [] $379=6\times 66-6-66/6.$
\item [] $380=6\times 66-6-(66-6)/6.$
\item [] $381=66\times 66/(6+6)+6+6+6.$
\item [] $382=6\times 66-6-6-(6+6)/6.$
\item [] $383=6\times 66-6-6-6/6.$
\item [] $384=6\times ((6+6)/6)^6.$
\item [] $385=6\times 66-66/6.$
\item [] $386=6\times 66-(66-6)/6.$
\item [] $387=6\times 66-(66-6-6)/6.$
\item [] $388=6\times 66-6-(6+6)/6.$
\item [] $389=6\times 66-6-6/6.$
\item [] $390=6\times 66-6.$
\item [] $391=6\times 66-6+6/6.$
\item [] $392=6\times 66-6+(6+6)/6.$
\item [] $393=6\times 66-6\times 6/(6+6).$
\item [] $394=6\times 66-(6+6)/6.$
\item [] $395=6\times 66-6/6.$
\item [] $396=6\times 66.$
\item [] $397=6\times 66+6/6.$
\item [] $398=6\times 66+(6+6)/6.$
\item [] $399=6\times 66+6\times 6/(6+6).$
\item [] $400=6\times 66+6-(6+6)/6.$
\item [] $401=6\times 66+6-6/6.$
\item [] $402=6\times 66+6.$
\item [] $403=6+6\times 66+6/6.$
\item [] $404=6+6\times 66+(6+6)/6.$
\item [] $405=6+6\times 66+6\times 6/(6+6).$
\item [] $406=6\times 66+(66-6)/6.$
\item [] $407=6\times 66+66/6.$
\item [] $408=6\times 66+6+6.$
\item [] $409=6\times 66+6+6+6/6.$
\item [] $410=6\times 66+6+6+(6+6)/6.$
\item [] $411=6\times 66+6+6+6\times 6/(6+6).$
\item [] $412=6\times 66+6+(66-6)/6.$
\item [] $413=6\times 66+6+66/6.$
\item [] $414=6\times 66+6+6+6.$
\item [] $415=6\times 66+6+6+6+6/6.$
\item [] $416=6\times 66+6+6+6+(6+6)/6.$
\item [] $417=6\times (66+(6\times 6+6)/(6+6)).$
\item [] $418=6\times 66+(66+66)/6.$
\item [] $419=6\times 66+6+6+66/6.$
\item [] $420=6\times (6+((6+6)/6)^6).$
\item [] $421=6\times (66+6)-66/6.$
\item [] $422=6\times (66+6)+(6-66)/6.$
\item [] $423=66-6+66\times 66/(6+6).$
\item [] $424=6\times (66-6)+((6+6)/6)^6.$
\item [] $425=6\times (66+6)-6-6/6.$
\item [] $426=6\times (66+6)-6.$
\item [] $427=6\times (66+6)-6+6/6.$
\item [] $428=6\times (66+6)-6+(6+6)/6.$
\item [] $429=66+66\times 66/(6+6).$
\item [] $430=6\times (66+6)-(6+6)/6.$
\item [] $431=6\times (66+6)-6/6.$
\item [] $432=6\times (66+6).$
\item [] $433=6\times (66+6)+6/6.$
\item [] $434=6\times (66+6)+(6+6)/6.$
\item [] $435=6\times (66+6)+6\times 6/(6+6).$
\item [] $436=6\times (66+6)+6-(6+6)/6.$
\item [] $437=6\times (66+6)+6-6/6.$
\item [] $438=6\times (66+6)+6.$
\item [] $439=6\times (66+6)+6+6/6.$
\item [] $440=6\times (66+6)+6+(6+6)/6.$
\item [] $441=6\times (66+6)+6+6\times 6/(6+6).$
\item [] $442=6\times (66+6)+(66-6)/6.$
\item [] $443=6\times (66+6)+66/6.$
\item [] $444=6\times (66+6)+6+6.$
\item [] $445=6\times (66+6)+6+6+6/6.$
\item [] $446=6\times (66+6)+6+6+(6+6)/6.$
\item [] $447=(6+6/6)\times ((6+6)/6)^6-6/6.$
\item [] $448=(6+6/6)\times ((6+6)/6)^6.$
\item [] $449=6\times (66+6)+6+66/6.$
\item [] $450=666-6\times 6\times 6.$
\item [] $451=66+6\times 66-66/6.$
\item [] $452=66+6\times 66-(66-6)/6.$
\item [] $453=66+6\times 66-(66-6-6)/6.$
\item [] $454=(6+6/6)\times ((6+6)/6)^6+6.$
\item [] $455=(6/6+6)\times (66-6/6).$
\item [] $456=66+6\times 66-6.$
\item [] $457=66+6\times 66-6+6/6.$
\item [] $458=66+6\times 66-6+(6+6)/6.$
\item [] $459=66+6\times 66-6\times 6/(6+6).$
\item [] $460=6\times 66+((6+6)/6)^6.$
\item [] $461=66+6\times 66-6/6.$
\item [] $462=66+6\times 66.$
\item [] $463=66+6\times 66+6/6.$
\item [] $464=66+6\times 66+(6+6)/6.$
\item [] $465=66+6\times 66+6\times 6/(6+6).$
\item [] $466=6\times 66+6+((6+6)/6)^6.$
\item [] $467=6\times (66+6+6)-6/6.$
\item [] $468=6\times (66+6+6).$
\item [] $469=6\times (66+6+6)+6/6.$
\item [] $470=6\times (66+6+6)+(6+6)/6.$
\item [] $471=6\times (66-6)+666/6.$
\item [] $472=6\times 66+6+6+((6+6)/6)^6.$
\item [] $473=6\times (66+6+6)+6-6/6.$
\item [] $474=6\times (66+6+6)+6.$
\item [] $475=6\times (66+6+6)+6+6/6.$
\item [] $476=(6+6/6)\times (66+(6+6)/6).$
\item [] $477=6\times (66-6)+6+666/6.$
\item [] $478=6\times (66+6+6)+(66-6)/6.$
\item [] $479=6\times (66+6+6)+66/6.$
\item [] $480=6\times (66+6+6)+6+6.$
\item [] $481=(6+6+6/6)\times (6\times 6+6/6).$
\item [] $482=(6+6/6)\times (66+(6+6)/6)+6.$
\item [] $483=(6+6/6)\times (66+6\times 6/(6+6)).$
\item [] $484=((66+66)/6)^{((6+6)/6)}.$
\item [] $485=6\times (66+6+6)+6+66/6.$
\item [] $486=666-6\times (6\times 6-6).$
\item [] $487=(6+6+6/6)\times (6\times 6+6/6)+6.$
\item [] $488=(6+(6+6)/6)\times (66-6+6/6).$
\item [] $489=666-66-666/6.$
\item [] $490=(6+6/6)\times (6+((6+6)/6)^6).$
\item [] $491=6\times 6+(6+6/6)\times (66-6/6).$
\item [] $492=66+6\times (66+6)-6.$
\item [] $493=(6+6)\times (6\times 6+6)-66/6.$
\item [] $494=(6+6)\times (6\times 6+6)-(66-6)/6.$
\item [] $495=6\times 66-6-6+666/6.$
\item [] $496=((6+6)/6)^6+6\times (66+6).$
\item [] $497=66+6\times (66+6)-6/6.$
\item [] $498=66+6\times (66+6).$
\item [] $499=66+6\times (66+6)+6/6.$
\item [] $500=66+6\times (66+6)+(6+6)/6.$
\item [] $501=6\times 66-6+666/6.$
\item [] $502=(6+6)\times (6\times 6+6)-(6+6)/6.$
\item [] $503=(6+6)\times (6\times 6+6)-6/6.$
\item [] $504=(6+6)\times (6\times 6+6).$
\item [] $505=(6+6)\times (6\times 6+6)+6/6.$
\item [] $506=(6+6)\times (6\times 6+6)+(6+6)/6.$
\item [] $507=6\times 66+666/6.$
\item [] $508=6\times 66+(666+6)/6.$
\item [] $509=(6+6)\times (6\times 6+6)+6-6/6.$
\item [] $510=(6+6)\times (6\times 6+6)+6.$
\item [] $511=(6+6)\times (6\times 6+6)+6+6/6.$
\item [] $512=((6+6)/6)^{(6+6\times 6/(6+6))}.$
\item [] $513=6\times 66+6+666/6.$
\item [] $514=6\times 66+6+(666+6)/6.$
\item [] $515=(6+6)\times (6\times 6+6)+66/6.$
\item [] $516=(6+6)\times (6\times 6+6)+6+6.$
\item [] $517=66\times (66+6\times 6/6)/6.$
\item [] $518=((6+6)/6)^{(6+6\times 6/(6+6))}+6.$
\item [] $519=6\times 66+6+6+666/6.$
\item [] $520=(6+(6+6)/6)\times (66-6/6).$
\item [] $521=(6+6)\times (6\times 6+6)+6+66/6.$
\item [] $522=666-(6+6)\times (6+6).$
\item [] $523=6+66\times (66+6\times 6/6)/6.$
\item [] $524=6\times 66+((6+6)/6)^{(6+6/6)}.$
\item [] $525=(6-6/6)\times (666/6-6).$
\item [] $526=(6+(6+6)/6)\times (66-6/6)+6.$
\item [] $527=66\times (6+(6+6)/6)-6/6.$
\item [] $528=66\times (6+(6+6)/6).$
\item [] $529=66\times (6+(6+6)/6)+6/6.$
\item [] $530=66\times (6+(6+6)/6)+(6+6)/6.$
\item [] $531=(6-6/6)\times (666/6-6)+6.$
\item [] $532=(6/6+6)\times (6+6+((6+6)/6)^6).$
\item [] $533=66\times (6+(6+6)/6)+6-6/6.$
\item [] $534=66\times (6+(6+6)/6)+6.$
\item [] $535=66\times (6+(6+6)/6)+6+6/6.$
\item [] $536=(6+(6+6)/6)\times (66+6/6).$
\item [] $537=6\times (66+6)-6+666/6.$
\item [] $538=666-((6+6)/6)^{(6+6/6)}.$
\item [] $539=(6+6/6)\times (66+66/6).$
\item [] $540=(6+6+6)\times (6\times 6-6).$
\item [] $541=(6+6+6)\times (6\times 6-6)+6/6.$
\item [] $542=(6+(6+6)/6)\times (66+6/6)+6.$
\item [] $543=6\times (66+6)+666/6.$
\item [] $544=666-(666+66)/6.$
\item [] $545=(6+6/6)\times (66+66/6)+6.$
\item [] $546=(6+6+6)\times (6\times 6-6)+6.$
\item [] $547=(6+6+6)\times (6\times 6-6)+6+6/6.$
\item [] $548=666-6-(666+6)/6.$
\item [] $549=666-6-666/6.$
\item [] $550=(6-6/6)\times (666-6)/6.$
\item [] $551=(6+6+6)\times (6\times 6-6)+66/6.$
\item [] $552=(6+6+6)\times (6\times 6-6)+6+6.$
\item [] $553=(6+6/6)\times (66+6+6+6/6).$
\item [] $554=666-(666+6)/6.$
\item [] $555=666-666/6.$
\item [] $556=(6666+6)/(6+6).$
\item [] $557=(6666+6)/(6+6)+6/6.$
\item [] $558=666-6\times (6+6+6).$
\item [] $559=666-6\times (6+6+6)+6/6.$
\item [] $560=(6-6/6)\times (666+6)/6.$
\item [] $561=666\times (6-6/6)/6+6.$
\item [] $562=(6666+6)/(6+6)+6.$
\item [] $563=666-66-6\times 6-6/6.$
\item [] $564=(6+6)\times (66+6\times 6/6).$
\item [] $565=(6+6)\times (66+6\times 6/6)+6/6.$
\item [] $566=(6-6/6)\times (666+6)/6+6.$
\item [] $567=666\times (6-6/6)/6+6+6.$
\item [] $568=(6666+6)/(6+6)+6+6.$
\item [] $569=66+(6+6)\times (6\times 6+6)-6/6.$
\item [] $570=66+(6+6)\times (6\times 6+6).$
\item [] $571=66+(6+6)\times (6\times 6+6)+6/6.$
\item [] $572=66\times (((6+6)/6)^6-6-6)/6.$
\item [] $573=66+6\times 66+666/6.$
\item [] $574=6\times (66+6\times 6-6)-(6+6)/6.$
\item [] $575=6\times (66+6\times 6-6)-6/6.$
\item [] $576=6\times (66+6\times 6-6).$
\item [] $577=6\times (66+6\times 6-6)+6/6.$
\item [] $578=6\times (66+6\times 6-6)+(6+6)/6.$
\item [] $579=6\times 66+6\times 6\times 66/(6+6).$
\item [] $580=(6-6/6)\times (6+(666-6)/6).$
\item [] $581=6\times (66+6\times 6-6)+6-6/6.$
\item [] $582=6\times (66+6\times 6-6)+6.$
\item [] $583=6\times (66+6\times 6-6)+6+6/6.$
\item [] $584=(6+(6+6)/6)\times (66+6+6/6).$
\item [] $585=(6-6/6)\times (6+666/6).$
\item [] $586=(((6+6)/6)^{(6+6)}+6)/(6+6/6).$
\item [] $587=666-66-6-6-6/6.$
\item [] $588=666-66-6-6.$
\item [] $589=666-66-66/6.$
\item [] $590=666-66-(66-6)/6.$
\item [] $591=(6-6/6)\times (6+666/6)+6.$
\item [] $592=6\times 6+(6666+6)/(6+6).$
\item [] $593=666-66-6-6/6.$
\item [] $594=666-66-6.$
\item [] $595=666-66-6+6/6.$
\item [] $596=666-6-((6+6)/6)^6.$
\item [] $597=666-66-6\times 6/(6+6).$
\item [] $598=666-66-(6+6)/6.$
\item [] $599=666-66-6/6.$
\item [] $600=666-66.$
\item [] $601=666-66+6/6.$
\item [] $602=666-((6+6)/6)^6.$
\item [] $603=666-((6+6)/6)^6+6/6.$
\item [] $604=666-66+6-(6+6)/6.$
\item [] $605=666-66+6-6/6.$
\item [] $606=666-66+6.$
\item [] $607=666-66+6+6/6.$
\item [] $608=666+6-((6+6)/6)^6.$
\item [] $609=6\times (66+6\times 6)-6\times 6/(6+6).$
\item [] $610=6\times (66+6\times 6)-(6+6)/6.$
\item [] $611=6\times (66+6\times 6)-6/6.$
\item [] $612=6\times (66+6\times 6).$
\item [] $613=6\times (66+6\times 6)+6/6.$
\item [] $614=6\times (66+6\times 6)+(6+6)/6.$
\item [] $615=(66\times 66-666)/6.$
\item [] $616=6\times (66+6\times 6)+6-(6+6)/6.$
\item [] $617=6\times (66+6\times 6)+6-6/6.$
\item [] $618=6\times (66+6\times 6)+6.$
\item [] $619=6\times (66+6\times 6)+6+6/6.$
\item [] $620=6\times (66+6\times 6)+6+(6+6)/6.$
\item [] $621=6+(66\times 66-666)/6.$
\item [] $622=66+(6666+6)/(6+6).$
\item [] $623=6\times (66+6\times 6)+66/6.$
\item [] $624=666-6\times 6-6.$
\item [] $625=(6-6/6)^{(6-(6+6)/6)}.$
\item [] $626=(6-6/6)^{(6-(6+6)/6)}+6/6.$
\item [] $627=666-6-6\times 66/(6+6).$
\item [] $628=666-6\times 6-(6+6)/6.$
\item [] $629=666-6\times 6-6/6.$
\item [] $630=666-6\times 6.$
\item [] $631=666-6\times 6+6/6.$
\item [] $632=666-6\times 6+(6+6)/6.$
\item [] $633=666-6\times 66/(6+6).$
\item [] $634=666-6\times 6+6-(6+6)/6.$
\item [] $635=666-6\times 6+6-6/6.$
\item [] $636=666-6\times 6+6.$
\item [] $637=666-6\times 6+6+6/6.$
\item [] $638=66\times (((6+6)/6)^6-6)/6.$
\item [] $639=666-6\times 66/(6+6)+6.$
\item [] $640=(66-6)\times ((6+6)/6)^6/6.$
\item [] $641=666-66+6\times 6/6.$
\item [] $642=6\times 6\times (6+6+6)-6.$
\item [] $643=6\times 6\times (6+6+6)-6+6/6.$
\item [] $644=666-(66+66)/6.$
\item [] $645=6\times (6\times 6\times 6\times 6-6)/(6+6).$
\item [] $646=6\times 6\times (6+6+6)-(6+6)/6.$
\item [] $647=6\times 6\times (6+6+6)-6/6.$
\item [] $648=6\times 6\times (6+6+6).$
\item [] $649=666-6-66/6.$
\item [] $650=666-6-(66-6)/6.$
\item [] $651=6\times (6\times 6\times 6\times 6+6)/(6+6).$
\item [] $652=666-6-6-(6+6)/6.$
\item [] $653=666-6-6-6/6.$
\item [] $654=666-6-6.$
\item [] $655=666-66/6.$
\item [] $656=666-(66-6)/6.$
\item [] $657=666-(66-6-6)/6.$
\item [] $658=666-6-(6+6)/6.$
\item [] $659=666-6-6/6.$
\item [] $660=666-6.$
\item [] $661=666-6+6/6.$
\item [] $662=666-6+(6+6)/6.$
\item [] $663=666-6\times 6/(6+6).$
\item [] $664=666-(6+6)/6.$
\item [] $665=666-6/6.$
\item [] $666=666.$
\item [] $667=666+6/6.$
\item [] $668=666+(6+6)/6.$
\item [] $669=666+6\times 6/(6+6).$
\item [] $670=666+6-(6+6)/6.$
\item [] $671=666+6-6/6.$
\item [] $672=666+6.$
\item [] $673=666+6+6/6.$
\item [] $674=666+6+(6+6)/6.$
\item [] $675=666+6+6\times 6/(6+6).$
\item [] $676=666+(66-6)/6.$
\item [] $677=666+66/6.$
\item [] $678=666+6+6.$
\item [] $679=666+6+6+6/6.$
\item [] $680=666+6+6+(6+6)/6.$
\item [] $681=66\times (6+6)-666/6.$
\item [] $682=666+6+(66-6)/6.$
\item [] $683=666+6+66/6.$
\item [] $684=666+6+6+6.$
\item [] $685=666+6+6+6+6/6.$
\item [] $686=666+6+6+6+(6+6)/6.$
\item [] $687=(6\times 6/(6+6))^6-6\times 6-6.$
\item [] $688=666+(66+66)/6.$
\item [] $689=666+6+6+66/6.$
\item [] $690=666+6+6+6+6.$
\item [] $691=66+6\times 66-66/6.$
\item [] $692=66+6\times 66-(66-6)/6.$
\item [] $693=(6\times 6/(6+6))^6-6\times 6.$
\item [] $694=666+6+(66+66)/6.$
\item [] $695=66+6\times 66-6-6/6.$
\item [] $696=66+6\times 66-6.$
\item [] $697=66+6\times 66-6+6/6.$
\item [] $698=66\times ((6+6)/6)^6/6-6.$
\item [] $699=(6\times 6/(6+6))^6-6\times 6+6.$
\item [] $700=66+6\times 66-(6+6)/6.$
\item [] $701=66+6\times 66-6/6.$
\item [] $702=66+6\times 66.$
\item [] $703=66+6\times 66+6/6.$
\item [] $704=66\times ((6+6)/6)^6/6.$
\item [] $705=(66\times ((6+6)/6)^6+6)/6.$
\item [] $706=666+6\times 6+6-(6+6)/6.$
\item [] $707=666+6\times 6+6-6/6.$
\item [] $708=666+6\times 6+6.$
\item [] $709=666+6\times 6+6+6/6.$
\item [] $710=66\times ((6+6)/6)^6/6+6.$
\item [] $711=(6\times 6/(6+6))^6-6-6-6.$
\item [] $712=66+6\times 66+(66-6)/6.$
\item [] $713=66+6\times 66+66/6.$
\item [] $714=(6+6)\times (66-6)-6.$
\item [] $715=66\times (66-6/6)/6.$
\item [] $716=(66\times (66-6/6)+6)/6.$
\item [] $717=(6\times 6/(6+6))^6-6-6.$
\item [] $718=(6+6)\times (66-6)-(6+6)/6.$
\item [] $719=(6+6)\times (66-6)-6/6.$
\item [] $720=(6+6)\times (66-6).$
\item [] $721=(6+6)\times (66-6)+6/6.$
\item [] $722=(6+6)\times (66-6)+(6+6)/6.$
\item [] $723=(6\times 6/(6+6))^6-6.$
\item [] $724=(66\times 66-6-6)/6.$
\item [] $725=(66\times 66-6)/6.$
\item [] $726=66\times 66/6.$
\item [] $727=(66\times 66+6)/6.$
\item [] $728=(6\times 6/(6+6))^6-6/6.$
\item [] $729=(6\times 6/(6+6))^6.$
\item [] $730=(6\times 6/(6+6))^6+6/6.$
\item [] $731=6+(66\times 66-6)/6.$
\item [] $732=66+666.$
\item [] $733=(66\times 66+6)/6+6.$
\item [] $734=(66\times 66+6+6)/6+6.$
\item [] $735=(6\times 6/(6+6))^6+6.$
\item [] $736=(6\times 6/(6+6))^6+6+6/6.$
\item [] $737=66\times (66+6/6)/6.$
\item [] $738=666+66+6.$
\item [] $739=(66\times 66+6)/6+6+6.$
\item [] $740=(6\times 6/(6+6))^6+66/6.$
\item [] $741=(6\times 6/(6+6))^6+6+6.$
\item [] $742=(6\times 6/(6+6))^6+6+6+6/6.$
\item [] $743=66\times (66+6/6)/6+6.$
\item [] $744=666+66+6+6.$
\item [] $745=(66\times 66+6)/6+6+6+6.$
\item [] $746=(6\times 6/(6+6))^6+6+66/6.$
\item [] $747=(6\times 6/(6+6))^6+6+6+6.$
\item [] $748=66\times (66+(6+6)/6)/6.$
\item [] $749=66\times (66+6/6)/6+6+6.$
\item [] $750=66\times (6+6)-6\times 6-6.$
\item [] $751=66+6\times 6\times (66-6/6)/6.$
\item [] $752=(6\times 6/(6+6))^6+6+6+66/6.$
\item [] $753=(6\times 6/(6+6))^6+6+6+6+6.$
\item [] $754=66\times (66+(6+6)/6)/6+6.$
\item [] $755=66\times (6+6)-6\times 6-6/6.$
\item [] $756=6\times (66+66-6).$
\item [] $757=6\times (66+66-6)+6/6.$
\item [] $758=6\times (66+66-6)+(6+6)/6.$
\item [] $759=6\times 6-6+(6\times 6/(6+6))^6.$
\item [] $760=6\times 6+(66\times 66-6-6)/6.$
\item [] $761=6\times 6+(66\times 66-6)/6.$
\item [] $762=6\times (66+66-6)+6.$
\item [] $763=6\times 6+(66\times 66+6)/6.$
\item [] $764=6\times 6+(6\times 6/(6+6))^6-6/6.$
\item [] $765=6\times 6+(6\times 6/(6+6))^6.$
\item [] $766=6\times 6+(6\times 6/(6+6))^6+6/6.$
\item [] $767=(6+6)\times ((6+6)/6)^6-6/6.$
\item [] $768=(6+6)\times ((6+6)/6)^6.$
\item [] $769=(6+6)\times ((6+6)/6)^6+6/6.$
\item [] $770=(6+6/6)\times (666-6)/6.$
\item [] $771=666-6+666/6.$
\item [] $772=666-6+(666+6)/6.$
\item [] $773=(6+6)\times ((6+6)/6)^6+6-6/6.$
\item [] $774=(6+6)\times ((6+6)/6)^6+6.$
\item [] $775=66\times (6+6)-6-66/6.$
\item [] $776=666+(666-6)/6.$
\item [] $777=666+666/6.$
\item [] $778=666+(666+6)/6.$
\item [] $779=66\times (6+6)-6-6-6/6.$
\item [] $780=(6+6)\times (66-6/6).$
\item [] $781=66\times (6+6)-66/6.$
\item [] $782=66\times (6+6)-(66-6)/6.$
\item [] $783=666+6+666/6.$
\item [] $784=(6+6/6)\times (666+6)/6.$
\item [] $785=66\times (6+6)-6-6/6.$
\item [] $786=66\times (6+6)-6.$
\item [] $787=66\times (6+6)-6+6/6.$
\item [] $788=66\times (6+6)-6+(6+6)/6.$
\item [] $789=66-6+(6\times 6/(6+6))^6.$
\item [] $790=66\times (6+6)-(6+6)/6.$
\item [] $791=66\times (6+6)-6/6.$
\item [] $792=66\times (6+6).$
\item [] $793=66\times (6+6)+6/6.$
\item [] $794=66\times (6+6)+(6+6)/6.$
\item [] $795=66+(6\times 6/(6+6))^6.$
\item [] $796=66\times (6+6)+6-(6+6)/6.$
\item [] $797=66\times (6+6)+6-6/6.$
\item [] $798=66\times (6+6)+6.$
\item [] $799=66\times (6+6)+6+6/6.$
\item [] $800=66\times (6+6)+6+(6+6)/6.$
\item [] $801=(6\times 6/(6+6))^6+66+6.$
\item [] $802=66\times (6+6)+(66-6)/6.$
\item [] $803=66\times (6+6)+66/6.$
\item [] $804=66\times (6+6)+6+6.$
\item [] $805=66\times (6+6)+6+6+6/6.$
\item [] $806=66\times (6+6)+6+6+(6+6)/6.$
\item [] $807=66+6+6+(6\times 6/(6+6))^6.$
\item [] $808=66\times (6+6)+6+(66-6)/6.$
\item [] $809=66\times (6+6)+6+66/6.$
\item [] $810=66\times (6+6)+6+6+6.$
\item [] $811=66\times (6+6)+6+6+6+6/6.$
\item [] $812=(6+6/6)\times (6+(666-6)/6).$
\item [] $813=66+6\times 66+666/6.$
\item [] $814=66\times (66+6+(6+6)/6)/6.$
\item [] $815=66\times (6+6)+6+6+66/6.$
\item [] $816=(6+6)\times (66+(6+6)/6).$
\item [] $817=66+6\times 6\times (6+6)-66/6.$
\item [] $818=(((6+6)/6)^{(6+6)}-6)/(6-6/6).$
\item [] $819=(6+6/6)\times (6+666/6).$
\item [] $820=66\times (6+6)+6+(66+66)/6.$
\item [] $821=66+6\times 6\times (6+6)-6-6/6.$
\item [] $822=66+6\times 6\times (6+6)-6.$
\item [] $823=66+6\times 6\times (6+6)-6+6/6.$
\item [] $824=6+(((6+6)/6)^{(6+6)}-6)/(6-6/6).$
\item [] $825=(6+6/6)\times (6+666/6)+6.$
\item [] $826=66+6\times 6\times (6+6)-(6+6)/6.$
\item [] $827=66+6\times 6\times (6+6)-6/6.$
\item [] $828=6\times (66+66+6).$
\item [] $829=66+6\times 6\times (6+6)+6/6.$
\item [] $830=66+6\times 6\times (6+6)+(6+6)/6.$
\item [] $831=66+6\times 6+(6\times 6/(6+6))^6.$
\item [] $832=(6+6+6/6)\times ((6+6)/6)^6.$
\item [] $833=66\times (6+6)+6\times 6+6-6/6.$
\item [] $834=66\times (6+6)+6\times 6+6.$
\item [] $835=66\times (6+6)+6\times 6+6+6/6.$
\item [] $836=(66\times 66+666-6)/6.$
\item [] $837=(66\times 66+666)/6.$
\item [] $838=(66\times 66+666+6)/6.$
\item [] $839=66+6\times 6\times (6+6)+66/6.$
\item [] $840=(6+6)\times (6+((6+6)/6)^6).$
\item [] $841=(6\times 6-6-6/6)^{((6+6)/6)}.$
\item [] $842=(66\times (66+66/6)+6)/6-6.$
\item [] $843=(66\times 66+666)/6+6.$
\item [] $844=66\times (6+6)-6-6+((6+6)/6)^6.$
\item [] $845=(6+6+6/6)\times (66-6/6).$
\item [] $846=666+6\times (6\times 6-6).$
\item [] $847=66\times (66+66/6)/6.$
\item [] $848=(66\times (66+66/6)+6)/6.$
\item [] $849=(66\times 66+666)/6+6+6.$
\item [] $850=66\times (6+6)-6+((6+6)/6)^6.$
\item [] $851=(6+6+6/6)\times (66-6/6)+6.$
\item [] $852=66+66\times (6+6)-6.$
\item [] $853=6\times (6+6)\times (6+6)-66/6.$
\item [] $854=(6+6/6)\times (666+66)/6.$
\item [] $855=6\times 6+(6+6/6)\times (6+666/6).$
\item [] $856=66\times (6+6)+((6+6)/6)^6.$
\item [] $857=66+66\times (6+6)-6/6.$
\item [] $858=66+66\times (6+6).$
\item [] $859=66+66\times (6+6)+6/6.$
\item [] $860=66+66\times (6+6)+(6+6)/6.$
\item [] $861=66+(((6/((6+6)/6))^6)+66).$
\item [] $862=6\times (6+6)\times (6+6)-(6+6)/6.$
\item [] $863=6\times (6+6)\times (6+6)-6/6.$
\item [] $864=6\times (6+6)\times (6+6).$
\item [] $865=6\times (6+6)\times (6+6)+6/6.$
\item [] $866=6\times (6+6)\times (6+6)+(6+6)/6.$
\item [] $867=6\times (6+6)\times (6+6)+6\times 6/(6+6).$
\item [] $868=6\times (6+6)\times (6+6)+6-(6+6)/6.$
\item [] $869=6\times (6+6)\times (6+6)+6-6/6.$
\item [] $870=6\times (6+6)\times (6+6)+6.$
\item [] $871=6\times (6+6)\times (6+6)+6+6/6.$
\item [] $872=6\times (6+6)\times (6+6)+6+(6+6)/6.$
\item [] $873=(6+6)\times (6+6)+(6\times 6/(6+6))^6.$
\item [] $874=6\times (6+6)\times (6+6)+(66-6)/6.$
\item [] $875=6\times (6+6)\times (6+6)+66/6.$
\item [] $876=6\times (6+6)\times (6+6)+6+6.$
\item [] $877=6\times (6+6)\times (6+6)+6+6+6/6.$
\item [] $878=6\times (6+6)\times (6+6)+6+6+(6+6)/6.$
\item [] $879=(6+6)\times (6+6)+6+(6\times 6/(6+6))^6.$
\item [] $880=666+6\times 6\times 6-(6+6)/6.$
\item [] $881=666+6\times 6\times 6-6/6.$
\item [] $882=666+6\times 6\times 6.$
\item [] $883=666+6\times 6\times 6+6/6.$
\item [] $884=666+6\times 6\times 6+(6+6)/6.$
\item [] $885=(6-6/6)\times (66+666/6).$
\item [] $886=666+6\times 6\times 6+6-(6+6)/6.$
\item [] $887=666+6\times 6\times 6+6-6/6.$
\item [] $888=666+6\times 6\times 6+6.$
\item [] $889=666+6\times 6\times 6+6+6/6.$
\item [] $890=666+6\times 6\times 6+6+(6+6)/6.$
\item [] $891=(6-6/6)\times (66+666/6)+6.$
\item [] $892=66+6\times 6\times (6+6)+((6+6)/6)^6.$
\item [] $893=666+6\times 66+6\times 6/6.$
\item [] $894=6\times ((6+6)\times (6+6)+6)-6.$
\item [] $895=6666/6-6\times 6\times 6.$
\item [] $896=(6/6+6)\times ((6+6)/6)^{(6+6/6)}.$
\item [] $897=66\times (6+6)-6+666/6.$
\item [] $898=6\times ((6+6)\times (6+6)+6)-(6+6)/6.$
\item [] $899=6\times ((6+6)\times (6+6)+6)-6/6.$
\item [] $900=6\times ((6+6)\times (6+6)+6).$
\item [] $901=6\times ((6+6)\times (6+6)+6)+6/6.$
\item [] $902=6\times ((6+6)\times (6+6)+6)+(6+6)/6.$
\item [] $903=66\times (6+6)+666/6.$
\item [] $904=66\times (6+6)+(666+6)/6.$
\item [] $905=6\times ((6+6)\times (6+6)+6)+6-6/6.$
\item [] $906=6\times ((6+6)\times (6+6)+6)+6.$
\item [] $907=6\times ((6+6)\times (6+6)+6)+6+6/6.$
\item [] $908=6\times ((6+6)\times (6+6)+6)+6+(6+6)/6.$
\item [] $909=66\times (6+6)+6+666/6.$
\item [] $910=(6+6/6)\times (66+((6+6)/6)^6).$
\item [] $911=6\times ((6+6)\times (6+6)+6)+66/6.$
\item [] $912=6\times ((6+6)\times (6+6)+6)+6+6.$
\item [] $913=66\times (66+6+66/6)/6.$
\item [] $914=66\times (6+6)+(666+66)/6.$
\item [] $915=66\times (6+6)+6+6+666/6.$
\item [] $916=(6+6/6)\times (66+((6+6)/6)^6)+6.$
\item [] $917=666+6\times (6\times 6+6)-6/6.$
\item [] $918=666+6\times (6\times 6+6).$
\item [] $919=666+6\times (6\times 6+6)+6/6.$
\item [] $920=666+6\times (6\times 6+6)+(6+6)/6.$
\item [] $921=66\times (6+6)+6+6+6+666/6.$
\item [] $922=66+66\times (6+6)+((6+6)/6)^6.$
\item [] $923=(6+6+6/6)\times (66+6-6/6).$
\item [] $924=66\times (6+6+(6+6)/6).$
\item [] $925=66\times (6+6+(6+6)/6)+6/6.$
\item [] $926=66\times (6+6+(6+6)/6)+(6+6)/6.$
\item [] $927=66+66+66+(6\times 6/(6+6))^6.$
\item [] $928=6\times (6+6)\times (6+6)+((6+6)/6)^6.$
\item [] $929=66+6\times (6+6)\times (6+6)-6/6.$
\item [] $930=66+6\times (6+6)\times (6+6).$
\item [] $931=66+6\times (6+6)\times (6+6)+6/6.$
\item [] $932=66+6\times (6+6)\times (6+6)+(6+6)/6.$
\item [] $933=6\times 6\times 6-6-6+(6\times 6/(6+6))^6.$
\item [] $934=(6+6)\times (66+6+6)-(6+6)/6.$
\item [] $935=(6+6)\times (66+6+6)-6/6.$
\item [] $936=(6+6)\times (66+6+6).$
\item [] $937=(6+6)\times (66+6+6)+6/6.$
\item [] $938=(6+6)\times (66+6+6)+(6+6)/6.$
\item [] $939=6\times 6\times 6+(6\times 6/(6+6))^6-6.$
\item [] $940=(6+6)\times (66+6+6)+6-(6+6)/6.$
\item [] $941=(6+6)\times (66+6+6)+6-6/6.$
\item [] $942=(6+6)\times (66+6+6)+6.$
\item [] $943=(6+6)\times (66+6+6)+6+6/6.$
\item [] $944=(6+6)\times (66+6+6)+6+(6+6)/6.$
\item [] $945=6\times 6\times 6+(6\times 6/(6+6))^6.$
\item [] $946=6\times 6\times 6+(6\times 6/(6+6))^6+6/6.$
\item [] $947=(6+6)\times (66+6+6)+66/6.$
\item [] $948=(6+6)\times (66+6+6)+6+6.$
\item [] $949=(6+6+6/6)\times (66+6+6/6).$
\item [] $950=(6\times 6-66/6)\times (6\times 6+(6+6)/6).$
\item [] $951=6\times 6\times 6+6+(6\times 6/(6+6))^6.$
\item [] $952=6\times 66+(6666+6)/(6+6).$
\item [] $953=(6+6)\times (66+6+6)+6+66/6.$
\item [] $954=666+6\times (6\times 6+6+6).$
\item [] $955=(6\times 6-6+6/6)^{((6+6)/6)}-6.$
\item [] $956=(6\times 6-66/6)\times (6\times 6+(6+6)/6)+6.$
\item [] $957=66\times (6\times (6\times 6-6)-6)/(6+6).$
\item [] $958=((6+6)/6)^{((66-6)/6)}-66.$
\item [] $959=(6+6/6)\times (66+66+6-6/6).$
\item [] $960=(66-6)\times (6+(66-6)/6).$
\item [] $961=(6\times 6-6+6/6)^{((6+6)/6)}.$
\item [] $962=(6\times 6-6+6/6)^{((6+6)/6)}+6/6.$
\item [] $963=66\times (6\times (6\times 6-6)-6)/(6+6)+6.$
\item [] $964=((6+6)/6)^{((66-6)/6)}-66+6.$
\item [] $965=66+6\times (6+6)\times (6+6)+6-6/6.$
\item [] $966=66+6\times (6+6)\times (6+6)+6.$
\item [] $967=(6\times 6-6+6/6)^{((6+6)/6)}+6.$
\item [] $968=66\times (66+(66+66)/6)/6.$
\item [] $969=6\times 6\times (6\times 6-6)-666/6.$
\item [] $970=6\times 6\times (6\times 6-6)-(666-6)/6.$
\item [] $971=(6+6+6)\times (66-6-6)-6/6.$
\item [] $972=(6+6+6)\times (66-6-6).$
\item [] $973=(6+6+6)\times (66-6-6)+6/6.$
\item [] $974=(6+6+6)\times (66-6-6)+(6+6)/6.$
\item [] $975=6\times (6+6)\times (6+6)+666/6.$
\item [] $976=(666+66)\times (6+(6+6)/6)/6.$
\item [] $977=(6+6+6)\times (66-6-6)+6-6/6.$
\item [] $978=(6+6+6)\times (66-6-6)+6.$
\item [] $979=6666/6-66-66.$
\item [] $980=(6\times 6-6/6)\times (6+(66+66)/6).$
\item [] $981=6\times (6\times 6+6)+(6\times 6/(6+6))^6.$
\item [] $982=((6+6)/6)^{((66-6)/6)}-6\times 6-6.$
\item [] $983=(6+6+6)\times (66-6-6)+66/6.$
\item [] $984=(6+6+6)\times (66-6-6)+6+6.$
\item [] $985=6-66-66+6666/6.$
\item [] $986=(6+66/6)\times (((6+6)/6)^6-6).$
\item [] $987=6\times (66\times (6\times 6-6)-6)/(6+6).$
\item [] $988=((6+6)/6)^{((66-6)/6)}-6\times 6.$
\item [] $989=666+6\times (66-6-6)-6/6.$
\item [] $990=666+6\times (66-6-6).$
\item [] $991=666+6\times (66-6-6)+6/6.$
\item [] $992=(6+66/6)\times (((6+6)/6)^6-6)+6.$
\item [] $993=6\times (66\times (6\times 6-6)+6)/(6+6).$
\item [] $994=((66-6)/6)^(6\times 6/(6+6))-6.$
\item [] $995=666+6\times 66-66-6/6.$
\item [] $996=(6+6)\times (66+6+66/6).$
\item [] $997=6\times 6+(6\times 6-6+6/6)^{((6+6)/6)}.$
\item [] $998=666+6\times 66-((6+6)/6)^6.$
\item [] $999=6\times 666/(6-(6+6)/6).$
\item [] $1000=((66-6)/6)^(6\times 6/(6+6)).$
\end{itemize}
\end{multicols}
}

\section{\textbf{Representations Using Number 7}}

{\footnotesize
\begin{multicols}{3}
\begin{itemize}
\item [] $101=7777/77.$
\item [] $102=77+7+7+77/7.$
\item [] $103=(777-7)/7-7.$
\item [] $104=777/7-7.$
\item [] $105=7\times (7+7)+7.$
\item [] $106=7\times (7+7)+7+7/7.$
\item [] $107=(777-77)/7+7.$
\item [] $108=7777/77+7.$
\item [] $109=(777-7-7)/7.$
\item [] $110=(777-7)/7.$
\item [] $111=777/7.$
\item [] $112=(777+7)/7.$
\item [] $113=(7+7+777)/7.$
\item [] $114=((7+7)/7)^7-7-7.$
\item [] $115=(777+77)/7-7.$
\item [] $116=7\times (7+7)+7+77/7.$
\item [] $117=(777-7)/7+7.$
\item [] $118=7+777/7.$
\item [] $119=7\times 7+77-7.$
\item [] $120=77+7\times 7-7+7/7.$
\item [] $121=((7+7)/7)^7-7.$
\item [] $122=(777+77)/7.$
\item [] $123=(777+77+7)/7.$
\item [] $124=7+7+(777-7)/7.$
\item [] $125=7+7+777/7.$
\item [] $126=7\times 7+77.$
\item [] $127=77+7\times 7+7/7.$
\item [] $128=((7+7)/7)^7.$
\item [] $129=7+(777+77)/7.$
\item [] $130=((7+7)/7)^7+(7+7)/7.$
\item [] $131=7+7+7+(777-7)/7.$
\item [] $132=7+7+7+777/7.$
\item [] $133=7\times 7+77+7.$
\item [] $134=77+7\times 7+7+7/7.$
\item [] $135=7+((7+7)/7)^7.$
\item [] $136=7+7/7+((7+7)/7)^7.$
\item [] $137=7\times 7+77+77/7.$
\item [] $138=7\times 7+77+(77+7)/7.$
\item [] $139=7\times (7+7+7)-7-7/7.$
\item [] $140=7\times (7+7+7)-7.$
\item [] $141=((7+7)\times (77-7)+7)/7.$
\item [] $142=7+7+((7+7)/7)^7.$
\item [] $143=77+77-77/7.$
\item [] $144=(7+7/7)\times (7+77/7).$
\item [] $145=7\times (7+7+7)-(7+7)/7.$
\item [] $146=7\times (7+7+7)-7/7.$
\item [] $147=7\times (7+7+7).$
\item [] $148=7\times (7+7+7)+7/7.$
\item [] $149=7\times (7+7+7)+(7+7)/7.$
\item [] $150=777+7\times 77/77.$
\item [] $151=7+(7+77/7)\times (7+7/7).$
\item [] $152=77+77-(7+7)/7.$
\item [] $153=77+77-7/7.$
\item [] $154=77+77.$
\item [] $155=77+77+7/7.$
\item [] $156=77+77+(7+7)/7.$
\item [] $157=77+77+(7+7+7)/7.$
\item [] $158=7\times (7+7+7)+77/7.$
\item [] $159=7\times 7+(777-7)/7.$
\item [] $160=777+7\times 7/7.$
\item [] $161=77+77+7.$
\item [] $162=77+77+7+7/7.$
\item [] $163=7+(7+7)\times (77+7/7)/7.$
\item [] $164=77+77+(77-7)/7.$
\item [] $165=77\times (7+7+7/7)/7.$
\item [] $166=7\times 7+7+(777-7)/7.$
\item [] $167=7\times 7+7+777/7.$
\item [] $168=77+77+7+7.$
\item [] $169=77+77+7+7+7/7.$
\item [] $170=7\times 7-7+((7+7)/7)^7.$
\item [] $171=7\times 7+(777+77)/7.$
\item [] $172=77+77+7+77/7.$
\item [] $173=77+7\times (7+7)+(7+7)/7.$
\item [] $174=77+7\times (7+7)-7/7.$
\item [] $175=77+7\times (7+7).$
\item [] $176=77+7\times (7+7)+7/7.$
\item [] $177=7\times 7+((7+7)/7)^7.$
\item [] $178=77+7777/77.$
\item [] $179=((7+7)/7)^{(7+7/7)}-77.$
\item [] $180=77-7+(777-7)/7.$
\item [] $181=77-7+777/7.$
\item [] $182=77+7+7\times (7+7).$
\item [] $183=7\times (7+7)+77+7+7/7.$
\item [] $184=7\times 7+7+((7+7)/7)^7.$
\item [] $185=(7+7)\times (7+7)-77/7.$
\item [] $186=(7+7)\times (7+7)-(77-7)/7.$
\item [] $187=77+(777-7)/7.$
\item [] $188=77+777/7.$
\item [] $189=(7+7)\times (7+7)-7.$
\item [] $190=(7+7)\times (7+7)-7+7/7.$
\item [] $191=77+((7+7)/7)^7-7-7.$
\item [] $192=(7\times 7-7/7)\times (77/7-7).$
\item [] $193=(7+7)\times (7+7)-(7+7+7)/7.$
\item [] $194=(7+7)\times (7+7)-(7+7)/7.$
\item [] $195=(7+7)\times (7+7)-7/7.$
\item [] $196=(7+7)\times (7+7).$
\item [] $197=(7+7)\times (7+7)+7/7.$
\item [] $198=77/7\times (7+77/7).$
\item [] $199=77+(777+77)/7.$
\item [] $200=(7+7)\times (7+7)+77/7-7.$
\item [] $201=(7+7)\times (7+7)+7-(7+7)/7.$
\item [] $202=(7+7)\times (7+7)+7-7/7.$
\item [] $203=(7+7)\times (7+7)+7.$
\item [] $204=(7+7)\times (7+7)+7+7/7.$
\item [] $205=77+((7+7)/7)^7.$
\item [] $206=(7+7)\times (7+7)+(77-7)/7.$
\item [] $207=(7+7)\times (7+7)+77/7.$
\item [] $208=7\times (7+7)+(777-7)/7.$
\item [] $209=7\times (7+7)+777/7.$
\item [] $210=(7+7)\times (7+7+7/7).$
\item [] $211=(7+7)\times (7+7)+7+7+7/7.$
\item [] $212=77+((7+7)/7)^7+7.$
\item [] $213=7\times (7+7+7)+77-77/7.$
\item [] $214=(7+7)\times (7+7)+7+77/7.$
\item [] $215=7\times 7\times 7-((7+7)/7)^7.$
\item [] $216=(7-7/7)^{((7+7+7)/7)}.$
\item [] $217=7\times (7\times 7-7)-77.$
\item [] $218=7\times (7\times 7-7)-77+7/7.$
\item [] $219=77+7+7+((7+7)/7)^7.$
\item [] $220=77\times (7+7+7-7/7)/7.$
\item [] $221=77+(7+77/7)\times (7+7/7).$
\item [] $222=777\times (7+7)/(7\times 7).$
\item [] $223=77+7\times (7+7+7)-7/7.$
\item [] $224=77+7\times (7+7+7).$
\item [] $225=(7+7+7/7)^{((7+7)/7)}.$
\item [] $226=7\times (7+7)+((7+7)/7)^7.$
\item [] $227=7\times (7\times 7-7-7)-7-77/7.$
\item [] $228=(7\times 7-77/7)\times (7-7/7).$
\item [] $229=(7+7)\times 777/(7\times 7)+7.$
\item [] $230=(77\times (7+7+7)-7)/7.$
\item [] $231=77+77+77.$
\item [] $232=7\times 7\times 7-777/7.$
\item [] $233=7\times 7\times 7-(777-7)/7.$
\item [] $234=(7+7+7)/7\times (77+7/7).$
\item [] $235=7\times (7+7+7)+77+77/7.$
\item [] $236=(7+7)\times (7+777/7)/7.$
\item [] $237=(777-7/7)-7\times 77.$
\item [] $238=777-77\times 7.$
\item [] $239=777-77\times 7+7/7.$
\item [] $240=(7-(7+7)/7)\times (7\times 7-7/7).$
\item [] $241=(7-7/7)\times (7\times 7-7)-77/7.$
\item [] $242=77\times (7+7+7+7/7)/7.$
\item [] $243=7\times (7\times 7-7-7)-(7+7)/7.$
\item [] $244=7\times (7\times 7-7-7)-7/7.$
\item [] $245=7\times (7\times 7-7-7).$
\item [] $246=7\times (7\times 7-7-7)+7/7.$
\item [] $247=7\times (7\times 7-7-7)+(7+7)/7.$
\item [] $248=77+(7\times 7\times 7\times 7-7)/(7+7).$
\item [] $249=((7+7)/7)^{(7+7/7)}-7.$
\item [] $250=(7\times 7+7/7)\times (7-(7+7)/7).$
\item [] $251=7\times (7\times 7-7-7)+7-7/7.$
\item [] $252=(7+7)\times (7+77/7).$
\item [] $253=7\times (7\times 7-7-7)+7+7/7.$
\item [] $254=7\times 7+77+((7+7)/7)^7.$
\item [] $255=7\times 7\times 7-77-77/7.$
\item [] $256=((7+7)/7)^{(7+7/7)}.$
\item [] $257=((7+7)/7)^{(7+7/7)}+7/7.$
\item [] $258=(7-7/7)\times (7\times 7-7+7/7).$
\item [] $259=7\times 7\times 7-77-7.$
\item [] $260=7\times 7\times 7-77-7+7/7.$
\item [] $261=7\times 7\times 7-77-7+(7+7)/7.$
\item [] $262=7\times 7\times 7-77+7-77/7.$
\item [] $263=7+((7+7)/7)^{(7+7/7)}.$
\item [] $264=7\times 7\times 7-77-(7+7)/7.$
\item [] $265=7\times 7\times 7-77-7/7.$
\item [] $266=7\times 7\times 7-77.$
\item [] $267=7\times 7\times 7-77+7/7.$
\item [] $268=7\times 7\times 7-77+(7+7)/7.$
\item [] $269=7\times (7\times 77-7/7)/(7+7).$
\item [] $270=(7+7+7/7)\times (7+77/7).$
\item [] $271=7\times 7\times 7+7-77-(7+7)/7.$
\item [] $272=7\times 7\times 7-77+7-7/7.$
\item [] $273=7\times 7\times 7-77+7.$
\item [] $274=7\times 7\times 7-77+7+7/7.$
\item [] $275=(7+7+77/7)\times 77/7.$
\item [] $276=7\times (7\times 7-7)-7-77/7.$
\item [] $277=7\times 7\times 7-77+77/7.$
\item [] $278=7\times 7\times 7-77+(77+7)/7.$
\item [] $279=7\times (7\times 7-7)-7-7-7/7.$
\item [] $280=7\times (7\times 7-7)-7-7.$
\item [] $281=7\times (7\times 7-7)-7-7+7/7.$
\item [] $282=(7\times 7-(7+7)/7)\times (7-7/7).$
\item [] $283=7\times (7\times 7-7)-77/7.$
\item [] $284=7\times (7\times 7-7)-(77-7)/7.$
\item [] $285=7\times (7\times 7-7)-7-(7+7)/7.$
\item [] $286=7\times (7\times 7-7)-7-7/7.$
\item [] $287=7\times (7\times 7-7)-7.$
\item [] $288=7\times (7\times 7-7)-7+7/7.$
\item [] $289=7\times (7\times 7-7)-7+(7+7)/7.$
\item [] $290=7\times (7\times 7-7)+7-77/7.$
\item [] $291=7\times (7\times 7-7)-(7+7+7)/7.$
\item [] $292=7\times (7\times 7-7)-(7+7)/7.$
\item [] $293=7\times (7\times 7-7)-7/7.$
\item [] $294=7\times (7\times 7-7).$
\item [] $295=7\times (7\times 7-7)+7/7.$
\item [] $296=7\times (7\times 7-7)+(7+7)/7.$
\item [] $297=7\times (7\times 7-7)+(7+7+7)/7.$
\item [] $298=7\times (7\times 7-7)-7+77/7.$
\item [] $299=7\times (7\times 7-7)+7-(7+7)/7.$
\item [] $300=7\times (7\times 7-7)+7-7/7.$
\item [] $301=7\times (7\times 7-7)+7.$
\item [] $302=7\times (7\times 7-7)+7+7/7.$
\item [] $303=7\times (7\times 7-7)+7+(7+7)/7.$
\item [] $304=(7\times 7-77/7)\times (7+7/7).$
\item [] $305=7\times (7\times 7-7)+77/7.$
\item [] $306=7\times (7\times 7-7)+(77+7)/7.$
\item [] $307=7\times (7\times 7-7)+7+7-7/7.$
\item [] $308=77\times (77/7-7).$
\item [] $309=7\times (7\times 7-7)+7+7+7/7.$
\item [] $310=7\times (7\times 7-7)+7+7+(7+7)/7.$
\item [] $311=7\times 7\times 7-7-7-7-77/7.$
\item [] $312=7\times (7\times 7-7)+7+77/7.$
\item [] $313=7\times (7\times 7-7)+7+(77+7)/7.$
\item [] $314=7\times (7\times 7+7)-77-7/7.$
\item [] $315=7\times (7\times 7+7)-77.$
\item [] $316=7\times (7\times 7+7)-77+7/7.$
\item [] $317=(7+77/7)^{((7+7)/7)}-7.$
\item [] $318=7\times 7\times 7-7-7-77/7.$
\item [] $319=7\times 7\times 7-7-7-(77-7)/7.$
\item [] $320=(7+7/7)\times (7\times 7-7-(7+7)/7).$
\item [] $321=7\times 7\times 7-(77+77)/7.$
\item [] $322=7\times 7\times 7-7-7-7.$
\item [] $323=7\times 7\times 7-7-7-7+7/7.$
\item [] $324=(77/7+7)^{((7+7)/7)}.$
\item [] $325=7\times 7\times 7-7-77/7.$
\item [] $326=7\times 7\times 7-7-(77-7)/7.$
\item [] $327=7\times 7\times 7-7-7-(7+7)/7.$
\item [] $328=7\times 7\times 7-7-7-7/7.$
\item [] $329=7\times 7\times 7-7-7.$
\item [] $330=7\times 7\times 7-7-7+7/7.$
\item [] $331=7\times 7\times 7-(77+7)/7.$
\item [] $332=7\times 7\times 7-77/7.$
\item [] $333=7\times 7\times 7-(77-7)/7.$
\item [] $334=7\times 7\times 7-7-(7+7)/7.$
\item [] $335=7\times 7\times 7-7-7/7.$
\item [] $336=7\times 7\times 7-7.$
\item [] $337=7\times 7\times 7-7+7/7.$
\item [] $338=7\times 7\times 7-7+(7+7)/7.$
\item [] $339=7\times 7\times 7+7-77/7.$
\item [] $340=7\times 7\times 7-(7+7+7)/7.$
\item [] $341=7\times 7\times 7-(7+7)/7.$
\item [] $342=7\times 7\times 7-7/7.$
\item [] $343=7\times 7\times 7.$
\item [] $344=7\times 7\times 7+7/7.$
\item [] $345=7\times 7\times 7+(7+7)/7.$
\item [] $346=7\times 7\times 7+(7+7+7)/7.$
\item [] $347=7\times 7\times 7-7+77/7.$
\item [] $348=7\times 7\times 7+7-(7+7)/7.$
\item [] $349=7\times 7\times 7+7-7/7.$
\item [] $350=7\times 7\times 7+7.$
\item [] $351=7\times 7\times 7+7+7/7.$
\item [] $352=7\times 7\times 7+7+(7+7)/7.$
\item [] $353=7\times 7\times 7+(77-7)/7.$
\item [] $354=7\times 7\times 7+77/7.$
\item [] $355=7\times 7\times 7+(77+7)/7.$
\item [] $356=7\times 7\times 7+7+7-7/7.$
\item [] $357=7\times 7\times 7+7+7.$
\item [] $358=7\times 7\times 7+7+7+7/7.$
\item [] $359=7\times 7\times 7+7+7+(7+7)/7.$
\item [] $360=(7\times 7+77/7)\times (7-7/7).$
\item [] $361=7\times 7\times 7+7+77/7.$
\item [] $362=7\times 7\times 7+7+(77+7)/7.$
\item [] $363=7\times 7\times 7+7+7+7-7/7.$
\item [] $364=7\times 7\times 7+7+7+7.$
\item [] $365=7\times 7\times 7+(77+77)/7.$
\item [] $366=7\times 7\times 7+7+7+7+(7+7)/7.$
\item [] $367=7\times 7\times 7+7+7+(77-7)/7.$
\item [] $368=7\times 7\times 7+7+7+77/7.$
\item [] $369=7\times 7\times 7+7+7+(77+7)/7.$
\item [] $370=7\times (7\times 7-7)+77-7/7.$
\item [] $371=7\times (7\times 7-7)+77.$
\item [] $372=7\times (7\times 7-7)+77+7/7.$
\item [] $373=7\times (7\times 7+7)-7-(77+7)/7.$
\item [] $374=7\times (7\times 7+7)-7-77/7.$
\item [] $375=7\times 7\times 7+7+7+7+77/7.$
\item [] $376=(7\times 7-(7+7)/7)\times (7+7/7).$
\item [] $377=7\times (7\times 7+7)-7-7-7/7.$
\item [] $378=7\times (7\times 7+7)-7-7.$
\item [] $379=7\times (7\times 7+7)-7-7+7/7.$
\item [] $380=(7-(7+7)/7)\times (77-7/7).$
\item [] $381=7\times (7\times 7+7)-77/7.$
\item [] $382=7\times (7\times 7+7)-(77-7)/7.$
\item [] $383=7\times (7\times 7+7)-7-(7+7)/7.$
\item [] $384=7\times (7\times 7+7)-7-7/7.$
\item [] $385=7\times (7\times 7+7)-7.$
\item [] $386=7\times (7\times 7+7)-7+7/7.$
\item [] $387=7\times (7\times 7+7)-7+(7+7)/7.$
\item [] $388=7\times (7\times 7+7)+7-77/7.$
\item [] $389=7\times (7\times 7+7)-(7+7+7)/7.$
\item [] $390=7\times (7\times 7+7)-(7+7)/7.$
\item [] $391=7\times (7\times 7+7)-7/7.$
\item [] $392=7\times (7\times 7+7).$
\item [] $393=7\times (7\times 7+7)+7/7.$
\item [] $394=7\times (7\times 7+7)+(7+7)/7.$
\item [] $395=7\times (7\times 7+7)+(7+7+7)/7.$
\item [] $396=7\times (7\times 7+7)-7+77/7.$
\item [] $397=7\times (7\times 7+7)+7-(7+7)/7.$
\item [] $398=7\times (7\times 7+7)+7-7/7.$
\item [] $399=7\times (7\times 7+7)+7.$
\item [] $400=7\times (7\times 7+7)+7+7/7.$
\item [] $401=7\times (7\times 7+7)+(7+7)/7+7.$
\item [] $402=7\times (7\times 7+7)+(77-7)/7.$
\item [] $403=7\times (7\times 7+7)+77/7.$
\item [] $404=7\times (7\times 7+7)+(77+7)/7.$
\item [] $405=7\times (7\times 7+7)+7+7-7/7.$
\item [] $406=7\times (7\times 7+7)+7+7.$
\item [] $407=7\times (7\times 7+7)+7+7+7/7.$
\item [] $408=(7\times 7+(7+7)/7)\times (7+7/7).$
\item [] $409=7\times 7\times 7+77-77/7.$
\item [] $410=7\times (7\times 7+7)+7+77/7.$
\item [] $411=7\times 77-((7+7)/7)^7.$
\item [] $412=7\times 7\times 7-7+77-7/7.$
\item [] $413=7\times (77-7)-77.$
\item [] $414=7\times (77-7)-77+7/7.$
\item [] $415=(77+7-7/7)\times (7-(7+7)/7).$
\item [] $416=(77\times 77-7)/(7+7)-7.$
\item [] $417=77\times 7-(777+77)/7.$
\item [] $418=77/7\times (7\times 7-77/7).$
\item [] $419=7\times 7\times 7+77-7/7.$
\item [] $420=7\times 7\times 7+77.$
\item [] $421=7\times 7\times 7+77+7/7.$
\item [] $422=7\times 7\times 7+77+(7+7)/7.$
\item [] $423=7\times 7\times 7+77+(7+7+7)/7.$
\item [] $424=7\times (77-7)-77+77/7.$
\item [] $425=7+7\times 7\times 7+77-(7+7)/7.$
\item [] $426=7\times 7\times 7+77+7-7/7.$
\item [] $427=7\times 7\times 7+77+7.$
\item [] $428=77\times 7-777/7.$
\item [] $429=77\times 7-(777-7)/7.$
\item [] $430=7\times (7\times 7+7+7)-77/7.$
\item [] $431=7\times 7\times 7+77+77/7.$
\item [] $432=(7\times 7-7/7)\times (7+(7+7)/7).$
\item [] $433=7\times (77-7-7)-7-7/7.$
\item [] $434=7\times (77-7-7)-7.$
\item [] $435=7\times 77-(777/7-7).$
\item [] $436=7\times 77+7-(777-7)/7.$
\item [] $437=7\times (77-7-7)+7-77/7.$
\item [] $438=7\times 7\times 7+77+7+77/7.$
\item [] $439=7\times (77-7-7)-(7+7)/7.$
\item [] $440=7\times (77-7-7)-7/7.$
\item [] $441=7\times (77-7-7).$
\item [] $442=7\times (77-7)-7\times 7+7/7.$
\item [] $443=7\times (77-7-7)+(7+7)/7.$
\item [] $444=777\times (77/7-7)/7.$
\item [] $445=7\times (77-7-7)-7+77/7.$
\item [] $446=(7\times 7+7)\times (7+7/7)-(7+7)/7.$
\item [] $447=7\times (77-7-7)+7-7/7.$
\item [] $448=7\times (77-7-7)+7.$
\item [] $449=7\times (77-7-7)+7+7/7.$
\item [] $450=(7\times 7+7/7)\times (7+(7+7)/7).$
\item [] $451=7\times 77-77-77/7.$
\item [] $452=7\times (7\times 7+7+7)+77/7.$
\item [] $453=7\times 7\times 7+(777-7)/7.$
\item [] $454=7\times 777+7\times 7/7.$
\item [] $455=7\times 77-77-7.$
\item [] $456=(7-7/7)\times (77-7/7).$
\item [] $457=7\times 77-77-7+(7+7)/7.$
\item [] $458=7\times 77-77+7-77/7.$
\item [] $459=7\times 77-77-(7+7+7)/7.$
\item [] $460=7\times 77-77-(7+7)/7.$
\item [] $461=7\times 77-77-7/7.$
\item [] $462=7\times 77-77.$
\item [] $463=7\times 77-77+7/7.$
\item [] $464=7\times 77-77+(7+7)/7.$
\item [] $465=(77-7)\times 7-7-7-77/7.$
\item [] $466=7\times 77+77/7-77-7.$
\item [] $467=7\times 77-77+7-(7+7)/7.$
\item [] $468=(77+7/7)\times (7-7/7).$
\item [] $469=7\times 77-77+7.$
\item [] $470=7\times 77-77+7+7/7.$
\item [] $471=7\times 7\times 7+((7+7)/7)^7.$
\item [] $472=7\times (77-7)-7-77/7.$
\item [] $473=7\times 77-77+77/7.$
\item [] $474=(77+(7+7)/7)\times (7-7/7).$
\item [] $475=7\times (77-7)-7-7-7/7.$
\item [] $476=7\times (77-7)-7-7.$
\item [] $477=7\times (77-7)-7-7+7/7.$
\item [] $478=7+((7+7)/7)^7+7\times 7\times 7.$
\item [] $479=7\times (77-7)-77/7.$
\item [] $480=(7\times 7-7/7)\times (77-7)/7.$
\item [] $481=7\times (77-7)-7-(7+7)/7.$
\item [] $482=7\times (77-7)-7-7/7.$
\item [] $483=7\times (77-7)-7.$
\item [] $484=7\times (77-7)-7+7/7.$
\item [] $485=7\times (77-7)-7+(7+7)/7.$
\item [] $486=7\times (77-7)+7-77/7.$
\item [] $487=7\times (77-7)-(7+7+7)/7.$
\item [] $488=7\times (77-7)-(7+7)/7.$
\item [] $489=7\times (77-7)-7/7.$
\item [] $490=7\times (77-7).$
\item [] $491=7\times (77-7)+7/7.$
\item [] $492=7\times (77-7)+(7+7)/7.$
\item [] $493=7\times (77-7)+(7+7+7)/7.$
\item [] $494=7\times (77-7)-7+77/7.$
\item [] $495=7\times (77-7)+7-(7+7)/7.$
\item [] $496=7\times (77-7)+7-7/7.$
\item [] $497=7\times (77-7)+7.$
\item [] $498=7\times (77-7)+7+7/7.$
\item [] $499=7\times (77-7)+7+(7+7)/7.$
\item [] $500=7\times (77-7)+(77-7)/7.$
\item [] $501=7\times (77-7)+77/7.$
\item [] $502=7\times (77-7)+(7+77)/7.$
\item [] $503=7\times (7\times 7+7)+777/7.$
\item [] $504=7\times (77-7)+7+7.$
\item [] $505=7\times (77-7)+7+7+7/7.$
\item [] $506=7\times (77-7)+7+7+(7+7)/7.$
\item [] $507=7\times 77-7-7-7-77/7.$
\item [] $508=7\times (77-7)+7+77/7.$
\item [] $509=7\times (77+7)-77-(7+7)/7.$
\item [] $510=7\times (77+7)-77-7/7.$
\item [] $511=7\times (77+7)-77.$
\item [] $512=((7+7)/7)^{(7+(7+7)/7)}.$
\item [] $513=7\times (77+7)-77+(7+7)/7.$
\item [] $514=7\times 77-7-7-77/7.$
\item [] $515=7\times (77-7)+7+7+77/7.$
\item [] $516=7\times 77-7-7-7-(7+7)/7.$
\item [] $517=7\times 77-7-7-7-7/7.$
\item [] $518=7\times 77-7-7-7.$
\item [] $519=7\times 77-7-7-7+7/7.$
\item [] $520=7\times 77-7-(77+7)/7.$
\item [] $521=7\times 77-7-77/7.$
\item [] $522=7\times 77-7-(77-7)/7.$
\item [] $523=7\times 77-7-7-(7+7)/7.$
\item [] $524=7\times 77-7-7-7/7.$
\item [] $525=7\times 77-7-7.$
\item [] $526=7\times 77-7-7+7/7.$
\item [] $527=77\times 7-(77+7)/7.$
\item [] $528=7\times 77-77/7.$
\item [] $529=7\times 77+(7-77)/7.$
\item [] $530=7\times 77-7-(7+7)/7.$
\item [] $531=7\times 77-7-7/7.$
\item [] $532=7\times 77-7.$
\item [] $533=7\times 77-7+7/7.$
\item [] $534=7\times 77-7+(7+7)/7.$
\item [] $535=7\times 77+7-77/7.$
\item [] $536=7\times 77-(7+7+7)/7.$
\item [] $537=7\times 77-(7+7)/7.$
\item [] $538=7\times 77-7/7.$
\item [] $539=7\times 77.$
\item [] $540=7\times 77+7/7.$
\item [] $541=7\times 77+(7+7)/7.$
\item [] $542=7\times 77+(7+7+7)/7.$
\item [] $543=7\times 77-7+77/7.$
\item [] $544=7\times 77+7-(7+7)/7.$
\item [] $545=7\times 77+7-7/7.$
\item [] $546=7\times 77+7.$
\item [] $547=7\times 77+7+7/7.$
\item [] $548=7\times 77+7+(7+7)/7.$
\item [] $549=7\times 77+(77-7)/7.$
\item [] $550=7\times 77+77/7.$
\item [] $551=7\times 77+(77+7)/7.$
\item [] $552=7\times 77+7+7-7/7.$
\item [] $553=7\times 77+7+7.$
\item [] $554=7\times 77+7+7+7/7.$
\item [] $555=7\times 77+7+7+(7+7)/7.$
\item [] $556=7\times 77+7+(77-7)/7.$
\item [] $557=7\times 77+7+77/7.$
\item [] $558=7\times 7\times (7+7)-((7+7)/7)^7.$
\item [] $559=7\times 77+7+7+7-7/7.$
\item [] $560=(77-7)\times (7+7/7).$
\item [] $561=7\times 77+7+7+7+7/7.$
\item [] $562=7+(7777-7)/(7+7).$
\item [] $563=7+(7777+7)/(7+7).$
\item [] $564=7\times 77+7+7+77/7.$
\item [] $565=7\times 77+7+7+(77+7)/7.$
\item [] $566=7\times (77-7)+77-7/7.$
\item [] $567=7\times (77-7)+77.$
\item [] $568=7\times (77-7)+77+7/7.$
\item [] $569=7\times (77+7)-7-(77+7)/7.$
\item [] $570=7\times (77+7)-7-77/7.$
\item [] $571=7\times 77+7+7+7+77/7.$
\item [] $572=7777/7-7\times 77.$
\item [] $573=7\times (77+7)-7-7-7/7.$
\item [] $574=7\times (77+7)-7-7.$
\item [] $575=7\times (77+7)-7-7+7/7.$
\item [] $576=7\times (77+7)-(77+7)/7.$
\item [] $577=7\times (77+7)-77/7.$
\item [] $578=7\times (77+7)-(77-7)/7.$
\item [] $579=7\times (77+7)-7-(7+7)/7.$
\item [] $580=7\times (77+7)-7-7/7.$
\item [] $581=7\times (77+7)-7.$
\item [] $582=7\times (77+7)-7+7/7.$
\item [] $583=7\times (77+7)-7+(7+7)/7.$
\item [] $584=7\times (77+7)+7-77/7.$
\item [] $585=7\times (77+7)-(7+7+7)/7.$
\item [] $586=7\times (77+7)-(7+7)/7.$
\item [] $587=7\times (77+7)-7/7.$
\item [] $588=7\times (77+7).$
\item [] $589=7\times (77+7)+7/7.$
\item [] $590=7\times (77+7)+(7+7)/7.$
\item [] $591=7\times (77+7)+(7+7+7)/7.$
\item [] $592=7\times (77+7)-7+77/7.$
\item [] $593=7\times (77+7)+7-(7+7)/7.$
\item [] $594=7\times (77+7)+7-7/7.$
\item [] $595=7\times (77+7)+7.$
\item [] $596=7\times (77+7)+7+7/7.$
\item [] $597=7\times (77+7)+7+(7+7)/7.$
\item [] $598=7\times (77+7)+(77-7)/7.$
\item [] $599=7\times (77+7)+77/7.$
\item [] $600=(7+7/7)\times (77-(7+7)/7).$
\item [] $601=7\times (77+7)+7+7-7/7.$
\item [] $602=7\times (77+7)+7+7.$
\item [] $603=7\times (77+7)+7+7+7/7.$
\item [] $604=7\times (77+7)+7+7+(7+7)/7.$
\item [] $605=7\times 77+77-77/7.$
\item [] $606=7\times (77+7)+7+77/7.$
\item [] $607=7\times (77+7)+7+(77+7)/7.$
\item [] $608=(7+7/7)\times (77-7/7).$
\item [] $609=7\times 77+77-7.$
\item [] $610=77\times (7+7/7)-7+7/7.$
\item [] $611=7\times 77+77-7+(7+7)/7.$
\item [] $612=7\times 77+77+7-77/7.$
\item [] $613=7\times (77+7)+7+7+77/7.$
\item [] $614=7\times 77+77-(7+7)/7.$
\item [] $615=7\times 77+77-7/7.$
\item [] $616=7\times 77+77.$
\item [] $617=7\times 77+77+7/7.$
\item [] $618=7\times 77+77+(7+7)/7.$
\item [] $619=7\times 77+77+(7+7+7)/7.$
\item [] $620=7\times 77+77-7+77/7.$
\item [] $621=(77-7-7/7)\times (7+(7+7)/7).$
\item [] $622=7\times 77+77+7-7/7.$
\item [] $623=7\times 77+77+7.$
\item [] $624=(77+7/7)\times (7+7/7).$
\item [] $625=(7-(7+7)/7)^(77/7-7).$
\item [] $626=7\times (77+7+7)-77/7.$
\item [] $627=7\times 77+77+77/7.$
\item [] $628=7\times (7+7+77)-7-(7+7)/7.$
\item [] $629=7\times (77+7+7)-7-7/7.$
\item [] $630=7\times (77+7+7)-7.$
\item [] $631=7\times (7\times (7+7)-7)-7+7/7.$
\item [] $632=(7+7/7)\times (77+(7+7)/7).$
\item [] $633=7\times (77+7+7)+7-77/7.$
\item [] $634=7\times 77+77+7+77/7.$
\item [] $635=7\times (77+7+7)-(7+7)/7.$
\item [] $636=7\times (77+7+7)-7/7.$
\item [] $637=7\times (77+7+7).$
\item [] $638=7\times (77+7+7)+7/7.$
\item [] $639=7\times (77+7+7)+(7+7)/7.$
\item [] $640=(7-(7+7)/7)\times ((7+7)/7)^7.$
\item [] $641=7\times (77+7+7)-7+77/7.$
\item [] $642=777-((7+7)/7)^7-7.$
\item [] $643=7\times 77-7+777/7.$
\item [] $644=7\times (77+7+7)+7.$
\item [] $645=7\times (77+7+7)+7+7/7.$
\item [] $646=7\times (77+7+7)+7+(7+7)/7.$
\item [] $647=7+(7-(7+7)/7)\times ((7+7)/7)^7.$
\item [] $648=7\times (77+7+7)+77/7.$
\item [] $649=777-((7+7)/7)^7.$
\item [] $650=7\times 77+777/7.$
\item [] $651=777-7\times 7-77.$
\item [] $652=777-7-7-777/7.$
\item [] $653=77\times (77/7+7\times 7)/7-7.$
\item [] $654=7\times (77+7)+77-77/7.$
\item [] $655=7+7\times (77+7+7)+77/7.$
\item [] $656=777+7-((7+7)/7)^7.$
\item [] $657=77\times 7+7+777/7.$
\item [] $658=7\times (77+7)-7+77.$
\item [] $659=777-7-777/7.$
\item [] $660=77/7\times (7\times 7+77/7).$
\item [] $661=77\times 7+(777+77)/7.$
\item [] $662=77+7\times (77+7)-(7+7+7)/7.$
\item [] $663=77+7\times (77+7)-(7+7)/7.$
\item [] $664=77+7\times (77+7)-7/7.$
\item [] $665=(77+((77+7)\times 7)).$
\item [] $666=777-777/7.$
\item [] $667=((7\times 77)+(((7+7)/7)^7)).$
\item [] $668=(777+(((7-777)+7)/7)).$
\item [] $669=7\times 7\times (7+7)-7-(77-7)/7.$
\item [] $670=7\times 7\times (7+7)-7-7-(7+7)/7.$
\item [] $671=7\times 7\times (7+7)-7-7-7/7.$
\item [] $672=(77+7)\times (7+7/7).$
\item [] $673=777+7-777/7.$
\item [] $674=7\times 77+7+((7+7)/7)^7.$
\item [] $675=7\times 7\times (7+7)-77/7.$
\item [] $676=7\times 7\times (7+7)-(77-7)/7.$
\item [] $677=7\times 7\times (7+7)-7-(7+7)/7.$
\item [] $678=7\times 7\times (7+7)-7-7/7.$
\item [] $679=7\times 7\times (7+7)-7.$
\item [] $680=7\times 7\times (7+7)-7+7/7.$
\item [] $681=7\times 7\times (7+7)-7+(7+7)/7.$
\item [] $682=7\times 7\times (7+7)+7-77/7.$
\item [] $683=7\times 7\times (7+7)-(7+7+7)/7.$
\item [] $684=7\times 7\times (7+7)-(7+7)/7.$
\item [] $685=7\times 7\times (7+7)-7/7.$
\item [] $686=7\times 7\times (7+7).$
\item [] $687=7\times 7\times (7+7)+7/7.$
\item [] $688=7\times 7\times (7+7)+(7+7)/7.$
\item [] $689=777-77-77/7.$
\item [] $690=7\times 7\times (7+7)-7+77/7.$
\item [] $691=7\times 7\times (7+7)+7-(7+7)/7.$
\item [] $692=7\times 7\times (7+7)+7-7/7.$
\item [] $693=7\times 7\times (7+7)+7.$
\item [] $694=7\times 7\times (7+7)+7+7/7.$
\item [] $695=7\times 7\times (7+7)+7+(7+7)/7.$
\item [] $696=7\times 7\times (7+7)+(77-7)/7.$
\item [] $697=7\times 7\times (7+7)+77/7.$
\item [] $698=777-77-(7+7)/7.$
\item [] $699=777-77-7/7.$
\item [] $700=777-77.$
\item [] $701=777-77+7/7.$
\item [] $702=(77+7/7)\times (7+(7+7)/7).$
\item [] $703=7\times 7\times (7+7)+7-(7-77)/7.$
\item [] $704=7\times 7\times (7+7)+7+77/7.$
\item [] $705=777-77+7-(7+7)/7.$
\item [] $706=777-77+7-7/7.$
\item [] $707=777-77+7.$
\item [] $708=777-77+7+7/7.$
\item [] $709=777-77+7+(7+7)/7.$
\item [] $710=777-7\times 7-7-77/7.$
\item [] $711=777-77+77/7.$
\item [] $713=777+7+7-77-7/7.$
\item [] $713=777-77+7+7-7/7.$
\item [] $714=777-77+7+7.$
\item [] $715=777-77+7+7+7/7.$
\item [] $716=7\times (77+7)+((7+7)/7)^7.$
\item [] $717=777-7\times 7-77/7.$
\item [] $718=777-77+7+77/7.$
\item [] $719=777-7-7\times 7-(7+7)/7.$
\item [] $720=777-7\times 7-7-7/7.$
\item [] $721=777-7\times 7-7.$
\item [] $722=777-7\times 7-7+7/7.$
\item [] $723=777-7\times 7-7+(7+7)/7.$
\item [] $724=7\times ((7+7)\times 7+7)-77/7.$
\item [] $725=777-7\times 7-(7+7+7)/7.$
\item [] $726=777-7\times 7-(7+7)/7.$
\item [] $727=777-7\times 7-7/7.$
\item [] $728=777-7\times 7.$
\item [] $729=777-7\times 7+7/7.$
\item [] $730=777-7\times 7+(7+7)/7.$
\item [] $731=777-7\times 7+(7+7+7)/7.$
\item [] $732=777-7\times 7-7+77/7.$
\item [] $733=7\times (7\times (7+7)+7)-(7+7)/7.$
\item [] $734=7\times (7\times (7+7)+7)-7/7.$
\item [] $735=7\times (7\times (7+7)+7).$
\item [] $736=7\times (7\times (7+7)+7)+7/7.$
\item [] $737=7\times (7\times (7+7)+7)+(7+7)/7.$
\item [] $738=777-7\times 7+(77-7)/7.$
\item [] $739=777-7\times 7+77/7.$
\item [] $740=7\times (7\times (7+7)+7)+7-(7+7)/7.$
\item [] $741=7\times (7\times (7+7)+7)+7-7/7.$
\item [] $742=7\times (7\times (7+7)+7)+7.$
\item [] $743=7\times (7\times (7+7)+7)+7+7/7.$
\item [] $744=7\times 77+77+((7+7)/7)^7.$
\item [] $745=777-7-7-7-77/7.$
\item [] $746=7\times (7\times (7+7)+7)+77/7.$
\item [] $747=7\times (7\times (7+7)+7)+(77+7)/7.$
\item [] $748=777-77+7\times 7-7/7.$
\item [] $749=777+7\times 7-77.$
\item [] $750=(7+7+7/7)\times (7\times 7+7/7).$
\item [] $751=777-7-7-(77+7)/7.$
\item [] $752=777-7-7-77/7.$
\item [] $753=777-(7\times 7\times 7-7)/(7+7).$
\item [] $754=777-7-7-7-(7+7)/7.$
\item [] $755=777-7-7-7-7/7.$
\item [] $756=777-7-7-7.$
\item [] $757=777-7-7-7+7/7.$
\item [] $758=777-7-(77+7)/7.$
\item [] $759=777-7-77/7.$
\item [] $760=777-7-(77-7)/7.$
\item [] $761=777-7-7-(7+7)/7.$
\item [] $762=777-7-7-7/7.$
\item [] $763=777-7-7.$
\item [] $764=777-7-7+7/7.$
\item [] $765=777-(77+7)/7.$
\item [] $766=777-77/7.$
\item [] $767=777-(77-7)/7.$
\item [] $768=777-7-(7+7)/7.$
\item [] $769=777-7-7/7.$
\item [] $770=777-7.$
\item [] $771=777-7+7/7.$
\item [] $772=777-7+(7+7)/7.$
\item [] $773=777+7-77/7.$
\item [] $774=777-(7+7+7)/7.$
\item [] $775=777-(7+7)/7.$
\item [] $776=777-7/7.$
\item [] $777=777.$
\item [] $778=777+7/7.$
\item [] $779=777+(7+7)/7.$
\item [] $780=777+7+7-77/7.$
\item [] $781=777-7+77/7.$
\item [] $782=777+7-(7+7)/7.$
\item [] $783=777+7-7/7.$
\item [] $784=777+7.$
\item [] $785=777+7+7/7.$
\item [] $786=777+7+(7+7)/7.$
\item [] $787=777+(77-7)/7.$
\item [] $788=777+77/7.$
\item [] $789=777+(77+7)/7.$
\item [] $790=777+7+7-7/7.$
\item [] $791=777+7+7.$
\item [] $792=777+7+7+7/7.$
\item [] $793=777+7+7+(7+7)/7.$
\item [] $794=777+7+(77-7)/7.$
\item [] $795=777+7+77/7.$
\item [] $796=777+7+(77+7)/7.$
\item [] $797=777+7+7+7-7/7.$
\item [] $798=777+7+7+7.$
\item [] $799=777+7+7+7+7/7.$
\item [] $800=(7\times 7+7/7)\times (7+7+(7+7)/7).$
\item [] $801=777+7+7+(77-7)/7.$
\item [] $802=777+7+7+77/7.$
\item [] $803=777+7+7+(77+7)/7.$
\item [] $804=777+77-7\times 7-7/7.$
\item [] $805=777+77-7\times 7.$
\item [] $806=777+77-7\times 7+7/7.$
\item [] $807=7\times 7\times (7+7)+((7+7)/7)^7-7.$
\item [] $808=777+7\times 7-7-77/7.$
\item [] $809=777+7+7+7+77/7.$
\item [] $810=(7-7/7)\times (7+((7+7)/7)^7).$
\item [] $811=777-77+777/7.$
\item [] $812=777+7\times 7-7-7.$
\item [] $813=777+7\times 7-7-7+7/7.$
\item [] $814=7\times 7\times (7+7)+((7+7)/7)^7.$
\item [] $815=777+7\times 7-77/7.$
\item [] $816=777+7\times 7-(77-7)/7.$
\item [] $817=777+7\times 7-7-(7+7)/7.$
\item [] $818=777+7\times 7-7-7/7.$
\item [] $819=777+7\times 7-7.$
\item [] $820=777+7\times 7-7+7/7.$
\item [] $821=777+7\times 7-7+(7+7)/7.$
\item [] $822=777+7\times 7+7-77/7.$
\item [] $823=777+7\times 7-(7+7+7)/7.$
\item [] $824=777+7\times 7-(7+7)/7.$
\item [] $825=777+7\times 7-7/7.$
\item [] $826=777+7\times 7.$
\item [] $827=777+7\times 7+7/7.$
\item [] $828=777+7\times 7+(7+7)/7.$
\item [] $829=(77-7/7)\times 77/7-7.$
\item [] $830=777+7\times 7-7+77/7.$
\item [] $831=777+7\times 7+7-(7+7)/7.$
\item [] $832=777+7\times 7+7-7/7.$
\item [] $833=777+7\times 7+7.$
\item [] $834=777+7\times 7+7+7/7.$
\item [] $835=77\times 77/7-(77+7)/7.$
\item [] $836=77\times (77-7/7)/7.$
\item [] $837=777+7\times 7+77/7.$
\item [] $838=(77\times 77-7-7)/7-7.$
\item [] $839=(77\times 77-7)/7-7.$
\item [] $840=77\times 77/7-7.$
\item [] $841=(77\times 77+7)/7-7.$
\item [] $842=(77\times 77+7+7)/7-7.$
\item [] $843=777+77-77/7.$
\item [] $844=777+7\times 7+7+77/7.$
\item [] $845=77\times 77/7-(7+7)/7.$
\item [] $846=77\times 77/7-7/7.$
\item [] $847=77\times 77/7.$
\item [] $848=77\times 77/7+7/7.$
\item [] $849=77\times 77/7+(7+7)/7.$
\item [] $850=(77-7)/7\times (77+7+7/7).$
\item [] $851=77\times (77+7/7)/7-7.$
\item [] $852=777+77-(7+7)/7.$
\item [] $853=777+77-7/7.$
\item [] $854=777+77.$
\item [] $855=777+77+7/7.$
\item [] $856=777+77+(7+7)/7.$
\item [] $857=77\times (77+7/7)/7-7/7.$
\item [] $858=77\times (77+7/7)/7.$
\item [] $859=777+77+7-(7+7)/7.$
\item [] $860=777+77+7-7/7.$
\item [] $861=777+77+7.$
\item [] $862=777+77+7+7/7.$
\item [] $863=777+77+7+(7+7)/7.$
\item [] $864=(7+77/7)\times (7\times 7-7/7).$
\item [] $865=777+77+77/7.$
\item [] $866=777+77+(77+7)/7.$
\item [] $867=7\times (77+7\times 7)-7-7-7/7.$
\item [] $868=777+77+7+7.$
\item [] $869=77\times (77+(7+7)/7)/7.$
\item [] $870=7\times (7\times 7+77)-(77+7)/7.$
\item [] $871=7\times (7\times 7+77)-77/7.$
\item [] $872=777+77+7+77/7.$
\item [] $873=7\times (7\times 7+77)-7-(7+7)/7.$
\item [] $874=777+7\times (7+7)-7/7.$
\item [] $875=777+7\times (7+7).$
\item [] $876=777+7\times (7+7)+7/7.$
\item [] $877=7\times (7\times 7+77)-7+(7+7)/7.$
\item [] $878=777+7777/77.$
\item [] $879=7\times (7\times 7+77)-(7+7+7)/7.$
\item [] $880=(7+7/7)\times (777-7)/7.$
\item [] $881=7\times (7\times 7+77)-7/7.$
\item [] $882=7\times (7\times 7+77).$
\item [] $883=7\times (7\times 7+77)+7/7.$
\item [] $884=7\times (7\times 7+77)+(7+7)/7.$
\item [] $885=7\times ((7+7)/7)^7-77/7.$
\item [] $886=777+(777-7-7)/7.$
\item [] $887=777+(777-7)/7.$
\item [] $888=777+777/7.$
\item [] $889=7\times ((7+7)/7)^7-7.$
\item [] $890=7\times (7\times 7+77)+7+7/7.$
\item [] $891=777-7-7+((7+7)/7)^7.$
\item [] $892=7-77/7+((7+7)/7)^7\times 7.$
\item [] $893=7\times (7\times 7+77)+77/7.$
\item [] $894=7\times ((7+7)/7)^7-(7+7)/7.$
\item [] $895=7\times ((7+7)/7)^7-7/7.$
\item [] $896=7\times ((7+7)/7)^7.$
\item [] $897=7\times ((7+7)/7)^7+7/7.$
\item [] $898=777-7+((7+7)/7)^7.$
\item [] $899=777+(777+77)/7.$
\item [] $900=(7\times 7+7/7)\times (77/7+7).$
\item [] $901=7\times ((7+7)/7)^7+7-(7+7)/7.$
\item [] $902=7\times ((7+7)/7)^7+7-7/7.$
\item [] $903=7\times ((7+7)/7)^7+7.$
\item [] $904=7\times ((7+7)/7)^7+7+7/7.$
\item [] $905=777+((7+7)/7)^7.$
\item [] $906=77\times (77+7-7/7)/7-7.$
\item [] $907=7\times ((7+7)/7)^7+77/7.$
\item [] $908=7\times ((7+7)/7)^7+(77+7)/7.$
\item [] $909=7\times ((7+7)/7)^7+7+7-7/7.$
\item [] $910=(77-7)\times (7+7-7/7).$
\item [] $911=7\times ((7+7)/7)^7+7+7+7/7.$
\item [] $912=(77+7)\times (77-7/7)/7.$
\item [] $913=77\times (77+7-7/7)/7.$
\item [] $914=7\times ((7+7)/7)^7+7+77/7.$
\item [] $915=7777/7-(7+7)\times (7+7).$
\item [] $916=(77\times (77+7)-7)/7-7.$
\item [] $917=77\times (77+7)/7-7.$
\item [] $918=(77\times (77+7)+7)/7-7.$
\item [] $919=777+7+7+((7+7)/7)^7.$
\item [] $920=77\times (77+7-7/7)/7+7.$
\item [] $921=77\times (77+7+7/7)/7-7-7.$
\item [] $922=77\times (77+7)/7-(7+7)/7.$
\item [] $923=77\times (77+7)/7-7/7.$
\item [] $924=77\times (77+7)/7.$
\item [] $925=(77\times (77+7)+7)/7.$
\item [] $926=77\times (77+7)/7+(7+7)/7.$
\item [] $927=77\times (77+7-7/7)/7+7+7.$
\item [] $928=77\times (77+7+7/7)/7-7.$
\item [] $929=7\times (7\times 7+77+7)-(7+7)/7.$
\item [] $930=7\times (7\times 7+77+7)-7/7.$
\item [] $931=7\times (7\times 7+77+7).$
\item [] $932=7\times (7\times 7+77+7)+7/7.$
\item [] $933=7\times (7\times 7+77+7)+(7+7)/7.$
\item [] $934=7\times (7+((7+7)/7)^7)-77/7.$
\item [] $935=77\times (77+7+7/7)/7.$
\item [] $936=(77+7)\times (77+7/7)/7.$
\item [] $937=7\times (7\times 7+7+77)+7-7/7.$
\item [] $938=7\times (7\times 7+7+77)+7.$
\item [] $939=7\times (7\times 7+7+77)+7+7/7.$
\item [] $940=7+7\times (7\times 7+77+7)+(7+7)/7.$
\item [] $941=(7\times 7+77+7)\times 7+(77-7)/7.$
\item [] $942=7\times (7\times 7+77+7)+77/7.$
\item [] $943=7\times (((7+7)/7)^7+7)-(7+7)/7.$
\item [] $944=(7+777/7)\times (7+7/7).$
\item [] $945=7\times (((7+7)/7)^7+7).$
\item [] $946=7\times (((7+7)/7)^7+7)+7/7.$
\item [] $947=7\times (((7+7)/7)^7+7)+(7+7)/7.$
\item [] $948=(77+7)\times (77+(7+7)/7)/7.$
\item [] $949=7\times (7\times 7+77+7)+7+77/7.$
\item [] $950=(7\times 7+7/7)\times (7+(77+7)/7).$
\item [] $951=7\times 7\times (7+7+7)-77-7/7.$
\item [] $952=7\times (7+((7+7)/7)^7)+7.$
\item [] $953=7\times 7\times (7+7+7)-77+7/7.$
\item [] $954=777+7\times 7+((7+7)/7)^7.$
\item [] $955=(7+7)\times (77-7-7/7)-77/7.$
\item [] $956=7\times (7+((7+7)/7)^7)+77/7.$
\item [] $957=(777+77\times 77-7)/7.$
\item [] $958=(777+77\times 77)/7.$
\item [] $959=7\times (7\times 7+77)+77.$
\item [] $960=7\times (77+7\times 7)+77+7/7.$
\item [] $961=77\times (77+77/7)/7-7.$
\item [] $962=(7+7)\times (77-7)-7-77/7.$
\item [] $963=(77-7)\times (7+7)-7+(7-77)/7.$
\item [] $964=7777/7-7\times (7+7+7).$
\item [] $965=777+77+777/7.$
\item [] $966=(7+7)\times (77-7-7/7).$
\item [] $967=77\times (7+7)-777/7.$
\item [] $968=77\times (77+77/7)/7.$
\item [] $969=(7+7)\times (77-7)-77/7.$
\item [] $970=(7+7)\times (77-7)-(77-7)/7.$
\item [] $971=(7+7)\times (77-7)-7-(7+7)/7.$
\item [] $972=(7+7)\times (77-7)-7-7/7.$
\item [] $973=(7+7)\times (77-7)-7.$
\item [] $974=(7+7)\times (77-7)-7+7/7.$
\item [] $975=(7+7)\times (77-7)-7+(7+7)/7.$
\item [] $976=(7+7)\times (77-7)+7-77/7.$
\item [] $977=(7+7)\times (77-7)-(7+7+7)/7.$
\item [] $978=(7+7)\times (77-7)-(7+7)/7.$
\item [] $979=(7+7)\times (77-7)-7/7.$
\item [] $980=(7+7)\times (77-7).$
\item [] $981=(7+7)\times (77-7)+7/7.$
\item [] $982=(7+7)\times (77-7)+(7+7)/7.$
\item [] $983=(7+7)\times (77-7)+(7+7+7)/7.$
\item [] $984=(7+7)\times (77-7)-7+77/7.$
\item [] $985=(7+7)\times (77-7)+7-(7+7)/7.$
\item [] $986=(7+7)\times (77-7)+7-7/7.$
\item [] $987=(7+7)\times (77-7)+7.$
\item [] $988=(7+7)\times (77-7)+7+7/7.$
\item [] $989=(7+7)\times (77-7)+7+(7+7)/7.$
\item [] $990=(7+7+7/7)\times (77-77/7).$
\item [] $991=(7+7)\times (77-7)+77/7.$
\item [] $992=(7+7)\times (77-7)+(77+7)/7.$
\item [] $993=(7+7)\times (77-7)+7+7-7/7.$
\item [] $994=(7+7)\times (77-7)+7+7.$
\item [] $995=(7+7)\times (77-7)+7+7+7/7.$
\item [] $996=(7+77)\times (7+77-7/7)/7.$
\item [] $997=(7+7)\times (77-7)+7+(77-7)/7.$
\item [] $998=(7+7)\times (77-7)+7+77/7.$
\item [] $999=777\times (7+(7+7)/7)/7.$
\item [] $1000=(7+7+7-7/7)\times (7\times 7+7/7).$
\end{itemize}
\end{multicols}
}

\section{\textbf{Representations Using Number 8}}

{\footnotesize
\begin{multicols}{3}
\begin{itemize}
\item [] $101=8888/88.$
\item [] $102=(888-8)/8-8.$
\item [] $103=(888-8\times 8)/8.$
\item [] $104=88+8+8.$
\item [] $105=88+8+8+8/8.$
\item [] $106=88+8+(88-8)/8.$
\item [] $107=88+8+88/8.$
\item [] $108=8+(888-88)/8.$
\item [] $109=(888-8-8)/8.$
\item [] $110=(888-8)/8.$
\item [] $111=888/8.$
\item [] $112=(888+8)/8.$
\item [] $113=(888+8+8)/8.$
\item [] $114=(888+88)/8-8.$
\item [] $115=888/8+8\times 8/(8+8).$
\item [] $116=8+8+(888-88)/8.$
\item [] $117=8\times (8+8)-88/8.$
\item [] $118=8+(888-8)/8.$
\item [] $119=8+888/8.$
\item [] $120=8\times (8+8)-8.$
\item [] $121=88\times 88/(8\times 8).$
\item [] $122=(888+88)/8.$
\item [] $123=88+8+8+8+88/8.$
\item [] $124=8\times (8+8)-8\times 8/(8+8).$
\item [] $125=8\times (8+8)+8-88/8.$
\item [] $126=8\times (8+8)-(8+8)/8.$
\item [] $127=8\times (8+8)-8/8.$
\item [] $128=8\times (8+8).$
\item [] $129=8\times (8+8)+8/8.$
\item [] $130=8\times (8+8)+(8+8)/8.$
\item [] $131=8\times (8+8)-8+88/8.$
\item [] $132=8\times (8+8)+8\times 8/(8+8).$
\item [] $133=(8+88/8)\times (8-8/8).$
\item [] $134=+8\times (8+8)+8-(8+8)/8.$
\item [] $135=8\times (8+8)+8-8/8.$
\item [] $136=8\times (8+8)+8.$
\item [] $137=8\times (8+8)+8+8/8.$
\item [] $138=8\times (8+8)+(88-8)/8.$
\item [] $139=8\times (8+8)+88/8.$
\item [] $140=8\times (8+8)+(88+8)/8.$
\item [] $141=88+8\times 8-88/8.$
\item [] $142=88+8\times 8-8-(8+8)/8.$
\item [] $143=88+8\times 8-8-8/8.$
\item [] $144=88+8\times 8-8.$
\item [] $145=88+8\times 8-8+8/8.$
\item [] $146=88+8\times 8-8+(8+8)/8.$
\item [] $147=8\times (8+8)+8+88/8.$
\item [] $148=888/(8-(8+8)/8).$
\item [] $149=88+8\times 8+8-88/8.$
\item [] $150=88+8\times 8-(8+8)/8.$
\item [] $151=88+8\times 8-8/8.$
\item [] $152=88+8\times 8.$
\item [] $153=88+8\times 8+8/8.$
\item [] $154=88+8\times 8+(8+8)/8.$
\item [] $155=88+8\times 8-8+88/8.$
\item [] $156=88+8\times 8+8\times 8/(8+8).$
\item [] $157=88+8\times 8+(88-8)/(8+8).$
\item [] $158=88+8\times 8+8-(8+8)/8.$
\item [] $159=88+8\times 8+8-8/8.$
\item [] $160=88+8\times 8+8.$
\item [] $161=88+8\times 8+8+8/8.$
\item [] $162=88+8\times 8+8+(8+8)/8.$
\item [] $163=88+8\times 8+88/8.$
\item [] $164=88+88-(88+8)/8.$
\item [] $165=88+88-88/8.$
\item [] $166=88+88-8-(8+8)/8.$
\item [] $167=88+88-8-8/8.$
\item [] $168=88+88-8.$
\item [] $169=88+88-8+8/8.$
\item [] $170=(8+8+8/8)\times (88-8)/8.$
\item [] $171=(8+88/8)\times (8+8/8).$
\item [] $172=88+88-8\times 8/(8+8).$
\item [] $173=88+88+8-88/8.$
\item [] $174=88+88-(8+8)/8.$
\item [] $175=88+88-8/8.$
\item [] $176=88+88.$
\item [] $177=88+88+8/8.$
\item [] $178=88+88+(8+8)/8.$
\item [] $179=88+88-8+88/8.$
\item [] $180=(88+(8+8)/8)\times (8+8)/8.$
\item [] $181=8\times (8+8+8)-88/8.$
\item [] $182=8\times (8+8+8)-8-(8+8)/8.$
\item [] $183=88+88+8-8/8.$
\item [] $184=88+88+8.$
\item [] $185=88+88+8+8/8.$
\item [] $186=88+88+(88-8)/8.$
\item [] $187=88+88+88/8.$
\item [] $188=88+(888-88)/8.$
\item [] $189=88+8888/88.$
\item [] $190=8\times (8+8+8)-(8+8)/8.$
\item [] $191=8\times (8+8+8)-8/8.$
\item [] $192=8\times (8+8+8).$
\item [] $193=8\times (8+8+8)+8/8.$
\item [] $194=8\times (8+8+8)+(8+8)/8.$
\item [] $195=88+88+8+88/8.$
\item [] $196=8\times (8+8+8)+8\times 8/(8+8).$
\item [] $197=88+(888-8-8)/8.$
\item [] $198=88+(888-8)/8.$
\item [] $199=88+888/8.$
\item [] $200=8\times (8+8+8)+8.$
\item [] $201=8\times (8+8+8)+8+8/8.$
\item [] $202=8\times (8+8+8)+(88-8)/8.$
\item [] $203=8\times (8+8+8)+88/8.$
\item [] $204=8\times (8+8+8)+(88+8)/8.$
\item [] $205=8\times (8+8)+88-88/8.$
\item [] $206=88+8+(888-8)/8.$
\item [] $207=88+8+888/8.$
\item [] $208=88+8\times (8+8)-8.$
\item [] $209=88+8\times (8+8)-8+8/8.$
\item [] $210=88+(888+88)/8.$
\item [] $211=8\times (8+8+8)+8+88/8.$
\item [] $212=8\times (8+8+8)+8+(88+8)/8.$
\item [] $213=88+8\times (8+8)+8-88/8.$
\item [] $214=88+8\times (8+8)-(8+8)/8.$
\item [] $215=88+8\times (8+8)-8/8.$
\item [] $216=88+8\times (8+8).$
\item [] $217=88+8\times (8+8)+8/8.$
\item [] $218=88+8\times (8+8)+(8+8)/8.$
\item [] $219=88+8\times (8+8)-8+88/8.$
\item [] $220=(888-8)\times (8+8)/(8\times 8).$
\item [] $221=888\times (8+8)/(8\times 8)-8/8.$
\item [] $222=888\times (8+8)/(8\times 8).$
\item [] $223=88+8\times (8+8)+8-8/8.$
\item [] $224=88+8\times (8+8)+8.$
\item [] $225=(8+8-8/8)^{((8+8)/8)}.$
\item [] $226=(888+8+8)\times (8+8)/(8\times 8).$
\item [] $227=88+8\times (8+8)+88/8.$
\item [] $228=(8+88/8)\times (88+8)/8.$
\item [] $229=8\times (8+8)+8888/88.$
\item [] $230=8+888\times (8+8)/(8\times 8).$
\item [] $231=8\times (8+8)-8+888/8.$
\item [] $232=88+88+8\times 8-8.$
\item [] $233=8+(8+8-8/8)^{((8+8)/8)}.$
\item [] $234=(8+8+8+(8+8)/8)\times (8+8/8).$
\item [] $235=8\times (8+8)+88+8+88/8.$
\item [] $236=((8+8)/8)^8-8-(88+8)/8.$
\item [] $237=((8+8)/8)^8-8-88/8.$
\item [] $238=8\times (8+8)+(888-8)/8.$
\item [] $239=8\times (8+8)+888/8.$
\item [] $240=88+88+8\times 8.$
\item [] $241=(8+8)\times (8+8-8/8)+8/8.$
\item [] $242=(88+88)\times 88/(8\times 8).$
\item [] $243=((8+8+8)/8)^{(8-(8+8+8)/8)}.$
\item [] $244=(8+8)\times (8+8)-(88+8)/8.$
\item [] $245=(8+8)\times (8+8)-88/8.$
\item [] $246=(8+8)\times (8+8)-(88-8)/8.$
\item [] $247=(8+8)\times (8+8)-8-8/8.$
\item [] $248=(8+8)\times (8+8)-8.$
\item [] $249=(8+8)\times (8+8)-8+8/8.$
\item [] $250=(8+8)\times (8+8)-8+(8+8)/8.$
\item [] $251=88+88+8\times 8+88/8.$
\item [] $252=8\times (8\times 8\times 8-8)/(8+8).$
\item [] $253=(8+8)\times (8+8)+8-88/8.$
\item [] $254=(8+8)\times (8+8)-(8+8)/8.$
\item [] $255=(8+8)\times (8+8)-8/8.$
\item [] $256=(8+8)\times (8+8).$
\item [] $257=(8+8)\times (8+8)+8/8.$
\item [] $258=(8+8)\times (8+8)+(8+8)/8.$
\item [] $259=(8+8)\times (8+8)+(8+8+8)/8.$
\item [] $260=8\times (8\times 8\times 8+8)/(8+8).$
\item [] $261=(8+8+8)\times (88-8/8)/8.$
\item [] $262=(8+8)\times (8+8)+8-(8+8)/8.$
\item [] $263=(8+8)\times (8+8)+8-8/8.$
\item [] $264=(8+8)\times (8+8)+8.$
\item [] $265=(8+8)\times (8+8)+8+8/8.$
\item [] $266=(8+8)\times (8+8)+8+(8+8)/8.$
\item [] $267=(8+8)\times (8+8)+88/8.$
\item [] $268=(8+8)\times (8+8)+(88+8)/8.$
\item [] $269=8\times (8+8+8)+88-88/8.$
\item [] $270=(8+8)\times (8+8)+8+8-(8+8)/8.$
\item [] $271=(8+8)\times (8+8+8/8)-8/8.$
\item [] $272=(8+8)\times (8+8+8/8).$
\item [] $273=(8+8)\times (8+8+8/8)+8/8.$
\item [] $274=(8+8)\times (8+8)+8+8+(8+8)/8.$
\item [] $275=(8+8)\times (8+8)+8+88/8.$
\item [] $276=(8+8)\times (8+8)+8+(88+8)/8.$
\item [] $277=88+8\times (8+8+8)-(8+8+8)/8.$
\item [] $278=(8\times (8+8)+88/8)\times (8+8)/8.$
\item [] $279=8\times (8+8+8)+88-8/8.$
\item [] $280=8\times (8+8+8)+88.$
\item [] $281=8\times (8+8+8)+88+8/8.$
\item [] $282=8\times (8+8+8)+88+(8+8)/8.$
\item [] $283=8\times (8+8+8)+88-8+88/8.$
\item [] $284=8\times 8\times (8\times 8+8-8/8)/(8+8).$
\item [] $285=((8+8+8)\times (88+8)-8-8-8)/8.$
\item [] $286=((8+8+8)\times (88+8)-8-8)/8.$
\item [] $287=((8+8+8)\times (88+8)-8)/8.$
\item [] $288=(8+8+8)\times (88+8)/8.$
\item [] $289=(8+8+8/8)^{((8+8)/8)}.$
\item [] $290=((8+8+8)\times (88+8)+8+8)/8.$
\item [] $291=8\times (8+8+8)+88+88/8.$
\item [] $292=8\times (8+8+8)+88+(88+8)/8.$
\item [] $293=((8+8)/8)^8+888/(8+8+8).$
\item [] $294=8+8\times 8+888\times (8+8)/(8\times 8).$
\item [] $295=8\times 8\times (8\times 8+8)/(8+8)+8-8/8.$
\item [] $296=8\times 888/(8+8+8).$
\item [] $297=(8+8+8/8)^{((8+8)/8)}+8.$
\item [] $298=8\times 888/(8+8+8)+(8+8)/8.$
\item [] $299=8\times 8\times (8\times 8+8)/(8+8)+88/8.$
\item [] $300=8\times 8\times (8\times 8+88/8)/(8+8).$
\item [] $301=(8+8)\times (8+8)+8\times 8-8-88/8.$
\item [] $302=(8+8)\times (8\times 8+88-8/8)/8.$
\item [] $303=(8+8)\times (8+88/8)-8/8.$
\item [] $304=(8+8)\times (8+88/8).$
\item [] $305=(8+8)\times (8+88/8)+8/8.$
\item [] $306=(8+8+(8+8)/8)\times (8+8+8/8).$
\item [] $307=8\times (88\times 8-88)/(8+8)-8/8.$
\item [] $308=8\times (88\times 8-88)/(8+8).$
\item [] $309=(8+8+8)\times (888/8-8)/8.$
\item [] $310=88+888\times (8+8)/(8\times 8).$
\item [] $311=8\times (8+8+8)+8+888/8.$
\item [] $312=(8+8)\times (8+8)+8\times 8-8.$
\item [] $313=(8+8)\times (8+8)+8\times 8-8+8/8.$
\item [] $314=(8+8)\times (8+8)+8\times 8-8+(8+8)/8.$
\item [] $315=(8\times 8-8/8)\times (8-(8+8+8)/8).$
\item [] $316=(8+8)\times (8+8)+8\times (8-8/(8+8)).$
\item [] $317=(8+8)\times (8+8)+8\times 8+8-88/8.$
\item [] $318=(8+8)\times (8+8)+8\times 8-(8+8)/8.$
\item [] $319=(8+8)\times (8+8)+8\times 8-8/8.$
\item [] $320=(8+8)\times (8+8)+8\times 8.$
\item [] $321=(8+8)\times (8+8)+8\times 8+8/8.$
\item [] $322=(8+8)\times (8+8)+8\times 8+(8+8)/8.$
\item [] $323=(8+8)\times (8+8)+8\times 8-8+88/8.$
\item [] $324=(8+8+(8+8)/8)^{((8+8)/8)}.$
\item [] $325=888\times (8+8+8)/(8\times 8)-8.$
\item [] $326=(8+8)\times (8+8)+8\times 8+8-(8+8)/8.$
\item [] $327=8\times (8+8)+88+888/8.$
\item [] $328=(8+8)\times (8+8)+8\times 8+8.$
\item [] $329=(8-8/8)\times (8\times 8-8-8-8/8).$
\item [] $330=(8-(8+8)/8)\times (8\times 8-8-8/8).$
\item [] $331=(8+8)\times (8+8)+88/8+8\times 8.$
\item [] $332=(888\times (88/8-8)-8)/8.$
\item [] $333=888\times (88/8-8)/8.$
\item [] $334=8\times 8\times 8-88-88-(8+8)/8.$
\item [] $335=8\times 8\times 8-88-88-8/8.$
\item [] $336=8\times 8\times 8-88-88.$
\item [] $337=8\times (8\times 8-8)-888/8.$
\item [] $338=(8+8)\times (8+8)+88-8+(8+8)/8.$
\item [] $339=(8+8)\times (8+8)+8\times 8+8+88/8.$
\item [] $340=8\times 8\times 88/(8+8)-(88+8)/8.$
\item [] $341=888\times (88/8-8)/8+8.$
\item [] $342=(8-(8+8)/8)\times (8\times 8-8+8/8).$
\item [] $343=(8-8/8)^{((8+8+8)/8)}.$
\item [] $344=(8+8)\times (8+8)+88.$
\item [] $345=(8+8)\times (8+8)+88+8/8.$
\item [] $346=(8+8)\times (8+8)+88+(8+8)/8.$
\item [] $347=8\times 8\times 88/(8+8)-8-8+88/8.$
\item [] $348=8\times (8\times 88-8)/(8+8).$
\item [] $349=8\times 8\times 88/(8+8)+8-88/8.$
\item [] $350=8\times 8\times 88/(8+8)-(8+8)/8.$
\item [] $351=8\times 8\times 88/(8+8)-8/8.$
\item [] $352=8\times 8\times 88/(8+8).$
\item [] $353=8\times 8\times 88/(8+8)+8/8.$
\item [] $354=8\times 8\times 88/(8+8)+(8+8)/8.$
\item [] $355=8\times 8\times 88/(8+8)+(8+8+8)/8.$
\item [] $356=8\times 8\times (88+8/8)/(8+8).$
\item [] $357=(8+8+8)\times (8+888/8)/8.$
\item [] $358=8\times (8\times 8-8)-88-(8+8)/8.$
\item [] $359=8\times (8\times 8-8)-88-8/8.$
\item [] $360=8\times (8\times 8-8)-88.$
\item [] $361=(8+88/8)^{((8+8)/8)}.$
\item [] $362=8\times (8\times 8-8)-88+(8+8)/8.$
\item [] $363=8\times 8\times 88/(8+8)+88/8.$
\item [] $364=8+8\times (8\times 88+8)/(8+8).$
\item [] $365=(8+8)\times (8+8+8)-8-88/8.$
\item [] $366=888-8\times 8\times 8-8-(8+8)/8.$
\item [] $367=(8+8)\times (8+8)+888/8.$
\item [] $368=888-8\times 8\times 8-8.$
\item [] $369=(8+88/8)^{((8+8)/8)}+8.$
\item [] $370=888-8\times 8\times 8-8+(8+8)/8.$
\item [] $371=8\times (8\times 8-8)-88+88/8.$
\item [] $372=(8+8)\times (8+8+8)-(88+8)/8.$
\item [] $373=(8+8)\times (8+8+8)-88/8.$
\item [] $374=888-8\times 8\times 8-(8+8)/8.$
\item [] $375=888-8\times 8\times 8-8/8.$
\item [] $376=888-8\times 8\times 8.$
\item [] $377=888-8\times 8\times 8+8/8.$
\item [] $378=888-8\times 8\times 8+(8+8)/8.$
\item [] $379=888-8\times 8\times 8+(8+8+8)/8.$
\item [] $380=(8+8)\times (8+8+8)-8\times 8/(8+8).$
\item [] $381=8\times (8\times 8-8-8)+8-88/8.$
\item [] $382=8\times (8\times 8-8-8)-(8+8)/8.$
\item [] $383=(8+8)\times (8+8+8)-8/8.$
\item [] $384=(8+8)\times (8+8+8).$
\item [] $385=(8+8)\times (8+8+8)+8/8.$
\item [] $386=(8+8)\times (8+8+8)+(8+8)/8.$
\item [] $387=(8+8)\times (8+8+8)+(8+8+8)/8.$
\item [] $388=(8+8)\times (8+8+8)+8\times 8/(8+8).$
\item [] $389=(8+8)\times (8+8+8)+8-(8+8+8)/8.$
\item [] $390=(8+8)\times (8+8+8)+8-(8+8)/8.$
\item [] $391=(8+8)\times (8+8+8)+8-8/8.$
\item [] $392=(8+8)\times (8+8+8)+8.$
\item [] $393=8\times 8\times 8-8-888/8.$
\item [] $394=(8+8)\times (8+8+8)+(88-8)/8.$
\item [] $395=(8+8)\times (8+8+8)+88/8.$
\item [] $396=(8+8)\times (8+8+8)+(88+8)/8.$
\item [] $397=8\times (8\times 88+88)/(8+8)+8/8.$
\item [] $398=(8+8)\times (88+888/8)/8.$
\item [] $399=(8+8)\times (8+8+8)+8+8-8/8.$
\item [] $400=(8+8)\times (8+8+8)+8+8.$
\item [] $401=8\times 8\times 8-888/8.$
\item [] $402=8\times 8\times 8-(888-8)/8.$
\item [] $403=8\times (8\times 8-8-8)+8+88/8.$
\item [] $404=8+88\times (8\times 8+8)/(8+8).$
\item [] $405=(8+8/8)\times (8\times 8-8-88/8).$
\item [] $406=(8+8+8)\times (8+8+8/8)-(8+8)/8.$
\item [] $407=(8+8+8)\times (8+8+8/8)-8/8.$
\item [] $408=(8+8+8)\times (8+8+8/8).$
\item [] $409=8\times 8\times 8+8-888/8.$
\item [] $410=8\times 8\times 8+8-(888-8)/8.$
\item [] $411=(8+8)\times (8+8+8)+8+8+88/8.$
\item [] $412=8\times 8\times (888/8-8)/(8+8).$
\item [] $413=8\times 8\times 8-88-88/8.$
\item [] $414=8\times 8\times 8-88-(88-8)/8.$
\item [] $415=8\times 8\times 8-88-8-8/8.$
\item [] $416=8\times 8\times 8-88-8.$
\item [] $417=8\times 8\times 8-88-8+8/8.$
\item [] $418=8\times 8\times 8-88-8+(8+8)/8.$
\item [] $419=8\times 8\times 8-88-8+(8+8+8)/8.$
\item [] $420=8\times 8\times 8-88-8\times 8/(8+8).$
\item [] $421=8\times 8\times 8-88-(8+8+8)/8.$
\item [] $422=8\times 8\times 8-88-(8+8)/8.$
\item [] $423=8\times 8\times 8-88-8/8.$
\item [] $424=8\times 8\times 8-88.$
\item [] $425=8\times 8\times 8-88+8/8.$
\item [] $426=8\times 8\times 8-88+(8+8)/8.$
\item [] $427=8\times 8\times 8-88+(8+8+8)/8.$
\item [] $428=8\times 8\times 8-88+8\times 8/(8+8).$
\item [] $429=8\times (8\times 8-8)-8-88/8.$
\item [] $430=8\times 8\times 8+8-88-(8+8)/8.$
\item [] $431=8\times 8\times 8+8-88-8/8.$
\item [] $432=8\times 8\times 8-88+8.$
\item [] $433=8\times 8\times 8-88+8+8/8.$
\item [] $434=(8-8/8)\times (8\times 8-(8+8)/8).$
\item [] $435=8\times 8\times 8-88+88/8.$
\item [] $436=8\times 888/(8+8)-8.$
\item [] $437=8\times (8\times 8-8)-88/8.$
\item [] $438=8\times (8\times 8-8)-(88-8)/8.$
\item [] $439=8\times (8\times 8-8)-8-8/8.$
\item [] $440=8\times (8\times 8-8)-8.$
\item [] $441=8\times (8\times 8-8)-8+8/8.$
\item [] $442=8\times (8\times 8-8)-8+(8+8)/8.$
\item [] $443=8\times 888/(8+8)-8/8.$
\item [] $444=8\times 888/(8+8).$
\item [] $445=8\times (8\times 8-8)+8-88/8.$
\item [] $446=8\times (8\times 8-8)-(8+8)/8.$
\item [] $447=8\times (8\times 8-8)-8/8.$
\item [] $448=8\times (8\times 8-8).$
\item [] $449=8\times (8\times 8-8)+8/8.$
\item [] $450=8\times (8\times 8-8)+(8+8)/8.$
\item [] $451=8\times (8\times 8-8)-8+88/8.$
\item [] $452=8\times 888/(8+8)+8.$
\item [] $453=8\times (8\times 8-8)+8-(8+8+8)/8.$
\item [] $454=8\times (8\times 8-8)+8-(8+8)/8.$
\item [] $455=8\times (8\times 8-8)+8-8/8.$
\item [] $456=8\times (8\times 8-8)+8.$
\item [] $457=8\times (8\times 8-8)+8+8/8.$
\item [] $458=8\times (8\times 8-8)+8+(8+8)/8.$
\item [] $459=8\times (8\times 8-8)+88/8.$
\item [] $460=8\times (8\times 8-8)+(88+8)/8.$
\item [] $461=8\times (8\times 8-8)+8+8-(8+8+8)/8.$
\item [] $462=(8\times 8+(8+8)/8)\times (8-8/8).$
\item [] $463=8\times (8\times 8-8)+8+8-8/8.$
\item [] $464=8\times (8\times 8-8)+8+8.$
\item [] $465=8\times (8\times 8-8)+8+8+8/8.$
\item [] $466=8\times (8\times 8-8)+8+8+(8+8)/8.$
\item [] $467=8\times (8\times 8-8)+8+88/8.$
\item [] $468=8\times 8\times 8-8\times 88/(8+8).$
\item [] $469=(8-8/8)\times (8\times 8+(8+8+8)/8).$
\item [] $470=8\times (8\times 8-8)+(88+88)/8.$
\item [] $471=(8+8)\times (8+8+8)+88-8/8.$
\item [] $472=(8+8)\times (8+8+8)+88.$
\item [] $473=(8+8)\times (8+8+8)+88+8/8.$
\item [] $474=(88-8-8/8)\times (8-(8+8)/8).$
\item [] $475=8\times 8\times 8-888/(8+8+8).$
\item [] $476=8\times (8+8/(8+8))\times (8-8/8).$
\item [] $477=8\times 8\times 8-8-8-8-88/8.$
\item [] $478=8\times 8\times (8-8/(8+8))-(8+8)/8.$
\item [] $479=8\times 8\times (8-8/(8+8))-8/8.$
\item [] $480=8\times 8\times (8-8/(8+8)).$
\item [] $481=8\times 8\times (8-8/(8+8))+8/8.$
\item [] $482=8\times 8\times (8-8/(8+8))+(8+8)/8.$
\item [] $483=8\times 8\times (8-8/(8+8))+(8+8+8)/8.$
\item [] $484=((88+88)/8)^{((8+8)/8)}.$
\item [] $485=8\times 8\times 8-8-8-88/8.$
\item [] $486=(88-8+8/8)\times (8-(8+8)/8).$
\item [] $487=8\times 8\times 8-8-8-8-8/8.$
\item [] $488=8\times 8\times 8-8-8-8.$
\item [] $489=8\times 8\times 8-8-8-8+8/8.$
\item [] $490=8\times 8\times 8-8-8-8+(8+8)/8.$
\item [] $491=8\times 8\times 8-8-8-8+(8+8+8)/8.$
\item [] $492=8\times (8\times 8-8/(8+8))-8-8.$
\item [] $493=8\times 8\times 8-8-88/8.$
\item [] $494=8\times 8\times 8-8-8-(8+8)/8.$
\item [] $495=8\times 8\times 8-8-8-8/8.$
\item [] $496=8\times 8\times 8-8-8.$
\item [] $497=8\times 8\times 8-8-8+8/8.$
\item [] $498=8\times 8\times 8-8-8+(8+8)/8.$
\item [] $499=8\times 8\times 8-(88+8+8)/8.$
\item [] $500=8\times 8\times 8-(88+8)/8.$
\item [] $501=8\times 8\times 8-88/8.$
\item [] $502=8\times 8\times 8-(88-8)/8.$
\item [] $503=8\times 8\times 8-8-8/8.$
\item [] $504=8\times 8\times 8-8.$
\item [] $505=8\times 8\times 8-8+8/8.$
\item [] $506=8\times 8\times 8-8+(8+8)/8.$
\item [] $507=8\times 8\times 8-8-8+88/8.$
\item [] $508=8\times (8\times 8-8/(8+8)).$
\item [] $509=8\times 8\times 8-(8+8+8)/8.$
\item [] $510=8\times 8\times 8-(8+8)/8.$
\item [] $511=8\times 8\times 8-8/8.$
\item [] $512=8\times 8\times 8.$
\item [] $513=8\times 8\times 8+8/8.$
\item [] $514=8\times 8\times 8+(8+8)/8.$
\item [] $515=8\times 8\times 8-8+88/8.$
\item [] $516=8\times (8\times 8+8/(8+8)).$
\item [] $517=8\times 8\times 8+8-(8+8+8)/8.$
\item [] $518=8\times 8\times 8+8-(8+8)/8.$
\item [] $519=8\times 8\times 8+8-8/8.$
\item [] $520=8\times 8\times 8+8.$
\item [] $521=8\times 8\times 8+8+8/8.$
\item [] $522=8\times 8\times 8+(88-8)/8.$
\item [] $523=8\times 8\times 8+88/8.$
\item [] $524=8\times 8\times 8+(88+8)/8.$
\item [] $525=8\times 8\times 8+8+8+8-88/8.$
\item [] $526=8\times 8\times 8+8+8-(8+8)/8.$
\item [] $527=8\times 8\times 8+8+8-8/8.$
\item [] $528=8\times 8\times 8+8+8.$
\item [] $529=8\times 8\times 8+8+8+8/8.$
\item [] $530=8\times 8\times 8+8+8+(8+8)/8.$
\item [] $531=8\times 8\times 8+8+88/8.$
\item [] $532=8\times 8\times 8+8+(88+8)/8.$
\item [] $533=8\times 8\times 8+8+(88+8+8)/8.$
\item [] $534=8\times 8\times 8+8+8+8-(8+8)/8.$
\item [] $535=8\times 8\times 8+8+8+8-8/8.$
\item [] $536=8\times (8\times 8-8)+88.$
\item [] $537=8\times (8\times 8-8)+88+8/8.$
\item [] $538=8\times (8\times 8-8)+88+(8+8)/8.$
\item [] $539=8\times 8\times 8+8+8+88/8.$
\item [] $540=8\times 8\times 8+8+8+(88+8)/8.$
\item [] $541=(88-8)\times 8-88-88/8.$
\item [] $542=8\times 8\times (8+8/(8+8))-(8+8)/8.$
\item [] $543=8\times 8\times (8+8/(8+8))-8/8.$
\item [] $544=8\times 8\times (8+8/(8+8)).$
\item [] $545=8\times 8\times (8+8/(8+8))+8/8.$
\item [] $546=8\times 8\times (8+8/(8+8))+(8+8)/8.$
\item [] $547=8\times 8\times 8+8+8+8+88/8.$
\item [] $548=8\times (88-8-8/(8+8))-88.$
\item [] $549=8\times (8\times 8+8)-8-8-88/8.$
\item [] $550=8\times (88-8)-88-(8+8)/8.$
\item [] $551=8\times (88-8)-88-8/8.$
\item [] $552=8\times (88-8)-88.$
\item [] $553=8\times (88-8)-88+8/8.$
\item [] $554=8\times (88-8)-88+(8+8)/8.$
\item [] $555=888\times (8-(8+8+8)/8)/8.$
\item [] $556=8\times 8\times 8+8\times 88/(8+8).$
\item [] $557=8\times (8\times 8+8)-8-88/8.$
\item [] $558=(8+8/8)\times (8\times 8-(8+8)/8).$
\item [] $559=(8\times 8-8)\times 8+888/8.$
\item [] $560=8\times 8\times 8+8\times 8-8-8.$
\item [] $561=8\times 8\times 8+8\times 8-8-8+8/8.$
\item [] $562=8\times (8+8\times 8)-8-8+(8+8)/8.$
\item [] $563=8\times (8\times 8+8)-(88+8+8)/8.$
\item [] $564=8\times (8\times 8+8-8/(8+8))-8.$
\item [] $565=8\times (8\times 8+8)-88/8.$
\item [] $566=8\times (8\times 8+8)-8-(8+8)/8.$
\item [] $567=(8+8/8)\times (8\times 8-8/8).$
\item [] $568=8\times (8\times 8+8)-8.$
\item [] $569=8\times (8\times 8+8)-8+8/8.$
\item [] $570=8\times (8\times 8+8)-8+(8+8)/8.$
\item [] $571=8\times (8+8\times 8)-8-8+88/8.$
\item [] $572=8\times (8\times 8+8-8/(8+8)).$
\item [] $573=8\times (8\times 8+8)-(8+8+8)/8.$
\item [] $574=8\times (8\times 8+8)-(8+8)/8.$
\item [] $575=8\times (8\times 8+8)-8/8.$
\item [] $576=8\times (8\times 8+8).$
\item [] $577=8\times (8\times 8+8)+8/8.$
\item [] $578=8\times (8\times 8+8)+(8+8)/8.$
\item [] $579=8\times (8\times 8+8)+(8+8+8)/8.$
\item [] $580=8\times (8\times 8+8+8/(8+8)).$
\item [] $581=8\times (8\times 8+8)+8+8-88/8.$
\item [] $582=8\times (8\times 8+8)+8-(8+8)/8.$
\item [] $583=8\times (8\times 8+8)+8-8/8.$
\item [] $584=8\times (8\times 8+8)+8.$
\item [] $585=8\times (8\times 8+8)+8+8/8.$
\item [] $586=8\times (8\times 8+8)+8+(8+8)/8.$
\item [] $587=8\times (8\times 8+8)+88/8.$
\item [] $588=8\times (8\times 8+8)+(88+8)/8.$
\item [] $589=8\times 8\times 8-88/8+88.$
\item [] $590=8\times 8\times 8+88-8-(8+8)/8.$
\item [] $591=8\times 8\times 8+88-8-8/8.$
\item [] $592=8\times 8\times 8+88-8.$
\item [] $593=8\times 88-888/8.$
\item [] $594=8\times 88-(888-8)/8.$
\item [] $595=(8+88-88/8)\times (8-8/8).$
\item [] $596=8\times (8\times 8-8/(8+8))+88.$
\item [] $597=88\times (8-8/8)-8-88/8.$
\item [] $598=8\times 8\times 8+88-(8+8)/8.$
\item [] $599=8\times 8\times 8+88-8/8.$
\item [] $600=8\times 8\times 8+88.$
\item [] $601=8\times 8\times 8+88+8/8.$
\item [] $602=8\times 8\times 8+88+(8+8)/8.$
\item [] $603=88\times 8-8888/88.$
\item [] $604=8\times 8\times 8+88+8\times 8/(8+8).$
\item [] $605=88\times 8-88-88/8.$
\item [] $606=8\times 88-88-8-(8+8)/8.$
\item [] $607=8\times 88-88-8-8/8.$
\item [] $608=8\times 88-88-8.$
\item [] $609=(8-8/8)\times (88-8/8).$
\item [] $610=88\times 8-88-8+(8+8)/8.$
\item [] $611=8\times 8\times 8+88+88/8.$
\item [] $612=8\times (88-8/(8+8))-88.$
\item [] $613=8\times (88-8)-8-8-88/8.$
\item [] $614=8\times 88-88-(8+8)/8.$
\item [] $615=8\times 88-88-8/8.$
\item [] $616=8\times 88-88.$
\item [] $617=8\times 88-88+8/8.$
\item [] $618=8\times 88-88+(8+8)/8.$
\item [] $619=8\times 88-88+(8+8+8)/8.$
\item [] $620=88\times (8-8/8)+8\times 8/(8+8).$
\item [] $621=8\times (88-8)-8-88/8.$
\item [] $622=88\times 8-88+8-(8+8)/8.$
\item [] $623=8\times 8\times 8+888/8.$
\item [] $624=8\times 88-88+8.$
\item [] $625=8\times 88-88+8+8/8.$
\item [] $626=8\times 88-88+8+(8+8)/8.$
\item [] $627=8\times 88-88+88/8.$
\item [] $628=8\times (88-8)-(88+8)/8.$
\item [] $629=8\times (88-8)-88/8.$
\item [] $630=8\times (88-8)-(88-8)/8.$
\item [] $631=8\times (88-8)-8-8/8.$
\item [] $632=8\times (88-8)-8.$
\item [] $633=8\times (88-8)-8+8/8.$
\item [] $634=8\times (88-8)-8+(8+8)/8.$
\item [] $635=8\times (88-8)-8-8+88/8.$
\item [] $636=8\times (88-8-8/(8+8)).$
\item [] $637=8\times (88-8)+8-88/8.$
\item [] $638=8\times (88-8)-(8+8)/8.$
\item [] $639=8\times (88-8)-8/8.$
\item [] $640=8\times (88-8).$
\item [] $641=8\times (88-8)+8/8.$
\item [] $642=8\times (88-8)+(8+8)/8.$
\item [] $643=8\times (88-8)-8+88/8.$
\item [] $644=8\times (88-8+8/(8+8)).$
\item [] $645=8\times (88-8)+8-(8+8+8)/8.$
\item [] $646=8\times (88-8)+8-(8+8)/8.$
\item [] $647=8\times (88-8)+8-8/8.$
\item [] $648=8\times (88-8)+8.$
\item [] $649=8\times (88-8)+8+8/8.$
\item [] $650=8\times (88-8)+(88-8)/8.$
\item [] $651=8\times (88-8)+88/8.$
\item [] $652=8\times (88-8)+(88+8)/8.$
\item [] $653=8\times (88-8)+(88+8+8)/8.$
\item [] $654=8\times (88-8)+8+8-(8+8)/8.$
\item [] $655=8\times (88-8)+8+8-8/8.$
\item [] $656=8\times (88-8)+8+8.$
\item [] $657=8\times (88-8)+8+8+8/8.$
\item [] $658=8\times (88-8)+8+8+(8+8)/8.$
\item [] $659=8\times (88-8)+8+88/8.$
\item [] $660=8\times 88\times (8-8/(8+8))/8.$
\item [] $661=8\times (8\times 8+8)+88-(8+8+8)/8.$
\item [] $662=8\times (8\times 8+8)+88-(8+8)/8.$
\item [] $663=8\times (8\times 8+8)+88-8/8.$
\item [] $664=8\times (8\times 8+8)+88.$
\item [] $665=8\times (8\times 8+8)+88+8/8.$
\item [] $666=888\times (8-(8+8)/8)/8.$
\item [] $667=8\times (88-8)+8+8+88/8.$
\item [] $668=8\times (88+8)-88-(88+8)/8.$
\item [] $669=8\times (88+8)-88-88/8.$
\item [] $670=(88+8)\times (8-8/8)-(8+8)/8.$
\item [] $671=(88+8)\times (8-8/8)-8/8.$
\item [] $672=(88+8)\times (8-8/8).$
\item [] $673=(8-8/8)\times (88+8)+8/8.$
\item [] $674=(8-8/8)\times (88+8)+(8+8)/8.$
\item [] $675=(8\times 8+88/8)\times (8+8/8).$
\item [] $676=8\times 88-8-8-(88+8)/8.$
\item [] $677=8\times 88-8-8-88/8.$
\item [] $678=8\times (88+8)-88-(8+8)/8.$
\item [] $679=8\times (88+8)-88-8/8.$
\item [] $680=8\times (88+8)-88.$
\item [] $681=8\times (88+8)-88+8/8.$
\item [] $682=8\times 88-(88+88)/8.$
\item [] $683=8\times 88-(88+88-8)/8.$
\item [] $684=8\times 88-8-(88+8)/8.$
\item [] $685=8\times 88-8-88/8.$
\item [] $686=8\times 88-8-8-(8+8)/8.$
\item [] $687=8\times 88-8-8-8/8.$
\item [] $688=8\times 88-8-8.$
\item [] $689=8\times 88-8-8+8/8.$
\item [] $690=8\times 88-8-8+(8+8)/8.$
\item [] $691=8\times 88-(88+8+8)/8.$
\item [] $692=8\times 88-(88+8)/8.$
\item [] $693=8\times 88-88/8.$
\item [] $694=8\times 88-(88-8)/8.$
\item [] $695=8\times 88-8-8/8.$
\item [] $696=8\times 88-8.$
\item [] $697=8\times 88-8+8/8.$
\item [] $698=8\times 88-8+(8+8)/8.$
\item [] $699=8\times 88-8+(8+8+8)/8.$
\item [] $700=8\times 88-8\times 8/(8+8).$
\item [] $701=8\times 88+8-88/8.$
\item [] $702=8\times 88-(8+8)/8.$
\item [] $703=8\times 88-8/8.$
\item [] $704=8\times 88.$
\item [] $705=8\times 88+8/8.$
\item [] $706=8\times 88+(8+8)/8.$
\item [] $707=8\times 88-8+88/8.$
\item [] $708=8\times (88+8/(8+8)).$
\item [] $709=8\times 88+8-(8+8+8)/8.$
\item [] $710=8\times 88+8-(8+8)/8.$
\item [] $711=8\times 88+8-8/8.$
\item [] $712=8\times 88+8.$
\item [] $713=8\times 88+8+8/8.$
\item [] $714=8\times 88+8+(8+8)/8.$
\item [] $715=8\times 88+88/8.$
\item [] $716=8\times 88+(88+8)/8.$
\item [] $717=8\times 88+(88+8+8)/8.$
\item [] $718=8\times 88+8+8-(8+8)/8.$
\item [] $719=8\times 88+8+8-8/8.$
\item [] $720=8\times 88+8+8.$
\item [] $721=8\times 88+8+8+8/8.$
\item [] $722=8\times 88+8+8+(8+8)/8.$
\item [] $723=88\times 8+8+88/8.$
\item [] $724=8\times (88+8/(8+8))+8+8.$
\item [] $725=8\times 88+(88+88-8)/8.$
\item [] $726=8\times 88+(88+88)/8.$
\item [] $727=8\times (88-8)+88-8/8.$
\item [] $728=8\times (88-8)+88.$
\item [] $729=8\times (88-8)+88+8/8.$
\item [] $730=8\times (88-8)+88+(8+8)/8.$
\item [] $731=8\times (88-8)+88-8+88/8.$
\item [] $732=8\times (88+8/(8+8))+8+8+8.$
\item [] $733=888-88-8\times 8-(8+8+8)/8.$
\item [] $734=888-88-8\times 8-(8+8)/8.$
\item [] $735=888-88-8\times 8-8/8.$
\item [] $736=888-88-8\times 8.$
\item [] $736=888-88-8\times 8.$
\item [] $737=888-88-8\times 8+8/8.$
\item [] $738=888-88-8\times 8+(8+8)/8.$
\item [] $739=8\times (88-8)+88+88/8.$
\item [] $740=8\times (88-8)+88+(88+8)/8.$
\item [] $741=8\times (88+8)-8-8-88/8.$
\item [] $742=8\times (88+8)-8-8-(88-8)/8.$
\item [] $743=8\times (88+8)-8-8-8-8/8.$
\item [] $744=8\times (88+8)-8-8-8.$
\item [] $745=8\times (88+8)-8-8-8+8/8.$
\item [] $746=8\times (88+8)-(88+88)/8.$
\item [] $747=888-88-8\times 8+88/8.$
\item [] $748=8\times (88+8)-8-(88+8)/8.$
\item [] $749=8\times (88+8)-8-88/8.$
\item [] $750=8\times (88+8)-8-(88-8)/8.$
\item [] $751=8\times (88+8)-8-8-8/8.$
\item [] $752=8\times (88+8)-8-8.$
\item [] $753=8\times (88+8)-8-8+8/8.$
\item [] $754=8\times (88+8)-8-8+(8+8)/8.$
\item [] $755=8\times (88+8)-(88+8+8)/8.$
\item [] $756=8\times (88+8)-(88+8)/8.$
\item [] $757=8\times (88+8)-88/8.$
\item [] $758=8\times (88+8)-(88-8)/8.$
\item [] $759=8\times (88+8)-8-8/8.$
\item [] $760=8\times (88+8)-8.$
\item [] $761=8\times (88+8)-8+8/8.$
\item [] $762=8\times (88+8)-8+(8+8)/8.$
\item [] $763=8\times (88+8)-8+(8+8+8)/8.$
\item [] $764=8\times (88+8-8/(8+8)).$
\item [] $765=8\times (88+8)+8-88/8.$
\item [] $766=(88+8)\times 8-(8+8)/8.$
\item [] $767=(88+8)\times 8-8/8.$
\item [] $768=(88+8)\times 8.$
\item [] $769=(88+8)\times 8+8/8.$
\item [] $770=(88+8)\times 8+(8+8)/8.$
\item [] $771=8\times (88+8)-8+88/8.$
\item [] $772=8\times (88+8)+8\times 8/(8+8).$
\item [] $773=8\times (88+8)+8+8-88/8.$
\item [] $774=8\times (88+8)+8-(8+8)/8.$
\item [] $775=8\times (88+8)+8-8/8.$
\item [] $776=8\times (88+8)+8.$
\item [] $777=888-888/8.$
\item [] $778=8\times (88+8)+8+(8+8)/8.$
\item [] $779=8\times (88+8)+88/8.$
\item [] $780=8\times (88+8+8/(8+8))+8.$
\item [] $781=8\times 88+88-88/8.$
\item [] $782=8\times (88+8)+8+8-(8+8)/8.$
\item [] $783=8\times (88+8)+8+8-8/8.$
\item [] $784=8\times 88+88-8.$
\item [] $785=888+8-888/8.$
\item [] $786=8\times (88+8)+8+(88-8)/8.$
\item [] $787=8\times (88+8)+8+88/8.$
\item [] $788=8\times (88-8/(8+8))+88.$
\item [] $789=888-88-88/8.$
\item [] $790=8\times 88+88-(8+8)/8.$
\item [] $791=8\times 88+88-8/8.$
\item [] $792=8\times 88+88.$
\item [] $793=8\times 88+88+8/8.$
\item [] $794=8\times 88+88+(8+8)/8.$
\item [] $795=8\times 88+88+(8+8+8)/8.$
\item [] $796=888-88-8\times 8/(8+8).$
\item [] $797=888-88-(8+8+8)/8.$
\item [] $798=888-88-(8+8)/8.$
\item [] $799=888-88-8/8.$
\item [] $800=888-88.$
\item [] $801=888-88+8/8.$
\item [] $802=888-88+(8+8)/8.$
\item [] $803=88\times (8\times 8+8+8/8)/8.$
\item [] $804=88+8+8\times (88+8/(8+8)).$
\item [] $805=8\times 88+8888/88.$
\item [] $806=8\times 88-8+(888-8)/8.$
\item [] $807=888-88+8-8/8.$
\item [] $808=888-88+8.$
\item [] $809=888-88+8+8/8.$
\item [] $810=888-88+(88-8)/8.$
\item [] $811=888-88+88/8.$
\item [] $812=888-8\times 8-(88+8)/8.$
\item [] $813=888-8\times 8-88/8.$
\item [] $814=8\times 88+(888-8)/8.$
\item [] $815=8\times 88+888/8.$
\item [] $816=888-8\times 8-8.$
\item [] $817=888-8\times 8-8+8/8.$
\item [] $818=888-8\times 8-8+(8+8)/8.$
\item [] $819=888-88+8+88/8.$
\item [] $820=888-8\times (8+8/(8+8)).$
\item [] $821=888-8\times 8-(8+8+8)/8.$
\item [] $822=888-8\times 8-(8+8)/8.$
\item [] $823=888-8\times 8-8/8.$
\item [] $824=888-8\times 8.$
\item [] $825=888-8\times 8+8/8.$
\item [] $826=888-8\times 8+(8+8)/8.$
\item [] $827=888-8\times 8+(8+8+8)/8.$
\item [] $828=888-8\times (8-8/(8+8)).$
\item [] $829=8\times (88+8+8)+8-88/8.$
\item [] $830=8\times (88+8+8)-(8+8)/8.$
\item [] $831=8\times (88+8+8)-8/8.$
\item [] $832=8\times (88+8+8).$
\item [] $833=8\times (88+8+8)+8/8.$
\item [] $834=8\times (88+8+8)+(8+8)/8.$
\item [] $835=888-8\times 8+88/8.$
\item [] $836=8\times (88+8+8+8/(8+8)).$
\item [] $837=8\times (88+8+8)+8-(8+8+8)/8.$
\item [] $838=8\times (88+8+8)+8-(8+8)/8.$
\item [] $839=8\times (88+8+8)+8-8/8.$
\item [] $840=8\times (88+8+8)+8.$
\item [] $841=8\times (88+8+8)+8+8/8.$
\item [] $842=8\times (88+8+8)+(88-8)/8.$
\item [] $843=8\times (88+8+8)+88/8.$
\item [] $844=888-8\times 88/(8+8).$
\item [] $845=8\times (88+8+8)+(88+8+8)/8.$
\item [] $846=(8+8/8)\times (88+8-(8+8)/8).$
\item [] $847=8\times (88+8+8)+8+8-8/8.$
\item [] $848=(8+8)\times (8\times 8-88/8).$
\item [] $849=8\times (88+8)+88-8+8/8.$
\item [] $850=8\times (88+8+8)+8+(88-8)/8.$
\item [] $851=8\times (88+8+8)+8+88/8.$
\item [] $852=8\times (88+8+8)+8+(88+8)/8.$
\item [] $853=888-8-8-8-88/8.$
\item [] $854=8\times (88+8)+88-(8+8)/8.$
\item [] $855=8\times (88+8)+88-8/8.$
\item [] $856=8\times (88+8)+88.$
\item [] $857=8\times (88+8)+88+8/8.$
\item [] $858=8\times (88+8)+88+(8+8)/8.$
\item [] $859=8\times (88+8)+88+(8+8+8)/8.$
\item [] $860=88+8\times (88+8+8/(8+8)).$
\item [] $861=888-8-8-88/8.$
\item [] $862=888-8-8-(88-8)/8.$
\item [] $863=888-8-8-8-8/8.$
\item [] $864=888-8-8-8.$
\item [] $865=888-8-8-8+8/8.$
\item [] $866=888-8-8-8+(8+8)/8.$
\item [] $867=8\times (88+8)+88+88/8.$
\item [] $868=888-8-(88+8)/8.$
\item [] $869=888-8-88/8.$
\item [] $870=888-8-(88-8)/8.$
\item [] $871=888-8-8-8/8.$
\item [] $872=888-8-8.$
\item [] $873=888-8-8+8/8.$
\item [] $874=888-8-8+(8+8)/8.$
\item [] $875=888-(88+8+8)/8.$
\item [] $876=888-(88+8)/8.$
\item [] $877=888-88/8.$
\item [] $878=888-(88-8)/8.$
\item [] $879=888-8-8/8.$
\item [] $880=888-8.$
\item [] $881=888-8+8/8.$
\item [] $882=888-8+(8+8)/8.$
\item [] $883=888-8-8+88/8.$
\item [] $884=888-8\times 8/(8+8).$
\item [] $885=888+8-88/8.$
\item [] $886=888-(8+8)/8.$
\item [] $887=888-8/8.$
\item [] $888=888.$
\item [] $890=888+(8+8)/8.$
\item [] $891=888-8+88/8.$
\item [] $892=888+8\times 8/(8+8).$
\item [] $893=888+8+8-88/8.$
\item [] $894=888+8-(8+8)/8.$
\item [] $895=888+8-8/8.$
\item [] $896=888+8.$
\item [] $897=888+8+8/8.$
\item [] $898=888+8+(8+8)/8.$
\item [] $899=888+88/8.$
\item [] $900=888+(88+8)/8.$
\item [] $901=888+(88+8+8)/8.$
\item [] $902=888+8+8-(8+8)/8.$
\item [] $903=888+8+8-8/8.$
\item [] $904=888+8+8.$
\item [] $905=888+8+8+8/8.$
\item [] $906=888+8+(88-8)/8.$
\item [] $907=888+8+88/8.$
\item [] $908=888+8+(88+8)/8.$
\item [] $909=888+8+(88+8+8)/8.$
\item [] $910=888+(88+88)/8.$
\item [] $911=888+8+8+8-8/8.$
\item [] $912=888+8+8+8.$
\item [] $913=888+8+8+8+8/8.$
\item [] $914=888+8+8+8+(8+8)/8.$
\item [] $915=888+8+8+88/8.$
\item [] $916=888+8+8+(88+8)/8.$
\item [] $917=888+8+(88+88-8)/8.$
\item [] $918=888+8+(88+88)/8.$
\item [] $919=888+8+8+8+8-8/8.$
\item [] $920=888+8+8+8+8.$
\item [] $921=888+8+8+8+8+8/8.$
\item [] $922=888+8+8+8+8+(8+8)/8.$
\item [] $923=888+8+8+8+88/8.$
\item [] $924=888+8+8+8+(88+8)/8.$
\item [] $925=888+888/(8+8+8).$
\item [] $926=8\times 88+888\times (8+8)/(8\times 8).$
\item [] $927=(8+8/8)\times (888/8-8).$
\item [] $928=888+8+8+8+8+8.$
\item [] $929=888+8+8+8+8+8+8/8.$
\item [] $930=(8\times 8-(8+8)/8)\times (8+8-8/8).$
\item [] $931=888+8+8+8+8+88/8.$
\item [] $932=888+8\times 88/(8+8).$
\item [] $933=888+8\times 8-8-88/8.$
\item [] $934=888+8\times 8-8-8-(8+8)/8.$
\item [] $935=888+8\times 8-8-8-8/8.$
\item [] $936=888+8\times 8-8-8.$
\item [] $937=888+8\times 8-8-8+8/8.$
\item [] $938=88\times (88-(8+8)/8)/8-8.$
\item [] $939=888+8\times 8-(88+8+8)/8.$
\item [] $940=888+8+88\times 8/(8+8).$
\item [] $941=888+8\times 8-88/8.$
\item [] $942=888+8\times 8-8-(8+8)/8.$
\item [] $943=888+8\times 8-8-8/8.$
\item [] $944=888+8\times 8-8.$
\item [] $945=888+8\times 8-8+8/8.$
\item [] $946=888+8\times 8-8+(8+8)/8.$
\item [] $947=888+8\times 8-8+(8+8+8)/8.$
\item [] $948=888+8\times (8-8/(8+8)).$
\item [] $949=8\times (8\times (8+8)-8)-88/8.$
\item [] $950=888+8\times 8-(8+8)/8.$
\item [] $951=888+8\times 8-8/8.$
\item [] $952=888+8\times 8.$
\item [] $953=888+8\times 8+8/8.$
\item [] $954=888+8\times 8+(8+8)/8.$
\item [] $955=888+8\times 8+(8+8+8)/8.$
\item [] $956=888+8\times (8+8/(8+8)).$
\item [] $957=88\times (88-8/8)/8.$
\item [] $958=8\times (8\times (8+8)-8)-(8+8)/8.$
\item [] $959=8\times (8\times (8+8)-8)-8/8.$
\item [] $960=8\times (8\times (8+8)-8).$
\item [] $961=8\times (8\times (8+8)-8)+8/8.$
\item [] $962=8\times (8\times (8+8)-8)+(8+8)/8.$
\item [] $963=888+8\times 8+88/8.$
\item [] $964=8\times (8\times (8+8)-8+8/(8+8)).$
\item [] $965=888+88-88/8.$
\item [] $966=(88\times 88-8-8)/8.$
\item [] $967=(88\times 88-8)/8.$
\item [] $968=88\times 88/8.$
\item [] $969=(88\times 88+8)/8.$
\item [] $970=(88\times 88+8+8)/8.$
\item [] $971=(8\times (8+8)-8)\times 8+88/8.$
\item [] $972=888+88-8\times 8/(8+8).$
\item [] $973=888+88-(8+8+8)/8.$
\item [] $974=888+88-(8+8)/8.$
\item [] $975=888+88-8/8.$
\item [] $976=888+88.$
\item [] $977=888+88+8/8.$
\item [] $978=888+88+(8+8)/8.$
\item [] $979=88\times (88+8/8)/8.$
\item [] $980=888+88+8\times 8/(8+8).$
\item [] $981=888+88+8-(8+8+8)/8.$
\item [] $982=888+88+8-(8+8)/8.$
\item [] $983=888+88+8-8/8.$
\item [] $984=888+88+8.$
\item [] $985=888+88+8+8/8.$
\item [] $986=888+88+8+(8+8)/8.$
\item [] $987=888+88+88/8.$
\item [] $988=888+88+(88+8)/8.$
\item [] $989=888+8888/88.$
\item [] $990=88\times (88+(8+8)/8)/8.$
\item [] $991=888+888/8-8.$
\item [] $992=888+88+8+8.$
\item [] $993=888+88+8+8+8/8.$
\item [] $994=888+88+8+8+(8+8)/8.$
\item [] $995=888+88+8+88/8.$
\item [] $996=888+88+8+(88+8)/8.$
\item [] $997=888+(888-8-8)/8.$
\item [] $998=888+(888-8)/8.$
\item [] $999=888+888/8.$
\item [] $1000=888+88+8+8+8.$
\end{itemize}
\end{multicols}
}

\section{\textbf{Representations Using Number 9}}

{\footnotesize
\begin{multicols}{3}
\begin{itemize}
\item [] $101=99+(9+9)/9.$
\item [] $102=999/9-9.$
\item [] $103=(999+9)/9-9.$
\item [] $104=(999+9+9)/9-9.$
\item [] $105=99+(99+9)/(9+9).$
\item [] $106=99+9-(9+9)/9.$
\item [] $107=99+9-9/9.$
\item [] $108=99+9.$
\item [] $109=(9+(9/9))+99.$
\item [] $110=(999-9)/9.$
\item [] $111=999/9.$
\item [] $112=(999+9)/9.$
\item [] $113=((9+9)+999)/9.$
\item [] $114=(999+9+9+9)/9.$
\item [] $115=99+9+9-(9+9)/9.$
\item [] $116=99+9+9-9/9.$
\item [] $117=99+9+9.$
\item [] $118=99+9+9+9/9.$
\item [] $119=9+(999-9)/9.$
\item [] $120=9+999/9.$
\item [] $121=9+(999+9)/9.$
\item [] $122=(999+99)/9.$
\item [] $123=(999+99+9)/9.$
\item [] $124=99+9+9+9-(9+9)/9.$
\item [] $125=99+9+9+9-9/9.$
\item [] $126=99+9+9+9.$
\item [] $127=99+9+9+9+9/9.$
\item [] $128=99+9+9+99/9.$
\item [] $129=9+9+999/9.$
\item [] $130=9+9+(999+9)/9.$
\item [] $131=9+(999+99)/9.$
\item [] $132=99\times(99+9)/(9\times9).$
\item [] $133=9\times(9+9)-9-9-99/9.$
\item [] $134=99+9+9+9+9-9/9.$
\item [] $135=99+9+9+9+9.$
\item [] $136=(9-9/9)\times(9+9-9/9).$
\item [] $137=99+9+9+9+99/9.$
\item [] $138=9+9+9+999/9.$
\item [] $139=9+9+9+(999+9)/9.$
\item [] $140=9+9+(999+99)/9.$
\item [] $141=9\times(9+9)-9-(99+9)/9.$
\item [] $142=9\times(9+9)-9-99/9.$
\item [] $143=9\times(9+9)-9-9-9/9.$
\item [] $144=(9+9)\times(9-9/9).$
\item [] $145=9\times(9+9)-9-9+9/9.$
\item [] $146=9\times(9+9)-9-9+(9+9)/9.$
\item [] $147=9+9+9+9+999/9.$
\item [] $148=99+(99\times9-9)/(9+9).$
\item [] $149=9\times(9+9)-(99+9+9)/9.$
\item [] $150=9\times(9+9)-(99+9)/9.$
\item [] $151=9\times(9+9)-99/9.$
\item [] $152=9\times(9+9)-9-9/9.$
\item [] $153=9\times(9+9)-9.$
\item [] $154=9\times(9+9)-9+9/9.$
\item [] $155=9\times(9+9)-9+(9+9)/9.$
\item [] $156=9\times(9+9)-9+(9+9+9)/9.$
\item [] $157=((9+9)/9)^{(9-9/9)}-99.$
\item [] $158=9\times(9+9)-(9\times9-9)/(9+9).$
\item [] $159=9\times(9+9)-(9+9+9)/9.$
\item [] $160=9\times(9+9)-(9+9)/9.$
\item [] $161=9\times(9+9)-9/9.$
\item [] $162=9\times(9+9).$
\item [] $163=9\times(9+9)+9/9.$
\item [] $164=9\times(9+9)+(9+9)/9.$
\item [] $165=9\times(9+9)+(9+9+9)/9.$
\item [] $166=(9+9)/9\times(9\times9+(9+9)/9).$
\item [] $167=9\times(9+9)+(99-9)/(9+9).$
\item [] $168=9\times(9+9)+9-(9+9+9)/9.$
\item [] $169=9\times9+99-99/9.$
\item [] $170=9\times(9+9)+9-9/9.$
\item [] $171=9\times(9+9)+9.$
\item [] $172=9\times(9+9)+9+9/9.$
\item [] $173=9\times(9+9)+99/9.$
\item [] $174=9\times(9+9)+(99+9)/9.$
\item [] $175=((9+9)/9)^{(9-9/9)}-9\times9.$
\item [] $176=(9+9)\times(99-99/9)/9.$
\item [] $177=99+9\times9-(9+9+9)/9.$
\item [] $178=99+9\times9-(9+9)/9.$
\item [] $179=99+9\times9-9/9.$
\item [] $180=99+9\times9.$
\item [] $181=99+9\times9+9/9.$
\item [] $182=9\times(9+9)+9+99/9.$
\item [] $183=9\times9-9+999/9.$
\item [] $184=9\times(9+9)+(99+99)/9.$
\item [] $185=99+99-(99+9+9)/9.$
\item [] $186=99+99-(99+9)/9.$
\item [] $187=99+99-99/9.$
\item [] $188=99+99-9-9/9.$
\item [] $189=99+99-9.$
\item [] $190=99+99-9+9/9.$
\item [] $191=99+9\times9+99/9.$
\item [] $192=9\times9+999/9.$
\item [] $193=9\times9+(999+9)/9.$
\item [] $194=9\times9+(999+9+9)/9.$
\item [] $195=99+99-(9+9+9)/9.$
\item [] $196=99+99-(9+9)/9.$
\item [] $197=99+99-9/9.$
\item [] $198=99+99.$
\item [] $199=99+99+9/9.$
\item [] $200=99+99+(9+9)/9.$
\item [] $201=99-9+999/9.$
\item [] $202=99-9+(999+9)/9.$
\item [] $203=9\times9+(999+99)/9.$
\item [] $204=(9+9)\times(999/9-9)/9.$
\item [] $205=99+99+9-(9+9)/9.$
\item [] $206=99+99+9-9/9.$
\item [] $207=99+99+9.$
\item [] $208=99+99+9+9/9.$
\item [] $209=99+(999-9)/9.$
\item [] $210=99+999/9.$
\item [] $211=99+(999+9)/9.$
\item [] $212=99+(999+9+9)/9.$
\item [] $213=(9+9)\times(9+9)-999/9.$
\item [] $214=99+(9+9)\times(9-9/9)/9.$
\item [] $215=99+99+9+9-9/9.$
\item [] $216=99+99+9+9.$
\item [] $217=9\times9\times9-((9+9)/9)^9.$
\item [] $218=99+9+(999-9)/9.$
\item [] $219=99+9+999/9.$
\item [] $220=99\times(9+99/9)/9.$
\item [] $221=99+(999+99)/9.$
\item [] $222=999\times(9+9)/(9\times9).$
\item [] $223=9\times(9+9+9)-9-99/9.$
\item [] $224=9\times(9+9+9)-9-9-9/9.$
\item [] $225=9\times(9+9+9)-9-9.$
\item [] $226=(9+9)\times(9+9)-99+9/9.$
\item [] $227=99+9+9+(999-9)/9.$
\item [] $228=99+9+9+999/9.$
\item [] $229=9+99\times(9+99/9)/9.$
\item [] $230=99+9+(999+99)/9.$
\item [] $231=9\times(9+9+9)-(99+9)/9.$
\item [] $232=9\times(9+9+9)-99/9.$
\item [] $233=9\times(9+9+9)-9-9/9.$
\item [] $234=9\times(9+9+9)-9.$
\item [] $235=9\times(9+9+9)-9+9/9.$
\item [] $236=9\times(9+9+9)-9+(9+9)/9.$
\item [] $237=9\times(9+9+9)-9+(9+9+9)/9.$
\item [] $238=((9+9)/9)^{(9-9/9)}-9-9.$
\item [] $239=9\times(9+9+9)-(9+9+9+9)/9.$
\item [] $240=9\times(9+9+9)-(9+9+9)/9.$
\item [] $241=9\times(9+9+9)-(9+9)/9.$
\item [] $242=9\times(9+9+9)-9/9.$
\item [] $243=9\times(9+9+9).$
\item [] $244=9\times(9+9+9)+9/9.$
\item [] $245=9\times(9+9+9)+(9+9)/9.$
\item [] $246=9\times(9\times9\times9+9)/(9+9+9).$
\item [] $247=((9+9)/9)^{(9-9/9)}-9.$
\item [] $248=(9\times9+9)/(9+9)+9\times(9+9+9).$
\item [] $249=9+(9+9)\times(9+999/9)/9.$
\item [] $250=9\times(9+9)+99-99/9.$
\item [] $251=9\times(9+9+9)+9-9/9.$
\item [] $252=9\times(9+9+9)+9.$
\item [] $253=9\times(9+9+9)+9+9/9.$
\item [] $254=9\times(9+9+9)+99/9.$
\item [] $255=((9+9)/9)^{(9-9/9)}-9/9.$
\item [] $256=((9+9)/9)^{(9-9/9)}.$
\item [] $257=((9+9)/9)^{(9-9/9)}+9/9.$
\item [] $258=9\times(9+9)+99-(9+9+9)/9.$
\item [] $259=9\times(9+9)+99-(9+9)/9.$
\item [] $260=9\times(9+9)+99-9/9.$
\item [] $261=9\times(9+9)+99.$
\item [] $262=9\times(9+9)+99+9/9.$
\item [] $263=9\times(9+9+9)+9+99/9.$
\item [] $264=9\times(9+9)-9+999/9.$
\item [] $265=((9+9)/9)^{(9-9/9)}+9.$
\item [] $266=((9+9)/9)^{(9-9/9)}+9+9/9.$
\item [] $267=((9+9)/9)^{(9-9/9)}+99/9.$
\item [] $268=9\times9+99+99-99/9.$
\item [] $269=999-9\times9\times9-9/9.$
\item [] $270=999-9\times9\times9.$
\item [] $271=999-9\times9\times9+9/9.$
\item [] $272=9\times(9+9)+99+99/9.$
\item [] $273=9\times(9+9)+999/9.$
\item [] $274=9\times(9+9)+(999+9)/9.$
\item [] $275=99\times(9+9+9-(9+9)/9)/9.$
\item [] $276=(9+9+9)\times(9\times9+99/9)/9.$
\item [] $277=99+99+9\times9-(9+9)/9.$
\item [] $278=99+99+9\times9-9/9.$
\item [] $279=99+99+9\times9.$
\item [] $280=99+99+9\times9+9/9.$
\item [] $281=999-9\times9\times9+99/9.$
\item [] $282=9\times(9+9)+9+999/9.$
\item [] $283=9\times(9+9)+9+(999+9)/9.$
\item [] $284=9\times(9+9)+(999+99)/9.$
\item [] $285=99+99+99-(99+9)/9.$
\item [] $286=99/9\times(9+9+9-9/9).$
\item [] $287=99+99+99-9-9/9.$
\item [] $288=99+99+99-9.$
\item [] $289=(9+9-9/9)^((9+9)/9).$
\item [] $290=(9+9/9)\times(9+9+99/9).$
\item [] $291=99+9\times9+999/9.$
\item [] $292=99+9\times9+(999+9)/9.$
\item [] $293=(9+9)\times(9+9)-9-(99+99)/9.$
\item [] $294=(9+9+9)\times(99-9/9)/9.$
\item [] $295=(9+9)\times(9+9)-9-9-99/9.$
\item [] $296=99+99+99-9/9.$
\item [] $297=99+99+99.$
\item [] $298=99+99+99+9/9.$
\item [] $299=(99\times(9+9+9)+9+9)/9.$
\item [] $300=(9+9+9)\times(99+9/9)/9.$
\item [] $301=(9+99/9)^((9+9)/9)-99.$
\item [] $302=(9+9)\times(9\times(9+9)-99/9)/9.$
\item [] $303=(9+9)\times(9+9)-9-(99+9)/9.$
\item [] $304=(9+9)\times(9+9)-9-99/9.$
\item [] $305=(9+9)\times(9+9)-9-9-9/9.$
\item [] $306=(9+9)\times(9+9-9/9).$
\item [] $307=(9+9)\times(9+9)-9-9+9/9.$
\item [] $308=99\times(9+9+9+9/9)/9.$
\item [] $309=99+99+999/9.$
\item [] $310=99+99+(999+9)/9.$
\item [] $311=(9+9)\times(9+9)-(99+9+9)/9.$
\item [] $312=(9+9)\times(9+9)-(99+9)/9.$
\item [] $313=(9+9)\times(9+9)-99/9.$
\item [] $314=(9+9)\times(9+9)-9-9/9.$
\item [] $315=(9+9)\times(9+9)-9.$
\item [] $316=(9+9)\times(9+9)-9+9/9.$
\item [] $317=(9+9)\times(9+9)-9+(9+9)/9.$
\item [] $318=99+99+9+999/9.$
\item [] $319=99\times(9+9+99/9)/9.$
\item [] $320=(9+9+9+9)\times(9\times9-9/9)/9.$
\item [] $321=(9+9)\times(9+9)-(9+9+9)/9.$
\item [] $322=(9+9)\times(9+9)-(9+9)/9.$
\item [] $323=(9+9)\times(9+9)-9/9.$
\item [] $324=(9+9)\times(9+9).$
\item [] $325=(9+9)\times(9+9)+9/9.$
\item [] $326=(9+9)\times(9+9)+(9+9)/9.$
\item [] $327=(9+9)\times(9+9)+(9+9+9)/9.$
\item [] $328=(9\times9\times9+9)\times(9-9/9)/(9+9).$
\item [] $329=(9+9)\times(9+9)+(99-9)/(9+9).$
\item [] $330=9\times99\times(9+9/9)/(9+9+9).$
\item [] $331=(9+9)\times(9+9)+9-(9+9)/9.$
\item [] $332=(9+9)\times(9+9)+9-9/9.$
\item [] $333=(9+9)\times(9+9)+9.$
\item [] $334=(9+9)\times(9+9)+9+9/9.$
\item [] $335=(9+9)\times(9+9)+99/9.$
\item [] $336=(9+9)\times(9+9)+(99+9)/9.$
\item [] $337=9\times9+((9+9)/9)^{(9-9/9)}.$
\item [] $338=9\times9+((9+9)/9)^{(9-9/9)}+9/9.$
\item [] $339=9\times(9+9+9)+99-(9+9+9)/9.$
\item [] $340=(9+9-9/9)\times(9+99/9).$
\item [] $341=9\times(9+9+9)+99-9/9.$
\item [] $342=9\times(9+9+9)+99.$
\item [] $343=9\times(9+9+9)+99+9/9.$
\item [] $344=(9+9)\times(9+9)+9+99/9.$
\item [] $345=9\times(9+9+9)-9+999/9.$
\item [] $346=((9+9)/9)\times((9+9)\times9+99/9).$
\item [] $347=9\times(9\times99+9/9)/(9+9)-99.$
\item [] $348=(9+9+99/9)\times(99+9)/9.$
\item [] $349=((9+9)/9)^9-9\times(9+9)-9/9.$
\item [] $350=((9+9)/9)^9-9\times(9+9).$
\item [] $351=9\times(9+9+9)+99+9.$
\item [] $352=(9+9+9/9)^((9+9)/9)-9.$
\item [] $353=(9+9)\times(9+9)+9+9+99/9.$
\item [] $354=9\times(9+9+9)+999/9.$
\item [] $355=99+((9+9)/9)^{(9-9/9)}.$
\item [] $356=9\times(9\times9\times9+9/9)/(9+9)-9.$
\item [] $357=(9+(9+99)/9)\times(9+9-9/9).$
\item [] $358=(9+9)\times(9+99/9)-(9+9)/9.$
\item [] $359=(9+9)\times(9+99/9)-9/9.$
\item [] $360=(9+9)\times(9+99/9).$
\item [] $361=(9+9+9/9)^((9+9)/9).$
\item [] $362=(9+9+9/9)^((9+9)/9)+9/9.$
\item [] $363=9\times(9+9+9)+9+999/9.$
\item [] $364=9\times(9\times9\times9-9/9)/(9+9).$
\item [] $365=9\times(9\times9\times9+9/9)/(9+9).$
\item [] $366=9\times(999+99)/(9+9+9).$
\item [] $367=9\times9+99\times(9+9+9-9/9)/9.$
\item [] $368=9+9+((9+9)/9)^9-9\times(9+9).$
\item [] $369=(9+9)\times(9+99/9)+9.$
\item [] $370=9\times99-9-((9+9)/9)^9.$
\item [] $371=(9+9)\times(99/9+9)+99/9.$
\item [] $372=9\times(9+9)+99+999/9.$
\item [] $373=9999/9-9\times9\times9-9.$
\item [] $374=(99+99)\times(9+9-9/9)/9.$
\item [] $375=(9+9)\times(9+9+9)-999/9.$
\item [] $376=(9+9+9+9+99/9)\times(9-9/9).$
\item [] $377=(9+9)\times(9+(9+99)/9)-9/9.$
\item [] $378=((99+9)/9+9)\times(9+9).$
\item [] $379=9\times99-((9+9)/9)^9.$
\item [] $380=(99/9+9)\times(9+9+9/9).$
\item [] $381=9\times(9+9)+9+99+999/9.$
\item [] $382=9999/9-9\times9\times9.$
\item [] $383=(9999+9)/9-9\times9\times9.$
\item [] $384=(9\times9+999/9)\times(9+9)/9.$
\item [] $385=99\times(9+9+9+9-9/9)/9.$
\item [] $386=(9+9)\times(9+9+9)-99-9/9.$
\item [] $387=(9+9)\times(9+9+9)-99.$
\item [] $388=9\times99-((9+9)/9)^9+9.$
\item [] $389=(9+9+9/9)\times(9+99/9)+9.$
\item [] $390=(999/9+9-9\times9)\times(9+9/9).$
\item [] $391=(9+99/9)^((9+9)/9)-9.$
\item [] $392=(9\times9-9)\times(99-9/9)/(9+9).$
\item [] $393=(9+9)\times(9+9)+9\times9-(99+9)/9.$
\item [] $394=(9+9)\times(9+9)+9\times9-99/9.$
\item [] $395=((9+9)/9)^9-9-99-9.$
\item [] $396=(9+9)\times(99+99)/9.$
\item [] $397=((9+9)\times(99+99)+9)/9.$
\item [] $398=(9+9)\times(99+99+9/9)/9.$
\item [] $399=(9+99/9)^((9+9)/9)-9/9.$
\item [] $400=(9+99/9)^((9+9)/9).$
\item [] $401=((9+9)/9)^9-999/9.$
\item [] $402=((9+9)/9)^9-99-99/9.$
\item [] $403=((9+9)/9)^9-99-9-9/9.$
\item [] $404=((9+9)/9)^9-99-9.$
\item [] $405=(9+9)\times(9+9)+9\times9.$
\item [] $406=(9+9)\times(9+9)+9\times9+9/9.$
\item [] $407=99\times(9+9+9+9+9/9)/9.$
\item [] $408=99+99+99+999/9.$
\item [] $409=(9+99/9)^((9+9)/9)+9.$
\item [] $410=((9+9)/9)^9+9-999/9.$
\item [] $411=((9+9)/9)^9-99-(9+9)/9.$
\item [] $412=((9+9)/9)^9-99-9/9.$
\item [] $413=((9+9)/9)^9-99.$
\item [] $414=(9+9)\times(9+9)+99-9.$
\item [] $415=(9+9)\times(9+9)+99-9+9/9.$
\item [] $416=(9+9)\times(9+9)+9\times9+99/9.$
\item [] $417=(9+9)\times(9+9)+9\times9+(99+9)/9.$
\item [] $418=(9+99/9)^((9+9)/9)+9+9.$
\item [] $419=((9+9)/9)^9-9\times9-(99+9)/9.$
\item [] $420=(9+9)\times(99+999/9)/9.$
\item [] $421=((9+9)/9)^9-99+9-9/9.$
\item [] $422=((9+9)/9)^9-99+9.$
\item [] $423=(9+9)\times(9+9)+99.$
\item [] $424=(9+9)\times(9+9)+99+9/9.$
\item [] $425=(9+9)\times(9+9)+99+(9+9)/9.$
\item [] $426=(9+9)\times(9+9)-9+999/9.$
\item [] $427=(9+9)\times(9+9)-9+(999+9)/9.$
\item [] $428=((9+9)/9)^9-9\times9-(9+9+9)/9.$
\item [] $429=((9+9)/9)^9-9\times9-(9+9)/9.$
\item [] $430=((9+9)/9)^9-9\times9-9/9.$
\item [] $431=((9+9)/9)^9-9\times9.$
\item [] $432=(9+9)\times(9+9)+99+9.$
\item [] $433=(9+9)\times(9+9)+99+9+9/9.$
\item [] $434=(9+9)\times(9+9)+(999-9)/9.$
\item [] $435=(9+9)\times(9+9)+999/9.$
\item [] $436=(9+9)\times(9+9)+(999+9)/9.$
\item [] $437=9\times(9\times99+9/9)/(9+9)-9.$
\item [] $438=9\times(9\times9-9)-99-999/9.$
\item [] $439=((9+9)/9)^9+9-9\times9-9/9.$
\item [] $440=((9+9)/9)^9+9-9\times9.$
\item [] $441=9\times(99\times9-9)/(9+9).$
\item [] $442=(9+9+9-9/9)\times(9+9-9/9).$
\item [] $443=(9+9)\times(9+9)+99+9+99/9.$
\item [] $444=999\times(9+9+9+9)/(9\times9).$
\item [] $445=9\times(9\times99-9/9)/(9+9).$
\item [] $446=9\times(9\times99+9/9)/(9+9).$
\item [] $447=9\times(9\times99+9/9)/(9+9)+9/9.$
\item [] $448=(9-9/9)\times(9+999)/(9+9).$
\item [] $449=((9+9)/9)^9+9+9-9\times9.$
\item [] $450=9\times(9\times99+9)/(9+9).$
\item [] $451=9\times(99\times9+99/9)/(9+9).$
\item [] $452=9\times(9\times99+9)/(9+9)+(9+9)/9.$
\item [] $453=9\times(9+9+9)+99+999/9.$
\item [] $454=9\times(99\times9-9/9)/(9+9)+9.$
\item [] $455=9\times(99\times9+9/9)/(9+9)+9.$
\item [] $456=(9+9+9+99/9)\times(9+99)/9.$
\item [] $457=(9+9+9)\times(9+9)-9-9-99/9.$
\item [] $458=((9+9)/9)^9+9+9+9-9\times9.$
\item [] $459=(9+9+9)\times(9+9-9/9).$
\item [] $460=9\times(99+9)-((9+9)/9)^9.$
\item [] $461=(9\times9-99/9)\times(9-9/9)-99.$
\item [] $462=9\times(9-9\times9)+(9999-9)/9.$
\item [] $463=9\times(9-9\times9)+9999/9.$
\item [] $464=(9+9+99/9)\times(9+9-(9+9)/9).$
\item [] $465=(9+9)\times(9+9+9)-9-(99+9)/9.$
\item [] $466=(9+9)\times(9+9+9)-9-99/9.$
\item [] $467=(9+9)\times(9+9+9-9/9)-9/9.$
\item [] $468=(9+9)\times(9+9+9-9/9).$
\item [] $469=(9+9)\times(9+9+9)-9-9+9/9.$
\item [] $470=(9+9)\times(9+9+9)-9-9+(9+9)/9.$
\item [] $471=(9+9)\times(9+9+9)-9-9+(9+9+9)/9.$
\item [] $472=(9+99/9)^((9+9)/9)+9\times9-9.$
\item [] $473=9\times9\times9-((9+9)/9)^{(9-9/9)}.$
\item [] $474=(9+9)\times(9+9+9)-(9+99)/9.$
\item [] $475=(9+9)\times(9+9+9)-99/9.$
\item [] $476=(9+9)\times(9+9+9)-9-9/9.$
\item [] $477=(9+9)\times(9+9+9)-9.$
\item [] $478=(9+9)\times(9+9+9)-9+9/9.$
\item [] $479=(9+9)\times(9+9+9)-9+(9+9)/9.$
\item [] $480=(9+9)\times(9+9+9)-9+(9+9+9)/9.$
\item [] $481=9\times9+(9+99/9)^((9+9)/9).$
\item [] $482=(9+9)\times(9+9+9)-(9\times9-9)/(9+9).$
\item [] $483=(9+9+9)\times((9+9)\times9-9/9)/9.$
\item [] $484=((99+99)/9)^((9+9)/9).$
\item [] $485=(9+9)\times(9+9+9)-9/9.$
\item [] $486=(9+9)\times(9+9+9).$
\item [] $487=(9+9)\times(9+9+9)+9/9.$
\item [] $488=(9+9)\times(9+9+9)+(9+9)/9.$
\item [] $489=(9+9+9)\times(9\times(9+9)+9/9)/9.$
\item [] $490=(99-9/9)\times(99-9)/(9+9).$
\item [] $491=((9+9)/9)^9-9-(9+99)/9.$
\item [] $492=((9+9)/9)^9-9-99/9.$
\item [] $493=((9+9)/9)^9-9-9-9/9.$
\item [] $494=((9+9)/9)^9-9-9.$
\item [] $495=(9+9)\times(9+9+9)+9.$
\item [] $496=(9+9)\times(9+9+9)+9+9/9.$
\item [] $497=(9+9)\times(9+9+9)+99/9.$
\item [] $498=(9+9)\times(9+9+9)+(99+9)/9.$
\item [] $499=(99/9+9)^((9+9)/9)+99.$
\item [] $500=((9+9)/9)^9-(9+99)/9.$
\item [] $501=((9+9)/9)^9-99/9.$
\item [] $502=((9+9)/9)^9-9-9/9.$
\item [] $503=((9+9)/9)^9-9.$
\item [] $504=((9+9)/9)^9-9+9/9.$
\item [] $505=((9+9)/9)^9-9+(9+9)/9.$
\item [] $506=((9+9)/9)^9-9+(9+9+9)/9.$
\item [] $507=((9+9)/9)^9-(9\times9+9)/(9+9).$
\item [] $508=((9+9)/9)^9-(9\times9-9)/(9+9).$
\item [] $509=((9+9)/9)^9-(9+9+9)/9.$
\item [] $510=((9+9)/9)^9-(9+9)/9.$
\item [] $511=((9+9)/9)^9-9/9.$
\item [] $512=((9+9)/9)^9.$
\item [] $513=((9+9)/9)^9+9/9.$
\item [] $514=((9+9)/9)^9+(9+9)/9.$
\item [] $515=((9+9)/9)^9+(9+9+9)/9.$
\item [] $516=((9+9)/9)^9+(9\times9-9)/(9+9).$
\item [] $517=((9+9)/9)^9+(9\times9+9)/(9+9).$
\item [] $518=((9+9)/9)^9+9-(9+9+9)/9.$
\item [] $519=((9+9)/9)^9+9-(9+9)/9.$
\item [] $520=((9+9)/9)^9+9-9/9.$
\item [] $521=((9+9)/9)^9+9.$
\item [] $522=((9+9)/9)^9+9+9/9.$
\item [] $523=((9+9)/9)^9+99/9.$
\item [] $524=((9+9)/9)^9+(99+9)/9.$
\item [] $525=((9+9)/9)^9+(99+9+9)/9.$
\item [] $526=((9+9)/9)^9+(99+9+9+9)/9.$
\item [] $528=((9+9)/9)^9+9+9-(9+9+9)/9.$
\item [] $528=((9+9)/9)^9+9+9-(9+9)/9.$
\item [] $529=((9+9)/9)^9+9+9-9/9.$
\item [] $530=((9+9)/9)^9+9+9.$
\item [] $531=9\times9\times9-99-99.$
\item [] $532=((9+9)/9)^9+9+99/9.$
\item [] $533=((9+9)/9)^9+9+(99+9)/9.$
\item [] $534=((9+9)/9)^9+9+(99+9+9)/9.$
\item [] $535=9\times(9\times9-9)-(999+9+9)/9.$
\item [] $536=9\times(9\times9-9)-(999+9)/9.$
\item [] $537=9\times(9\times9-9)-999/9.$
\item [] $538=((9+9)/9)^9+9+9+9-9/9.$
\item [] $539=((9+9)/9)^9+9+9+9.$
\item [] $540=(9+9+9)\times(9+99/9).$
\item [] $541=9\times(9\times9-9)-9-99+9/9.$
\item [] $542=9\times(9\times9-9)-9-99+(9+9)/9.$
\item [] $543=(9+9+9)\times(99+9/9+9\times9)/9.$
\item [] $554=9\times(9\times9-9-9)-(99+9+9)/9.$
\item [] $545=99+9\times(9\times99+9/9)/(9+9).$
\item [] $546=9\times(9\times9-9)-(999/9-9).$
\item [] $547=9\times(9\times9-9-9)-9-99/9.$
\item [] $548=9\times(9\times9-9)-99-9/9.$
\item [] $549=9\times(9\times9-9)-99.$
\item [] $550=9\times(9\times9-9)-99+9/9.$
\item [] $551=9\times(9\times9-9)-99+(9+9)/9.$
\item [] $552=(9\times9-(99+9)/9)\times(9-9/9).$
\item [] $553=(9\times9-(9+9)/9)\times(9-(9+9)/9).$
\item [] $554=9\times(9\times9-9-9)-(99+9+9)/9.$
\item [] $555=9\times(9\times9-9-9)-(99+9)/9.$
\item [] $556=9\times(9\times9-9-9)-99/9.$
\item [] $557=9\times(9\times9-9-9)-9-9/9.$
\item [] $558=9\times(9\times9-9-9)-9.$
\item [] $559=9\times(9\times9-9-9)-9+9/9.$
\item [] $560=(9\times9-99/9)\times(9-9/9).$
\item [] $561=9\times(9\times9-9-9)-9+(9+9+9)/9.$
\item [] $562=9\times(9+9)+(9+99/9)^((9+9)/9).$
\item [] $563=9\times(9\times9-9-9)-(9+9+9+9)/9.$
\item [] $564=9\times(9\times9-9-9)-(9+9+9)/9.$
\item [] $565=9\times(9\times9-9-9)-(9+9)/9.$
\item [] $566=9\times(9\times9-9-9)-9/9.$
\item [] $567=9\times(9\times9-9-9).$
\item [] $568=9\times(9\times9-9-9)+9/9.$
\item [] $569=9\times(9\times9-9-9)+(9+9)/9.$
\item [] $570=9\times(9\times9-9-9)+(9+9+9)/9.$
\item [] $571=9\times(9\times9-9-9)+(9+9+9+9)/9.$
\item [] $572=(9+9+9-9/9)\times(99+99)/9.$
\item [] $573=9\times9+((9+9)/9)^9-9-99/9.$
\item [] $574=(9-(9+9)/9)\times(9\times9+9/9).$
\item [] $575=(9\times9-9)\times(9-9/9)-9/9.$
\item [] $576=(9\times9-9)\times(9-9/9).$
\item [] $577=(9\times9-9)\times(9-9/9)+9/9.$
\item [] $578=9\times(9\times9-9-9)+99/9.$
\item [] $579=9\times(9\times9-9-9)+(99+9)/9.$
\item [] $580=(9+9+99/9)\times(9+99/9).$
\item [] $581=(9\times9+(9+9)/9)\times(9-(9+9)/9).$
\item [] $582=9\times9+((9+9)/9)^9-99/9.$
\item [] $583=9\times9+((9+9)/9)^9-9-9/9.$
\item [] $584=9\times9+((9+9)/9)^9-9.$
\item [] $585=9\times(9\times9-9-9)+9+9.$
\item [] $586=9\times(9\times9-9-9)+9+9+9/9.$
\item [] $587=9\times(9\times9-9-9)+9+99/9.$
\item [] $588=(99-9/9)\times(9-(9+9+9)/9).$
\item [] $589=9\times9+9+(9\times999-9)/(9+9).$
\item [] $590=9\times9+((9+9)/9)^9-(9+9+9)/9.$
\item [] $591=9\times9+((9+9)/9)^9-(9+9)/9.$
\item [] $592=9\times9+((9+9)/9)^9-9/9.$
\item [] $593=9\times9+((9+9)/9)^9.$
\item [] $594=99\times(99+9)/(9+9).$
\item [] $595=99\times(99+9)/(9+9)+9/9.$
\item [] $596=99\times(99+9)/(9+9)+(9+9)/9.$
\item [] $597=(9+9)\times(9+9+9)+999/9.$
\item [] $598=9\times9\times9-9-(99+999)/9.$
\item [] $599=9999/9-((9+9)/9)^9.$
\item [] $600=(99+9)\times(99+9/9)/(9+9).$
\item [] $601=((9+9)/9)^9+99-9-9/9.$
\item [] $602=((9+9)/9)^9+99-9.$
\item [] $603=((9+9)/9)^9+99-9+9/9.$
\item [] $604=((9+9)/9)^9+99-9+(9+9)/9.$
\item [] $605=99\times(999-9)/(9\times(9+9)).$
\item [] $606=9\times9\times9-(999+99+9)/9.$
\item [] $607=9\times9\times9-(999+99)/9.$
\item [] $608=9\times9\times9-99\times99/(9\times9).$
\item [] $609=9\times9\times9-9-999/9.$
\item [] $610=((9+9)/9)^9+99-9/9.$
\item [] $611=((9+9)/9)^9+99.$
\item [] $612=9\times9\times9-9-9-99.$
\item [] $613=9\times9\times9-9-9-99+9/9.$
\item [] $614=((9+9)/9)^9-9+999/9.$
\item [] $615=9\times9\times9-(999+9+9+9)/9.$
\item [] $616=9\times9\times9-(999+9+9)/9.$
\item [] $617=9\times9\times9-(999+9)/9.$
\item [] $618=9\times9\times9-999/9.$
\item [] $619=9\times9\times9-99-99/9.$
\item [] $620=9\times9\times9-99-9-9/9.$
\item [] $621=9\times9\times9-99-9.$
\item [] $622=9\times9\times9-9-99+9/9.$
\item [] $623=((9+9)/9)^9+999/9.$
\item [] $624=(9-9/9)\times(9\times9-(9+9+9)/9).$
\item [] $625=(9+9+9-(9+9)/9)^((9+9)/9).$
\item [] $626=9+9\times9\times9-(999+9)/9.$
\item [] $627=9\times9\times9+9-999/9.$
\item [] $628=9\times9\times9-99-(9+9)/9.$
\item [] $629=9\times9\times9-99-9/9.$
\item [] $630=9\times9\times9-99.$
\item [] $631=9\times9\times9-99+9/9.$
\item [] $632=9\times9\times9-99+(9+9)/9.$
\item [] $633=9\times9\times9-99+(9+9+9)/9.$
\item [] $634=9\times9\times9-99+(9+9+9+9)/9.$
\item [] $635=9\times99-((9+9)/9)^{(9-9/9)}.$
\item [] $636=9\times(9\times9-9)-(99+9)/9.$
\item [] $637=9\times(9\times9-9)-99/9.$
\item [] $638=9\times(9\times9-9)-9-9/9.$
\item [] $639=9\times(9\times9-9)-9.$
\item [] $640=9\times(9\times9-9)-9+9/9.$
\item [] $641=9\times(9\times9-9)-9+(9+9)/9.$
\item [] $642=9\times(9\times9-9)-9+(9+9+9)/9.$
\item [] $643=9\times(9\times9-9)-(99-9)/(9+9).$
\item [] $644=(9\times9+99/9)\times(9-(9+9)/9).$
\item [] $645=9\times(9\times9-9)-(9+9+9)/9.$
\item [] $646=9\times(9\times9-9)-(9+9)/9.$
\item [] $647=9\times(9\times9-9)-9/9.$
\item [] $648=9\times(9\times9-9).$
\item [] $649=9\times(9\times9-9)+9/9.$
\item [] $650=9\times(9\times9-9)+(9+9)/9.$
\item [] $651=9\times(9\times9-9)+(9+9+9)/9.$
\item [] $652=9\times(9\times9-9)+(9\times9-9)/(9+9).$
\item [] $653=9\times(9\times9-9)+(99-9)/(9+9).$
\item [] $654=9\times(9\times9-9)+9-(9+9+9)/9.$
\item [] $655=9\times(9\times9-9)+9-(9+9)/9.$
\item [] $656=9\times(9\times9-9)+9-9/9.$
\item [] $657=9\times(9\times9-9)+9.$
\item [] $658=9\times(9\times9-9)+9+9/9.$
\item [] $659=9\times(9\times9-9)+99/9.$
\item [] $660=9\times(9\times9-9)+(99+9)/9.$
\item [] $661=9\times(9\times9-9)+(99+9+9)/9.$
\item [] $662=9\times(9\times9-9)+9+(99-9)/(9+9).$
\item [] $663=9\times(9\times9-9)+9+9-(9+9+9)/9.$
\item [] $664=(9\times9+(9+9)/9)\times(9-9/9).$
\item [] $665=9\times(9\times9-9)+9+9-9/9.$
\item [] $666=9\times(9\times9-9)+9+9.$
\item [] $667=9\times(9\times9-9)+9+9+9/9.$
\item [] $668=9\times(9\times9-9)+9+99/9.$
\item [] $669=9\times(9\times9-9)+9+(99+9)/9.$
\item [] $670=9\times(9\times9-9)+(99+99)/9.$
\item [] $671=99\times(9\times9-9-99/9)/9.$
\item [] $672=(9-9/9)\times(9\times9+(9+9+9)/9).$
\item [] $673=9\times(9+9)+((9+9)/9)^9-9/9.$
\item [] $674=9\times(9+9)+((9+9)/9)^9.$
\item [] $675=9\times(9\times9-9)+9+9+9.$
\item [] $676=(9+9+9-9/9)^{((9+9)/9)}.$
\item [] $677=9\times(9\times9-9)+9+9+99/9.$
\item [] $678=999/9+9\times(9\times9-9-9).$
\item [] $679=9\times9\times9-(9\times99+9)/(9+9).$
\item [] $680=9\times9\times9-(9\times99-9)/(9+9).$
\item [] $681=-999/9+99\times(9-9/9).$
\item [] $682=(9-(9+9)/9)\times99-99/9.$
\item [] $683=((9+9)/9)^9+9+9\times(9+9).$
\item [] $684=99\times(9-(9+9)/9)-9.$
\item [] $685=(9-9/9)\times(99-9/9)-99.$
\item [] $686=(9-(9+9)/9)\times(99-9/9).$
\item [] $687=(99-(99+9)/9)\times(9-9/9)-9.$
\item [] $688=9\times9\times9-(9\times9\times9+9)/(9+9).$
\item [] $689=9\times9\times9-(9\times9\times9-9)/(9+9).$
\item [] $690=9\times(99-9)-9-999/9.$
\item [] $691=9\times9\times9-9-9-9-99/9.$
\item [] $692=99\times(9-(9+9)/9)-9/9.$
\item [] $693=99\times(9-(9+9)/9).$
\item [] $694=99\times(9-(9+9)/9)+9/9.$
\item [] $695=99\times(9-(9+9)/9)+(9+9)/9.$
\item [] $696=(9-9/9)\times(99-(99+9)/9).$
\item [] $697=9\times(99-9)-(999+9+9)/9.$
\item [] $698=9\times(99-9)-(999+9)/9.$
\item [] $699=9\times(9\times9+9)-999/9.$
\item [] $700=9\times9\times9-9-9-99/9.$
\item [] $701=9\times9\times9-9-9-9-9/9.$
\item [] $702=9\times(99-9)-99-9.$
\item [] $703=9\times9\times9-9-9-9+9/9.$
\item [] $704=(99-99/9)\times(9-9/9).$
\item [] $705=9\times9\times9-(99+99+9+9)/9.$
\item [] $706=9\times9\times9-(99+99+9)/9.$
\item [] $707=9\times9\times9-(99+99)/9.$
\item [] $708=9\times9\times9-9-(99+9)/9.$
\item [] $709=9\times9\times9-9-99/9.$
\item [] $710=9\times9\times9-9-9-9/9.$
\item [] $711=9\times9\times9-9-9.$
\item [] $712=9\times9\times9-9-9+9/9.$
\item [] $713=9\times9\times9-9-9+(9+9)/9.$
\item [] $714=9\times9\times9-9-9+(9+9+9)/9.$
\item [] $715=9\times9\times9-9-(99-9)/(9+9).$
\item [] $716=9\times9\times9-(99+9+9)/9.$
\item [] $717=9\times9\times9-(99+9)/9.$
\item [] $718=9\times9\times9-9-(9+9)/9.$
\item [] $719=9\times9\times9-9-9/9.$
\item [] $720=9\times9\times9-9.$
\item [] $721=9\times9\times9-9+9/9.$
\item [] $722=9\times9\times9-9+(9+9)/9.$
\item [] $723=9\times9\times9-9+(9+9+9)/9.$
\item [] $724=9\times9\times9-(99-9)/(9+9).$
\item [] $725=9\times9\times9-(9+9+9+9)/9.$
\item [] $726=9\times9\times9-(9+9+9)/9.$
\item [] $727=9\times9\times9-(9+9)/9.$
\item [] $728=9\times9\times9-9/9.$
\item [] $729=9\times9\times9.$
\item [] $730=9\times9\times9+9/9.$
\item [] $731=9\times9\times9+(9+9)/9.$
\item [] $732=9\times9\times9+(9+9+9)/9.$
\item [] $733=9\times9\times9+(9\times9-9)/(9+9).$
\item [] $734=9\times9\times9+(99-9)/(9+9).$
\item [] $735=9\times9\times9+(99+9)/(9+9).$
\item [] $736=9\times9\times9+9-(9+9)/9.$
\item [] $737=9\times9\times9+9-9/9.$
\item [] $738=9\times9\times9+9.$
\item [] $739=9\times9\times9+9+9/9.$
\item [] $740=9\times9\times9+99/9.$
\item [] $741=9\times9\times9+(99+9)/9.$
\item [] $742=9\times9\times9+(99+9+9)/9.$
\item [] $743=9\times9\times9+9+(99-9)/(9+9).$
\item [] $744=9\times9\times9+9+9-(9+9+9)/9.$
\item [] $745=9\times9\times9+9+9-(9+9)/9.$
\item [] $746=9\times9\times9+9+9-9/9.$
\item [] $747=9\times9\times9+9+9.$
\item [] $748=9\times9\times9+9+9+9/9.$
\item [] $749=9\times9\times9+9+99/9.$
\item [] $750=9\times9\times9+9+(99+9)/9.$
\item [] $751=9\times9\times9+(99+99)/9.$
\item [] $752=9\times9\times9+(99+99+9)/9.$
\item [] $753=9\times99-999/9-9-9-9.$
\item [] $754=9\times9\times9+9+9+9-(9+9)/9.$
\item [] $755=9\times9\times9+9+9+9-9/9.$
\item [] $756=9\times9\times9+9+9+9.$
\item [] $757=9\times9\times9+9+9+9+9/9.$
\item [] $758=9\times9\times9+9+9+99/9.$
\item [] $759=9\times(9\times9-9)+999/9.$
\item [] $760=9\times9\times9+9+(99+99)/9.$
\item [] $761=99\times(9\times9-99/9)/9-9.$
\item [] $762=9\times99-9-9-999/9.$
\item [] $763=(9/9+99+9)\times(9-(9+9)/9).$
\item [] $764=9\times9\times9+9+9+9+9-9/9.$
\item [] $765=9\times9\times9+9+9+9+9.$
\item [] $766=9\times9\times9+9+9+9+9+9/9.$
\item [] $767=9\times9\times9+9+9+9+99/9.$
\item [] $768=(9-9/9)\times(99-(9+9+9)/9).$
\item [] $769=9\times99-(99+999)/9.$
\item [] $770=99\times(9\times9-99/9)/9.$
\item [] $771=9\times99-9-999/9.$
\item [] $772=9\times99-9-99-99/9.$
\item [] $773=9\times99-99-9-9-9/9.$
\item [] $774=9\times99-99-9-9.$
\item [] $775=(99-9/9)\times(9-9/9)-9.$
\item [] $776=(9-9/9)\times(99-(9+9)/9).$
\item [] $777=999\times(9-(9+9)/9)/9.$
\item [] $778=9\times9\times9+(9\times99-9)/(9+9).$
\item [] $779=9\times99-(999+9)/9.$
\item [] $780=9\times99-999/9.$
\item [] $781=9\times99-99-99/9.$
\item [] $782=9\times99-99-9-9/9.$
\item [] $783=9\times99-99-9.$
\item [] $784=9\times99-99-9+9/9.$
\item [] $785=9\times99-99-9+(9+9)/9.$
\item [] $786=9\times99-99-9+(9+9+9)/9.$
\item [] $787=9999/9-(9+9)\times(9+9).$
\item [] $788=9\times99+9-(999+9)/9.$
\item [] $789=9\times99+9-999/9.$
\item [] $790=9\times(9\times9+9)-9-99/9.$
\item [] $791=9\times99-99-9/9.$
\item [] $792=9\times99-99.$
\item [] $793=9\times99-99+9/9.$
\item [] $794=9\times99-99+(9+9)/9.$
\item [] $795=9\times99-99+(9+9+9)/9.$
\item [] $796=9\times99-99+(9+9+9+9)/9.$
\item [] $797=9\times(99-9)-(99+9+9)/9.$
\item [] $798=9\times(99-9)-(99+9)/9.$
\item [] $799=9\times(99-9)-99/9.$
\item [] $800=9\times(99-9)-9-9/9.$
\item [] $801=9\times(99-9)-9.$
\item [] $802=9\times(99-9)-9+9/9.$
\item [] $803=9\times99-99+99/9.$
\item [] $804=9\times(99-9)-9+(9+9+9)/9.$
\item [] $805=9\times(99-9)-(99-9)/(9+9).$
\item [] $806=9\times(99-9)-(9+9+9+9)/9.$
\item [] $807=9\times(99-9)-(9+9+9)/9.$
\item [] $808=9\times(99-9)-(9+9)/9.$
\item [] $809=9\times(99-9)-9/9.$
\item [] $810=9\times(99-9).$
\item [] $811=9\times(99-9)+9/9.$
\item [] $812=9\times(99-9)+(9+9)/9.$
\item [] $813=9\times(99-9)+(9+9+9)/9.$
\item [] $814=9\times(99-9)+(9\times9-9)/(9+9).$
\item [] $815=9\times(99-9)+(99-9)/(9+9).$
\item [] $816=9\times9\times9+99-(99+9)/9.$
\item [] $817=9\times9\times9+99-99/9.$
\item [] $818=9\times(99-9)+9-9/9.$
\item [] $819=9\times9\times9+99-9.$
\item [] $820=9\times9\times9+99-9+9/9.$
\item [] $821=9\times(99-9)+99/9.$
\item [] $822=9\times(99-9)+(99+9)/9.$
\item [] $823=9\times(99-9)+(99+9+9)/9.$
\item [] $824=((999+9)/9-9)\times(9-9/9).$
\item [] $825=(999/9-9)\times(9-9/9)+9.$
\item [] $826=9\times9\times9+99-(9+9)/9.$
\item [] $827=9\times9\times9+99-9/9.$
\item [] $828=9\times9\times9+99.$
\item [] $829=9\times9\times9+99+9/9.$
\item [] $830=9\times(9\times9+9)+9+99/9.$
\item [] $831=9\times9\times9-9+999/9.$
\item [] $832=9\times9\times9-9+(999+9)/9.$
\item [] $833=9\times9\times9-9+(999+9+9)/9.$
\item [] $834=9\times9\times9+99+9-(9+9+9)/9.$
\item [] $835=9\times9\times9+99+9-(9+9)/9.$
\item [] $836=999-9\times(9+9)-9/9.$
\item [] $837=999-9\times(9+9).$
\item [] $838=9\times9\times9+99+9+9/9.$
\item [] $839=9\times9\times9+99+99/9.$
\item [] $840=9\times9\times9+999/9.$
\item [] $841=(9+9+99/9)^((9+9)/9).$
\item [] $842=9\times99-(99\times9-9)/(9+9).$
\item [] $843=9\times(99-9)+9\times99/(9+9+9).$
\item [] $844=9\times99-9-9-9-9-99/9.$
\item [] $845=9\times9\times9+99+9+9-9/9.$
\item [] $846=9\times9\times9+99+9+9.$
\item [] $847=9\times9\times9+99+9+9+9/9.$
\item [] $848=9\times9\times9+9+(999-9)/9.$
\item [] $849=9\times9\times9+9+999/9.$
\item [] $850=9\times9\times9+99\times99/(9\times9).$
\item [] $851=9\times9\times9+(99+999)/9.$
\item [] $852=9\times(99+9)-9-999/9.$
\item [] $853=9\times99-9-9-9-99/9.$
\item [] $854=9\times99-9-9-9-9-9/9.$
\item [] $855=9\times99-9-9-9-9.$
\item [] $856=9\times99-9-9-9-9+9/9.$
\item [] $857=9\times9\times9+99+9+9+99/9.$
\item [] $858=9\times9\times9+9+9+999/9.$
\item [] $859=9\times99-9-(99+99+9)/9.$
\item [] $860=9\times99-9-(99+99)/9.$
\item [] $861=9\times(99+9)-999/9.$
\item [] $862=9\times99-9-9-99/9.$
\item [] $863=9\times99-9-9-9-9/9.$
\item [] $864=(99+9)\times(9-9/9).$
\item [] $865=9\times(99+9)-99-9+9/9.$
\item [] $866=9\times99-9-9-9-9+99/9.$
\item [] $867=9\times9\times9+9+9+9+999/9.$
\item [] $868=9\times99-(99+99+9)/9.$
\item [] $869=9\times99-(99+99)/9.$
\item [] $870=9\times99-9-(99+9)/9.$
\item [] $871=9\times99-9-99/9.$
\item [] $872=9\times99-9-9-9/9.$
\item [] $873=9\times(99+9)-99.$
\item [] $874=9\times99-9-9+9/9.$
\item [] $875=9\times99-9-9+(9+9)/9.$
\item [] $876=9\times99-9-9+(9+9+9)/9.$
\item [] $877=9\times99-9-(99-9)/(9+9).$
\item [] $878=9\times99-(99+9+9)/9.$
\item [] $879=9\times99-(99+9)/9.$
\item [] $880=9\times99-99/9.$
\item [] $881=9\times99-9-9/9.$
\item [] $882=9\times99-9.$
\item [] $883=9\times99-9+9/9.$
\item [] $884=9\times99-9+(9+9)/9.$
\item [] $885=9\times99-9+(9+9+9)/9.$
\item [] $886=9\times99-(99-9)/(9+9).$
\item [] $887=9\times99-(9\times9-9)/(9+9).$
\item [] $888=999-999/9.$
\item [] $889=9\times99-(9+9)/9.$
\item [] $890=9\times99-9/9.$
\item [] $891=9\times99.$
\item [] $892=9\times99+9/9.$
\item [] $893=9\times99+(9+9)/9.$
\item [] $894=9\times99+(9+9+9)/9.$
\item [] $895=9\times99+(9\times9-9)/(9+9).$
\item [] $896=(999+9)\times(9-9/9)/9.$
\item [] $897=9\times99+9-(9+9+9)/9.$
\item [] $898=9\times99+9-(9+9)/9.$
\item [] $899=9\times99+9-9/9.$
\item [] $900=9\times99+9.$
\item [] $901=9\times99+9+9/9.$
\item [] $902=9\times99+99/9.$
\item [] $903=9\times99+(99+9)/9.$
\item [] $904=9\times99+(99+9+9)/9.$
\item [] $905=9\times99+9+(99-9)/(9+9).$
\item [] $906=9\times99+9+9-(9+9+9)/9.$
\item [] $907=9\times99+9+9-(9+9)/9.$
\item [] $908=9\times99+9+9-9/9.$
\item [] $909=9\times99+9+9.$
\item [] $910=9\times99+9+9+9/9.$
\item [] $911=9\times99+9+99/9.$
\item [] $912=9\times99+9+(99+9)/9.$
\item [] $913=9\times99+(99+99)/9.$
\item [] $914=9\times99+(99+99+9)/9.$
\item [] $915=999-9\times9-(9+9+9)/9.$
\item [] $916=999-9\times9-(9+9)/9.$
\item [] $917=999-9\times9-9/9.$
\item [] $918=999-9\times9.$
\item [] $919=999-9\times9+9/9.$
\item [] $920=999-9\times9+(9+9)/9.$
\item [] $921=9\times(9\times9+9)+999/9.$
\item [] $922=9\times(99-9)+(999+9)/9.$
\item [] $923=9\times99+9+(99+99+9)/9.$
\item [] $924=9\times(99+99/(9+9+9)).$
\item [] $925=((9+9)/9)^(9+9/9)-99.$
\item [] $926=999-9\times9+9-9/9.$
\item [] $927=999-9\times9+9.$
\item [] $928=999-9\times9+9+9/9.$
\item [] $929=999-9\times9+99/9.$
\item [] $930=(999/9-9-9)\times(9/9+9).$
\item [] $931=9999/9-9\times9-99.$
\item [] $932=9\times99+(9\times9\times9+9)/(9+9).$
\item [] $933=999-9\times9+9+9-(9+9+9)/9.$
\item [] $934=((9+9)/9)^(9+9/9)-9\times9-9.$
\item [] $935=999-9\times9+9+9-9/9.$
\item [] $936=999-9\times9+9+9.$
\item [] $937=999-9\times9+9+9+9/9.$
\item [] $938=999-9\times9+9+99/9.$
\item [] $939=9\times9\times9+99+999/9.$
\item [] $940=9\times99+(9\times99-9)/(9+9).$
\item [] $941=9\times99+(9\times99+9)/(9+9).$
\item [] $942=9\times(99+9+9)-999/9.$
\item [] $943=((9+9)/9)^(9+9/9)-9\times9.$
\item [] $944=9\times(9+99)-9-9-9-9/9.$
\item [] $945=9\times(9+99)-9-9-9.$
\item [] $946=9\times(9+99)-9-9-9+9/9.$
\item [] $947=9\times99+(999+9)/(9+9).$
\item [] $948=(99+9)\times(9\times9-(9+9)/9)/9.$
\item [] $949=9999/9-9\times(9+9).$
\item [] $950=9\times(99+9)-(99+99)/9.$
\item [] $951=9\times(99+9)-9-(99+9)/9.$
\item [] $952=9\times(9+99)-9-99/9.$
\item [] $953=9\times(9+99)-9-9-9/9.$
\item [] $954=9\times(9+99)-9-9.$
\item [] $955=9\times(9+99)-9-9+9/9.$
\item [] $956=9\times(9+99)-9-9+(9+9)/9.$
\item [] $957=99\times(99-(99+9)/9)/9.$
\item [] $958=9+9999/9-9\times(9+9).$
\item [] $959=9\times(99+9)-(99+9+9)/9.$
\item [] $960=(99+9)\times(9\times9-9/9)/9.$
\item [] $961=9\times(99+9)-99/9.$
\item [] $962=9\times(99+9)-9-9/9.$
\item [] $963=9\times(99+9)-9.$
\item [] $964=9\times(99+9)-9+9/9.$
\item [] $965=9\times(99+9)-9+(9+9)/9.$
\item [] $966=9\times(99+9)-9+(9+9+9)/9.$
\item [] $967=9\times(99+9)-(99-9)/(9+9).$
\item [] $968=99\times(99-99/9)/9.$
\item [] $969=9\times(99+9)-(9+9+9)/9.$
\item [] $970=9\times(99+9)-(9+9)/9.$
\item [] $971=9\times(99+9)-9/9.$
\item [] $972=9\times(99+9).$
\item [] $973=9\times(99+9)+9/9.$
\item [] $974=9\times(99+9)+(9+9)/9.$
\item [] $975=9\times(99+9)+(9+9+9)/9.$
\item [] $976=999-(99+99+9)/9.$
\item [] $977=999-(99+99)/9.$
\item [] $978=999-9-(99+9)/9.$
\item [] $979=999-9-99/9.$
\item [] $980=(9+9/9)\times(99-9/9).$
\item [] $981=999-9-9.$
\item [] $982=999-9-9+9/9.$
\item [] $983=9\times(99+9)+99/9.$
\item [] $984=(99+9)\times(9\times9+9/9)/9.$
\item [] $985=999-(99+9+9+9)/9.$
\item [] $986=999-(99+9+9)/9.$
\item [] $987=999-(99+9)/9.$
\item [] $988=999-99/9.$
\item [] $989=999-9-9/9.$
\item [] $990=999-9.$
\item [] $991=999-9+9/9.$
\item [] $992=999-9+(9+9)/9.$
\item [] $993=999-9+(9+9+9)/9.$
\item [] $994=999-(99-9)/(9+9).$
\item [] $995=999-(9+9+9+9)/9.$
\item [] $996=999-(9+9+9)/9.$
\item [] $997=999-(9+9)/9.$
\item [] $998=999-9/9.$
\item [] $999=999.$
\item [] $1000=999+9/9.$
\end{itemize}
\end{multicols}
}

\bigskip
\section{\textbf{Acknowledgement}}

\bigskip
The author is thankful to T.J. Eckman, Georgia, USA (email: jeek@jeek.net) in programming the script to develop these representations.

\begin{center}
-------------------------------------------------
\end{center}
\end{document}